\documentclass{amsart}

\usepackage{mathrsfs}
\usepackage{amsmath,amssymb,enumerate,tikz}
\usepackage{amsfonts,amssymb}
\usepackage{pdfpages}
\usepackage{booktabs}
\usepackage[colorlinks,linkcolor=red,anchorcolor=blue,
citecolor=green,CJKbookmarks=True]{hyperref}
\usepackage{xypic}%

\usepackage{lipsum} %number by letters in the appendix
\usepackage[T1]{fontenc}
\usepackage{csquotes}
\usepackage{soul, xcolor, todonotes}

\newif\ifistoreview
\istoreviewtrue

\theoremstyle{plain}
\newtheorem{theorem}{Theorem}[section]
\newtheorem{lemma}[theorem]{Lemma}
\newtheorem{corollary}[theorem]{Corollary}

\newtheorem{proposition}[theorem]{Proposition}

\theoremstyle{definition}
\newtheorem{definition}[theorem]{Definition}

\newtheorem{remark}[theorem]{Remark}
\newtheorem{remarks}[theorem]{Remarks}

\numberwithin{equation}{section}
\newcommand{\R}{\mathbb R}
\newcommand{\C}{\mathbb C}

\newcommand{\Z}{\mathbb Z}
\renewcommand{\H}{\mathbb{H}}
\newcommand{\PSL}{\mathrm{PSL}}
\newcommand{\SL}{\mathrm{\SL}}
\renewcommand{\d}{\textup{d}}

\newcommand{\MF}{\mathcal{MF}}
\newcommand{\PMF}{\mathcal{PMF}}
\newcommand{\ML}{\mathcal{ML}}
\newcommand{\GL}{\mathcal{GL}}
\newcommand{\PML}{\mathcal{PML}}
\newcommand{\T}{\mathcal{T}}
\newcommand{\M}{\mathcal{M}}

\renewcommand{\i}{\mathrm{i}}

\renewcommand{\ss}{stretch }
\newcommand{\TR}{\mathbf{TR}}
\newcommand{\HR}{\mathbf{HR}}

\newcommand{\SR}{\mathbf{SR}}
\renewcommand{\SL}{\mathbf{SL}}
\newcommand{\HSR}{\mathbf{HSR}}
\newcommand{\HSL}{\mathbf{HSL}}
\newcommand{\ESR}{\mathbf{ESR}}
\newcommand{\PSR}{\mathbf{PSR}}
\renewcommand{\PSL}{\mathbf{PSL}}

\newcommand{\HDR}{{\mathbf{hr}}}
\newcommand{\HDD}{\mathbf{HD}}
\newcommand{\hd}{\mathbf{hd}}
\newcommand{\TD}{\mathbf{TD}}

\newcommand{\Hopf}{\mathrm{Hopf}}

\newcommand{\G}{\mathscr{G}}

\newcommand{\Ext}{\mathrm{Ext}}

\newcommand{\F}{\mathcal{F}}
\renewcommand{\Vert}{\mathrm{Vert}}
\newcommand{\Hor}{\mathrm{Hor}}
\newcommand{\Lip}{\mathrm{Lip}}

\newcommand{\tec}{Teichm\"uller }

\newcommand{\Sing}{\mathrm{Sing}}
\newcommand{\Reg}{\mathrm{Reg}}
\newcommand{\Crit}{\mathrm{Crit}}
\newcommand{\Singt}{\widetilde{\mathrm{Sing}}}
\newcommand{\Regt}{\widetilde{\mathrm{Reg}}}
\newcommand{\Critt}{\widetilde{\mathrm{Crit}}}

\newcommand{\dist}{\mathrm{dist}}
\newcommand{\distfunction}{\mathrm{dist}(u(\cdot),v(\cdot))}

\newcommand{\mh}{\mathcal{H}}
\newcommand{\ml}{\mathcal{L}}

\newcommand{\mw}{\mathcal{W}}
\newcommand{\Hull}{\mathrm{hull}}
\newcommand{\dth}{d_{\mathrm{Th}}}
\renewcommand{\mod}{\mathrm{Mod}}
\newcommand{\Hess}{\mathrm{Hess}}
\title{Ray structures on Teichm\"uller Space}
\author{Huiping Pan}
\address{Huiping Pan,
School of Mathematics, South China  University of Technology\\
Wushan Rd 381, Tianhe, Guangzhou, China, 510641}
\email{panhp@scut.edu.cn}

\author{Michael Wolf}
\address{Michael Wolf,
School of Mathematics\\ 
Georgia Institute of Technology\\ 
Atlanta, GA USA 30332}
\email{mwolf40@gatech.edu}
\date{\today}
\begin{document}

\begin{abstract}
 While there may be  many Thurston metric geodesics between a pair of points in \tec space, we find that by imposing an additional energy minimization constraint on the geodesics, thought of as limits of harmonic map rays, we select a unique Thurston geodesic through those points. Extending the target surface to the Thurston boundary yields, for each point $Y$ in \tec space, an \enquote{exponential map} of rays from that point $Y$ onto \tec space with visual boundary the Thurston boundary of \tec space.

We first depict harmonic map ray structures on \tec space as a geometric transition between \tec ray structures and Thurston geodesic ray structures. In particular, by appropriately degenerating the source of a harmonic map between hyperbolic surfaces (along \enquote{harmonic map dual rays}), the harmonic map rays through the target converge to a Thurston geodesic; by appropriately degenerating the target of the harmonic map, those harmonic map dual rays
through the domain converge to \tec geodesics. We then extend this transition to one from \tec disks through Hopf differential disks to stretch-earthquake disks. These results apply to surfaces with boundary, resolving a question on stretch maps between such surfaces. 
\end{abstract}

%\includepdf[pages={1-2}]{coverletter}
\maketitle
\tableofcontents

%==========================
\section{Introduction}

\subsection{Overview} \tec space admits several ray structures, each arrived at through a similar process.  One begins with a pair of points in \tec space and decides whether to construe those points as complex structures or hyperbolic structures on that surface. Then one chooses a variational problem to solve, the solution to which defines an auxiliary object on the surface, for example a holomorphic differential or a geodesic lamination and a scaling constant. A family of pairs where the auxiliary object is fixed projectively then defines a path in \tec space.

Probably the most famous example of this process are the \tec geodesics.  Here we attempt to minimize the quasiconformal dilatation between a pair of Riemann surfaces, finding a solution that may be described in terms of holomorphic quadratic differentials on each surface; looking for families in which the projective class of such a differential on one of the surfaces is fixed defines a \tec ray.

Thurston \cite{Thurston1998}  developed an analogue of this perspective in the 1980's where he viewed points of \tec space as hyperbolic structures on a surface.  Then one sought minimizers of the Lipschitz stretch of maps between the surfaces, obtained by a map with maximal stretch locus a chain-recurrent geodesic lamination.  When, for example, that lamination was maximal, one could consider the family of target surfaces for which the lamination was  fixed, and this family defined a Thurston geodesic for his (asymmetric) metric.

There are several other families of minimizers, for example Kerckhoff's \enquote{lines of minima} \cite{Kerckhoff1992}, but here we focus on ones related to minimization of $L^2$ Energy. Because of a conformal invariance in two dimensions, this energy depends on the hyperbolic structure of the target but only the conformal (hence complex) structure of the domain, and so occupies a middle ground between the two variational problems we opened with. Here also the solution, a harmonic map, determines and may be understood through the use of a holomorphic quadratic (Hopf) differential, and once again the families of targets who projectively share such a differential foliate \tec space.

As a summary comment, though these ray structures are developed through similar processes, the \tec geodesic rays, the Thurston geodesic rays and the harmonic map Hopf differential rays seem unrelated.

We need one more construction. A somewhat more obscure ray structure comes from seeking a family of Riemann surfaces whose Hopf differentials share one particular (horizontal) projective measured foliation. (See \cite{Tabak1985}.) We can imagine these as defining a type of \enquote{dual} harmonic map rays which we will clarify a bit later.

With all this context in place, we may explain our goals in this paper.  We show a collection of results that together portray the harmonic map rays structures as providing a transition between the \tec ray structures and the Thurston ray structures.  In particular, we show the following, stated roughly in this overview, and in terms of an initial harmonic map 
$u:X \to Y$ and of course in terms of the Hopf differential $\Hopf(X,Y)$.  

First, if we allow the domain $X$ to degenerate along harmonic map dual rays, then the harmonic map rays through $Y$ converge to Thurston geodesics. (See Theorems~\ref{thm:HR:SR} and \ref{thm:HR:limit:geodesic}.)

Second, if we allow the target $Y$ to degenerate along harmonic map rays, the harmonic map dual rays through $X$ converge to \tec geodesics. (See Theorem~\ref{thm:HR:AHR:limit}.)

Third, there are natural ways that the \tec and Thurston rays may be arranged into disks, either as \tec disks in the complex case or into stretch-earthquake disks in the hyperbolic case. Of course, the Hopf differentials also admit a $\C^{\ast}$ action, and the results above on harmonic map rays defining a transition between \tec and Thurston geodesics extend to the disk setting.

Finally, we take up some applications of these results.  We were slightly careful in our description of the Thurston Lipschitz problem to distinguish solutions where the maximally stretched lamination was maximal (so that the complementary regions were ideal triangles) from the more general solutions. Thurston noted that in the non-maximal case, there was no canonical (Thurston) geodesic between the pair of surfaces. There have been a number of proposals for some more canonical geodesics between some pairs of surfaces in some settings (see \cite{PapadopoulosYamada2017,HuangPapadopoulos2019,CalderonFarre2021}). Here we observe that if we additionally require the geodesic from $Y$ to $Z$ to solve an energy-minimization 
problem in the sense of being a limit of harmonic map rays that proceed from $Y$ through $Z$, then this results in a uniquely defined geodesic, a \enquote{harmonic stretch line}. (See Theorem~\ref{thm:unique:geodesic}.)

Moreover, a simple extension of that technique produces a well-defined Thurston geodesic proceeding from a hyperbolic surface to a point on the Thurston boundary of \tec space, represented by a projective measured foliation. Thus in Theorem~\ref{thm:exponentialmap} we show that the harmonic stretch rays from a point $Y$ foliate \tec space and accumulate only at their endpoint on the Thurston boundary: thus they provide an \enquote{exponential} map for the Thurston metric with visual boundary the Thurston boundary. 

 This \enquote{exponential} map for the Thurston metric defines one of two distinct versions we may define for the Thurston geodesic flow on the bundle  over the \tec space with fibers 
 the space $\PMF$ of projective measured foliations{: here} the orbit of a hyperbolic surface $X$ and a projective measured foliation $[\eta]$ is  the harmonic stretch line which passes through $X$ and converges to $[\eta]$ in the Thurston boundary in the forward direction.  The other version of the Thurston geodesic flow is defined such that the orbit of $(X,[\eta])$ is  the harmonic stretch line arising as the limit of harmonic map rays through $X$ when the domain degenerates along the harmonic map dual ray determined by $X$ and $\eta$. The first version of the Thurston geodesic flow describes the contracting foliation along its orbits (which is also the endpoint on the Thurston boundary) while the second version describes the stretching lamination. Both of these two flows commute with the mapping class group action, and hence descend to the bundle over moduli space with fiber $\PMF$.

Along the way, we show some other results, for example an extension (Theorem \ref{thm:arc:Lip}) of this theory to surfaces with geodesic boundary, proving that the optimal Lipschitz constants between hyperbolic surfaces with boundaries are always realized by some surjective Lipschitz homeomorphisms. This verifies a conjecture of Alessandrini and Disarlo \cite{AlessandriniDisarlo2019} (in the case of unpunctured surfaces).

In several places, we prove not only subconvergence of the approximating rays but actual convergence of the families.  The technique here often comes down to regarding a harmonic map as an equivariant minimal graph over a hyperbolic surface with values in an $\R$-tree, and we prove a uniqueness theorem in that case. (See Theorem~\ref{thm:minimalgraph:uniqueness}.) The difficulties here are both that the tree typically does not admit an orientation and the graphs we imagine have infinite diameter, but our result might still be regarded as in the spirit of the Jenkins-Serrin uniqueness theorem (see \cite{JS66}). Having broached this topic, in the appendix we provide some completeness to the discussion by also proving a corresponding existence theorem for graphs with values in a real tree and asymptotic boundary values over hyperbolic domains with geodesic boundary  (that result will also play a role in the proof of Lemma~\ref{lem:analyticity1}).  (Here, by extending the range to trees, we extend the possible domains to surfaces whose polygonal ends are not limited to an even number of sides, as in \cite{JS66} and elsewhere.)

In the next section, we define our terms and state our results somewhat more carefully.

 \subsection{Harmonic map rays and harmonic stretch rays.}
 The remainder of this introductory section is devoted to stating our results. We quickly declare some notation:  a fuller treatment of the constructions behind those definitions is given in Section~\ref{sec:definitions} and Section~\ref{sec:harmonic:maps}.
 
 \subsubsection{Four old families of rays and one new one.}\label{subsec:fourrays} Let $S$ be an orientable closed surface of genus at least two and $\T(S)$ be the \tec space of $S$ (see Section \ref{subsec:teich:space} for the definition of $\T(S)$). A \tec ray $\TR_{X,\Phi}:[1,\infty)\to\T(S)$ is defined by a base Riemann surface $X$ and a quadratic differential $\Phi$ which is holomorphic on $X$ (cf. Section \ref{subsec:teich}).  It is convenient to also denote that ray by $\TR_{X,\lambda}:[1,\infty)\to\T(S)$, where $\lambda$ is the horizontal measured (or projective measured) lamination of $\Phi$.
 
 A Thurston stretch ray $\SR_{Y,\lambda}:[1,\infty)\to\T(S)$ is defined with base point a hyperbolic surface $Y$ and a maximal geodesic lamination $\lambda$ (cf. Section \ref{subsec:stretch:maps}).
 
 A harmonic map ray $\HR_{X,\Phi}:[1,\infty)\to\T(S)$ is the family of surfaces so that the harmonic map from $X$ to any member of that family has Hopf differential proportional to $\Phi$, a holomorphic quadratic differential on $X$ (cf. Definition \ref{def:HR}).  We will also have need of the notation $\HR_{X,Y}:[1,\infty)\to\T(S)$ which indicates the restriction of that harmonic map ray to begin at the hyperbolic surface $Y$.
 
 The most obscure previously defined ray we will consider initially is the harmonic map dual ray $\HDR_{Y,\lambda}=\HDR_{Y,\lambda}(t)$, defined by the condition that the horizontal measured foliation of the Hopf differential $\Hopf(\HDR_{Y,\lambda}(t)\to Y)$ is measure equivalent to the measured lamination $t\lambda$ (cf. Definition \ref{def:HDR}). Alternatively, $\HDR_{Y,\lambda}(t)$ is the conformal structure underlying the projection of the minimal surface which is a graph over $\widetilde Y$ in the product $\widetilde Y \times T_{t\lambda}$, where $T_{t\lambda}$ is the tree dual to the lift of the lamination $t\lambda$.
 
 Eventually we will define a family of \enquote{harmonic  stretch lines} (see Definition \ref{defn:harmonic:stretch:line}), a distinguished family of Thurston geodesics arising as limits of harmonic map rays.
 
For convenience, we collect all types of rays and lines used in this paper in Table \ref{table}, with names, notations, and references.

 \begin{table}
 \begin{tabular}{lcr}
  \toprule
  Names of rays & Notations & Reference
  \\ \midrule
  \tec ray & $\TR_{X,\Phi}$, $\TR_{X,\lambda}$, or $\TR_{X,Y}$ & Section \ref{subsec:teich} 
  \\\midrule    Thurston stretch ray or line & $\SR_{X,\lambda}$ or $\SL_{X,\lambda}$ & Section \ref{subsec:stretch:maps}  
  \\\midrule    Harmonic map ray & $\HR_{X,\Phi}$, $\HR_{X,\lambda}$, or $\HR_{X,Y}$ & Definition \ref{def:HR}
  \\\midrule    Harmonic map dual  ray & $\HDR_{Y,\lambda}$ or $\HDR_{X,Y}$ & Definition  \ref{def:HDR}
  \\\midrule    Harmonic stretch ray or line & $\HSR_{X,Y}$ or $\HSL_{X,Y}$ & Definition  \ref{defn:harmonic:stretch:line} 
  \\\midrule Piecewise harmonic stretch ray or line &$\PSR_{Y,\lambda,f}$ or $\PSL_{Y,\lambda,f}$& Theorem \ref{thm:generalized:stretchmap}
  \\ \bottomrule
\end{tabular}
  \caption{Various rays and lines in this paper.}\label{table}
 \end{table}

 \subsubsection{Foundational convergence results} We now state our first results on asymptotic relationships between these ray families.
 
 \begin{theorem}\label{thm:HR:SR}
   Let $Y\in\T(S)$ be a hyperbolic surface and $\lambda$ a measured lamination.  Then the harmonic map rays $\HR_{X_t,Y}:[1,\infty)\to\T(S)$ converge to a Thurston geodesic locally uniformly as $X_t$ diverges along the harmonic map dual ray $\HDR_{Y,\lambda}:[1,\infty)\to\infty$. 

   Moreover, if $\lambda$ is maximal, then the limit Thurston geodesic is exactly the Thurston stretch ray $\SR_{Y,\lambda}:[1,\infty)\to\T(S)$ defined by $Y$ and $\lambda$. 
 \end{theorem}

\begin{remark}\label{rmk:convergence:TR}
   It is natural to ask whether the result still holds if $X_t$ degenerates along the \tec ray $\HR_{X_t, Y}:[1,\infty)\to\T(S)$. The answer is affirmative if $\lambda$ is maximal (see Theorem \ref{thm:convergence:TR}). (We believe that answer holds in general, and we hope to take up this issue in a future paper.)
\end{remark}

 We also have the following compactness result.
 \begin{theorem}\label{thm:HR:limit:geodesic}
  For any fixed $Y\in\T(S)$, let $X_n\in\T(S)$ be any divergent sequence. Then the sequence of harmonic map rays $\HR_{X_n,Y}:[1,\infty)\to\T(S)$ contains a subsequence which converges to some  Thurston geodesic   locally uniformly.
 \end{theorem}
 
 %Geodesics obtained from above  are called \enquote{harmonic-stretch} lines.

 The limit Thurston geodesics in Theorem \ref{thm:HR:SR} for non-maximal measured lamination will be characterized later as instances of a special class -- called the \enquote{harmonic stretch rays} (Definition \ref{defn:harmonic:stretch:line}) -- of \enquote{piecewise harmonic stretch rays}   (see Theorem~\ref{thm:generalized:stretchmap}, Section \ref{sec:convergence:stretchrays}).
 
 \subsubsection{Harmonic map dual rays and Teichm\"uller rays}

 Informally, these results we just presented assert that if we look at a family of harmonic map rays $\HR_{X_s,Y}$ through a hyperbolic structure $Y$ and then degenerate the domains $X_s$ appropriately, say along  a harmonic map dual ray, then the rays limit on a Thurston geodesic.  
 
 It is natural to ask, what happens to the limits of harmonic map rays if we degenerate in the other direction, say by letting the target surface, now $Y_s$, tend to infinity? We should expect some object to emerge that is defined in terms of the complex structure $X$. Roughly, the answer is that the relevant rays converge to \tec geodesics through $X$. More formally, the rays through $X$ to focus on are the \enquote{harmonic map dual rays} defined in Subsection~\ref{subsec:fourrays}. We then show 
 that if we degenerate appropriately, this time letting the target $Y_s$ degenerate along a harmonic map ray, then the harmonic map dual rays through $X$ increasingly approximate \tec geodesics. 
 
 Precisely, fix a complex structure $X\in \T(S)$. Let $\Phi$ be a holomorphic quadratic differential on $X$.  Let $Y_s:=\HR_{X,\Phi}(s)$; we will be degenerating the harmonic map dual rays defined by $X$ and $Y_s$.

 To that end, let $\lambda\in\ML(S)$ be the measured lamination which is measure equivalent to the horizontal foliation of $\Phi$. This data defines a family of harmonic map dual rays through $X$ as follows.

 We say that a {\it harmonic map dual ray} $\HDR_{Y, \lambda}$ is the set of points $X_t$ so that $\Hor (\Hopf(X_{t}\to Y))=t\lambda$ (cf. \cite{Tabak1985}). An alternative, perhaps more geometric, description (\cite{Wolf1998}) is that the universal cover $\widetilde{X_{t}}$ is the minimal surface graph over the universal cover $\widetilde{Y}$ in the product $\widetilde{Y} \times tT_{\lambda}$, where $T_{\lambda}$ is the real tree dual to $\lambda$ (here of course lifting equivariantly to universal covers).

 Extending this definition to encompass a family of such harmonic map dual rays, indexed by $s$, we set $X_{s,t}\in \HDR_{Y_s,\sqrt{s}\lambda}$ to be such that $\Hor (\Hopf(X_{s,t}\to Y_s))=t\sqrt{s}\lambda$.  Then $X_{s,1}\equiv X$  for all $s>0$.  (Again the more geometric description is that $\widetilde{X_{s,t}}$ is the minimal surface graph over $\widetilde{Y_s}$ in the product $\widetilde{Y_s} \times t\sqrt{s} T_{\lambda}$, where $T_{\lambda}$ is the real tree dual to $\lambda$.)
 
 Then in that setting, we prove our convergence result.
 
 \begin{theorem}\label{thm:HR:AHR:limit}
    The family of harmonic map dual rays $$\HDR_{Y_s,\sqrt{s}\lambda}:[1,\infty)\to\T(S)$$ converges locally uniformly to the Teichm\"uller geodesic ray $\TR_{X,\Phi}:[1,\infty)\to\T(S)$,  as $s\to+\infty$. 
 \end{theorem}

%==========================

%==========================

 %=====================================
 \subsection{Piecewise harmonic stretch lines and harmonic stretch lines}

 Theorems~\ref{thm:HR:SR}  and \ref{thm:HR:limit:geodesic} provide for a collection of distinguished Thurston geodesics, those that arise as limits of harmonic map rays.  We turn next to describing the locations of these special geodesics in \tec space (and its Thurston compactification), finding existence and uniqueness theorems for such geodesics between given points in those spaces.

 To begin, Theorem \ref{thm:HR:SR} enables us to generalize Thurston's construction of stretch lines (maps) from maximal laminations to non-maximal ones using harmonic maps. 
 
A basic tool for us will be \enquote{piecewise harmonic stretch maps}. These maps will have the same structure of the limits of the harmonic maps defined by a harmonic map ray -- and indeed they will often arise as such limits -- but they are not required {\it a priori} to have any compatibility between components.

Later we will combine their properties with a uniqueness result (cf. Theorem~\ref{thm:minimalgraph:general} below) to argue for Theorem~\ref{thm:HR:SR}, with the convergence property in that theorem relying on a compatibility across the components inherited from the limiting process.
  
  \begin{theorem}[piecewise harmonic \ss map]\label{thm:generalized:stretchmap}
    Let  $Y\in\T(S)$ be any closed hyperbolic surface, and let  $\lambda$ be any geodesic lamination. Then for any harmonic diffeomorphism $f:X\to Y\setminus \lambda$ from some (possibly disconnected) punctured surface $X$,  there is a new hyperbolic surface
    \begin{equation*}
     Y_t:= \PSL_{Y,\lambda,f}(t)\in\T(S)
    \end{equation*}
    depending analytically on $\{t>0\}$ such that
    \begin{enumerate}[(a)]
    \item the induced map $f_t:X\to Y_t\setminus\lambda$ is a   harmonic diffeomorphism with  Hopf differential  $t\Hopf(f)$;
      \item for any $0< s<t$, the map ($f_t\circ f_s^{-1}$)  extends to a homeomorphism from $Y_s$ to $Y_t$ that is $\sqrt{t/s}$-Lipschitz with (pointwise) Lipschitz constant strictly less than  $\sqrt{t/s}$ in $Y_s\backslash\lambda$, but
       exactly expands arc length on $\lambda$ by the constant factor $\sqrt{t/s}$.
    \end{enumerate}
  \end{theorem}
  The family of hyperbolic structures $\PSL_{Y,\lambda,f}(t)$  constructed above is called a {\it piecewise harmonic \ss  line}. It admits a canonical orientation coming from the orientation of the  positive real ray $\{t>0\}$. In that orientation, a piecewise harmonic stretch line is a (reparametrized) geodesic in the Thurston metric.  Whenever we say a piecewise harmonic stretch line, we mean a directed line.
  Similarly to Thurston's construction of concatenation of stretch lines, one can construct a geodesic between any two points in Teichm\"uller space by a concatenation of piecewise harmonic stretch line segments.

Now, in addition to the piecewise harmonic stretch lines described in Theorem~\ref{thm:generalized:stretchmap} above, an important focus for us will be families in a subclass we call {\it harmonic stretch lines}.  A piecewise harmonic stretch line is called a {\it harmonic stretch line} if it is the limit of a sequence of harmonic map rays.
 Given hyperbolic surfaces $Y,Z \in \T(S)$, a harmonic stretch line passing through $Y$ to $Z$ will be a limit of a family of harmonic map rays from base points $X_n \in \T(S)$ that all proceed through   $Y$ to $Z$,  as the base points $X_n$ degenerate in $\T(S)$. Of course, as in the case of the piecewise harmonic stretch lines, these harmonic stretch lines are also directed.
 
 Our basic theorem on these harmonic stretch lines is the following.

  \begin{theorem}[Uniqueness of harmonic stretch lines]\label{thm:unique:geodesic}
    For any two distinct hyperbolic surfaces $Y, Z\in \T(S)$, there exists a unique harmonic stretch line proceeding from $Y$ through $Z$. 
  \end{theorem}

\begin{remark}
 From the construction of Thurston stretch lines and piecewise harmonic stretch lines, we see that every Thurston stretch line is a piecewise harmonic stretch line (see Lemma \ref{lem:id:PSR:SR}). However, Thurston stretch lines are not necessarily harmonic stretch lines. In fact, a Thurston stretch line is a harmonic stretch line if and only if the defining maximal lamination is {\it chain-recurrent}, see Corollary \ref{cor:Thurston:harmonic:stretch}. 
\end{remark}

\begin{remark}\label{rmk:boundary:puncture}
 The constructions above extend to the case when $Y$ is a hyperbolic surface with geodesic boundary components. Using the uniqueness result Theorem~\ref{thm:unique:geodesic}, we verify a conjecture (cf. Section~\ref{sec:concluding}) of Alessandrini-Disarlo  (\cite[Conjecture 1.8]{AlessandriniDisarlo2019}) in the case of unpunctured surfaces on the existence of stretch homeomorphisms with non-maximal stretch laminations. 
\end{remark}

\begin{remark}
    In a subsequent paper \cite{PW2024}, we use the existence and uniqueness of harmonic stretch lines to describe the \enquote{envelope} $\mathrm{Env}(X,Y)$ of Thurston geodesics from $X\in \T(S)$ to $Y\in \T(S)$.
\end{remark} 

  \subsection{Visual boundary of the Thurston metric and an exponential map}
  
  Finally, we extend the existence/uniqueness theory of harmonic stretch lines to rays whose terminal point is a projective measured lamination, representing a point on the boundary of the Thurston compactification of the \tec space $\T(S)$.  Properties of these rays, which are also Thurston geodesics, allow us to construct an \enquote{exponential} map on \tec space from any base point $Y\in \T(S)$, with the boundary $\PML(S)$ of the Thurston compactification appearing as the visual boundary for this family of rays.

   \begin{theorem} \label{thm:exponentialmap}
   For any $Y\in \T(S)$ and  any $[\eta]\in\PML(S)$, there exists a unique harmonic stretch ray  starting at $Y$, which converges to $[\eta]\in\PML(S)$ in the Thurston compactification.
   
   Moreover, these rays foliate $\T(S)$ if we fix $Y$ and let $[\eta]$ vary in $\PML(S)$, or if we fix $[\eta]$ and let $Y$ vary in $\T(S)$. 
 \end{theorem}
 
 \begin{remark}\label{rmk:exp}
    For any $X\in\T(S)$, the set of unit vectors tangent to Thurston stretch lines (whose stretch lamination is maximal) has  Hausdorff measure zero in $T^1_X\T(S)$ (\cite[Theorem 10.5]{Thurston1998}), while the set of unit vectors tangent to harmonic stretch lines is exactly $T^1_X\T(S)$ (see Proposition \ref{prop:existence-exponential-lines}).
 \end{remark}

\begin{remarks}\label{rmk:clarification}
   So far, we have introduced several types of geodesics of the Thurston metric. (i) Harmonic stretch rays (lines) are special cases of piecewise harmonic stretch rays (lines).    (ii) Those Thurston geodesics considered in Theorem \ref{thm:HR:SR}, Theorem \ref{thm:HR:limit:geodesic},  Theorem \ref{thm:unique:geodesic}, and Theorem \ref{thm:exponentialmap} are harmonic stretch rays (lines) (see the second statement of Proposition \ref{prop:limit:HR:PSL} and Definition \ref{defn:harmonic:stretch:line}). (iii) For maximal geodesic laminations, piecewise harmonic stretch lines and Thurston stretch lines coincide (Lemma \ref{lem:id:PSR:SR}).
\end{remarks}

 \subsection{Geodesic flow for the Thurston metric}
 The Teichm\"uller geodesic flow on the moduli space
 has been extensively studied in the literature. However, the notion of geodesic flow does not exist naturally for the Thurston metric.  A natural challenge by Rafi (\cite[Problem 3.10]{Su2016}) is to introduce a suitable notion of geodesic flow for the Thurston metric.  Here we respond to this question by defining two versions of the geodesic flow for the Thurston metric. 
 
 Theorem \ref{thm:exponentialmap} allows us to define  the Thurston geodesic flow  $$\psi_t: \T(S)\times \PML(S)\longrightarrow \T(S)\times\PML(S)$$
 such that the orbit through $(Y,[\eta])\in \T(S)\times \PML(S)$ is the harmonic stretch line determined by $Y$ and $[\eta]$ via Theorem \ref{thm:exponentialmap}. Moreover, every harmonic stretch line appears as a (forward) orbit.

 There is another version of the Thurston geodesic flow  
 $$\phi_t:\T(S)\times\PML(S) \to \T(S)\times\PML(S)$$
 such that the orbit through $(Y,[\eta])\in \T(S)\times\PML(S)$ is the harmonic stretch line  obtained from Theorem \ref{thm:HR:SR} as the limit of harmonic map rays $\HR_{X_t,Y}$ when $X_t$ diverges along the harmonic map dual ray $\HDR_{Y,\eta}$.  

 Both of these flows are mapping class group equivariant, and so they descend to $\M(S)\times \PML(S)$, where $\M(S)$ is the moduli space. Moreover, the earthquake flow and $\phi_t$ are compatible in the sense that earthquake translates of $\phi_t$-orbits are still $\phi_t$-orbits, and hence define an action of the upper-triangular subgroup of $SL(2,\R)$, see Proposition \ref{prop:compatibility}. We imagine it could be interesting to study the dynamical properties (invariant measures, ergodicity, and mixing) of these two versions of the Thurston geodesic flow  over the moduli space.

 %====================
 \subsection{Minimal graphs}

 An important tool for us will be conformal harmonic graphs over a hyperbolic domain with values in a real tree. Such maps will enjoy uniqueness properties which we can leverage to prove that the families of harmonic map rays, such as those in Theorem~\ref{thm:HR:SR}, have a unique limit for any chosen subsequence.  The graphs we study have infinite diameter and are extensions to the singular setting of the properly embedded classical minimal surfaces studied by Jenkins-Serrin \cite{JS66} (see also \cite{CollinRosenberg2010}). A complication is that because the trees do not, in general, fold, some of the usual arguments only partially generalize. Our principal result is the following theorem.
 
 \begin{theorem}\label{thm:minimalgraph:general}
Let $Y$ be a crowned hyperbolic surface.  Let $T$ be an admissible dual tree of $Y$. Then for any prescribed admissible partial boundary map there exists a unique  $\pi_1(Y)$-equivariant minimal graph over $\widehat{Y}$ in $\widetilde{Y}\times {T}$, where $\widetilde{Y}$ is the universal cover of $Y$. 
\end{theorem}

We define \emph{admissble partial boundary map} just after Definition~\ref{def:admissibleTree:general}.

We apply the uniqueness portion of this result in our central results and so prove it within the body of the paper.  Having opened the discussion of these sorts of graphs, we treat the existence portion in an appendix.

\subsection{Extension of results from Rays to Disks.} 
 So far, we have focused on rays, i.e. sets defined in terms of a single real parameter. Yet an important topic in \tec theory are \tec disks, i.e. images of \tec maps from $X$ whose holomorphic quadratic differentials (on $X$) are proportional by a complex number. 
 
 On the hyperbolic geometric side, there are also distinguished disks, defined in terms of a measured geodesic lamination $\lambda$ on a surface $Y$, and two operations that deform a hyperbolic structure that each use $\lambda$.  Naturally, one may stretch the structure along the lamination, as we have described throughout this section.  One can also perform an earthquake along this lamination. The two operations, defined with data $X$ and a measured lamination $\lambda$, together define a real two-dimensional family of hyperbolic surfaces through $X$, called a {\it stretch-earthquake disk}.
 
 We state results that assert that the Hopf differential disks well-approximate stretch-earthquake disks for nearly degenerate domains $X$, and that there is a reasonable notion of dual Hopf differential disks that well-approximate \tec disks for nearly degenerate ranges $Y$. 
 
  \subsubsection{Convergence to Thurston stretch-earthquake disks} 
  
  We begin with the approximation of stretch-earthquake disks.
 
 Let $Y\in\T(S)$  be a hyperbolic surface and $\lambda\in \ML(S)$ a measured foliation/lamination.  Let $\PSL_{Y,\lambda}$ be the piecewise harmonic stretch line obtained from Theorem \ref{thm:HR:SR} as the limit  of harmonic map rays $\HR_{X_t,Y}$ where $X_t$ diverges along the harmonic map dual ray $\HDR_{Y,\lambda}$. Let $\mathcal{E}^s_{\lambda}(Y)$ be the surface obtained from $Y$ by a time $s$ earthquake along $\lambda$. Define the {\it stretch-earthquake disk} $\mathbf{SED}(Y,\lambda)$ of $(Y,\lambda)$ to be the set:
   \begin{equation*}
     \bigcup_{-\infty<s<+\infty}\PSL_{\mathcal{E}^s_{\lambda}(Y),\lambda}({0},+\infty).
   \end{equation*}

 We wish to say, roughly, that the images of Hopf differential disks through a hyperbolic surface $Y$ converge to a stretch-earthquake disk through $Y$ -- as we let the center $X$ of the disks degenerate appropriately.  To state this properly, we quickly introduce a bit more notation. Let $X_t\in\HDR_{Y,\lambda}$ be the Riemann surface such that the horizontal foliation of  $\Hopf(X_t,Y)=\Phi_t$ is $t\lambda$.
 Let $Y(t,r,s)$ be the hyperbolic surface such that $\Hopf(X_t,Y(t,r,s))=re^{i\frac{s}{2t}}\Phi_t$ and $Y_r=\PSL_{Y,\lambda}(r)\in \PSL_{Y,\lambda}$.
 For such data $(Y, X_t, \Phi_t)$, the {\it Hopf differential disk} $(\mathbf{HDD}(X_t,\Phi_t),Y)$ comprises surfaces $Z$ satisfying $\Hopf(X_t, Z)=\zeta\Phi_t$ for some $\zeta \in \C$. (Naturally, $Y\in \mathbf{HDD}(X_t,\Phi_t),Y)$ for the choice of $\zeta = 1$ in the above.)

 In this language, we may state our result on the convergence of disks in the harmonic setting to disks in the hyperbolic geometric setting.

 \begin{theorem}\label{thm:Hopfdisk:limit1}
    Let $Y\in\T(S)$  be a hyperbolic surface and $\lambda\in \ML(S)$ a measured foliation/lamination. Then, for $X_t$ as above, the family of Hopf differential disks $(\mathbf{HDD}(X_t,\Phi_t),Y)$ with base point $Y$ locally uniformly converge to the stretch-earthquake disk $(\mathbf{SED}(Y,\lambda),Y)$ with base point $Y$. Namely, for any prescribed $\mathbf{s}>0$ and $0<\mathbf{r}<\mathbf{r}'$, the family
    $ Y(t,r,s)$ of surfaces converges to $\mathcal{E}^s_{\lambda}(Y_r)$ uniformly in $(r,s)\in[\mathbf{r},\mathbf{r}']\times [-\mathbf{s},\mathbf{s}]$ as $t\to\infty$.
 \end{theorem}
 
 %=========================
  \subsubsection{Convergence to \tec disks}
  
  On the other hand, we may look for a family of disks defined via harmonic maps that well-approximate \tec disks. There are two possibilities for approximates, but we focus in this introduction on just the first: here we consider an $S^1$ family of harmonic map dual rays defined by a point $Y$ very far from $X$ in \tec space, and then show that this \enquote{disk} of dual rays converges to a classical \tec disk as $Y$ diverges.

 More precisely, begin with a surface $Y$ and, fixing a domain $X$, a Hopf differential $\Phi= \Hopf(X,Y)$. Let $Y_{s,\theta}$ be the hyperbolic surface such that $\Hopf(X,Y_{s,\theta})=se^{2i\theta}\Phi$.  In particular, $Y_{1,0}=Y$ and $Y_{s,\theta}$ diverges as $s \to \infty$ (for any choice of $\theta$). Let
      $$\hd_{X,\Phi,s}=\bigcup_{0\leq\theta\leq\pi} \HDR_{Y_{s,\theta},\lambda_{\theta}}$$
      denote the (variable target) harmonic map dual disk, where $\lambda_\theta$ is the horizontal foliation of $e^{2i\theta}\Phi$.
 
Let $\TD_{X,\Phi}$ be the Teich\"uller disk determined by $X$ and $\Phi$. 

\begin{theorem}\label{thm:convergence:disks:v1}
   $\hd_{X,\Phi,s}$ locally uniformly converges to the Teichm\"uller disk $\TD_{X,\Phi}$, as $s\to+\infty$. 
\end{theorem}

In Section~\ref{sec:convergence:TeichDisks}, we also provide for a version of convergence to \tec horodisks.
  
%======================
\subsection{Previous results}  A number of authors have studied how harmonic maps might approximate \tec maps, or how Thurston stretch maps might be approached from an analytic perspective.

Before Bers' exposition \cite{Bers1960:TeichmullersTheorem} of a proof of Teichm\"uller's theorem, Gerstenhaber and Rauch (\cite{GerstenhaberRauchI1954}, \cite{GerstenhaberRauchII1954}) began a program to find \tec maps as limits of energy-minimizing maps, optimized over conformal metrics.  This program was completed by Mese in \cite{Mese2004}.

 In the other direction, Daskalopoulos and Uhlenbeck \cite{DaskalopolousUhlenbeck,DaskalopolousUhlenbeck2022} focused on developing an analytic approach to finding least stretch maps: they build on literature on least gradient maps (see for example \cite{SWZ1992,MRL2014,SZ1993}), but they also develop the analytic foundations of what they term \enquote{$J_p$ harmonic maps}.

 In \cite{BonMonSch2013} and \cite{BonMonSch2015}, Bonsante, Mondello and Schlenker develop an $S^1$ action on $\T(S)\times \T(S)$, the \enquote{landslide flow},  which limits to the earthquake flow when one of the parameters goes to a measured lamination in the Thurston boundary of $\T(S)$. Certainly, these landslides, as smooth approximates to earthquakes, have many interesting properties. Similarly to our analysis, these landslides $L(h, h^{\ast})$ are defined in terms of the circle action on the Hopf differential disks, here centered at a Riemann surface $c$ that is the \enquote{midpoint} of the ray between $h$ and $h^{\ast}$; the limit to the earthquake flow on the second variable in $\T(S)\times \T(S)$ occurs as the first variable (not the midpoint) decays appropriately in \tec space.  In contrast, described in terms of the language of \cite{BonMonSch2013} and \cite{BonMonSch2015}, here we study the limits of the Hopf differential disks as the center $c$ degenerates along a dual ray; these are more specific paths than those studied by these authors, but allow us to show the convergence we require in our setting. The precise relationship between the two constructions is not clear. 

There have been a number of studies comparing various rays in \tec space. Choi, Rafi and Series show that short curves along lines of minima coincide with short curves along \tec geodesics (with the same defining laminations) (\cite{CRS2008b}), and  they prove that every line of minima is a  \tec quasi-geodesic (\cite{CRS2008a}); Choi, Dumas and Rafi prove that every grafting ray is a \tec quasi-geodesic which stays within a bounded distance of some \tec geodesic (\cite{CDR2012}); Gupta proves that every grafting ray is asymptotic to a \tec geodesic (\cite{Gupta2014,Gupta2015}); Lenzhen, Rafi and Tao compare the \tec geodesics and Thurston geodesics in \cite{LRT2012}.

\subsection{Organization of the paper} We provide a basic background and define our notation in Sections 2 and 3. Section~\ref{sec:compactness} is devoted to finding some subsequential limits of a family of harmonic maps $f_n:X_n \to Y$ from a degenerating sequence $X_n$ of Riemann surfaces: crucial to our studies here and elsewhere is Minsky's \cite{Minsky1992} convex regions $\mathscr{P}_R$, and we include some first results here.  In Section~\ref{sec:GeneralizedJenkinsSerrin} we develop a key tool for our uniqueness and convergence claims:  a type of Jenkins-Serrin uniqueness theorem for graphs over hyperbolic domains with values in a real tree. Sections~\ref{sec:subconvergence:stretchrays}, \ref{sec:stretchline:construction}, and \ref{sec:convergence:stretchrays} provide a proof of Theorem~\ref{thm:HR:SR}:  subconvergence is proved in Section~\ref{sec:subconvergence:stretchrays}, we prove Theorem~\ref{thm:generalized:stretchmap} in Section~\ref{sec:stretchline:construction} and use that result and the uniqueness theorem for minimal graphs (from section~\ref{sec:GeneralizedJenkinsSerrin}) to reach the final convergence result in Section~\ref{sec:convergence:stretchrays}. 

We change foci in the next two sections. In  Sections~\ref{sec:convergence:TeichRays} and \ref{sec:convergence:TeichDisks}, we consider limits of harmonic map dual rays when the range of the harmonic map degenerates, and find that \tec rays and disks emerge:  we prove Theorems~\ref{thm:HR:AHR:limit} and \ref{thm:convergence:disks:v1} (as well as a second version of this convergence). 

We then return to limits of the harmonic map rays and disks as the domain surface degenerates for the rest of the paper.  In section~\ref{sec:convergence:Stretch-EarthquakeDisk}, we prove the convergence result Theorem~\ref{thm:Hopfdisk:limit} on limits of families of harmonic map rays which define an earthquake-stretch disk. Then in Section~\ref{sec:stretchray:uniqueness}, we prove our basic existence and uniqueness result for harmonic stretch lines, our refinement of Thurston geodesics for non-maximal laminations. Sections~\ref{sec:exponential} and \ref{sec:concluding} offer consequences of our study, including a proof of Theorem~\ref{thm:exponentialmap} in Section~\ref{sec:exponential}.

The paper concludes by resolving a dangling issue:  we prove in the Appendix an existence result for our Jenkins-Serrin problem to complement the uniqueness result Theorem~\ref{thm:minimalgraph:uniqueness} that we used throughout.
%=====================

\subsection{Acknowledgements} The first author is supported by National Natural Science Foundation of China NSFC 12371073 and NSFC 11901241. The second author gratefully acknowledges support from the Simons Foundation and the U.S. NSF grant DMS-2005551 and DMS-2429005. He also benefited from a casual insightful remark of Yair Minsky and a conversation with J\'er\'emy Toulisse about \cite{Tholozan2017}. 

The authors wish to express their deep gratitude to the referee for the very careful reading and thoughtful remarks.
\section{Definitions and background stretch maps and Teichm\"uller maps} \label{sec:definitions}

%====================
%====================
\subsection{\tec space}\label{subsec:teich:space} 
 Let $S$ be an oriented closed surface of genus $g\geq2$. A \textit{marked Riemann surface} is a pair $(X,f)$ where $X$ is a Riemann surface and $f:S\to X$ is an orientation preserving homeomorphism. Two marked Riemann surfaces $(X,f)$ and $(X',f')$ are called \textit{equivalent} if there exists a conformal map in the homotopy class of $f'\circ f^{-1}$.  The \tec space $\T(S)$ is then defined as the space of equivalence classes of marked Riemann surfaces. Topologically, $\T(S)$ is homeomorphic to $\mathbb{R}^{6g-6}$.

 The \tec space can also be defined via hyperbolic surfaces.
 By the uniformization theorem, every Riemann surface of genus at least two admits a unique conformal metric of  constant Gaussian curvature $-1$, the hyperbolic metric.  A \textit{marked hyperbolic surface} is a pair $(X,f)$ where $X$ is a hyperbolic surface and $f:S\to X$ is an orientation preserving homeomorphism. Two marked {hyperbolic surfaces} $(X,f)$ and $(X',f')$ are called \textit{equivalent} if there exists an isometry in the homotopy class of $f'\circ f^{-1}$.  The \tec space $\T(S)$ is then defined as the space of equivalence classes of marked hyperbolic surfaces.
 
For simplicity, we denote the (equivalence class of the) marked Riemann/hyperbolic surface $(X,f)$ by $X$. Whenever we say a map from $X=(X,f)$ to $X'=(X',f')$, we mean it is a map homotopic to $f'\circ f^{-1}$, the change of marking.  (An alternative perspective on equivalence is acquired by pulling $X$ back to $S$ by the homeomorphism $f$: then any equivalent pair $(X,f)$ and $(X',f')$ may be compared via the map $f^{-1} \circ f'$, which is homotopic to the identity map on $S$. We will often compare structures on $S$ via maps of $S$ that are homotopic to the identity.)

 \subsection{Quadratic differentials, measured foliations and measured laminations} \label{subsec:foliation:lamination}
 Given a closed Riemann surface $X$ of genus $g$, a  holomorphic quadratic differential $\Phi$ is a holomorphic section of $K_X^2$, where $K_X$ is the holomorphic cotangent bundle over $X$ and $K_X^2$ is a tensor product of $K_X$.  The space of holomorphic quadratic differentials on $X$, denoted by {$H^0(X, K_X^2)$}, is a vector space of real dimension $6g-6$.
 
 Every holomorphic quadratic differential $\Phi$ defines two measured foliations, the \textit{horizontal foliation} $\mathrm{Hor}(\Phi)$ by the imaginary  part of $\sqrt{\Phi}$, and the \textit{vertical foliation} $\mathrm{Vert}(\Phi)$ by the real part of $\sqrt{\Phi}$. Conversely, every measured foliation on $X$ is realized as the horizontal/vertical foliation of a unique holomorphic quadratic differential on $X$ (\cite{HM1979,Wolf1996}). The space of measured foliations on $S$, denoted by $\MF(S)$, is homeomorphic to $\R^{6g-6}$ (see for instance \cite[Corollary 7.8]{FLP}). 
   
   The hyperbolic counterpart of foliations are geodesic laminations. A \textit{geodesic lamination} on a hyperbolic surface $X$ is a closed subset which is a disjoint union of complete simple geodesics, called \textit{leaves}. Typical examples are simple closed geodesics. A geodesic lamination $\lambda$ is said be \textit{maximal} if the complementary regions on $X$ are ideal triangles. Consider the space of geodesic laminations on $S$ equipped with the Hausdorff topology of the space of closed subsets on $X$. This is a compact metric space with the Hausdorff distance.  Since there is a canonical correspondence between the geodesics of any two marked hyperbolic metrics on a closed surface $S$,  the space of geodesic laminations, denoted by $\GL(S)$,  depends only on the topology of $S$. A collection of disjoint simple closed geodesics is called a \textit{multicurve}.
   A geodesic lamination is called \textit{chain-recurrent} if it can be approximated by multicurves in $\GL(S)$ in the Hausdorff metric on $\GL(S)$. 
  
  A \textit{measured geodesic lamination} (or \textit{measured lamination} for short) is a geodesic lamination equipped with a \textit{transverse invariant measure}, which associates to every arc transverse to the lamination a {Radon measure} (see  \cite[Page 11-13]{Bonahon2001} for more details). 
  Typical examples are simple closed geodesics with the Dirac measures. The \textit{intersection number} between simple closed geodesics extends naturally to the setting of measured geodesic laminations as follows. 
   Given a measured geodesic lamination $\lambda$ and a simple closed geodesic $\alpha$, the intersection number, denoted by $i(\lambda,\alpha)$, is defined to be the mass of $\alpha$ given by the transverse measure of $\lambda$.   This associates to every  measured geodesic lamination  a function over the set of simple closed geodesics.    Let $\ML(S)$ be the space of measured geodesic laminations, equipped with the weak-* topology of the space of functions over   the set of simple closed geodesics.  Topologically, $\ML(S)$ is homeomorphic to $\mathbb{R}^{6g-6}$ (see \cite[Theorem 6.1]{FLP} or \cite[Theorem 3.1.1]{PennerHarer} ).   
   There is a natural correspondence between the space of measure foliations and the space of measured geodesic laminations by straightening leaves of foliations using hyperbolic geodesics ( \cite{Levitt1983}).

 %====================
 %====================
 \subsection{Extremal length and hyperbolic length}\label{subsec:ExtremalHyperbolicLength}

 Given the Riemann surface $X$ and a simple closed curve $\alpha$, the \textit{extremal length} $\Ext_X(\alpha)$ is defined as:
 $$ \Ext_X(\alpha)=\sup_{\rho}
 \frac{\ell^2_\rho(\alpha)}{\mathrm{Area}(\rho)}, $$
 where  the supremum  ranges over all conformal metrics $\rho$ on $X$, and where $\ell_\rho(\alpha)$ is the infimum of the length of curves (free) homotopic to $\alpha$. Alternatively, we can also define the extremal length via the conformal modulus of cylinders. Any embedded cylinder $C$ in $X$ inherits a conformal structure from $X$ and is conformal to a unique flat cylinder, up to a scaling.  The \emph{conformal modulus} (or \emph{modulus} for short) $\mod(C)$ of $C$ is defined as the height of flat cylinder divided by the circumference. Accordingly, we have $$\Ext_C(\alpha)=\frac{1}{\mod(C)},$$
 where $\alpha$ is the core curve of $C$. For Riemann surface $X$ and any essential simple closed curve $\alpha$ on $X$, we have
 \begin{equation}\label{eq:ext:mod}
 	\Ext_X(\alpha)=\inf_C\frac{1}{ \mod (C)},
 \end{equation}
 where the infimum ranges over all embedded cylinders with core curves homotopic to $\alpha$  (see \cite[Section 3]{Kerckhoff1980} for an explanation of \eqref{eq:ext:mod}). In practice, one obtains a lower bound on  the  extremal length via the first formulation and an upper bound via the second formulation.

 From the perspective of hyperbolic geometry, one can define the \textit{hyperbolic length function} $\ell_X(\alpha)$ by associating to $\alpha$ the  length of the unique geodesic representative with respect to the hyperbolic metric of $X$.  Both the extremal length function and the hyperbolic length function can be extended to measured laminations and measured foliations (\cite{Kerckhoff1980,Kerckhoff1983}).

%=====================
%=====================
\subsection{\tec maps and \tec rays} \label{subsec:teich}

 The \tec map between {a pair of Riemann surfaces} is the solution to the variational problem of finding the minimizing quasiconformal constant within a given homotopy class. Let $X\in\T(S)$ be a Riemann surface and $\Phi$ a holomorphic quadratic differential on $X$. Near the regular points of $\Phi$, there are natural coordinates in which one  represents $\Phi$ as $dz^2=(dx+idy)^2$. For any $k\geq1$, consider   the quadratic differential $\Phi_k$ locally defined by $(k^{1/2} dx+i k^{-{1/2}}dy)^2.$   It defines a unique Riemann surface  $X_k\in\T(S)$ on which $\Phi_k$ is holomorphic.  The ray $\TR_{X,\Phi}:[1,\infty)\to\T(S)$ sending $k\geq 1$ to $X_k$ is  called a \textit{\tec ray}.  (Note that the parametrization of the Teichm\"uller ray (as well as the Thurston stretch ray later in Section \ref{subsec:stretch:maps}) is not the commonly used geodesic parametrization. The reason we use the current parametrization is to make parametrizations of various rays consistent with each other.) A {\it Teichm\"uller line} arises from allowing $k$ in the construction of $X_k$ to take on all positive values.

 With respect to the natural coordinates of $\Phi$ and $\Phi_t$, the identity map $X\to X_t$ is a \tec map. Conversely, all \tec maps arise in this way. Namely, for any two distinct Riemann surfaces $X,Y\in\T(S)$, there exist a unique holomorphic quadratic differential $\Phi$ on $X$ (the \textit{initial quadratic differential}) and  a unique holomorphic quadratic differential $\Psi$ on $Y$ (the \textit{terminal quadratic differential} ),  such that the (unique) \tec map $X\to Y$ is locally defined by $(x,y)\mapsto (k^{1/2}x,k^{-1/2}y)$ with respect to the natural coordinates of $\Phi$ and $\Psi$, where $k>1$ is a positive constant (see \cite{Bers1960:TeichmullersTheorem}, \cite[Chapter 6]{Gardiner1987} or \cite[Chapter 5]{IT1992} for the existence and uniqueness of \tec maps between Riemann surfaces of finite type). Consider the \tec ray determined by $X$ and the initial quadratic differential $\Phi$. We denote it by $\TR_{X,Y}$, indicating that its initial point is at $X$ and it passes through $Y$.

 Given $X,Y\in\T(S)$, the \tec distance $d_T$ on $\T(S)$ is defined as
 $$ d_T(X,Y)=\frac{1}{2}\log K, $$
 where $K$ is the infimum of quasiconformal constants among quasiconformal maps from $X$ to $Y$ in the homotopy class determined by the markings of $X$ and $Y$. Another formulation is given by Kerckhoff \cite{Kerckhoff1980} in terms of extremal lengths:
 $$  d_T(X,Y)=\frac{1}{2}\log \sup_{\alpha}
 \frac{\Ext_Y(\alpha)}{\Ext_X(\alpha)}, $$
 where the supremum ranges over all simple closed curves. In this regard, the projective class of the  horizontal foliation of the initial quadratic differential of the Teichm\"uller map from $X$ to $Y$ is characterized as the unique one realizing the maximum $\max_{\mu\in\MF(S)} \frac{\Ext_Y(\mu)}{\Ext_X(\mu)}$.

 Teichm\"uller lines are geodesics under the Teichm\"uller metric. Conversely, every geodesic under the Teichm\"uller metric is a Teichm\"uller line. Given a  \tec line $\TR_{X,\Phi}$, for any two points $X_t,X_s$ in this line,  we have
 $ d_T(X_t,X_s)=|\log \sqrt{s/t} |$.

%==================
%==================

\subsection{Remark on ray and disk structures on \tec space.}\label{subsec:rays} The paper concerns \enquote{ray structures} on \tec space, where we imagine a number of different ways of defining the \enquote{rays}.  Here, by ray structure, we will mean that for every point $Y \in \T(S)$,
there is a family of parametrized halflines with initial point at $Y$ which together foliate \tec space. Probably the most famous such ray structure is that of \tec rays from $Y$ parametrized by rays in the vector space $H^0(Y, K_Y^2)$.  That complex vector space $H^0(Y, K_Y^2)$ also has complex lines in it, and we imagine a \enquote{disk structure} on \tec space $\T(S)$ to be an organization of a ray structure with some naturally defined free $S^1$ action on the rays.

 %=====================
 %=====================

\subsection{Stretch maps}\label{subsec:stretch:maps}
  A seminal paper \cite{Thurston1998} defined a new (Finsler) metric (\enquote{Thurston's asymmetric metric}) on \tec space, characterized it in terms of the variational problem of finding the minimizing Lipschitz stretch between two hyperbolic surfaces, described a solution to that variational problem, and used the solution in describing shortest paths in \tec space for the metric.  In addition to Thurston's original paper, other  sources for exposition of these topics, and further refinements of the theory, are \cite{PapadopoulosTheret2007, ChoiRafi2007, LRT2012,Walsh2014, Pan2020,  DLRT2020,HOP2021,Papadopoulos2021,PanSu2024}. 
  
 A stretch map relies for its definition upon   a maximal geodesic lamination $\lambda$ on a hyperbolic surface.   Recall that a maximal geodesic lamination has complementary regions which are ideal triangles.  An ideal triangle admits a partial foliation by horocyclic arcs centered at the ideal points which both respects the symmetries of the ideal triangle and also meets the complete boundary of the ideal triangle.  One defines a $K$-stretch map of the ideal triangle as a $K$-Lipschitz map which preserves the horocyclic foliation, fixes the three boundary points common to the foliations about each ideal point and stretches distance by $K$ along the bounding geodesics (see \cite[Proposition 3.2]{Thurston1998}). 

 Remarkably, Thurston's careful analysis shows that these stretch maps of complementary domains of the maximal geodesic lamination extend to  a $K$-Lipschitz {map}  between hyperbolic surfaces, and moreover, this extended map has the least possible Lipschitz constant among all maps in its homotopy class (\cite{Thurston1998}, Corollary~4.2).

 We introduce some notation to summarize these constructions.  If $Y$ is a hyperbolic surface and $\lambda$ is a maximal geodesic lamination, then we set $\SR_{Y,\lambda}(t)$ to define the surface constructed from $Y$ and $\lambda$ by a stretch map of Lipschitz constant $\sqrt{t}\geq1$  along the lamination.  We set $\SR_{Y,\lambda}(1, \infty)$ (resp. $\SL_{Y,\lambda}(0,\infty)$) to be the family of hyperbolic surfaces obtained by allowing $t$ in $\SR_{Y,\lambda}(t)$ to range over all values in $[1, \infty)$ (resp. $(0,\infty)$), and we  refer to that set as a {\it Thurston stretch ray} (resp. \emph{Thurston stretch line}). 
  
  Given $X,Y\in \T(S)$, the  Thurston (asymmetric) distance  is defined by 
 $$ \dth(X,Y):=\log L, $$
 where $L$ is the infimum of Lipschitz constants of Lipschitz homeomorphisms from $X$ to $Y$ in the homotopy class determined by the markings of $X$ and $Y$. In particular, for $1\leq s\leq t$, we see that $\dth(\SR_{Y,\lambda}(s),\SR_{Y,\lambda}(t))=\log\sqrt{t/s}$. Thurston also characterized this distance in terms of ratios of length functions of simple closed curves:
 \begin{equation}
 \label{eq:thurstonmetric}
   \dth(X,Y)=\log\sup_\alpha
 \frac{\ell_Y(\alpha)}{\ell_X(\alpha)},  
 \end{equation} 
 where the supremum ranges over all simple closed curves on $X$. 

\subsubsection{Maximally stretched laminations} 
\label{subsec:maximally:stretched:lamination}
 As we mentioned earlier, for any two different Riemann surfaces $X,Y\in\T(S)$,  the maximal ratio of extremal lengths $\max_{\mu\in\MF(S)}\frac{\Ext_Y(\mu)}{\Ext_X(\mu)}$ is uniquely realized by the projective class of the horizontal measured foliation of the initial quadratic differential of the Teichm\"uller map from $X$ to $Y$. For the ratio of hyperbolic length functions, the uniqueness is no longer true  in the setting of measured laminations.   If a non-uniquely ergodic measured lamination $\mu$ realizes the maximal ratio $\max_{\mu\in\ML(S)}\frac{\ell_Y(\mu)}{\ell_X(\mu)}$, then any measured lamination with the same support as $\mu$ also realizes the maximal ratio. Nevertheless, Thurston showed that there is still a well-defined object in this setting, the \textit{maximal maximally stretched  chain-recurrent lamination} (or the \textit{maximally stretched lamination} for short).  For every two different hyperbolic surfaces $X$ and $Y$, the maximally stretched lamination, denoted by $\Lambda(X,Y)$, is defined as the union of   chain-recurrent laminations to which the restriction of every optimal Lipschitz map from $X$ to $Y$ is an affine map with the optimal Lipschitz constant.
  
\begin{theorem}[\cite{Thurston1998}, Theorem 8.4] \label{thm:maximally:stretched:lamination}
 Let $g$ and $h$ be any two distinct hyperbolic structures on $S.$ If $g_{i}$ and $h_{i}$ are sequences of hyperbolic structures converging to $g$ and $h$ respectively, then $\Lambda(g, h)$ contains any lamination in the limit set of $\Lambda\left(g_{i}, h_{i}\right)$ in the Hausdorff topology.
\end{theorem}
 In particular, if $\Lambda(g_i,h_i)$ contains a simple closed curve $\alpha_i$ which converges to a maximal lamination $\lambda$ in the Hausdorff topology, then $\Lambda(g,h)=\lambda$ because the only  chain recurrent lamination containing $\lambda$ is $\lambda$ itself.

\subsubsection{Generalized stretch maps and their rays} Not all of the laminations that will arise in our discussions will be maximal;  indeed, a focus of this paper is the non-maximal case. Corresponding to the case when a quadratic differential will have higher order zeroes, we will find ourselves considering versions of stretch maps where the lamination will have complementary domains on a surface which are hyperbolic ideal polygons, or more generally, the so-called crowned hyperbolic surfaces.    We extend the definition of stretch maps to this case, called {\it piecewise harmonic stretch maps}, see Section \ref{sec:stretchline:construction}; we will also define in Section~\ref{sec:stretchray:uniqueness} a special class of Thurston geodesics which we call {\it harmonic stretch rays}.

%===================
%===================
\subsection{Crowned hyperbolic surfaces}
\begin{definition}
[Crown end] A \emph{crown} with $m \geq 1$ {boundary cusps} is an open and incomplete hyperbolic surface bounded by a closed geodesic boundary $c$, and a \emph{crown end} comprising bi-infinite geodesics $\left\{\gamma_{i}\right\}_{1 \leq i \leq m}$ arranged in a cyclic order, such that the right half-line of the geodesic $\gamma_{i}$ is asymptotic to the left half-line of geodesic $\gamma_{i+1}$, where $\gamma_{m+1}=\gamma_{1} .$ 
\end{definition}

\begin{definition}
 A \emph{crowned hyperbolic surface} is obtained by attaching crowns to a compact hyperbolic surface with geodesic boundaries by isometries along some of their closed boundaries and removing the remaining geodesic boundaries. This results in an open and incomplete hyperbolic metric of finite area on the surface. Topologically, the underlying surface is the interior of a compact surface. \end{definition}

\remark The definition of crowned hyperbolic surfaces here is slightly different from the definition in \cite{Gupta2017}. The definition in \cite{Gupta2017} requires all ends to be crown ends.

 A {\it truncation} of a crown $\mathcal C$ with $m$ boundary cusps is obtained from $\mathcal C$ by removing a collection of disjoint horodisk neighborhoods  $U_1,U_2,...,U_m$, one at each  ideal vertex of  $\mathcal C$.
\begin{definition}\label{def:metric:residue}
  The {\it metric residue} of the hyperbolic crown $\mathcal C$ with $m$ boundary cusps is defined to be zero when $m$ is odd, and equal to {the absolute value of} the alternating sum of lengths of geodesic sides of a truncation when $m$ is even. 
\end{definition}
\begin{remark}
	The metric residue does not depend on the choice of truncation.
\end{remark} 
%====================
%====================
 \subsection{Horizontal foliations of meromorphic quadratic differentials} 
 \label{subsec:meromorphic:leaves}
 Let $\Phi$ be a meromorphic quadratic differential on a compact Riemann surface $X$.

\begin{definition}\label{def:critical:graph}
    Let $X$ and $\Phi$ be as above. A \emph{saddle connection} is a finite $|\Phi|$-geodesic with endpoints at zeroes or simple poles of $\Phi$ and whose interior points are regular points of $\Phi$. A \emph{horizontal saddle connection} of $\Phi$ is a segment of a horizontal leaf with endpoints on the zeroes or simple poles of $\Phi$. The \emph{critical graph} of the horizontal foliation of $\Phi$ is the closure of the union of horizontal saddle connections of $\Phi$ and critical half-infinite leaves entering poles of order at least two of $\Phi$. 
\end{definition}
In particular, the critical graph of the horizontal foliation of $\Phi$ does not contain any recurrent (half-infinite) critical leaf of $\Phi$. Consider a component $X_0$ complementary to the critical graph of $\Phi$.  Then there are five possibilities (\cite[Section 11.4]{Strebel1984}):  
  \begin{enumerate}[(i)]
      \item $X_0$ is a horizontal cylinder of finite height foliated by closed horizontal leaves, or  
       \item  $X_0$ is precompact inside the complement in $X$ of poles of $\Phi$ of order at least two but is not a horizontal cylinder, or 
      \item $X_0$ is a half-infinite cylinder foliated by closed horizontal leaves, or
      \item $X_0$ is an infinite strip of finite height foliated  by bi-infinite horizontal leaves, or  
      \item $X_0$ is a half-plane foliated by bi-infinite horizontal leaves.
  \end{enumerate} 
 Items (i) and (iii) correspond to  \cite[Section 11.4, type (1)]{Strebel1984}; item (ii) corresponds to \cite[Section 11.4, type (2c)]{Strebel1984};  items (iv) and (v) correspond to  \cite[Section 11.4, type (2b)]{Strebel1984}.
  Finite cylinders in the first item and the closures of leaves in the second item are called \textit{compact components} of $\Hor(\Phi)$.
  Half-planes, strips and half-infinite cylinders are called \textit{non-compact components} of  $\Hor(\Phi)$.
  Notice that a strip may spiral to a non-horizontal half-infinite cylinder, or be parallel to some half-plane, or both.

\section{Harmonic maps}\label{sec:harmonic:maps}
In this section, we provide a brief overview of harmonic maps between surfaces. 

\subsection{Harmonic maps}\label{subsec:harmonic:map}
Let $(M,\sigma(z)|dz|^2)$ and $(N,\rho(w)|dw|^2)$ be two Riemannian surfaces. A differentiable map $f:M\to N$ is said be \textit{harmonic} if it satisfies the \textit{Euler-Lagrange Equation}
$$ f_{z\bar{z}}+ (\log \rho)_w f_zf_{\bar{z}}=0.  $$
If $M$ and $N$ are compact, then $f$ is harmonic if and only if it is a critical point of the energy functional: 
 $$ E(f):=\int_M e_f(z)\sigma(z) dzd\bar{z} ,$$
 where $e_f(z):=\frac{\rho(f(z))}{\sigma(z)}(|f_z|^2+|f_{\bar{z}}|^2).$ 
 Notice that the energy depends on the conformal structure on $M$ and the metric on $N$. The energy of the map $f:X\to Y$ is denoted  by $E(f)=E(f:X\to Y)=E(X,Y)$ depending on the context that is required.

The basic existence result  of harmonic maps in a homotopy class was established by Eells and Sampson in  \cite{EellsSampson1964} and by  Hamilton in \cite{Hamilton1975},    if the target manifold has nonpositive sectional curvature.   The uniqueness of harmonic maps in a homotopy class was obtained by Al'ber \cite{alber1964} and Hartman \cite{Hartman1967} if the target manifold has negative sectional curvature and if the image is not contractible to a point or a geodesic.  Moreover, Sampson \cite{Sampson1978} and Schoen-Yau \cite{SchoenYau1978} proved that any harmonic map between compact surfaces which is homotopic to a diffeomorphism is a diffeomorphism, provided that the target surface has nonpositive curvature. 

\subsection{Hopf differentials} \label{subsec:Hopf} 
Let $X$ and $Y$ be two hyperbolic surfaces. 
Let 
\begin{equation}\label{eq:pullback:metric:0}
    f:(X,\sigma(z)|dz|^2)\to (Y,\rho(w)|dw|^2)
\end{equation} be a harmonic diffeomorphism.   Consider the pullback  of $\rho$ by $f$:
$$ f^*(\rho)(z)= \rho(z) f_z \overline{f_{\bar{z}}}dz^2+e_f(z)\sigma(z) dz d\bar{z} +\rho(z) \overline{f_z}f_{\bar{z}} d\bar{z}^2.$$
 The $(2,0)$-part of $f^*(\rho)$ is called the \textit{Hopf differential} of $f$. The harmonicity of  $f$ implies that the Hopf differential of $f$ is holomorphic (see \cite{Schoen1984,Jost1984}).
 The Hopf differential for the map $f: X\to Y$ is denoted by $\Hopf(f) = \Hopf(f:X\to Y) = \Hopf(X, Y) = \Phi_X(Y) = \Phi$, depending on the context that is required. 
 
 Recall that every holomorphic quadratic differential defines two measured foliations on $X$, the horizontal foliation and the vertical foliation.
 The leaves of the horizontal foliation  (resp. vertical foliation) of $\Phi:=\Hopf(f)$ are exactly the maximally stretched  (resp. minimally stretched) directions of $f$. Choosing a local coordinate  $z=x+iy$ on $M$ such that the leaves of the horizontal foliation (resp. vertical foliation) of $\Hopf(f)$ are tangent to the $x$-axis (resp. $y$-axis). Then, in this coordinate, $\Phi=|\Phi|dz^2$. Hence,  
 \begin{equation}\label{eq:pullback:metric}
     f^*\rho~=|\Phi|dz^2+e_f\sigma dzd\bar{z}+|\Phi|d\bar{z}^2=(e_f\sigma+2|\Phi|) dx^2+(e_f\sigma-2|\Phi|)dy^2.
 \end{equation} 
 In particular, if we choose the coordinate $z=x+iy$ such that $\Phi=dz^2$ and choose $\sigma$ to be the singular flat metric induced by $|\Phi|$, then $f^*(\rho)$ can be simply expressed as 
 $$ f^*\rho=(e_f+2)dx^2+(e_f-2)dy^2. $$ 
 Let $\nu(z):=\frac{f_{\bar{z}}d\bar{z}}{f_zdz}$ be the Beltrami differential of $f$. Set $\mathscr{G}(z)=\log(1/|\nu(z)|)$.
 By calculation, we see that
 $\cosh\mathscr{G}=\frac{\sigma e_f}{2|\Phi|}$.
 Substituting this into  (\ref{eq:pullback:metric}) yields
 \begin{equation}\label{eq:pullback:metric:2}
     (f^*\rho)(z)=2|\Phi(z)|(\cosh\mathscr{G}(z)+1)dx^2+2|\Phi(z)|(\cosh\mathscr{G}(z)-1)dy^2. 
 \end{equation}

 %=================
 %=================
\subsection{Harmonic maps to trees and minimal suspensions} \label{subsec:minimalsuspension}
We begin with a harmonic map $u: X \to Y$ with Hopf differential $\Phi = \Hopf(u)$. We lift the setting to the universal cover with a map $\tilde{u}:\tilde{X} \to \tilde{Y}$ and a Hopf differential $\tilde{\Phi} = \Hopf(\tilde{u})$. We consider the projection $p: \tilde{X} \to T_{h}$ from $\tilde{X}$ to the leaf space of the horizontal foliation $\Hor(\tilde{\Phi})$.  That leaf space $T= T_{h}$, consisting of equivalence classes of points on the universal cover $\widetilde S$ of $S$ where two points are equivalent if they are contained in a connected leaf  of $\Hor(\tilde{\Phi})$ (including leaves which branch at zeroes of $\widetilde{\Phi}$), has the structure of a tree, with topology induced from $\Hor(\tilde{\Phi})$.  The tree acquires a well-defined distance $d= d_{T_{h}}$ from the push-forward $p_*\mu_{\tilde{\Phi}, h}$ of the lift $\mu_{\tilde{\Phi}, h}$ of the measure $\mu_{\Phi, h}$ on arcs transverse to the horizontal foliation of $\Phi$. The metric tree $(T,d)$ is not locally compact \cite[Section 2.4]{Wolf1996} when the genus of $X$ is at least two, even if the vertices are typically of finite valence.

By construction, the map $p: \tilde{X} \to (T_{h}, 2d)$ is harmonic in the senses of \cite{Wolf1996} (Definition 2.1) and \cite{KorevaarSchoen1993} (p. 643, p. 656). 
 In the former case, it is a reflection that distances on a tree are submean near a vertex and proven (Proposition~3.1) within the paper, while in the latter case, it is an easy computation. In that case, the energy of the map to the tree is stationary, as a domain variation leads to a variation of energy expressed as the integrated product of the resulting infinitesimally trivial Beltrami differential for the domain variation against the Hopf differential of the projection, a pairing known to vanish by an integration by parts (cf. \cite{Ahlfors1961}); then, further, because the map $p$ is a locally a projection, a local variation in the target along the map induces a domain variation, here using that because the tree is NPC, we may restrict to local variations that stay within the convex finite subtree that is the image of a small disk. Hence as both $\tilde{u}: \tilde{X} \to \tilde{Y}$ and $p: \tilde{X} \to T_{h}$ are harmonic, so is the product map $(\tilde{u},p): \tilde{X} \to \tilde{Y} \times T_{h}$.  We call the latter product map the {\it minimal suspension}  of the harmonic map $u$ because, by construction, the product map $(\tilde{u},p)$ is conformal, and since it is also harmonic, it is therefore minimal. The minimal suspension is stable in an appropriate sense (see \cite{Wolf1998} for this definition and further details), and indeed, given a tree $(T,d)$ dual to a measured foliation $\F$, there is a unique minimal suspension $(\tilde{u},p): \tilde{X} \to \tilde{Y} \times (T, 2d)$ so that both $\tilde{u}$ and $p$ are harmonic with Hopf differentials that are additive inverses (so that $(\tilde{u},p)$ is conformal and harmonic).  A key portion of the present work, principally in Section~\ref{sec:GeneralizedJenkinsSerrin} but also in the Appendix, may be seen as a less technically restrictive approach to the case where $X$ is compact and also an extension to the case where $X$ is complete but with punctures. Not formally related to any of these results but partially aligned in spirit is recent work of Markovic \cite{Markovic2021} in higher codimension.

To the equivariant harmonic map $p: \tilde{X} \to (T_h, 2d)$ we associate the \emph{equviariant energy} $E(X,T_h)$, defined as the integral of energy density over a fundamental domain on $\tilde X$ of the fundamental group $\pi_1(X)$. Choose local coordinate $z=x+iy$ on $\tilde X$ such that $\tilde \Phi=dz^2$. With respect to this coordinate, the map $p$ is represented as $(x,y)\mapsto 2y$. Accordingly, the equivariant energy $E(X,T_h)$ is:
\begin{equation}\label{eq:energy:tree}
	E(X,T_h)=\int_{\tilde X/\pi_1(S)} 2dxdy= 2\|\Phi\| =2\Ext_X(\Hor(\Phi)),
\end{equation} 
where $\Hor(\Phi)$ is the horizontal measured foliation of $\Phi$.

\subsection{Harmonic map rays and dual rays}\label{subsec:HR:HDR}
This paper centers on harmonic map rays in \tec space as an interpolating structure between Thurston stretch rays and \tec rays.  In this subsection, we define these rays in two ways: the first is in terms of the Hopf differential, and the second is in terms of the minimal suspension in the previous subsection.  The latter definition then suggests a dual construction of a different sort of ray structure which we call harmonic map dual rays.

\subsubsection{Harmonic map rays.} If $\Phi \in H^0(X, K_X^2)$ is a holomorphic quadratic differential on the Riemann surface $X$, then we are led to study the ray $s\Phi \subset H^0(X, K_X^2)$ for $s>0$. By the identification of \tec space $\T(S)$ with $H^0(X, K_X^2)$ via $Y \in \T(S) \mapsto \Hopf(X,Y)$ (cf. \cite{Wolf1989}, \cite{Hit87}), we then obtain a family $Y_s \subset \T(S)$ of surfaces for which $\Hopf(X,Y_s) = s\Phi$ for some non-trivial element $\Phi \in H^0(X, K_X^2)$.
\begin{definition}\label{def:HR}
A family of hyperbolic surfaces $Y_s$ for which 	$\Hopf(X,Y_s) = s\Phi$ for some non-trivial element $\Phi \in H^0(X, K_X^2)$ and $s\geq 0$ is called a \emph{harmonic map ray}  and is denoted $\HR_{X,\Phi}(s)$ or sometimes $\HR_{X,Y}(s)$ when we want to describe the ray passing between $X$ and $Y$. In particular, $Y_0=X$.
\end{definition}

Note that the ray $Y_s = \HR_{X,\Phi}(s)$  describes a family of minimal suspensions in $\widetilde{Y_s} \times (T_\Phi, 2s^{\frac{1}{2}}d)$  all with the same second factor $T_\Phi$  up to a scaling of the metric on the tree by $s^{\frac{1}{2}}$. In particular, we can imagine the ray as determining a change in the  first factor $\widetilde{Y_s}$ of the two factors $\widetilde Y$ and $T$  so that the minimal surface's conformal structure $\widetilde{X}$ is held constant while the second factor is scaled.

Regarding a minimal suspension $X \subset Y \times (T,2d)$ as having the three variables of the two factors $Y$ and $T$ as well as the conformal structure $X$, we see that a new ray structure naturally presents itself: we may fix the first factor and the projective class of the second factor and then let the conformal structure vary as the scale of the tree varies.

\begin{definition}\label{def:HDR}
The \emph{harmonic map dual ray} $\HDR_{Y, \Phi}(t)$ determined by a hyperbolic surface $Y\in\T(S)$ and holomorphic quadratic differential $\Phi$ is  the family of conformal structures $X_t$  so that $\Hor (\Hopf(X_{t} \to Y)) = t\Hor({\Phi})$, i.e. the horizontal foliation of the Hopf differential of the harmonic map from $X_t$ to $Y$ is proportional to that of $\Phi$ by a factor of $t\geq0$. In particular, $X_0=Y$. 
\end{definition}  
Indeed, the harmonic map dual ray has its parametrization defined only by the scaling of the dual tree to the horizontal measured foliation, say $\lambda$, of $\Phi$. Thus, we often also use the notation $\HDR_{Y,\lambda}(t)$ where $\lambda$ is the horizontal measured foliation of $\Phi$.

Equivalently, a harmonic map dual ray is a family of conformal structures defined via the minimal suspensions in a way dual to that of the harmonic map rays.  The family $\HDR_{Y, \Phi}(t)$ is the family of underlying conformal structures to the minimal graph in $Y \times (T_{\Phi}, 2td)$ parametrized by $t$ (we refer to \cite[Theorem 3.1]{Wolf1998} to the existence and uniqueness of minimal graphs).  Here the duality is expressed as follows: for the harmonic map rays $\HR_{X,\Phi}(t)$, we fix the conformal structure of the minimal surface $X$ and vary the surface $Y \in \HR_{X,\Phi}$, while for the dual ray $\HDR_{Y, \Phi}(t)$, we fix the surface $Y$ and let the conformal structure $X$ vary in $\HDR_{Y, \Phi}$. 

To see that these rays provide a ray structure (cf. Section \ref{subsec:rays})  for \tec space $\T(S)$, see \cite[Theorem 3.1]{Wolf1998} and \cite[Theorem 1.3]{Tabak1985}.

Our first focus in this paper will be the effect on the rays $\HR_{X,\Phi}(s)$ and $\HDR_{Y, \Phi}(t)$ through a point on the ray if we allow a defining surface to degenerate along the other ray. More precisely, we fix a point $Y \in \T(S)$ and consider all the harmonic map rays $\HR_{X,\Phi}(s)$ which pass through $Y$ and then study the limits of these as $X$ tends to infinity along a harmonic map dual ray.  Dually, we consider a fixed point $X$ on a family of harmonic map dual rays and study the limits of those dual rays $\HDR_{Y, \Phi}(t)$ as $Y$ diverges along a harmonic map ray.
%==============
%==============
\subsection{Analyticity of harmonic map rays}\label{sec:analyticity:HR} We recall the analyticity result of harmonic map rays in \cite{Wolf1991b}. Before that we need to introduce some notation. Let $X\in\T(S)$ and $\HR_{X,\Phi}:[0,\infty)\to\T(S)$ be the harmonic map ray defined by a holomorphic differential $\Phi$ on $X$. Consider the harmonic diffeomorphism $f_t:X\to \HR_{X,\Phi}(t)$ homotopic to the identity. In particular, the Hopf differential satisfies $\Hopf(f_t)=t\Phi$.  Define the conformal energy density and anti-conformal energy density  $$\mathcal{H}_t:=\frac{\rho(f_t)}{\sigma}|\frac{\partial f_t}{\partial z}|^2,\qquad\mathcal{L}_t:=\frac{\rho(f_t)}{\sigma}|\frac{\partial f_t}{\partial \bar{z}}|^2,$$
 and the Laplacian 
  \begin{equation*}
    \Delta_\sigma:=
    \frac{4}{\sigma}\frac{\partial^2}
    {\partial z\partial \bar{z}}.
  \end{equation*}
  Then (see \cite{Wolf1989}(p. 453); see also \cite{Sampson1978}, \cite{SchoenYau1978}) 
  \begin{eqnarray}
    e_t&=&\mathcal{H}_t+\mathcal{L}_t\label{eq:density:H:L}\\{t^2|\Phi|^2}/{\sigma^2}&=&\mh_t \ml_t \label{eq:quadratic:norm}\\
    \nu_t&=&{\partial{f_t}{\partial\bar{z}}}/{\partial{f_t}{\partial{z}}}={\overline{t\Phi}}/
    {\sigma\mh_t}\label{eq:Beltrami}\\
    |\nu_t|^2&=&{\ml_t}/{\mh_t}\label{eq:Beltrami:square}\\
    \Delta_\sigma \log \mh_t&=& 2\mh_t -2\ml_t +2K(\sigma), ~~(\text{where } \mh\neq0) \label{eq:Bochner:H}\\
    \Delta_\sigma \log \ml_t&=& 2\ml_t -2\mh_t +2K(\sigma), 
    ~~(\text{where } \ml\neq0)\label{eq:Bochner:L}\\
    t\Phi&=&\sigma \mh_t \overline{\nu_t}\label{eq:Phi}.
  \end{eqnarray}

 \begin{lemma}[\cite{Wolf1991b}, Theorem 2.2]\label{lem:analyticity2}
 With notations as above,   the conformal energy density $\mathcal{H}_t$ is real analytic in $t>0$. 
 \end{lemma}
\begin{proof}
    It is proved in \cite[Theorem 2.2]{Wolf1991b} that $\mathcal{H}_{t}$ real analytic near $t=0$. With slight modification, that proof works for all $t\in(-\infty,+\infty)$. Indeed, for any $\alpha>0$, consider the functions space $C^{2,\alpha}(X)$ and the operator:
    \begin{equation*}
        F: C^{2,\alpha}(X)\times\C \longmapsto C^{2,\alpha}(X)
    \end{equation*}
    given by 
    \begin{equation*}
        F(\mathcal{H},t):=\Delta _\sigma \log \mathcal{ H}-2\mathcal{H}+2t^2\frac{|\Phi|^2}{\sigma^2\mathcal{H}}+2.
    \end{equation*}
    For any $\tau\in(-\infty,+\infty)$, we have $F(\mathcal{H}_{\tau},\tau)=0$ and that operator $F$ is complex analytic near $(\mathcal{H}_{\tau},\tau)$. Furthermore, by elementary computation,  the linear operator $dF_{\mathcal H}$ at $(\mathcal{H}_{\tau},\tau)$ is given by:
    \begin{eqnarray*}
        dF_{\mathcal H}(\mathcal{H}_{\tau},\tau)[\psi]&=&\Delta_\sigma \frac{\psi}{\mathcal{H}_{\tau}}-2\psi -2\tau^2\frac{|\Phi|^2}{\sigma^2\mathcal{H}^2_\tau} \psi \\
        &=& \left\{ \Delta_\sigma- 2\mathcal{H}_{\tau}-2\tau^2\frac{|\Phi|^2}{\sigma^2\mathcal{H}_\tau}\right\} [\frac{\psi}{\mathcal{H}_{\tau}}] \\
        &=& \left\{ \Delta_\sigma- 2\mathcal{H}_{\tau}-2\mathcal{L}_\tau\right\} [\frac{\psi}{\mathcal{H}_{\tau}}]\qquad\quad \text{ by \eqref{eq:quadratic:norm}}\\
        &=& \left\{ \Delta_\sigma- 2e_{\tau}\right\} [\frac{\psi}{\mathcal{H}_{\tau}}] \qquad\quad \text{ by \eqref{eq:density:H:L}}.
    \end{eqnarray*}
Since $\Delta_\sigma$ is a negative operator, $\mathcal{H}_t\geq1$ (\cite[Lemma 5.1]{Wolf1989}), and hence $e_{\tau}\geq 1$, we see, by the maximum principle as well as the bound $1\leq \mathcal{H}_{\tau}\leq C<\infty$ (for some bound $C$ due to the compactness of $X$), that $dF_{\mathcal H}$ is invertible at $(\mathcal{H}_{\tau},\tau)$. It then follows from the analytic implicit function theorem \cite[Theorem 3.3.2, page 134]{Berger} that there exists a unique solution $\mathcal{H}=\mathcal{H}(t)$ of $F(\mathcal H,t)$ near $(\mathcal{H}_{\tau},\tau)$ and that $\mathcal{H}(t)$ is complex analytic near $t=\tau$.   
On the other hand, when $t$ is real, the uniqueness of the solution implies that $\mathcal{H}(t)=\mathcal H_{t}$. In particular, we see that $\mathcal H_{t}$ is real analytic in $t$ near $t=\tau$. The arbitrariness of $\tau$ then implies that $\mathcal H_{t}$ is real analytic in $t>0$. This completes the proof. 
\end{proof}

 As a direct consequence, we have:
\begin{lemma}\label{lem:analyticity3} Let $X\in\T(S)$ be a Riemann surface and $\Phi$ be a holomorphic quadratic differential on $X$. Then the harmonic map ray 
\begin{equation*}
    \HR_{X,\Phi}:[0,\infty)\to \T(S)
\end{equation*}
is real analytic in $t>0$. 
\end{lemma}
\begin{proof}
     Note that, by \eqref{eq:pullback:metric:0}, the pullback on $X$ of the hyperbolic metric of $\HR_{X,\Phi}(t)$ by $f_t$ is  
     \begin{equation}\label{eq:pullback:t}
        t\Phi(z) + e_{t}(z) \sigma (z) dz d\bar{z} +\overline{t\Phi(z)},
     \end{equation}
     where $e_{t}(z)=e_{f_t}(z)$ is the energy density of the harmonic map $f_t$ and $\sigma(z) dzd\bar{z}$ is the hyperbolic metric of $X$.
By \eqref{eq:density:H:L} and   \eqref{eq:quadratic:norm}, 
     \begin{equation*}
         e_t=\mathcal{H}_t+\frac{t^2|\Phi|^2}{\sigma^2\mathcal{H}_t}.
     \end{equation*}
     By Lemma \ref{lem:analyticity2}, we see that $\mathcal{H}_t$ is real analytic in $t>0$. 
     Therefore, the energy density $e_t$ is also real analytic in $t>0$, so is the pullback metric on $X$ of the hyperbolic metric on $\HR_{X,\Phi}(t)$. Accordingly, the harmonic map ray $\HR_{X,\Phi}(t)$ is real analytic in $t>0$.     
\end{proof}
 
%==================
%==================
 \subsection{Minsky's estimates} 
 The function $\mathscr{G}$  appearing in (\ref{eq:pullback:metric:2}) is nearly zero at points which are far away from the zeros of $\Phi$. More precisely, 
 \begin{lemma}[\cite{Minsky1992},Lemma 3.2 and Lemma 3.3]\label{lem:Minsky:decay}
 Let $X,\Phi,\mathscr{G}$ be as defined in Section \ref{subsec:Hopf}.  Let $p\in X$ be at a $|\Phi|$-distance at least $d$ from any zero of $\Phi$. Then 
     $$ \mathscr{G}(p)\leq 
   \frac{\sinh^{-1}(|\chi(M)|/d^2)}{\exp(d)}.$$
 \end{lemma}

 We also need the following estimate about the gradient of $\mathscr G$, with respect to the singular flat metric induced from $|\Phi|$: 
 \begin{lemma}
    \label{lem:grad:decay} 
  Let $X,\Phi,\mathscr{G}$ be as defined in Section \ref{subsec:Hopf}. Let $p\in X$ be at a $|\Phi|$-distance at least $d$ from any zero of $\Phi$. Then there exists a constant $c$ depending only on $d$ and the topology of $X$ such that 
     $$ |\nabla \mathscr{G}(p)|\leq c \cdot \mathscr{G}(p).$$ 
 \end{lemma}
 \begin{proof}
     Recall that 
     \begin{equation}\label{eq:G:laplacian}
         \Delta \mathscr{ G}=4\sinh \mathscr{G}
     \end{equation} (see for instance \cite[Equation (3.2)]{Minsky1992}, where we choose $\sigma=|\Phi|$ and $K(\rho)=-1$).
     Combined with Lemma \ref{lem:Minsky:decay}, this implies that 
     \begin{equation*}
         0\leq \Delta \mathscr{G}\leq c_1\cdot  \mathscr{G}
     \end{equation*}
     for some constant depending only on $d$ and the topology of $X$. By \cite[Theorem 1]{Heinz1956}, we see that there exists a constant $c_2$ depending only on $d$ and $c_1$, hence on $d$ and the topology of $X$, such that
     \begin{equation}\label{eq:G:Harnack}
         \mathscr{G}(q)\leq c_2\cdot \mathscr G(p)
     \end{equation}
     for any $q$ in the $d/2$ ball $B_{d/2}(p)$ centered at $p$ under the $|\Phi|$-metric. 

     Next, we claim that there exists a constant $c_3$ depending only $d$ such that
     \begin{equation}\label{eq:G:grad}
		\sup_{q\in B_{d/4}(p)}|\nabla \mathscr{G}(q)|\leq c_3 \sup_{q\in B_{d/2}(p)}\mathscr{G}(q).
	\end{equation}
    The proof of the claim follows almost line-by-line from the proof of \cite[Theorem 2.1]{Colding2005} for harmonic functions via the Bochner formula and testing functions. Equation \eqref{eq:G:laplacian} together with the fact that $\mathscr{G}>0$ 
		ensure that
		\begin{equation*}
			<\nabla\Delta \mathscr{G}, \nabla \mathscr{G}>~\geq 0 \text{ and } \mathscr{G}\Delta \mathscr{G}\geq 0.
		\end{equation*}
		This is where the assumption of harmonicity is used in the proof of \cite[Theorem 2.1]{Colding2005}. 
The lemma then follows from \eqref{eq:G:Harnack} and \eqref{eq:G:grad}.
 \end{proof}

 Minsky defined a family $\mathscr{P}_{R}$ of regions whose geometry is controlled and in whose complement the harmonic map is nearly a projection.  We now summarize what we need of this work. Before stating the result, we need the notion of \emph{boundary-convex}. Given a closed Riemann surface $X$ and a holomorphic differential $\Phi$ on $X$, a subset $C$ of $X$ is \emph{boundary-convex} if any geodesic arc, with respect to the singular flat metric induced from $\Phi$, deformable rel endpoints into $C$, in fact lies in $C$. Equivalently, $C$ is boundary-convex if the lift of each component of $C$ to the universal cover $\widetilde X$  of $X$ is convex under the singular flat metric induced from $\Phi$. (We refer to \cite[page 172]{Minsky1992} for further explanation about boundary-convex.)

 \begin{theorem}[\cite{Minsky1992}, Theorem 5.1]
 \label{thm:Minsky:polygon}
Let $X$ be a closed surface of genus $g$ carrying a flat metric induced by a holomorphic quadratic differential $\Phi$, and let $s>0$ and $c_1,\cdots,c_{3g-3} > 0$ be chosen constants. Then there exist positive constants $A_1, K_1, K_2$,  and $K_3$ depending only on $s$ and $\chi(X)$, and  a nested family $\{\mathscr P_R\}_{R>0}$ of boundary-convex set $\mathscr P_R$  with the following properties: 
 \begin{enumerate}[(i)]
     \item  $ \mathscr{P}_{R}$ contains the $R$ neighbourhood of zeros of $\Phi$. 
     \item Every component of $\partial \mathscr{P}_{R}$ is either polygonal comprising alternatively  horizontal geodesic segments and vertical geodesic segments, or regular horizontal geodesics. The regular geodesic components occur in pairs bounding homotopically distinct flat cylinders.
     \item If $\mathscr{F}_{k}$ is the $k$-th maximal horizontal cylinder whose subcylinders occur as components of $X\backslash\mathscr{P}_{r}$ for some $r \leq R$, and $\mathscr{F}_{k}$ is partially contained in  $\mathscr{P}_{R}$ (that is, $\mathscr{F}_{k}$ itself is not contained in $\mathscr P_R$ but a subcylinder of $\mathscr{F}_{k}$ is contained in $\mathscr P_R$), then $\mathscr{F}_{k} \cap \mathscr{P}_{R}$ is a pair of flat cylinders with height at least $r_{W}+R+c_{k} R^{2} / W$, where $W=W\left(\mathscr{F}_{k}\right)$ is the circumference of $\mathscr{F}_{k}$  and $r_W=W+\log\frac{1}{2}\sinh ^{-1}(2|\chi(X)|/W^2)$.
     \item  $\ell\left(\partial \mathscr{P}_{R}\right) \leq K_{1} R$.
     \item $\operatorname{Area}\left(\mathscr{P}_{R}\right) \leq A_{1}+\left(K_{2}+2 \sum_{i=1}^{k} c_{i}\right) R^{2}$, where $k$ is the number of flat cylinder components of $M\backslash\mathscr{P}_{r}$ that have occurred for $r \leq R$. 
     \item  Each edge of a polygonal boundary component has length at least $K_{3} R$.
     \item The polygonal components of $\partial \mathscr{P}_{R}$ are $s$-separated. Namely, for any two components $\gamma_1,\gamma_2$ of $\partial \mathscr{P}_{R}$ (where possibly $\gamma_1=\gamma_2$), any arc in $X$ with endpoints in $\gamma_1 $ and $\gamma_2$ that can not be deformed (rel endpoints) into $ \mathscr{P}_{R}$ has length greater than $s\max\{\ell(\gamma_1),\ell(\gamma_2)\}$.
   \end{enumerate}
 \end{theorem}

Before going further, let us make some comments about  $\mathscr P_R$. 
\begin{itemize}
	\item  Very roughly speaking, the basic strategy of the construction/proof of Theorem \ref{thm:Minsky:polygon} is to start with the $R$-neighborhood of the zeros of $\Phi$, boundary-convexify it, \enquote{square out the corners}, and control the size of the resulting set. (Flat cylinders need separate treatment to satisfy property (iii).)  In particular, we may assume that each component of  $\mathscr P_R$ contains at least one zero of $\Phi$. 
	\item When $R$ is small, say much less than half of the shortest $|\Phi|$-distance between zeros of  $\Phi$, the region $\mathscr P_R$ is a union of polygons, each surrounds a zero of $\Phi$. As $R$ increase, those polygons may merge, possibly resulting in new regular boundary geodesics. Eventually, when $R$ is bigger than  the $|\Phi|$-diameter of $X$, the region $\mathscr P_R$ is the whole surface $X$. 
	\item From item (vii) of Theorem \ref{thm:Minsky:polygon}, for each polygonal boundary component $\partial_i \mathscr P_R$, one can append an embedded annulus of radius $\frac{s}{2}\cdot \ell(\partial_i \mathscr P_R)$ outside $\mathscr P_R$.  Furthermore, those appended annuli are pairwise disjoint.  This is not necessary true for regular geodesic boundary components. On the other hand, both polygonal boundary components and regular geodesic boundary components admit annuli inside $\mathscr P_R$.  Indeed, for regular geodesic boundary components, this is the content of item (iii) in Theorem \ref{thm:Minsky:polygon}; for  polygonal boundary components, this is the content of Lemma \ref{lem:in:annulus} below. 
\end{itemize}   

\begin{lemma}\label{lem:in:annulus} 	Let $X,\Phi$, $s$, $c_1,\cdots, c_{3g-3}$ $K_3$, and $\mathscr P_R$  be as in Theorem \ref{thm:Minsky:polygon}. Then each polygonal boundary component $\partial_i\mathscr P_R$  of $\mathscr P_R$ admits an embedded annulus $\mathscr A_i$ in $\mathscr P_R$ such that 
	\begin{enumerate}[(i)]
		\item $\partial_i\mathscr P_R$ is one of the two boundary components of $\mathscr A_i$, 
		\item the other component $\alpha'$ of $\partial\mathscr A_i$ is polygonal comprising alternatively horizontal geodesic segments and vertical geodesic segments parallel to the corresponding geodesic segments of $\partial_i\mathscr P_R$ 
		\item each point in $\alpha'$ is (exactly) of distance $K_3'R/4$ from $\partial _i\mathscr P_R$, where $K_3'=\min\{K_3,1\}$,
		\item the angle at each corner of $\alpha'$ measured within $\mathscr A_i$ is $3\pi/2$.
	\end{enumerate}
\end{lemma}
\begin{proof}

Let $D$ be the infimum of the lengths of arcs in $\mathscr P_R$ that connect one point in $ \partial_i\mathscr P_R$ to $\partial  \mathscr P_R$ and are not homotopic rel endpoints to an arc in $\partial_i \mathscr P_R$. Since $\mathscr P_R$ is boundary convex, it follows that the minimum $D$ is realized by (at least) one geodesic arc, say $\gamma$. Note that $\partial \mathscr P_R$ contains no zeros of $\Phi$ and comprises  horizontal geodesic segments and vertical geodesic segments. Hence the angle at each corner of $\partial\mathscr P_R$ is either $\pi/2$ or $3\pi/2$.  Being boundary convex implies that the angle at each corner of $\partial \mathscr P_R$, from the pointview of the interior of $\mathscr{P}_R  $, is exactly $\pi/2$. Hence, as a minimizing path in $\mathscr P_R$, the arc $\gamma$ has to avoid the corners of $\partial \mathscr P_R$ and meet $\partial \mathscr P_R$ perpendicularly.

Next, we claim that
\begin{quote}
	Claim 1: there is a minimizing arc $\gamma'$ realizing $D$  that contains at least one zero of $\Phi$.
\end{quote}
   To see this, we start with the minimizing arc $\gamma$ mentioned above. If $\gamma$ itself contains a zero of $\Phi$, then we are done. Suppose that $\gamma$ does not contain any zero of $\Phi$. Recall that $\gamma$ meets $\partial \mathscr P_R$ perpendicularly. Therefore, the arc $\gamma$ has to connect a horizontal segment of $\partial _i \mathscr P_R$ to a horizontal segment or connect a vertical segment to a vertical segment, as $\gamma$ is a flat non-singular geodesic in a flat metric whose initial tangent is vertical, resp. horizontal, and thus is always vertical, resp. horizontal. From this we see that there exists an embedded rectangle near $\gamma$ which contains $\gamma$ as one of its boundary sides, with $\gamma$ connecting two other parallel boundary sides of $\partial\mathscr P_R$. We enlarge this rectangle by pushing outward the other boundary side $\hat\gamma$ (that is parallel to $\gamma$) until it hits a corner of $\partial \mathscr P_R$ or a zero of $\Phi$. If $\hat\gamma$ hits a corner of $\partial \mathscr P_R$, then $\hat\gamma$ contains a segment of $\partial\mathscr P_R$, then removing this subsegment yields a strictly shorter subarc that is homotopic to $\hat\gamma$ and hence also homotopic to $\gamma$ rel   $\partial \mathscr P_R$. Since $\hat\gamma$ and $\gamma$ minimize distance, this contradicts the assumption that $\gamma$ realizes the minimum $D$.  Therefore, the arc $\hat\gamma$ must contain a zero of $\Phi$. This establishes Claim 1 by taking $\hat\gamma$ as the required $\gamma'$.

From Claim 1 and item (i) of Theorem \ref{thm:Minsky:polygon}, we see that the length of $\gamma'$ (or equivalently $\hat \gamma$) is at least $2R$. Hence, the minimizing distance satisfies 
\begin{equation}\label{eq:D:min}
	D\geq 2R.
\end{equation}
  Let $I_1, I_2, \cdots, I_{2m}$ be the edges of the polygonal boundary component $\partial_i\mathscr P_R$, denoted in such a way that they appear consecutively.  For each edge $I_j$ of $\partial _i\mathscr P_R$, let $\Box_j$ be the rectangle in $\mathscr P_R$ that contains $I_j$ as one of its boundary sides and whose second pair of parallel sides (those on $I_{j-1}$ and $I_{j+1}$ in the usual notation) is also contained in $\partial_i\mathscr P_R$ and is of length $K_3'R/2$ in total (i.e. of length $K_3'R/4$ for each of the second pair of sides on $I_{j-1}$ and $I_{j+1}$), where $K_3'=\min\{K_3,1\}$ and $K_3$ is the constant from Theorem \ref{thm:Minsky:polygon}.    We claim that
\begin{quote}
	Claim 2: $\Box_j\cap\Box_k\neq \emptyset$ if and only if  $I_j$ and $I_k$ are the two sides of a corner of $\partial_i\mathscr{P}_R$.
	\end{quote}
	To see this, for one direction, suppose that $I_j$ and $I_k$ are the two sides of a corner, then $\Box_j\cap\Box_k$ contains that corner, and hence is not an emptyset. For the converse direction, suppose that $\Box_j\cap\Box_k\neq \emptyset$, let $p$ be a point in $\Box_j\cap\Box_k\neq \emptyset$. Since three of the four sides of $\Box_j$ are contained in $\partial_i\mathscr P_R$, it follows that $\Box_j$ and $\Box_k$ either coincide or share a common corner. The first possibility happens only if the component of $\mathscr P_R$ containing $\partial_i\mathscr P_R$ is itself a rectangle without any zero of $\Phi$, which is excluded from our consideration (see the first comment following Theorem \ref{thm:Minsky:polygon}). Therefore, the two rectangles $\Box_j$ and $\Box_k$ must share a corner. Accordingly, the two edges $I_j$ and $I_k$ are the two sides of that corner in $\partial _i\mathscr P_R$. This proves Claim 2.
	
	From Claim 2, we see that the union $\cup_j \Box_j$ is a topological annulus with $\partial _i\mathscr P_R$ as one of its boundary components and that the angle at each corner of the boundary component other than $\partial _i\mathscr P_R$ is exactly $3\pi/2$ measured within the annulus $\cup_j \Box_j$. The lemma now follows by taking the union  $\cup_j \Box_j$ as the required annulus $\mathscr A_i$.
\end{proof}

 \begin{theorem}[\cite{Minsky1992}, Theorem 7.1]
 \label{thm:Minsky:traintrack} 
 Let $f:X\to Y$ be a harmonic diffeomorphism between closed hyperbolic surfaces of genus $g$ with Hopf differential $\Phi$. There are choices of constants $s>0$ and $c_{1}, \cdots$, $c_{3g-3}>0$  for the construction of the polygonal region $\mathscr{P}_R$ and an $R_0>0$, such that in the complement of $\mathscr{P}_{R_0}$ there is  a map $\pi$ from the leaves of $\mathrm{Hor}(\Phi)$ to the lamination  corresponding to $f_*(\mathrm{Hor}(\Phi))$ that factors through $f$, and is a local diffeomorphism on each leaf of $\mathrm{Hor}(\Phi)$, mapping it to the corresponding geodesic representative of its image. For any point $p$ on a leaf in $X-\mathscr{P}_{R_0}$,
 $$ d_Y(f(p),\pi(p))<a\exp \left(-b ~d_{|\Phi|}(p,\mathscr{P}_{R_0})\right), $$
where $d_Y$  is the hyperbolic distance on $Y$, and the derivative of $\pi$ along leaves with respect to the $|\Phi|$-metric satisfies
 $$ | |d\pi|-2|\leq a\exp \left(-b ~d_{|\Phi|}(p,\mathscr{P}_{R_0})\right), $$
 where $a$ and $b$ are positive constants depending only on $\chi(X)$.
 \end{theorem}

 \begin{theorem}[\cite{Minsky1992}, Theorem 7.2]\label{thm:Minsky:energy}
   There exists a constant $C_0$ depending on the topology of $S$, such that for every $X, Y \in \T (S)$
 \begin{equation}\label{eq:energy:length:ratio}
   \sup_{\mu\in\ML(S)} \frac{1}{2}\frac{\ell^2_Y(\mu)}{\Ext_{X}(\mu)}\leq E(X,Y)\leq  \frac{1}{2}\frac{\ell^2_Y(\gamma)}{\Ext_{X}(\gamma)}+C_0
 \end{equation}
 where $\gamma$ is the horizontal measured foliation of the Hopf differential of the harmonic diffeomorphism $f:X\to Y$ homotopic to the change of marking. 
 \end{theorem}

 \remark The method of Minsky actually proves a slightly stronger estimate in the following sense. Let $\Phi$ be the Hopf differential of the harmonic diffeomorphism $f:X\to Y$ with $\gamma$ being the horizontal measured foliation/lamination.    Let $\{X_i\}$ be the components complementary to the critical graph  (see Definition \ref{def:critical:graph}) of (the horizontal foliation $\gamma$ of) $\Phi$ and let $\gamma_i$ be the restriction of $\gamma$ to $X_i$. 
 Since each $X_i$ is locally convex and contains every leaf of $\gamma_i$, the construction of Minsky's polygonal region in Theorem \ref{thm:Minsky:polygon} and train-track approximate in the proof of Theorem \ref{thm:Minsky:traintrack}    occur within each $X_i$ (cf. especially sections 5 and 6 of \cite{Minsky1992}). Then Theorem \ref{thm:Minsky:energy} holds on $X_i$. Namely,
 \begin{equation}\label{eq:length:ratio:component}
   \sup_{\mu\in\ML(X_i)} \frac{1}{2}\frac{\ell^2_Y(\mu)}{\Ext_{X_i}(\mu)}\leq E(f|_{X_i})\leq  \frac{1}{2}\frac{\ell^2_Y(\gamma_i)}
   {\Ext_{X_i}(\gamma_i)}+C_0,
 \end{equation}
 where $\ML(X_i)$ represents the measured foliations/laminations on $X$ that can be homotoped into $X_i$ 
 
Recall \cite[Page 165, Equation (3.4)]{Minsky1992} (see also \cite[Lemma 3.2]{Wolf1989}) that for any subsurface $U\subset X$,  the energy and norm of Hopf differential $\Phi$ satisifes
\begin{equation}\label{eq:energy:Hopfdifferential}
   2\|\Phi|_U\|\leq E(f|_U)\leq 2\|\Phi|_U\|+2\pi |\chi(S)|. 
 \end{equation}

  Notice that for the horizontal measured foliation/lamination $\gamma$ of $\Phi$,
 \begin{equation}\label{eq:ext:quadraticnorm}
   \Ext_X(\gamma)=\|\Phi\|.
 \end{equation}
 Combining (\ref{eq:energy:length:ratio}), (\ref{eq:energy:Hopfdifferential}), and (\ref{eq:ext:quadraticnorm}), we obtain some basic first estimates for this paper.

 \begin{lemma}\label{lem:hyplength:quadraticnorm} 
 For any closed orientable surface $S$ of genus at least two there exists a constant $C$ depending only on $S$ such that the following holds.
  Let $X,Y\in\T(S)$ be two hyperbolic surfaces  and $f:X\to Y$ the unique harmonic map homotopic to the change of marking. Let $\Phi$ be the Hopf differential of $f$. Let $\gamma$  be the measured lamination corresponding to the horizontal measured foliation of $\Phi$. Then 
\begin{equation}\label{eq:phi:length}
  \ell_Y(\gamma)-C\leq 2\|\Phi\|\leq \ell_Y(\gamma)+C
 \end{equation}
 and
 \begin{equation}\label{eq:energy:length}
    \ell_Y(\gamma)-C\leq E(f)\leq \ell_Y(\gamma)+C.
 \end{equation}

 Moreover, let $\gamma=\sum_{i=1}^{k}\gamma_i$  be the component decomposition of $\gamma$.  Let $X_i$ be the component supporting $\gamma_i$ in the complementary region of the critical graph (see Definition \ref{def:critical:graph}) of  $\Phi$.  Then
    \begin{equation}\label{eq:phi:length:component}
  \ell_Y(\gamma_i)-C\leq 2\|\Phi|_{X_i}\|\leq \ell_Y(\gamma_i)+C
 \end{equation}
 and
 \begin{equation}\label{eq:energy:length:component}
    \ell_Y(\gamma_i)-C\leq E(f|_{X_i})\leq \ell_Y(\gamma_i)+C.
 \end{equation}
 \end{lemma}
\begin{proof}
 It follows from (\ref{eq:energy:length:ratio}), (\ref{eq:energy:Hopfdifferential}), and (\ref{eq:ext:quadraticnorm}) that
 \begin{equation*}
  4\|\Phi\|(\|\Phi\|-C_0/2)\leq \ell^2_Y (\gamma)\leq 4\|\Phi\|(\|\Phi\|+\pi|\chi(S)|)
 \end{equation*}
 which implies (\ref{eq:phi:length}). Inequality (\ref{eq:energy:length}) then follows from (\ref{eq:energy:Hopfdifferential}) and (\ref{eq:phi:length}).
 Upon  replacing (\ref{eq:energy:length:ratio}) by (\ref{eq:length:ratio:component}) in the above reasoning, we obtain (\ref{eq:phi:length:component}) and (\ref{eq:energy:length:component}).  
 \end{proof}

%=================
%=================
 
 \section{Compactness of harmonic maps to a fixed target} \label{sec:compactness}
 
 The goal of this section is to prove a compactness result for  harmonic maps to a fixed target (Lemma \ref{lem:precompactness}).  We continue with the notation introduced in the previous section.
 
 \subsection{Extension of estimates on Minsky's polygonal regions}
 We begin with an extension of Minsky's analysis of the polygonal region that we will need in our discussion of the compactness of a family of Riemann surface domains for the harmonic maps.

 \begin{lemma}\label{lem:inj:lowerbound}
   Let $Y\in\T(S)$ be a hyperbolic surface. For any $R$ sufficiently large, there exist two positive constants $\delta$ and $\rho$ depending on $R$ and $Y$ such that for any harmonic map $f:X\to Y$ with $X\in\T(S)$ homotopic to the change of marking, each component of the polygonal region $\mathscr{P}_R$ has injectivity radius at least $\delta$ and diameter at most $\rho$ with respect to the Hopf differential metric.
 \end{lemma}

\begin{remark}
    This lemma will be important in controlling the shape of limits of harmonic maps.  As one simple example, note that this lemma means that one cannot pinch the domain along short curves containing a zero of the Hopf differential:  the zeroes are always in the interior of the Minsky regions, and the Lemma states that these regions have injectivity radii bounded away from zero.
\end{remark}
 
 \begin{proof}[Proof of Lemma \ref{lem:inj:lowerbound}]
 
   If $\mathscr{P}_R$ is a topological disk, then the conclusion follows directly from Minsky's construction. Now let us assume that $\mathscr{P}_R$ is not a topological disk.
 
   First, we claim that 
   \begin{quote}
   \textbf{Claim 1:} 
   \textit{there exists a positive constant $\delta_1$ depending on $R$ and $Y$ such that the extremal length of any  essential simple closed curve $\alpha$ of $\mathscr{P}_R$ is at least $\delta_1$:
   $ \Ext_{\mathscr{P}_R}(\alpha)\geq \delta_1. $}
   \end{quote}
    To see this first claim, note that by (\ref{eq:energy:Hopfdifferential}), we have $E(f|_{\mathscr{P}_R})\leq 2\|\Hopf(f)\|_{\mathscr{P}_R}+\mathrm{Area}(Y)\leq c R^2+\mathrm{Area}(Y)$. It then follows from (\ref{eq:energy:length:ratio})  that
   $$ \Ext_{\mathscr{P}_R}(\alpha) \geq \frac{\ell_{Y}^2(\alpha)}{2 E(f|_{\mathscr{P}_R})}\geq \frac{\mathrm{Syst}(Y)^2}{2c R^2+2\mathrm{Area}(Y)}.$$
   This proves Claim 1.

 \vskip5pt
   Next, we claim that
   \begin{quote} \textbf{Claim 2:}
    \textit{
     there exists a positive constant $\delta_2$ depending on the topology of $S$ such that,  for each  boundary component $\alpha$ of $\mathscr{P}_R$, we have $\Ext_{\mathscr{P}_R}(\alpha)\leq \delta_2.$}
   \end{quote}
  To see this second claim, note that there are two cases here, depending on whether the given boundary component $\alpha$ is polygonal or regular.  We start with the case where $\alpha$ is a regular geodesic. By item (iii) and item (iv) of Theorem \ref{thm:Minsky:polygon},  $\alpha$ admits a horizontal cylinder  $\mathscr C$ inside $\mathscr{P}_R$ of height at least $R$ and of circumference $\ell(\alpha)$ at most $K_1R$ for some positive constant $K_1$ depending the topology of $S$. Hence,  the modulus $\mod(\mathscr C)$ satisfies 
   $$\mod(\mathscr C)\geq  \frac{R}{\ell(\alpha)} \geq \frac{1}{K_1}. $$
   According to \eqref{eq:ext:mod}, we see that 
\begin{equation}\label{eq:clm:ext:bnd}
	\Ext_{\mathscr{P}_R}(\alpha)\leq 
   \frac{1}{\mod(\mathscr C)}\leq K_1.
\end{equation}
  Now we turn to the second case
   where  $\alpha$ is a polygonal boundary component. Let $\mathscr A\subset\mathscr P_R$ be the annulus attached to $\alpha$ obtained from Lemma 
  \ref{lem:in:annulus}. In particular, the annulus $\mathscr A$ contains $\alpha$ as one of its boundary components. Suppose  that $\alpha$ has $m$ corners. Then it has $m$ horizontal segments $\{h_1,\cdots,h_m\}$ and $m$ vertical segments $\{v_1,\cdots,v_m\}$. By item (ii) and (iii) of Lemma \ref{lem:in:annulus}, we infer that the other boundary component $\alpha'$ of the annulus $\mathscr A$ consists of  $m$ horizontal segments  $\{h'_1,\cdots,h'_m\}$  with lengths  $\ell(h_i')=\ell(h_i)-K_3'R/2$, $m$ vertical segments $\{v'_1,\cdots,v'_m\}$ with lengths $\ell(v_i')=\ell(v_i)-K_3'R/2$. Therefore, the lengths of $\alpha'$ and $\alpha$ satisfy:
  \begin{equation*}
  	\ell(\alpha')=\sum_i (\ell(h_i')+\ell(v_i'))=\sum_i (\ell(h_i)-K_3'R/2+\ell(v_i)-K_3'R/2)=\ell(\alpha)-mK_3'R.
  \end{equation*}
  From item (iii) of Lemma \ref{lem:in:annulus}, it follows that the annulus $\mathscr A$ contains the regular annulus $\mathscr A'$ of $\alpha'$ of radius $K_3'R/4$:
  \begin{equation*}
  	\mathscr A':=\{p\in \mathscr A: d(p,\alpha')=K_3'R/4\}
  \end{equation*}
  where $d$ refers to the singular $|\Phi|$-metric.  By item (iv) of Lemma \ref{lem:in:annulus}, the angle at each corner of $\alpha'$ with respect to $\mathscr A$ is $3\pi/2$. So the total curvature $\mathbb K_{\alpha'}$ of $\alpha'$ with respect to both $\mathscr A$ and $\mathscr A'$ is $\mathbb K_{\alpha'}=-m\pi /2$.  The  boundary component $\alpha''$ of  $\mathscr A'$ other than $\alpha'$ consists of  $m$ horizontal segments  $\{h''_1,\cdots,h''_m\}$  with lengths  $\ell(h_i'')=\ell(h_i')$, and  $m$ vertical segments $\{v''_1,\cdots,v''_m\}$ with lengths $\ell(v_i'')=\ell(v_i')$, and $m$ circular arcs centered at each corner of $\alpha'$ of radius $K_3'R/4$ and of angle $\pi/2$. Hence the length $\ell(\alpha'')$ satisfies:
  \begin{equation} \label{eq:length:alpha:dprime}	\ell(\alpha'')=\sum_i(\ell(h_i'')+\ell(v_i''))+m\cdot \frac{\pi}{2}\cdot \frac{K_3'R}{4}=\ell(\alpha')+ \frac{mK_3'R\pi}{8}.
  \end{equation}
   Combined with \cite[Theorem 4.5]{Minsky1992}, this implies  that the modulus of $\mathscr A'$ satisfies:
  \begin{eqnarray*}
  \mod (\mathscr A') &\geq& \frac{1}{|\mathbb K_{\alpha'}|}\log \frac{\ell(\alpha'')}  {\ell(\alpha')} \qquad(\text{\cite[Theorem 4.5]{Minsky1992} }) \\&=& \frac{1}{m\pi/2}\log \frac{\ell(\alpha')+mK_3'R\pi/8}{\ell(\alpha')}\\&\geq & 
  	 \frac{1}{m\pi/2}\log \frac{mK_3'R\pi}{8\ell(\alpha')}\\
  	 &= & \frac{1}{m\pi/2}\log \frac{mK_3'R\pi}{8(\ell(\alpha)-mK_3'R)} \text{ by \eqref{eq:length:alpha:dprime} }\\&\geq&  \frac{1}{m\pi/2}\log \frac{mK_3'R\pi}{8\ell(\alpha)}\\
  	 &\geq & \frac{1}{m\pi/2}\log \frac{mK_3'\pi}{8K_1}\qquad (\text{ item (iv) of Theorem \ref{thm:Minsky:polygon}}). 
  \end{eqnarray*}
  Recall that the angle at each corner of each polygonal boundary component of $\mathscr P_R$ with respect to the interior of $\mathscr P_R$ is $\pi/2$. Applying the Gauss-Bonnet formula to $\mathscr P_R$ with the singular $|\Phi|$-metric, we see that total number of corners of $\partial \mathscr P_R$, hence also the number $m$ of corners of $\alpha$, is bounded from above by some constant $\mathbf{k}$ depending only the topology of the underlying surface $X$, that is,
  \begin{equation}\label{eq:corner}
     m=\#\{\text{corners of }\partial \mathscr P_R \}\leq \mathbf{k}.
  \end{equation}
  Inserting this fact into the above-displayed equation about moduli, we see that $\mod(\mathscr A')\geq C$ for some positive constant $C$ depending on the topology of $X$.  Therefore, by \eqref{eq:ext:mod}, the extremal length of $\alpha$ on $\mathscr P_R$ satisfies:
  \begin{equation*}
  	\Ext_{\mathscr P_R}(\alpha)\leq \frac{1}{\mod(\mathscr A')} \leq \frac{1}{C}
  \end{equation*}
  for some constant $C$ depending only on the topology of $X$. Combined with \eqref{eq:clm:ext:bnd}, this establishes Claim 2.

   \vskip5pt
   Thirdly, let us consider the double of $\mathscr{P}_R$, denoted by $\mathscr{P}_R^d$ obtained by gluing a copy of $\mathscr{P}_R$ to (the oppositely oriented) $\mathscr{P}_R$ along the corresponding boundary edges.    We claim that
    \begin{quote}
    \textbf{Claim 3:}
     \textit{there exists a positive constant $\delta_3$ such that, for any essential simple closed curve $\alpha$ on $\mathscr{P}_R^d$, we have $\Ext_{\mathscr{P}_R^d}(\alpha)\geq \delta_3$.} 
    \end{quote}
   If $\alpha$ is a boundary component of $\mathscr{P}_R$ or an interior curve of $\mathscr{P}_R$, then
   $$ \Ext_{\mathscr{P}_R^d}(\alpha)\geq \frac{1}{2}\Ext_{\mathscr{P}_R}(\alpha) \geq \frac{\delta_1}{2}. $$
   If $\alpha$ intersects one of the boundary components of $\mathscr{P}_R$ say $\beta$, then
   $$ \Ext_{\mathscr{P}_R^d}(\alpha)\geq
   \frac{i(\alpha,\beta)^2}{\Ext_{\mathscr{P}_R^d}(\beta)}
   \geq \frac{1}{\Ext_{\mathscr{P}_R^d}(\beta)}
   \geq \frac{1}{\Ext_{\mathscr{P}_R}(\beta)}\geq\frac{1}{\delta_2},  $$
   where the last inequality follows from the second claim above. This proves Claim 3.
  
  \vskip5pt
   Finally, we claim that there exist $\delta_4>0$ and $\rho>0$ depending only on $R$ and the topology of $S$ such that
   \begin{quote}
   \textbf{ Claim 4:} \textit{
   there exist constants $\delta_4$ and $\rho$ so that the polygonal region $\mathscr{P}_R$ has injectivity radius at least $\delta_4$ and diameter at most $\rho$ with respect to the Hopf differential metric.}
   \end{quote}
   Notice that the restriction of the Hopf differential to $ \mathscr{P}_R$ and its copy gives a meromorphic quadratic differential $q^d$ on $\mathscr{P}_R^d$. This differential has simple poles at the corners of the boundary components of $\mathscr{P}_R$; the number of these points is bounded from above by some constant $\mathbf{k}$ depending only on the topology of $S$ (see \eqref{eq:corner}).   Notice that the genus of $\mathscr{P}_R^d$ is at most $\mathrm{genus}(X)+ (3\mathrm{genus}(X)-2)=4\mathrm{genus}(X)-2$.
   Set  $\mathbf{g}=4\mathrm{genus}(X)-2$.
    Then $$(\mathscr{P}_R^d,\frac{q^d}{2\|\mathscr{P}_R\|})\in \cup _{g\leq \mathbf{g}, \kappa\leq \mathbf{k} } \mathcal{Q}^\kappa{\mathcal M}_g$$ 
    where $\mathcal{Q}^\kappa{\mathcal M}_g^\kappa$ represents the bundle over the moduli space $\mathcal{M}_g$ whose fiber, over a Riemann surface  $M\in \mathcal{M}_g$, is the space $Q^\kappa(M)$ of area one meromorphic quadratic differentials having $\kappa$ simple poles.  Moreover, combined with Mumford's compactness theorem \cite[Corollary 3]{Mumford1971} and an inequality between hyperbolic length and extremal length \cite[Corollary 3]{Maskit1985}, the third claim above implies that for each $g$ there exists a compact subset $K_{g}\subset \mathcal{M}_g$ such that $\mathscr{P}_R^d\subset K_{g}$,   i.e. $(\mathscr{P}_R^d,\frac{q^d}{2\|\mathscr{P}_R\|})\in \cup _{g\leq \mathbf{g}, \kappa\leq \mathbf{k} } \mathcal{Q}^\kappa K_g$. Note that, in general, simple poles may collide to form less singular points, so for any $\kappa$, neither $\mathcal{Q}^\kappa(M)$ nor $\mathcal{Q}^\kappa K_g$ is compact. On the other hand,  simple poles cannot collide to form poles of higher order because of the assumption of area one, so they could only collide to form regular points or zeros. Therefore,  both $\cup_{0\leq \kappa \leq \mathbf{k}}\mathcal{Q}^\kappa (M)$ and $\cup_{0\leq \kappa \leq \mathbf{k}}\mathcal Q^\kappa K_g$ are compact
    (as a subcomplex of the stratified space of quadratic differentials of unit area and with at most simple poles).  Then for each  $g\leq\mathbf{g}$, there exist two positive constants $\delta_g$ and $\rho_g$ such that 
   $$ \inf_{\substack{
        0\leq \kappa \leq \mathbf{k}\\(M,q')\in \mathcal Q^\kappa K_g(M)}}\frac{\mathrm{inj}(q') }{\mathrm{inj}(M)}\geq \delta_{g} \quad\text{and}\quad \ \sup_{\substack{
      {  0\leq \kappa \leq \mathbf{k}}\\(M,q')\in \mathcal Q^\kappa K_g(M)}}\frac{\mathrm{diam}(q') }{\mathrm{diam}(M)}\leq  \rho_{g},$$
    where $\mathrm{inj}(M)$ and $\mathrm{diam}(M)$ represent respectively the injectivity radius and diameter of the hyperbolic metric of $M$, and where  $\mathrm{inj}(q')$ and $\mathrm{diam}(q')$ represent respectively the injectivity radius and diameter of the singular flat metric induced by $q'$.
    Hence for any $M\in K_g$ and any $q'\in \cup_{0\leq\kappa\leq\mathbf{k}}Q^\kappa(M)$, we have
    $$\mathrm{inj}(q')\geq \delta_{g} \min_{M\in K_g}\mathrm{inj}(M)\quad {\text{ and } \quad \mathrm{diam}(q')\leq \rho_{g} \max_{M\in K_g}\mathrm{inj}(M) } $$
   That $0\leq g\leq \mathbf{g}<\infty$ then implies that
    $$ \delta_4':=\min_{ 0<g\leq \mathbf{g}}  \delta_{g} \min_{M\in K_g}\mathrm{inj}(M)>0, $$
    and 
    $$ 
    \rho_4':=\max_{0<g\leq \mathbf{g}}  \rho_{g} \max_{M\in K_g}\mathrm{inj}(M)<\infty.
    $$
   In particular, the injectivity radius of $\frac{q^d}{2\|\mathscr{P}_R\|}$ with respect to the singular flat metric $\frac{q^d}{2\|\mathscr{P}_R\|}$ is at least $\delta_4'$.   Hence the injectivity radius of $q^d$ is at least ${2\|\mathscr{P}_R\|\delta_4'}$, which is at least $3\pi R^2\delta_4'$ because $\mathscr{P}_R$ contains the $R$ neighbourhood of its zeros. Recall that $q^d$ is obtained as a double of $\mathscr{P}_R$ along its boundary edges. Therefore, the injectivity radius of $\mathscr{P}_R$ is at least  $\frac{3}{2}\pi R^2\delta_4'$. Similarly, the diameter of $\mathscr P_R$ is at most $2\|\mathscr P_R\|\rho_4'$. Together with item (v) of Theorem \ref{thm:Minsky:polygon}, this implies that the diameter of $\mathscr P_R$ is bounded from above by some constant depending on $R$ and the topology of $S$. This finishes the proof of Claim 4 and hence the proof of the lemma.
 \end{proof}

We will need an estimate on the total error that accumulates in using that the harmonic map outside $\mathscr{P}_R$ is a projection.

 \begin{lemma} \label{lem:energy:difference:infinity}
   Let $R_0$ be the constant from Theorem \ref{thm:Minsky:traintrack}. Let $f:X\to Y$ be a harmonic diffeomorphism
   between closed hyperbolic surfaces $X,Y\in\T(S)$. Let $e(f)$ be the energy density with respect to the metric induced by the Hopf differential $\Hopf(f)$. Then for any $R>R_0$,
   $$ \int_{X\backslash\mathscr{P}_R} |e(f)-2| dA\leq C|\chi(S)|e^{-R/2}$$
   where  $C$ is a constant depending on $\chi(X)$,  $\mathscr{P}_R$ is  Minsky's polygonal region and $dA$ is the area measure induced by the Hopf differential $\Hopf(f)$. 
  \end{lemma}
  \begin{proof}
   Let  $R_0$ be the constant from Minsky's estimate. Then by (\ref{eq:pullback:metric:2}) and Lemma \ref{lem:Minsky:decay},  for any $R>R_0$, and any $p\in X\backslash \mathscr{P}_R$,
   \begin{equation}
       |e(f)-2|\leq Ce^{-\mathrm{dist}(p, \mathscr Z)},
   \end{equation}
   where $C$ is a constant depending on $\chi(X)$, and $\mathrm{dist}(p, \mathscr Z)$ represents the distance from $p$ to the zero set $\mathscr Z$ of $\Hopf(f)$ with respect to the flat metric induced by $\Hopf(f)$.

   Consider the  Voronoi decomposition of $X$ with respect to $\Hopf(f)$, where the 2-cells  are the path components of the set of points which have unique length-minimizing paths to the zero set of $\Hopf(f)$. The number of such 2-cells is exactly the number of zeros of $\Hopf(f)$ counted without multiplicity. Within each 2-cell, consider the horizontal critical segments initiating from the corresponding zero of $\Hopf(f)$. Let $\kappa$ be the order of this zero.  These segments cut the underlying 2-cell into $\kappa+2$ sub-cells, each of which can be identified with a subset of the upper or lower half plane in $\mathbb C$. The total number of such sub-cells is at most  $3(4g-4)$, which corresponds to the case where all zeros of $\Hopf(f)$ are simple.  Integrating $|e(f)-2)|$ over each sub-cell, we obtain
    \begin{eqnarray}\label{eq:energy:difference:upper}
        &&  \int_{X\backslash\mathscr{P}_R}|e(f)-2|dA\\
         \nonumber  &\leq & 3(4g-4)C\int_{\{z\in\mathbb C: |z|\geq R, \mathrm{Im}z\geq0\}} e^{-|z|}dA\\
         \nonumber &\leq& 6(2g-2)C\pi e^{-R/2}.
    \end{eqnarray}

  \end{proof}

   The proof  above also proves the following extension of Lemma~\ref{lem:energy:difference:infinity} from the setting of a closed target $Y$ to the setting of a crowned surface $Y$:
  \begin{lemma} \label{lem:energy:difference:infinity:2}  
   Let $R_0$ be the constant from Theorem \ref{thm:Minsky:traintrack}. Let $f:X\to Y$ be a harmonic diffeomorphism
   from a punctured Riemann surface (possibly disconnected) to a crowned hyperbolic surface $Y$. Let $e(f)$ be the energy density with respect to the metric induced by the Hopf differential $\Hopf(f)$.  Then for any $R>R_0$,
   $$ \int_{X\backslash\mathscr{P}_R} |e(f)-2| dA\leq C|\chi(X)|e^{-R/2}$$
   where $C$ is a  constant depending on $\chi(X)$, $\mathscr{P}_R$ is the Minsky's polygonal region and $dA$ is the area measure induced by the Hopf differential $\Hopf(f)$.
  \end{lemma}

 %%%%%%%%%%%%%%%%%%
  \subsection{Harmonic maps with varying domains}\label{sec:limit:varingdomains} In this subsection, we shall prove a compactness result for harmonic diffeomorphisms with a fixed target but with varying domains. To start, we briefly explain the notion of convergence of harmonic maps with varying domains. For more detailed description of this notion, we refer to \cite[Section 3 and Section 4]{Gupta2019}.
    Let $Y\in\T(S)$ be a fixed hyperbolic surface.
    Let $X_n\in\T(S)$ be an arbitrary divergent sequence of Riemann surfaces, and let $\Phi_n$ be the Hopf differential of the harmonic map $f_n:X_n\to Y$ homotopic to the change of marking. Let $R_m>0$ be a sequence of divergent positive real numbers. Since the number of zeros of $\Phi_n$ is at most $2|\chi(S)|$, it follows that there exists a positive integer $k\leq 2|\chi(S)|$  such that for sufficiently large $m$, the polygonal region $\mathscr{P}_{R_m}(\Phi_n)$ contains exactly $k$ components,  up to a subsequence of $(X_n,\Phi_n)$. Choose a zero for each component of $\mathscr{P}_{R_m}(\Phi_n)$. Let $\mathfrak{p}_n=\{p_{1,n},\cdots, p_{k,n}\}$ be the choice of zeros of $\Phi_n$.
    Consider the family of pointed singular flat surfaces  $(\mathscr P_{R_m}(\Phi_n);\mathfrak{p}_n)$.  By Lemma \ref{lem:inj:lowerbound}, these pointed surfaces have bounded injective radius and bounded diameter, hence (sub)converge to some pointed singular flat surface $(Z_m;\mathfrak{q})$ in the Gromov-Haursdorff topology. Letting $m\to\infty$ and applying a diagonal argument, we get a nested sequence of singular flat surfaces
    $$ (Z_1;\mathfrak{p})\subsetneq (Z_2;\mathfrak{p})\subsetneq \cdots \subsetneq (Z_m;\mathfrak{p})\subsetneq\cdots $$
    such that  there exists a subsequence of $(\Phi_n;\mathfrak{p}_n)$ converging to the pointed singular flat surface $(\cup_{m\geq1}Z_m;\mathfrak{p})$. For simplicity, we still denote this subsequence by $(\Phi_n;\mathfrak{p}_n)$.
    Let $X$ be the Riemann surface underlying $\cup_{m\geq1}Z_m$. Then there exists  a family of quasiconformal embeddings $ \iota_{m,n}: \mathscr{P}_{R_{m}}(\Phi_n)\to X $ with quasiconformal constant uniformly converging to 1 as $n\to1$, whose images exhaust $X$, that is, $\cup_{m\geq1}\lim_{n\to
    \infty}\iota_{m,n}(\mathscr{P}_{R_{m}}(\Phi_n))=X$.

     Applying a standard energy estimate (see for instance \cite[Lemma 3]{Cheng}, \cite[Proposition 3.1]{Wolf1991b}, or \cite[Proposition 1.3]{Schoen1993}), we see that, up to a subsequence if necessary,  the sequence $f_{n}:(X_{n};\mathfrak{p}_n)\to Y$ converges to a harmonic map  $f:(X;\mathfrak{p})\to Y$ with Hopf differential $\Phi$ in the following sense.
     Take an arbitrary compact exhaustion $\{\mathcal{K}_j\}$ of $X$. For each $j$,
    the sequence of composition maps $f_{n}\circ (\iota_{m,n})^{-1}|_{\mathcal{K}_j}:\mathcal{K}_j\to Y$ converges to $f|_{\mathcal{K}_j}:{\mathcal{K}_j}\to Y$ uniformly for sufficiently large $m$.  (We also consider the situation where $Y_n\to Y$. Let $\xi _n:Y_n\to Y$ be a diffeomorphism with Lipschitz constant converging to 1 as $n\to\infty$.  Consider the harmonic diffeomorphism $h_{n}:(X_{n};\mathfrak{p}_n) \to Y_n$. 
    An analogous argument shows 
    that $h_{n}:(X_{n};\mathfrak{p}_n)\to Y_n$ converges to a harmonic map $h:(X;\mathfrak{p})\to Y$ with Hopf differential $\Phi$: for each $j$,  the sequence of composition maps $\xi_n\circ h_{n}\circ (\iota_{m,n})^{-1}|_{\mathcal{K}_j}:\mathcal{K}_j\to Y$ converges to $f|_{\mathcal{K}_j}:{\mathcal{K}_j}\to Y$ uniformly for sufficiently large $m$.)

    In the rest of this paper, for the sake of simplicity, we identify $U\subset X$ with its images $\iota_{m,n}^{-1}(U)\subset X_n$ without mentioning the map $\iota_{m,n}^{-1}$.

   We next state a result regarding limits of sequences of harmonic maps and the limit of the associated Hopf differentials. In this direction, we define our notion of limit in this setting. 

\begin{definition} \label{defn:limitharmonicmap}
Let $X_n, Y_n,Y\in\T(S)$. Let  $f_n: X_n \to Y_n$ be a harmonic diffeomorphism with Hopf differential $\Phi_n$. Let $f_\infty:X_\infty \to Y$ be a harmonic map (not necessarily surjective) from a (possibly disconnected) punctured Riemann surface $X_\infty$ to $Y$ with Hopf differential $\Phi_\infty$.  
We say  $f_n: X_n \to Y_n$ {\it converges to}  $f: X_\infty \to Y$ if (i) each component of $X_\infty$ equipped with the $|\Phi_\infty|$-metric is a Gromov-Hausdorff limit of $X_n$ under the $|\Phi_n|$-metric for some choice of base points $p_n\in \mathscr Z_n$, where $\mathscr Z_n$ is the set of zeros of $\Phi_n$; (ii) the surface $X_\infty$ equipped with the $|\Phi_\infty|$-metric contains all pointed Gromov-Hausdorff limits of $(X_n,p_n)$  under the $|\Phi_n|$-metric with base point $p_n\in \mathscr Z_n$; and (iii) on each Gromov-Hausdorff component $U_\infty$ of $X_\infty$, we have that the harmonic map $f_\infty$ is the limit of the harmonic maps $f_n$ on a domain $U_n \subset X_n$, where $U_n$ converges to $U_\infty$. 
\end{definition}

  \begin{lemma}[compactness]\label{lem:precompactness}
  Let $Y\in\T(S)$ be a hyperbolic surface.
    Let $X_n\in\T(S)$ be an arbitrary divergent sequence of Riemann surfaces with $f_n:X_n\to Y$ being the corresponding harmonic diffeomorphism.  Then there exist a chain-recurrent geodesic lamination $\lambda$ on $Y$,  and  a subsequence  $f_{n_m}:X_n\to Y$ which converges to a surjective harmonic diffeomorphism $f:X\to Y\backslash\lambda$ from some (possibly disconnected) punctured surface $X$. Moreover,
    \begin{itemize}
      \item the quadratic differential $\Hopf(f)$ has a pole of order at least two at each puncture; and
      \item $\lambda =\lim\limits _{R\to\infty}\lim\limits_{m\to\infty}\overline{f_{n} ( X_{n_m}\backslash\mathscr{P}_{R}(\Hopf(f_{n_m})))}.$
    \end{itemize}
  \end{lemma}

 \begin{proof}
    We continue using the notations introduced in the beginning of this subsection. In particular,  $X$ is the Riemann surface underlying the singular flat metric $\cup_{m\geq1} Z_m$.
 
Notice that on each connected component of $X$, the limiting harmonic map $f$ is a harmonic diffeomorphism onto its image (\cite{Wolf1991b}, proof of Proposition 3.4,  following \cite{SchoenYau1978}).
  By Theorem \ref{thm:Minsky:polygon}, for sufficiently large $n$, the image $f_n(X_n\backslash \mathscr{P}_{R_n}(\Phi_n))$ is contained in some $\epsilon_n$ neighbourhood of the geodesic lamination corresponding to the horizontal measured foliation of $\Phi_n$,  with $\epsilon_n\to 0$ as $n\to\infty$.   This implies that   $f_{n}(X_{n}\backslash \mathscr{P}_{R_n}(\Phi_n))$ converges to some geodesic lamination $\lambda$ as  $n\to\infty$ in the Hausdorff topology. Recall that every measured lamination is chain-recurrent. Combining with the fact that the Hausdorff limit of a sequence of chain-recurrent laminations is again chain-recurrent (\cite[Proposition 6.2]{Thurston1998}), we see that
    $\lambda$ is a chain-recurrent.

   We claim that  $f(X)\subset Y\backslash \lambda$.  Suppose to the contrary that there exists some $x\in X$ such that $f(x)\in \lambda$. Then there exists a neighbourhood $V$ of $f(x)$ with $V\subset f(X)$. 
   Since $f_n:X_n\to Y$ converges to $f:X\to Y$, it follows that there exists small neighbourhoods $U_n\subset X_n$ with $f_n(U_n)\subset V$ and also with $U_n$ approximating some fixed region in $X$ and hence at a uniformly bounded distance from the zeroes of $\Phi_n$.   On the other hand, the assumption that $\lambda$ is the Hausdorff limit of $f_{n}(X_{n}\backslash \mathscr{P}_{R_n})$ implies that  there exists a sequence of points $p_n\in X_{n}\backslash \mathscr{P}_{R_n}$  whose distance to the zeros of $\Phi_n$ diverges such that $f_n(p_n)\to f(x)\in \lambda$. In particular, $f_n(p_n)\in f_n(U_n)$ but $p_n\notin U_n$.  This contradicts the fact that $f_n$ is a homeomorphism.

  Recall that (from the second paragraph of the proof)  $f_n(X_n\backslash\mathscr P_{R_n}(\Phi_n))$ converges to $\lambda$ as $n\to\infty$. Combined with the fact $f_n(X_n)\equiv Y$, this implies that $f(X)$, being the limit of $f_n(\mathscr P_{R_n}(\Phi_n))$ as $n\to\infty$, contains $Y\backslash\lambda$. On the other hand, the claim in the previous paragraph implies that $f(X)\subset Y\backslash\lambda$. Hence, we have $f(X)=Y\backslash\lambda$.

  Finally, we show that $f:X\to Y$ is globally injective. Recall that for each component $X^i$ of $X$, the restriction $f|_{X^i}$ is injective.  In particular, $f$ is open. To prove that  $f:X\to Y$ is globally injective, it suffices to show that  for $i\neq j$, we have $f(X^i)\cap f(X^j)= \emptyset$.  Suppose to the contrary that there exist $x^i\in X^i$ and $x^j\in X^j$ such that $f(x^i)=f(x^j)$. Then there exist disjoint neighbourhoods $U^i$ of $x^i$ and $U^j$ of $x^j$ such that $f(U^i)=f(U^j)$. It follows that for   $n$ sufficiently large,
 $f_{n}(U^i)\cap f_{n}(U^j)\neq \emptyset$. Again this contradicts the fact that $f_{n}$ is a homeomorphism for every $n$. 
 \end{proof}

\section{The generalized Jenkins-Serrin problem}\label{sec:GeneralizedJenkinsSerrin}
In this section, we consider minimal surfaces in ${M}\times T$ where $M$ is a hyperbolic surface and $T$ is a tree satisfying some conditions. The story begins with the Scherk's example (\cite{Scherk1835,JS66}), which is a minimal graph over a square in $\mathbb R^2$ with boundary values $\pm \infty$ alternately, \enquote{the Dirichlet problem with infinite boundary values}. This was generalized to minimal graphs over $2n$-gons in $\mathbb R^2$ in \cite{JS66} and over ideal $2n$-gons in $\mathbb H^2$ in \cite{NR02}. (In \cite{JS66} and \cite{NR02,NR07}, the domains also allow strictly convex arcs as part of the boundary.) In these cases, the surface $M$ is a polygon with an even number of edges and the tree $T$ is simply the real axis $\mathbb R$. For our purpose, we need to consider the case where $M$ is the universal cover of a \enquote{hyperbolic crowned surface} and $T$ is a tree dual to some measured foliation. In particular, by extending to trees (which may not admit a folding to a real line), we extend the scope of the results to include the case where $M$ is a hyperbolic ideal polygon with an odd number of edges. This section concerns the uniqueness problem of minimal graphs (see Theorem \ref{thm:minimalgraph:uniqueness}). The existence problem will be addressed in Appendix~\ref{sec:appendix} (see Theorem \ref{thm:minimalgraph:existence}).

\subsection{Minimal graphs over domains in $\H^2$}\label{subsec:minimal:graph}
In this subsection, we collect some results about the minimal graphs over domains in $\H^2$. For more details, we refer to \cite{NR02}.
Consider the unit disk model of $\mathbb{H}^2$. Denote by $(x_1,x_2,x_3)$ the coordinates on the product $\H^2\times \R$. The metric on $\H^2\times \R$ is
$$ d\sigma^2=\frac{dx_1^2+dx_2^2}{F} +dx_3^2$$
where $$F=\left(\frac{1-x_1^2-x_2^2}{2}\right)^2.$$
Let $D\subset \H^2$ be a hyperbolic domain. The graph of a function $u:D\to \R$ is \textit{minimal} if and only if satisfies the \textit{minimal surface equation}:
  \begin{equation}\label{eq:MSE}
    \mathrm{div}\left(\frac{\nabla u}{\tau_u}\right) =0,
  \end{equation}
 where $\tau_u=\sqrt{1+|\nabla u|^2}$, and $\nabla u$ and $\mathrm{div}$ are the  gradient and divergence with respect to $\H^2$.

 Let $u$ be a solution of the minimal surface equation. It is clear that the differential
  \begin{equation}\label{eq:dual:differential}
    du^*:=\frac{F( u_1dx_2-u_2 dx_1)}{\tau_u}
  \end{equation}  is closed
  on $D$, where
  $u_i:=\frac{\partial u}{\partial x_i}$. Locally, we may then define a function $u^*$ on $D$, uniquely up to an additive constant, the \textit{conjugate function} of $u$.
 Geometrically, we can interpret $du^*$ as follows. Let $U\subset D$ be a subdomain and $\alpha\subset \partial U$ a boundary arc with arc length parametrization $s$ such that the domain is on the left. Then
 $$ \int_\alpha du^*=\int_\alpha\left<\frac{\nabla u}{\tau_u},\nu\right> ds$$
 where $\nu$ is the outward unit conormal to $U$ along $\alpha$. (The integral $\int_\alpha du^*$ is called the \textit{flux} of $u$ across $\alpha$.) In particular,
  \begin{equation}\label{eq:u-integral}
    \int_\alpha \left| du^*\right|\leq |\alpha|,
  \end{equation} where $|\alpha|$ is the hyperbolic length of $\alpha$. 

  Suppose $\gamma\subset\H^2$ is a geodesic segment. Consider the infinite strip bounded by the two geodesics which are orthogonal to $\gamma$ and which pass through the endpoints of $\gamma$. For each $\epsilon>0$, there are two level curves contained in this strip each of which consists of points of distance $\epsilon$ to $\gamma$.  Each of these two level curves is said to be an $\epsilon$-translate of $\gamma$. 
  We need the following estimate from \cite[Lemma 1]{NR02}, which we state slightly differently here. Geometrically, this lemma says that the tangent planes near \enquote{divergence points} are almost vertical.

 \begin{theorem}[\cite{NR02} Lemma 1, see also \cite{CollinRosenberg2010} Flux theorem]\label{thm:flux:theorem}
 Let $D\subset \H^2$ be a convex domain. Let $\alpha\subset\partial D$ be a compact geodesic arc. Let  $\alpha_\epsilon\subset D$ be an $\epsilon$-translate of $\alpha$. Let $u:D\to \R$ be a solution of the minimal surface equation (\ref{eq:MSE}).
 \begin{itemize}
   \item  If $u|_\alpha=+\infty$, then  \begin{enumerate}
        \item
         $\left<\frac{\nabla u}{\tau_u},\nu\right>$  converges uniformly to 1 as $\epsilon\to0$, where $\nu$ is the unit  field normal to $\alpha_\epsilon$ and pointing toward $\alpha$;
         \item $\lim\limits_{\epsilon\to0} \int_{\alpha_\epsilon} du^*=|\alpha|$,    where $|\alpha|$ is the hyperbolic length of $\alpha$.
         \end{enumerate}
   \item If $u|_\alpha=-\infty$, then \begin{enumerate}
        \item[(3)]
         $\left<\frac{\nabla u}{\tau_u},\nu\right>$  converges uniformly to $-1$ as $\epsilon\to0$, where $\nu$ is the unit  field normal to $\alpha_\epsilon$ and pointing toward $\alpha$;
         \item[(4)] $\lim\limits_{\epsilon\to0} \int_{\alpha_\epsilon} du^*=-|\alpha|$,    where $|\alpha|$ is the hyperbolic length of $\alpha$.
         \end{enumerate}
 \end{itemize}
 \end{theorem}

\subsection{Minimal graphs in $M\times T$}
\label{subsec:minimal:graph2}
Given a measured foliation $F$ on a crowned hyperbolic surface $Y$, the \emph{critical graph} of $F$ is the union of critical leaves that either connect singularities of $F$ or leave all compacta of $F$. In particular, the critical graph does not include half-infinite but precompact critical leaves. (Compare Definition \ref{def:critical:graph}.)
\begin{definition}[Admissible foliations]\label{def:admissibility:general} 
   A measured foliation $F$ on $Y$ is said to be \textit{admissible} if it satisfies the following properties. 
  \begin{enumerate}[(I)]
      \item   Each component in the complement of the critical graph of $F$ is either (i) precompact, or (ii) a half-infinite cylinder foliated by closed leaves homotopic to a closed boundary geodesic, or (iii) an infinite strip of finite height foliated by homotopically nontrivial bi-infinite leaves such that each end of the strip either spirals around a closed boundary geodesic or is asymptotic to an ideal point of a crown end, or (iv) a half-plane foliated by bi-infinite leaves parallel to a single ideal geodesic boundary arc. 
      \item Each closed boundary geodesic $Y$ is either parallel to the closed leaves of one half-infinite cylinder component of $F$ or is 
      the (spiral) limit of (at least) one infinite strip component of $F$, in the sense that each leaf in the strip limits on the closed boundary geodesic.
      \item Each ideal geodesic boundary arc of $Y$ is parallel to the bi-infinite leaves of one half-plane component of $F$.
      \item Each ideal point of a crown end of $Y$ is asymptotic to (at least) one half-infinite critical leaf.
  \end{enumerate}   
\end{definition}
\begin{remark}
 An admissible measured foliation $F$ on a crowned hyperbolic surface $Y$ has finitely many precompact components, finitely many half-infinite cylinders (at most the number of closed boundary geodesics of $Y$), finitely many infinite strips (since strips are pairwise disjoint), and finitely many half-planes (the same number as that of ideal boundary geodesic arcs of crown ends).  Accordingly, an admissible measured foliation has only a finite number of vertices/singularities, each of finite valence.
\end{remark}

\begin{definition}[Admissible dual trees]\label{def:admissibleTree:general}
Let $Y$ be a crowned hyperbolic surface.
The (metric)  tree dual to the lift of an admissible measured foliation on $Y$ to the universal cover $\widetilde{Y}$  is said to be an \textit{admissible dual tree}.
\end{definition}

\begin{remark}
 Note that the lift of a proper path to a boundary curve or a crown projects to a half-infinite path in the tree.
\end{remark}

   Let $T$ be the dual tree of the lift of some admissible measured foliation $F$ on $Y$.  Let $\iota:\widetilde{Y}\to T$ be the projection map along leaves of  $\widetilde{F}$.
  Consider the boundary behaviour of $\iota$. If $p\in \widetilde{Y}$ approaches some point on a lift of an ideal boundary geodesic arc, then $\iota(p)$ goes to infinity along a half-infinite ray dual to the lift of the half-infinite plane of $F$ given by item (III) of Definition \ref{def:admissibility:general}. If $p$ approaches some point on a lift of a closed boundary geodesic that is parallel to the closed leaves of a half-infinite cylinder of $F$, then $\iota(p)$ goes to infinity along a half-infinite ray dual to the lift of that half-infinite cylinder. If $p$ approaches some point on a lift of a closed boundary geodesic that is spiralled towards by (finitely many) infinite strips of $F$, then $\iota(p)$ goes to infinity along a half-infinite ray dual to the lift of those strips.  Here, distinct lifts of these strips induce adjoining finite segments in the tree which collect to form a single half-infinite ray; this ray is the image in the tree of the lift of a neighborhood of a finite subsegment of the boundary (compare Lemma~\ref{lem:Hopf:differential}(iv).). Let $\partial _{hi}T\subset \partial T$ be the subset of the Gromov boundary $\partial T$ of $T$, which consists of ends determined by the three types of half-infinite rays mentioned above.  In summary, the map $\iota$ induces a \emph{partial boundary map} $\partial \iota:\partial_{ibg} \widetilde{Y}\to \partial_{hi} T$, where $\partial_{ibg}\widetilde{Y}$ is the union of the lift of closed boundary geodesics and ideal boundary geodesic arcs of $Y$. Note that any other admissible measured foliation $F'$ equivalent to $F$ differs from $F$ by an isotopy followed by a sequence of Whitehead moves, that is, splitting or collapsing singularities of admissible measured foliations. Hence the partial boundary map $\partial \iota$ is independent of the choice of $F$ in its equivalence class.    A partial boundary map $\partial_{ibg}\widetilde{Y}\to \partial _{hi}T$ is said to be an \textit{admissible partial boundary map} if it is the partial boundary map of some projection map $\widetilde{Y}\to T$ along the leaves of an admissible measured foliation. 
  
    \begin{remark}\label{rmk:boundary:map}
     Here the map $\iota$ is not defined on the Gromov boundary of $\widetilde{Y}$. For instance, it is not defined at any boundary cusp of $\widetilde{Y}$. In fact, for any sequence of points ${p_i}\in \widetilde{Y}$ converging to a boundary cusp bounded by two ideal geodesic boundary arcs $\gamma^\pm$, the accumulation set of the image sequence $\iota(p_i)$ could be the bi-infinite geodesic on $T$ containing the two half-infinite edges corresponding to the two ideal geodesic boundary arcs $\gamma^\pm$. What we actually need (and describe) is the behavior of $\iota$ as one approaches an ideal geodesic boundary arc of $\widetilde{Y}$. 
 \end{remark} 
   
   Our main result in this section is the following uniqueness result, which we will rely on in places in order to show the well-definedness of some limits of sequences.

 \begin{theorem}\label{thm:minimalgraph:uniqueness}
 Let $Y$ be a crowned hyperbolic surface.  Let $T$ be an admissible dual tree. Then there exists at most one  $\pi_1(Y)$-equivariant minimal graph in $\widetilde{Y}\times T$ over $\widetilde{Y}$ with a prescribed admissible partial boundary map. 
 \end{theorem}

Here a \emph{minimal graph} in $\widetilde{Y}\times T$ over $\widetilde{Y}$ is a map $u:\widetilde{Y}\to T$ that satisfies the following  properties.
\begin{itemize}
    \item There is a discrete set $\{p_i\}\subset \widetilde Y$ of singular points such that for any $p\neq p_i$, there exists a neighbourhood $U\subset \widetilde{Y}$ of $p$, an interval $I\subset \R$, and an isometric embedding $j:I\to T$ such that: (i) $u(U)\subset j(I)$,  and (ii) the composition map $j^{-1}\circ u: U\to I$ satisfies the \emph{minimal surface equation} \eqref{eq:MSE}.
\item  At each singular point $p_i$, 
   there exists a neighborhood $U$ of $p_i$ disjoint from any other singular point, a $k$-valence star $V$ of finite total length ($3\leq k<\infty$), and an isometric embedding $j:V\to T$ that sends the center of $V$ to  $u(p_i)$ such that:  (i) $u(U)\subset j(V)$, and (ii)  the composition map $j^{-1}\circ u:U\to V$ is a projection map along leaves of a (singular) measured foliation on $U$ whose leaf space is isometric to $V$.

\item Moreover, the projection maps $(z, u(z)) \mapsto z$ and $(z, u(z)) \mapsto u(z)$, from the graph of $u$ in $\widetilde{Y} \times T$ to either factor, are harmonic. 
\end{itemize}

In light of the discussion in subsection  \ref{subsec:minimalsuspension}, such a minimal graph induces an (equivariant) conformal and harmonic product map from $\widetilde{X}$ to the product $\widetilde{Y} \times T$, where $\widetilde{X}$ is the conformal structure induced on the image $\{(p,u(p))\in\widetilde{Y} \times T\}$.

We briefly describe the organization of the proof. Of course, we want to compare the minimal graphs of two maps from $\tilde{Y}$ to $T$, say $u:\tilde{Y}\to T$ and $v:\tilde{Y}\to T$. The proof is divided into two steps. In the first step, we  prove that the distance function $\distfunction$ is bounded and the supremum is realized at some point. The idea of this step is to fix a point $p \in \tilde{Y}$ which is not a zero of any Hopf differential, and then, for any $q$, consider the distances
 \[\mathbf{u}(q):=u(q)-u(p),~~\mathbf{v}(q):=v(q)-v(p).\]
These distances are not well-defined in a tree, but we are principally interested in a form 
$\widetilde\Psi=
 (\mathbf{u}-\mathbf{v})d(\mathbf{u}^*-\mathbf{v}^*)$, and we find that this form is well-defined on a neighbourhood of $\widetilde{\partial {Y}}$ and descends to a neighborhood of $\partial Y$.
 Analyzing the level sets of $\distfunction$ near $\partial Y$, we also find that   $d\mathbf{u}^*-d\mathbf{v}^*$ is \enquote{nearly well defined}  on simply connected domains near $\partial Y$ (regardless of the number of ideal geodesics of each crown end), and on cylindrical domains near crown ends which have an even number of ideal geodesics, which will then turn out to be the only case which remains. 
 Estimating these forms,  and applying Stokes' theorem, we prove that $\distfunction$ is bounded and the supremum is realizable. In the second step, we prove that $\distfunction$ is identically zero, using the maximum principle and the fact that harmonic maps are conformal at zeros of Hopf differentials.
 This concludes the outline.
We now describe the argument more precisely.

\begin{remark}\label{rmk:polygon:ADT}
 The existence problem about the equivariant minimal graphs in $\widetilde{Y}\times T$ will be addressed in  the Appendix, see Theorem \ref{thm:minimalgraph:existence}.  As a direct consequence, we are able to parametrize harmonic diffeomorphisms from the complex plane to any ideal hyperbolic $n$-gon, using trees with $n$ half-infinite prongs (and no vertices of valence one).  Let $P$ be an ideal hyperbolic $n$-gon.   Let $\mathrm{ADT}(P)$ be the set of admissible dual trees of $P$. Then there is a bijection between $\mathrm{ADT}(P)$ and the set of harmonic diffeomorphisms from $\C$ to $P$. From \cite[Theorem 3.3]{MulasePenkava} or \cite[Proposition 3.5]{GuptaWolf2016}, it follows that $\mathrm{ADT}(P)$ is homeomorphic to $\R^{n-3}$. (Consider the admissible tree $T_0$ with only one vertex and $n$ half-infinite prongs. Then every other admissible dual tree of $P$ is obtained from $T_0$ by replacing the vertex with a metric tree,  a \enquote{metric expansion} of $T_0$. In other words, the set $\mathrm{ADT}(P)$ is the set of metric expansions of $T_0$.)   
\end{remark}

 For $u: \tilde{Y} \to T$ the equivariant map from $\tilde{Y}$ to $T$, we let $\tilde{X}$ be the Riemann surface underlying the graph $(z, u(z))$, and we let $\tilde{f}: \tilde{X} \to \tilde{Y}$ be the equivariant projection map $(z, u(z)) \mapsto z$ from $\tilde{X}$ to $\tilde{Y}$, and we set $f:X \to Y$ its descent. 

  Let $\{a_i,c_j\}$ be the set of punctures of $X$ labeled in such a way that $a_i$  corresponds to the closed geodesic boundary $\alpha_i$ of $Y$ while  $c_j$  corresponds to the crown end $C_j$ of $Y$. Then we have the following.
\begin{lemma}\label{lem:Hopf:differential}
   \begin{enumerate}[(i)]
     \item Each puncture $a_i$ is a pole of order two, with residue of non-vanishing real part, of $\Hopf(f)$ and $c_j$ is a pole of order $k_j$ of $\Hopf(f)$, where $k_j\geq 3$ is the number of ideal geodesics contained in the boundary of the crown end $C_j$.
     \item The tree $T$ is the dual tree to the horizontal measured foliation of $\Hopf(\tilde{f})$.
     \item Let $\gamma$ be an arbitrary ideal geodesic in $\partial Y$. Then for any lift $\tilde\gamma$ and any compact subsegment $I\subset \tilde{\gamma}$,   there is a convex domain $U\subset \widetilde{Y}$ with $I \subset\partial U$ such that $u(U)$ is
      contained in a single half-infinite edge of $T$ and that $u(p)\to\infty $ as $p\to I$.
     \item Let $\alpha$ be an arbitrary closed geodesic loop in $\partial Y$. Then for any lift $\tilde\alpha$ and any compact subsegment $I\subset \tilde{\alpha}$,   there is a convex domain $U\subset \widetilde{Y}$ with $I \subset\partial U$ and a geodesic ray $r:[0,+\infty)\to T$ such that $u(U)$ is
      contained in the image of $r$ and that $u(p)\to\infty $ as $p\to I$.
        \end{enumerate}
 \end{lemma}
 Here the phrasing that $u(p)\to\infty $ as $p\to I$ means that $u(p)$ leaves all compact sets in the half-infinite (half-closed) edge of $T$ as $p$ approaches the segment $I$. 
  \begin{proof}
 The first two items follow from \cite{Gupta2017} except for the statement that the residue at $a_i$ has non-vanishing real part.
 Suppose that the residue at some $a_i$ is purely imaginary then the  half-infinite cylinder corresponding to $a_i$ is vertical. Let $\omega_d$ be the core curve whose distance to the compact boundary of this half-infinite cylinder is exactly $d$.  By Theorem \ref{thm:Minsky:traintrack}, the length of image $f(\omega_d)$ under the harmonic map $f:X\to Y$ converges to zero as $d\to\infty$. This contradicts the fact that the length of non-trivial simple closed curves on  $Y$ have a uniform lower bound away from zero. Hence the residue can not be purely imaginary.

 For the third item, let $U' \subset \widetilde{X}$ be the half-plane (cf. also \cite{Gupta2017}) corresponding to $\widetilde{\gamma}$. Then $\widetilde{\gamma}\subset \partial \widetilde{f}(U') $ and  $u(\widetilde{f}(U'))$ is contained in the half-infinite edge of $T$ corresponding to $\widetilde{\gamma}$. We may choose $U$ to be a convex domain of $\widetilde{f}(U')$ with $I\subset \partial U$.

 It remains to show the fourth item. By the first item of this lemma, we know that the puncture corresponding to $\alpha$ is a second order pole of $\Phi$. (Here, from the first statement, the residue of this second order pole has non-vanishing real part. In the following argument, the non-vanishing of the real part of the residue will play no role; note that the assumption that $\alpha$ has positive length precludes the case of a purely imaginary residue.) This gives a half-infinite cylinder $C$ in the flat metric $|\Phi|$ (which is horizontal if and only $\Phi$ has purely real residue at this puncture).  Let $\widetilde{C}$ be a lift of $C$ corresponding to $\widetilde{\alpha}$. Then $\tilde{f}(\widetilde{C})$ is a simply connected domain of $\widetilde{Y}$ with $\tilde{\alpha}\subset \partial \tilde{f}(\widetilde{C})$. Moreover, $u(\tilde{f}(\widetilde{C}))$ is a geodesic ray with $u(p)\to\infty$ as $p\in \tilde{f}(\widetilde{C}) $ approaches $\tilde{\alpha}$. We may choose $U$ to be a convex domain  of $\widetilde{f}(\widetilde{C})$ with $I\subset \partial U$.
\end{proof}

\subsection{Two differentials}
\label{subsec:two:differentials}

Recall that $u:\widetilde{Y}\to T$ and $v:\widetilde{Y}\to T$ are two $\pi_1(Y)$-equivariant minimal graphs in $\widetilde{Y}\times T$ with the same admissible partial boundary map.

 Let $\widetilde{X}_u$ be the graph of $u:\widetilde{Y}\to T$ in $\widetilde{Y}\times T$, and $\widetilde{f}_u:
    \widetilde{X}\to \widetilde{Y}$ the equivariant projection map which is harmonic. Then $T$ is dual to the horizontal measured foliation of $\Hopf(\widetilde f_u)$. Let  $f_u$ be  the projection map descended from  $\widetilde{f}_u$.  Consider the map $v:\widetilde{Y}\to T$. Define $\widetilde{X}_v$ and $\widetilde{f}_v$ similarly. Then $T$ is also dual to the horizontal measured foliation of $\Hopf(\widetilde{f}_v)$. Let $X_u=\widetilde{X}_u/\pi_1(Y)$ $X_v=\widetilde{X}_v/\pi_1(Y)$ be the quotient surfaces. Then the horizontal measured foliations $\Hor(\Hopf(f_u))$ and $\Hor(\Hopf(f_v))$ of $\Hopf(f_u)$ and $\Hopf(f_v)$ respectively are topologically equivalent. Let $F_u=f_u(\Hor(\Hopf(f_u)))$ and $F_v=f_v(\Hor(\Hopf(f_v)))$ be the associated measured foliations on $Y$. Since $f_u$ and $f_v$ are both homotopic to the identity, it follows  that $F_u$ and $F_v$ differ by an isotopy and Whitehead moves.  Let $\Sing(u)$ be the set of singular points of $F_u$ and $\Crit(u)$ be the union of critical leaves of $F_u$.   Let $\Reg(u)$ be the complement of $\Crit(u)\cup\Sing(u)$ in ${Y}$. Let $\Regt(u),\Critt(u),\Singt(u)$ be respectively the lifts to $\widetilde{Y}$ of $\Reg(u),\Crit(u),\Sing(u)$. 
    Then 
    \begin{itemize}
        \item 
     $\Regt(u)$ is the subset of points $p\in \widetilde{Y}$ such that $u(p)$ is not a vertex of $T$;
    \item $\Critt(u)$ is the subset of points $p\in \widetilde{Y}$ such that $u(p)$ is  a vertex of $T$ but $p$ admits a neighbourhood whose image under $u$ is a geodesic segment;
    \item $\Singt(u)$ is the subset of points $p\in \widetilde{Y}$ such that $u(p)$ is vertex and $p$ admits a neighbourhood  whose image under $u$ is a star with $m\geq 3$ edges/prongs. 
     \end{itemize}
     Let $\Sing(v)$, $\Crit(v)$,  $\Reg(v)$, $\Singt(v)$, $\Critt(v)$, and $\Regt(v)$ be similarly defined.
 
  We now begin our analysis of how the horizontal measured foliations $F_u$ and $F_v$ align, and the implications for the maps $u$ and $v$.
  
 \begin{lemma}\label{lem:reg:geodesic}
     Let $p\in\Regt(u)\cap \Regt(v)$. Then there exists a neighbourhood $U\subset\widetilde{Y}$ of $p$ such that the convex hull of $u(U)\cup v(U)$ is a geodesic. 
     \end{lemma}
     \begin{proof}
     Let $U$ be a neighbourhood of $p$ such that $u(U)$ and $v(U)$ are both geodesic segments. Let $\Hull(U)$ be the convex hull of $u(U)\cup v(U)$ in $T$. If $\Hull(U)$ is a geodesic segment, then we are done. Otherwise, the assumption that neither $u(p)$ nor $v(p)$ is a vertex of $T$ implies that there exists a subdomain $U'\ni p$ such that $u(U')\cup v(U')$ avoids the (discrete set of) points of $\Hull(U)$ which have valence at least three. It then follows that the convex hull of $u(U')\cup v(U')$ is a geodesic segment in $T$.
     \end{proof}

 Consider the minimal graphs $u: \widetilde{Y}\to T$ and $v: \widetilde{Y}\to T$. If there exists a folding $\xi:T\to \R$ (cf. \cite{Farb-Wolf} or \cite{Morgan-Otal}) then we can compare $u$ and $v$ by considering the difference $\xi\circ u-\xi\circ v$. Theorem \ref{thm:minimalgraph:uniqueness} then follows directly from  the argument  in the proof of step 6 (uniqueness) of \cite[Theorem 3]{NR02}. But such a folding does not exist in general. Nevertheless, the observation below (Lemma~\ref{lem:Psi:definition}) allows us to get past this obstruction. 

  Let $p\in\Regt(u)\cap \Regt(v)$. Let  $U\subset\widetilde{Y}$  be a neighbourhood of $p$ such that the convex hull of $u(U)\cup v(U)$ is a geodesic, which is denoted by $\mathbf{L}$, (Lemma \ref{lem:reg:geodesic}). 
  We notice that this geodesic admits two orientations; we pick one of these. Now, let us define two functions $\mathbf{u}$ and $\mathbf{v}$ on $U$ as follows.
  For each $q\in U$, let
 \[\mathbf{u}(q):=u(q)-u(p),~~\mathbf{v}(q):=v(q)-v(p)\]
 be respectively the oriented distance from $u(p)$ to $u(q)$ and from $v(p)$ to $v(q)$; here we may measure distance along the geodesic we have constructed. Let $ \mathbf u^*$, $ \mathbf v^*$, $ d\mathbf{u}^*$  and  $ d\mathbf v^*$ be  defined as in (\ref{eq:dual:differential}). Notice that both $\mathbf{u}-\mathbf{v}$ and $d(\mathbf{u}^*-\mathbf{v}^*)$ are independent on the particular choice of $p\in U$. If we reverse the orientation of the geodesic $\mathbf{L}$, then  $\mathbf{u}-\mathbf{v}$ differs by a negative sign, so does $d(\mathbf{u}^*-\mathbf{v}^*)$. This implies that the differential $(\mathbf{u}-\mathbf{v})d(\mathbf{u}^*-\mathbf{v}^*)$
 is independent on the choice of the orientation of $\mathbf{L}$. Therefore, we obtain a well defined differential $\widetilde\Psi$  on $\Regt(u)\cap \Regt(v)= \widetilde{Y}-u^{-1}(u(\widetilde{\Sing}(u)))\cup v^{-1}(v(\widetilde{\Sing}(v)))$, the complementary region of the union of the preimages of the vertices of $T$ by $u$ and $v$. Locally, the differential $\widetilde{\Psi}$  can be represented as
 $(\mathbf{u}-\mathbf{v})d(\mathbf{u}^*-\mathbf{v}^*)$.

 Notice that $\widetilde{\Psi}$ is not closed in general. In fact, we have
  \begin{eqnarray}
  \nonumber d \widetilde{\Psi}&=&
    \left[(\mathbf{u}_1-\mathbf{v}_1) \left( \frac{\mathbf{u}_1}{\tau_{\mathbf{u}}}-
     \frac{\mathbf{v}_1}{\tau_{\mathbf{v}}} \right)+
     (\mathbf{u}_2-\mathbf{v}_2) \left( \frac{\mathbf{u}_2}{\tau_{\mathbf{u}}}-
     \frac{\mathbf{v}_2}{\tau_{\mathbf{v}}} \right)\right]dx_1dx_2\\
     &=&\left(\frac{\tau_{\mathbf{u}}
     +\tau_{\mathbf{v}}}{2}\right)\left[
     \left( \frac{\mathbf{u}_1}{\tau_{\mathbf{u}}}-
     \frac{\mathbf{v}_1}{\tau_{\mathbf{v}}} \right)^2
     +\left( \frac{\mathbf{u}_2}{\tau_{\mathbf{u}}}-
     \frac{\mathbf{v}_2}{\tau_{\mathbf{v}}} \right)^2
     + \frac{1}{F}\left( \frac{1}{\tau_{\mathbf{u}}}-
     \frac{1}{\tau_{\mathbf{v}}} \right)^2
     \right]dx_1dx_2 \label{eq:dPsi}
  \end{eqnarray}
  where $\mathbf{u}_i=\frac{\partial\mathbf{u}}{\partial x_i}$, $\mathbf{v}_i=\frac{\partial\mathbf{v}}{\partial x_i}$, $\tau_\mathbf{u}$ and $\tau_\mathbf{v}$ are defined in \eqref{eq:MSE}; the second identity is taken from \cite[Page 280]{NR02}.

 We first extend the domain of definition of $\widetilde{\Psi}$ to the lift of a cylindrical neighbourhood of $\partial Y$.

  \begin{lemma}\label{lem:Psi:definition}
    Let $C$ be a  closed geodesic loop boundary component or crown end of $Y$. There exists a cylindrical neighbourhood $U$ of $C$ such that the following holds. 
    \begin{enumerate}[(i)]
        \item The differential  $\widetilde{\Psi}$ is well defined on the lift of $U$ to the universal cover and descends to a differential $\Psi$ on $U$.
        \item For any simple arc $\eta$ of $U$ which cuts $U$ into a simply connected domain, the restriction of $u$ and $v$ to $U\setminus \eta$ are minimal graphs in $\H^2\times \R$.
    \end{enumerate}
  \end{lemma}
  \begin{proof} 
    Notice that $Y$ may have geodesic boundary components or crown ends.     To prove the Lemma, it suffices to find a cylindrical neighbourhood for each boundary component which satisfies the two mentioned properties. 
    
    Let $\alpha$ be an arbitrary geodesic boundary loop of  $Y$. Let $I\subset \widetilde{\alpha}$ be a segment of a lift $\widetilde{\alpha}$ of $\alpha$ to $\widetilde{Y}$ with $\ell(I)>\ell(\alpha)$. Then by item (iv) in Lemma \ref{lem:Hopf:differential}, we see that there exists a convex neighbourhood $U_I\subset \widetilde{Y}$ of $I$ with $I\subset\partial U_I$ and a geodesic ray $r:[0,+\infty)\to T$ such that 
    \begin{itemize}
        \item both $u(U_I)$ and $v(U_I)$ are contained in the image of $r$,
        \item  $u(p)\to\infty$ and $v(p)\to\infty$, as $p\to I$.
    \end{itemize}
     By the definition of $\widetilde{\Psi}$, we see that $\widetilde{\Psi}$ is well defined on any point in $U_I$. Any cylindrical subdomain of the quotient of $U_I$, which admits $\alpha$ as a boundary component, satisfies the conditions (i,ii).

    Now we move on to consider the case of crown ends. Let $C$ be a crown end with ideal geodesic boundary arcs $\gamma_1,\cdots,\gamma_k$ labelled cyclically.  
    Let $\widetilde{X}_u$ and $\widetilde{X}_v$ be respectively the minimal graphs of $u:\widetilde{Y}\to T$ and $v:\widetilde{Y}\to T$. Let $f_u:\widetilde{X}_u\to \widetilde{Y}$ and $f_v:\widetilde{X}_u\to \widetilde{Y}$ be the harmonic projection maps described at the outset of this subsection. Consider the horizontal foliations $\Hor(\Hopf(f_u))$ and $\Hor(\Hopf(f_v))$ of the Hopf differentials of $f_u$ and $f_v$, respectively.  Recall that $\Crit(u)\subset Y$ is the image under $f_u$ of the critical leaves of $\Hor(\Hopf(f_u))$, which is also the preimage of the vertices of $T$ under $u$. Of course, $\Crit(v)$ is defined similarly. For the cusp $P_i$ bounded  by $\gamma_i$ and $\gamma_{i+1}$, there exists a cusp neighbourhood $W_i$ such that $W_i\cap (\Crit(u)\cup\Crit(v))$
    consists of half-infinite simple arcs approaching the cusp $P_i$. Observe that these arcs are contained in the interior of $W_i$. We may take smaller cusp neighbourhoods so that for any $i\neq j$, we have that $W_i$ and $W_j$ are disjoint.  We may then take $U$ to be  a cylindrical neighbourhood of the crown $C$ such that
    \begin{itemize}
      \item[(a)] $\partial {U}\cap \mathrm{int}({W_i})$ is connected and nonempty,
      \item[(b)] $(\Crit(u)\cap U) \subset \cup_i {W_i}$ and $(\Crit(v)\cap {U}) \subset \cup_i {W_i}$.
    \end{itemize}
    In particular, this proves  item (ii). It remains to prove (i).    Let $\widetilde{W_i}$ be a connected component of the lift of $W_i$ to the universal cover $\widetilde{Y}$. Then 
    for each $W_i$, we have that the images $u(W_i)=v(W_i)$, which is isometric to the real line $\R$ (by applying Lemma~\ref{lem:Hopf:differential}(iii) twice). By the definition of $\widetilde{\Psi}$, we see that it is well defined on $\widetilde{W_i}$, and so $\tilde{\Psi}$ descends to a differential $\Psi$ on $W_i$. 
    
    Consider the complementary components of $\Crit(u)\cup \Crit(v)$ in $U$. For each $\gamma_i$, there exists exactly one of these components, say $R_i$, which contains $\gamma_i$ in the boundary. Let $\widetilde{R_i}$ be a connected component of the lift of $R_i$ to the universal cover $\widetilde{Y}$. Let $\widetilde{\gamma_i}$ be a lift of $\gamma_i$ which is contained in $\partial \widetilde{R_i}$. 
    Let $e_i$ be the half-infinite edge,   corresponding to the leaves of $\Hor(\Hopf(f_u))$, or equivalently $\Hor(\Hopf(f_v))$, in the half-plane parallel to $\widetilde{\gamma_i}$, of $T$ (here the finite endpoint of $e_i$ corresponds to the critical leaf of $\Hor(\Hopf(f_u))$, or equivalently $\Hor(\Hopf(f_v))$, bounding the half plane parallel to $\widetilde{\gamma_i}$). Then both  $u(\widetilde{R_i})$ and   $v(\widetilde{R_i})$ are contained in $e_i$. By the definition of $\widetilde{\Psi}$, we see that it is well defined on $\widetilde{R_i}$, which descends to a differential $\Psi$ on $R_i$.
    Since $U=\cup_i(R_i\cup W_i)$, it follows that $\Psi$ is well defined on $U$, and that $\widetilde{\Psi}$ is well defined on $\widetilde{U}$. This finishes the proof.
  \end{proof}

  If a crown end $C$ consists of an even number of ideal geodesic arcs $\gamma_1,\cdots,\gamma_{2n}$, then we can say a bit more about  $u$ and $v$ as follows. Let $\Hor(\Hopf(f_u))$ and $\Hor(\Hopf(f_v))$ be as defined in the above proof (see also the beginning of Section \ref{subsec:two:differentials}). Let $U$ and $R_i$ also be as defined in the above proof. Then the restriction to $U$ of both $F_u:=f_u(\Hor(\Hopf(f_u)))$ and $F_v:=f_v(\Hor(\Hopf(f_v)))$ are orientable. We choose one orientation which induces the same orientation on $R_i\subset U$.  Each choice of orientation yields  two 1-forms $\omega_u$ and $\omega_v$ on $U$, defined by $F_u|_U$ and $F_v|_U$, respectively. For each cusp $P_i$ of $C$, recalling section~\ref{subsec:meromorphic:leaves}, 
  let $\mathrm{Strip}_i(u)$ (resp. $\mathrm{Strip}_i(v)$) be the  restriction to the cusp region $W_i$ (mentioned in the proof of Lemma \ref{lem:Psi:definition}) of the union of those (possibly degenerate) closed infinite-strips of $F_u$ (resp. $F_v$) that approach $P_i$.  The heights of $\mathrm{Strip}_i(u)$ and $\mathrm{Strip}_i(v)$ are determined by $T$, hence are equal. Let $h_i\geq0$ be the height. Then for any oriented simple closed curve $\alpha$ in $U$, we have 
\begin{equation}\label{eq:int:omega:alpha}
      \int_\alpha \omega_u=\int_\alpha \omega_v= a\sum_{i=1}^{2n} (-1)^{i}h_i, 
  \end{equation}
  where $a\in\{1,-1\}$ depends on the orientation we select for the  foliations on $U$ as well as the orientation of $\alpha$. In other words, both $u$ and $v$ induce functions from $U$ to $\R/ (h\Z)$, where $h=\sum_{i=1}^{2n} (-1)^{i}h_i$. In particular, this implies that both $u-v$ and $du^*-dv^*$ are well defined on $U$. We summarize the discussion below.  
  \begin{lemma}
  \label{lem:conjugate:differential}
    Suppose that $C$ is  a crown end of $Y$ which  consists of an even number of ideal geodesic arcs. Let $U$ be the cylindrical neighbourhood of $C$ obtained from Lemma \ref{lem:Psi:definition}. Let $F_u$ and $F_v$ be defined as above.   Then there are compatible choices of orientations of the restriction $F_u|_U$ and $F_v|_U$ so that both $u-v$ and $du^*-dv^*$ are well defined on $U$. 
  \end{lemma}
  
   For any crown end with $k$ ideal geodesic arcs $\gamma_1,\cdots,\gamma_k$, where $k$ is not necessarily an even number, we have the following more general result. We continue with the notations as above. Recall that $h_i$ is the height of the strip of $F_u|_U$  approaching the cusp $P_i$, which is also the height of the strip of $F_v|_U$  approaching the cusp $P_i$. 
     Consider the measured foliations $F_u$. Let $\mathrm{HP}_i(u)$ be the half plane corresponding to the ideal geodesic arc $\gamma_i$.  Consider the half-infinite strip $\mathrm{Strip}_i(u)$ of $F_u$ approaching $P_i$ (constructed in the paragraph above \eqref{eq:int:omega:alpha}). Let $\mathbf{G}$ the leaf space of $F_u|_{\cup_i(\mathrm{HP}_i(u)\cup \mathrm{Strip}_i(u))}$, i.e. the restriction of $F_u$ to $\cup_i(\mathrm{HP}_i\cup \mathrm{Strip}_i(u))$. (Informally, then, the graph $\mathbf{G}$ is a collection of half-infinite prongs attached to a circle.) 
Let $E_i$ be the half-infinite edge of $\mathbf{G}$ dual to $\mathrm{HP}_i(u)$. Let $Q_i$ be the vertex of $E_i$, and $\overline{Q_iQ_{i+1}}$ be the finite edge of $\mathbf{G}$ dual to $\mathrm{Strip}_i(u)$. [Here our notation allows $\overline{Q_iQ_{i+1}}$ to be degenerate, i.e. of zero length, so that some successive such edges $\overline{Q_iQ_{i+1}}$ might coincide (when all the strips incident to some successive cusps are degenerate).] Then the orientation of $Y$ induces a cyclic order of all edges of $\mathbf{G}$ as follows:
   \begin{equation}
    \label{eq:edge:cyclicorder1}
     E_1,\overline{Q_1Q_2}, E_2, \overline{Q_2Q_3},\cdots, E_k, \overline{Q_kQ_1}, E_1,
  \end{equation}
   In particular, the orientation of $Y$ induces not only a cyclic order of the half-infinite edges corresponding to the successive finite edges of the graph $\mathbf{G}$ but also a cyclic order of edges incident at any common vertex of $\mathbf{G}$. 
  Moreover, the map $u:\widetilde{Y}\to T$ descends to a map $u_U:U\to \mathbf{G}$.   Considering the foliation $F_v$ similarly, we get a similarly defined map $v_U:U\to \mathbf{G}$.
  In summary, we have the following.
  
  \begin{lemma}
  \label{lem:maps:graph}
    Let  $C$ be a crown end with $k$ ideal geodesic arcs $\gamma_1,\cdots,\gamma_k$. Let $\mathbf{G}$ be the graph defined as above. Then $u$ and $v$ descend to maps $u_U:U\to \mathbf{G}$ and $v_U:U\to \mathbf{G}$, respectively. 
  \end{lemma}

\subsection{Boundedness of $\distfunction$}
\label{subsec:bounded:dist}

We next apply the structure theory of the previous subsection to find that, for the Jenkins-Serrin maps we are studying, any pair of them have their maximum distance realized at (the lift of) an interior point of $Y$. Here we see the importance of ensuring that $\Psi$ is well-defined in Lemma~\ref{lem:Psi:definition}, as we are allowed to then use Stokes theorem, and especially the sign of the expression for $d\Psi$ in \eqref{eq:dPsi} in our computations. The subsection is devoted to a proof of the following Lemma.

  \begin{lemma}\label{lem:uv:localconstant}
  
  Let $u: \widetilde{Y}\to T$ and $v: \widetilde{Y}\to T$ be two minimal graphs sharing the same admissible partial boundary map. Then there exists a point  $p\in \widetilde Y$  which realizes the supremum of the distance function on $\widetilde Y$, i.e. 
   $$ \dist(u(p),v(p))=
   \sup_{q\in \widetilde Y}
   \dist(u(q),v(q)). $$
  \end{lemma}
  \begin{proof}
  Let $U$ be the union of the cylinder neighbourhoods of $\partial Y$ obtained in Lemma \ref{lem:Psi:definition}. Set $K:=Y\setminus U$. Then $K$ is a compact subset. Notice that the distance function $\distfunction$ descends to $Y$.  Suppose to the contrary that the conclusion (of the lemma) does not hold.  Then we have 
    \begin{equation}
        \label{eq:strictlyless}
        \max\limits_{p\in  K} \mathrm{dist}(u(p),v(p))< \sup\limits_{p\in  Y} \mathrm{dist}(u(p),v(p)).
    \end{equation}
    Let $M$ be an arbitrary  constant such that $$\max_{p\in K} \mathrm{dist}(u(p),v(p))< M < \sup_{p\in Y} \mathrm{dist}(u(p),v(p)).$$
    In particular, $M>0$. Consider a component $\Omega$ of  
    $$\{p\in Y:\mathrm{dist}(u(p),v(p))>M\}.
    $$ 
    It is clear that $\Omega\subset U$, as the points in the complement of $U$ in $Y$ have images at distance less than $M$. In the following, we shall show that  $\mathrm{dist}(u(p),v(p))\equiv M$ on $\Omega$. The arbitrariness of $M$ then implies that 
    $$ \max\limits_{p\in K} \mathrm{dist}(u(p),v(p))= \sup\limits_{p\in Y} \mathrm{dist}(u(p),v(p)),$$
    which contradicts (\ref{eq:strictlyless}).

    Consider the topology of $\Omega$. There are three possibilities:
    \begin{enumerate}[(i)]
        \item  $\Omega$ is simply connected and $\overline{\Omega}$ is contained in the interior of $U$,
        \item  $\Omega$ is simply connected but $\overline{\Omega}$ is not contained in the interior of $U$, 
     \item $\Omega$ is multiconnected.
    \end{enumerate}
       
       \vskip 5pt
      \textbf{Case (i)}: $\Omega$ is simply connected and $\overline{\Omega}$ is contained in the interior of $U$. Then we may take $u$ and $v$ as functions valued in $\R$ with signs chosen so that $u-v=M$ on $\partial\Omega$. Therefore,  since both $du^*$ and $dv^*$ are locally well-defined closed differentials, we have
         $\int_{\partial \Omega}\Psi=M\int_{\partial \Omega}(du^*-dv^*)=0.$
    Combining (\ref{eq:dPsi}) and the Stokes theorem, we see that $\nabla u=\nabla v$ on $\Omega$. Hence $u-v$ is a constant function over $\Omega$, which is exactly $M$  since $u-v=M$ on $\partial \Omega$. (Note that the proof here relies on the assumption that the distance function is identically $M$ on $\partial \Omega$ and avoids the value $M$ in $\Omega$, but does not rely on the specific assumption that $\dist(u(p),v(p))>M$ for $p\in \Omega$. So the conclusion still holds if $\dist(u(p),v(p))<M$ for $p\in \Omega$ and $\dist(u(p),v(p))=M$ for $p\in \partial\Omega$. This observation will be used in case (iii-b).)

    \vskip 5pt
      \textbf{Case (ii)}:
    $\Omega$ is simply connected but $\overline{\Omega}$ is not contained in the interior of $U$. 
   We approximate $\Omega$ by a domain $\Omega^{\delta,\epsilon}$ which is slightly separated from $\partial Y$, and we then prove that the error used in the approximation is negligible.
    Notice that  any consecutive pair of ideal geodesics $\gamma_i,\gamma_{i+1}$  in the crown end $\mathcal{C}$ are asymptotic, and hence determines a point at infinity which we denoted by $P_i$. Let $\delta$ be a  small positive constant to be determined. Let $F^\delta_{i}$ be a neighbourhood of $P_i$ which is bounded by segments in $\gamma_i$, $\gamma_{i+1}$, and a horocyclic arc centered at $P_i$ with length $\delta$. Let $\gamma_i^\delta:=\gamma_i\cap  \partial ( \Omega-(\cup_i F_i^\delta))$. Let $\epsilon$ be another small constant.  Let $\Omega^{\delta,\epsilon}\subset \Omega-(\cup_i F_i^\delta)$ be the subsurface in $\Omega$ consisting of points whose distance to the boundary of $Y$ is at least $\epsilon$. The boundary $\partial \Omega^{\delta,\epsilon}$ consists of  arcs $\gamma_{i}^{\delta,\epsilon}$ corresponding to $\gamma_i$,  arcs $h_{i}^{\delta,\epsilon}$ contained in the horocycle boundary of $F_i^\delta$, and a subarc $\omega^{\delta,\epsilon}$ of $\partial\Omega\cap U$. 
    In particular, \begin{equation}\label{eq:horocycle}
     |h_{i}^{\delta,\epsilon}|<\delta
   \end{equation} for each $i=1,2,\cdots, k$.   Since both $u$ and $v$ are well-defined functions over $\Omega$ with values in $\R$,
   \begin{equation} \label{eq:Psi:duv}
      \int_{\partial \Omega\cap U} \Psi =M\int_{\partial \Omega\cap U}(du^*-dv^*), ~ ~     \int_{\omega^{\delta,\epsilon}} \Psi =M\int_{\omega^{\delta,\epsilon}}(du^*-dv^*).
   \end{equation}
 Combining with the fact that both $du^*$ and $dv^*$ are closed differentials on $\Omega$, we see that
   \begin{eqnarray}
   \label{eq:Psi:Omega1}
     &&\int_{\omega^{\delta,\epsilon}} (du^*-dv^*) +\sum_{i=1}^{k}
   \int_{\gamma_i^{\delta,\epsilon}} (du^*-dv^*)
   +\sum_{i=1}^{k}
   \int_{h_i^{\delta,\epsilon}}(du^*-dv^*)\\
   \nonumber
   &=&\int_{\Omega^{\delta,\epsilon}}d(du^*-dv^*)\\
   &=&0. \nonumber
   \end{eqnarray}

    First, we consider the integration over $h_i^{\delta,\epsilon}$, 
   \begin{eqnarray}\label{eq:flux:horocycle}
     \nonumber\left|\int_{h_i^{\delta,\epsilon}}(du^*-dv^*)\right|&\leq & \int_{h_i^{\delta,\epsilon}} (
   |d\mathbf{u}^*|+|d\mathbf{v}^*|)\\
   \nonumber&\leq & \int_{h_i^{\delta,\epsilon}} 2 ds \qquad\qquad(\text{by (\ref{eq:u-integral}) })
   \\
   \nonumber&=&2 |h_i^{\delta,\epsilon}|\\
   & <&2\delta, \qquad\qquad\text{ (by (\ref{eq:horocycle})) }
   \end{eqnarray}

  We move on to consider $\int_{\gamma_i^{\delta,\epsilon}}(du^*-dv^*)$.
   Notice that $u$ and $v$ can be viewed as minimal graphs in $\H^2\times \R$ with $u|_{\gamma_i}=+\infty$ and $v|_{\gamma_i}=+\infty$.   Combining this with Theorem \ref{thm:flux:theorem}, we see  that for any $\eta>0$, we may choose $\epsilon$ small enough so that
   \begin{equation}\label{eq:grad}
     \left|\left<\frac{\nabla u}{\tau_u},\nu\right>-1\right|<\eta~~\text{ and }
     \left|\left<\frac{\nabla v}{\tau_v},\nu\right>-1\right|<\eta
   \end{equation}
   on the lift $\gamma_i^{\delta,\epsilon}$, where $\nu$ is the unit normal field of $\gamma_i^{\delta,\epsilon}$ pointing toward $\gamma_i$.   Hence
   \begin{eqnarray}
   \nonumber
     \left|\int_{\gamma_i^{\delta,\epsilon}}(du^*-dv^*)\right|
     &=&
      \left|\int_{\gamma_i^{\delta,\epsilon}}  \left(\left<\frac{\nabla u}{\tau_u},\nu\right>-\left<\frac{\nabla v}{\tau_v},\nu\right>\right)ds \right|\\
    \nonumber &=&\left|\int_{\gamma_i^{\delta,\epsilon}} \left((\left<\frac{\nabla u}{\tau_u},\nu\right>-1)-(\left<\frac{\nabla v}{\tau_v},\nu\right>-1)\right)ds \right|\\
   \nonumber  &\leq & \int_{\gamma_i^{\delta,\epsilon}} \left(\left|\left<\frac{\nabla u}{\tau_u},\nu\right>-1\right|+
     \left|\left<\frac{\nabla v}{\tau_v},\nu\right>-1\right|\right)ds\\
   \nonumber  &\leq & \int_{\gamma_i^{\delta,\epsilon}} 2\eta \qquad\qquad\text{by (\ref{eq:grad}) }\\
     &\leq & 2\eta (|\gamma_i^\delta|+2\epsilon). \label{eq:flux:gamma}
   \end{eqnarray}

   Combining  (\ref{eq:Psi:duv}), (\ref{eq:Psi:Omega1}), (\ref{eq:flux:horocycle}) (\ref{eq:flux:gamma}), we get
    \begin{equation*}
    \label{eq:Psi:Omega}
     \left|
         \int_{\omega^{\delta,\epsilon}} \Psi\right| \leq \sum_{i=1}^k2M\delta+\sum_{i=1}^k 2M\eta(|\gamma_i^\delta|+2\epsilon).
   \end{equation*}
   Letting  $\epsilon\to 0$, we have
   \begin{equation*}
    \lim_{\epsilon\to0} \left|
     \int_{\omega^{\delta,\epsilon}} \Psi\right|\leq  \sum_{i=1}^{k}2M\eta |\gamma_i^\delta|+ \sum_{i=1}^{k}2M\delta.
   \end{equation*}
   The arbitrariness of $\eta$ then implies that
    \begin{equation*}\label{eq:limit:integral:Psi:M}
    \lim_{\epsilon\to0} \left|\int_{\omega^{\delta,\epsilon}}\Psi\right|\leq 2kM\delta.
   \end{equation*}
 Therefore,
 \begin{equation}
     \label{eq:int:psi:ii}
     \left|\int_{\partial \Omega\cap U}\Psi\right|
 =\lim_{\delta\to0}\lim_{\epsilon\to0} \left|\int_{\omega^{\delta,\epsilon}}\Psi\right|=0.
 \end{equation}

\medskip
 It remains to show that the distance function is constant. Suppose to the contrary that the distance function is not constant. Then $M<\sup_{p\in \Omega}\dist(u(p),v(p))$. Let $M'$ be a constant with $M<M'<\sup_{p\in \Omega}\dist(u(p),v(p))$. Consider a component $\Omega'$ of $\{p\in \Omega:\dist(u(p),v(p))>M'\}$. If $\Omega'$ admits a boundary component that encloses a disc in $\Omega$, then by case (i), we see that the distance function is a constant on $\Omega'$. This contradicts the definition of $\Omega'$. Hence, the distance function is constant. If  $\Omega'$ admits a boundary component $\partial \Omega'\cap U$ which intersects $\cup_i\gamma_i$, then by the discussion as for  \eqref{eq:int:psi:ii}, we know that $\int_{\partial\Omega'\cap U}\Psi=0$. Consider the region $A$ in $\Omega$ bounded by $\partial \Omega \cap U$ and $\partial \Omega'\cap U$. Approximating $\partial A\cap (\cup_i \gamma_i)$ similarly as above and applying the Stokes theorem to the approximating regions using (\ref{eq:dPsi}), we see that the distance function on $A$ is a constant. This contradicts the assumption that $\dist(u(p),v(p))=M$ for $p\in \partial \Omega\cap U$ while $\dist(u(p),v(p))=M'>M$ for $p\in \partial \Omega'\cap U$. Hence, the distance function is constant on $\Omega$.  (Note that the proof here relies on the assumption that the distance function is identically $M$ on $\partial \Omega\cap U$ and avoids the value $M$ in $\Omega$, but does not rely on the specific assumption that $\dist(u(p),v(p))>M$ for $p\in \Omega$. So the conclusion still holds if $\dist(u(p),v(p))<M$ for $p\in \Omega$ and $\dist(u(p),v(p))=M$ for $p\in \partial\Omega\cap U$. This observation, along with the analogous observation we made at the end of the discussion of case (i), will be used in case (iii-b).)

  \vskip 5pt
      \textbf{Case (iii)}:
  $\Omega$ is multiconnected. 
  Let $\zeta$ be a boundary component of $\Omega$. Since $U$ is a cylinder, there are two subcases: either $\zeta$ is homotopic to the core curve of $U$ or it encloses a simply connected domain in $U$. 
  
  \textbf{Case (iii-a)}:  $\zeta$ is homotopic to the core curve of $U$.
 By Lemma \ref{lem:maps:graph}, we see that $u$ and $v$ descend to two maps, say  $u_U:U\to \mathbf{G}$
  and $v_U:U\to \mathbf{G}$, from $U$ to the graph $\mathbf{G}$. Notice that the cyclic order of the vertices $\{Q_i\}$ of $\mathbf{G}$ induces a cyclic order of all edges of $\mathbf{G}$ as follows: 
  \begin{equation}
     \label{eq:edge:cyclicorder}
     E_1,\overline{Q_1Q_2}, E_2, \overline{Q_2Q_3},\cdots, E_k, \overline{Q_kQ_1}, E_1,
  \end{equation}
  where $E_i$ is the half-infinite edge attached to $Q_i$.
   If some edge $\overline{Q_iQ_{i+1}}$ has length zero, then we remove it from the above list, meaning that $E_i$ and $E_{i+1}$ are now consecutive.   Consider the restrictions $u_U:\zeta\to \mathbf{G}$ and $v_U:\zeta\to \mathbf{G}$.  
   
    For any point $p$ in the cusp neighbourhood $W_i$ (defined in the proof of Lemma \ref{lem:Psi:definition}), we see that both $u(p)$ and $v(p)$ are contained in the union $E_i\cup \overline{Q_iQ_{i+1}}\cup E_{i+1}$: this union is isometric to $\R$. Therefore, $d_{\mathbf{G}}(u_U(p),v_U(p))=\mathrm{dist}(u(p),v(p))$ for any $p\in W_i$. 
   Recall from the last paragraph of the proof of Lemma \ref{lem:Psi:definition} that there is a component $R_i$ of  $U\setminus(\Crit(u)\cup\Crit(v))$ which contains $\gamma_i$ in the boundary. In particular, for any point $p\in R_i$, we see that both $u(p)$ and $v(p)$ are contained in the half-infinite edge $E_i$. This implies that $d_{\mathbf{G}}(u_U(p),v_U(p))=\mathrm{dist}(u(p),v(p))$ for any $p\in R_i$. Hence  $d_{\mathbf{G}}(u_U(p),v_U(p))=\mathrm{dist}(u(p),v(p))$ for any $p\in U $ because $U=\cup_i( R_i\cup W_i)$. In particular, this equation of distances holds on $\zeta$.  
   For each $i$, let $p_i$ be a point in $R_i\cap \zeta$. Combining the discussion above and the  assumption that $\distfunction$ is identically $M>0$  on $\zeta$ and the fact that $d_{\mathbf{G}}(u_U(\cdot),v_U(\cdot))=\mathrm{dist}(u(\cdot),v(\cdot))$ on $\zeta$, 
   we see that 
   $d_{\mathbf{G}}(u_U(p_i),Q_i)-d_{\mathbf{G}}(v_U(p_i),Q_i)$ is well-defined and has signs which are alternatively $+$ and $-$ with $i$. 
   
   In particular, this implies that $\mathbf{G}$ has even number of edges. Correspondingly, the crown end $C$ in consideration also has an even number of ideal geodesic arcs.

    It then follows from Lemma \ref{lem:conjugate:differential} that there are compatible choices of orientations of the restrictions  $F_u|_U$ and $F_v|_U$ so that both $u-v$ and $du^*-dv^*$ are well-defined on $U$. We choose an orientation of $F_u|_U$ so that $u-v=M$ on $\zeta$.   Applying a similar argument as in case (ii) we see that $\int_{\zeta}(du^*-dv^*)=0$. Hence $\int_\zeta\Psi=M\int_{\zeta}(du^*-dv^*)=0.$

  To prove that the distance function is a constant over $\Omega$, let $M''>M$ be a constant which is close enough to $M$ such that $\Omega'':=\{p\in\Omega:\dist(u(p),v(p))>M''\}$ admits a boundary component say $\partial \Omega''\cap \Omega$ which is homotopic to $\partial \Omega\cap U$. The discussion in the previous paragraphs yields that $\int_{\partial\Omega''\cap \Omega}\Psi=0$. Applying  Stokes theorem to the cylinder region bounded by $\partial \Omega\cap U$ and $\partial \Omega''\cap \Omega$ and using (5.4) we see that the distance function is a constant on this cylinder region. This contradicts the assumption that the distance function is $M$ on  $\partial\Omega\cap U$ but is $M''>M$ on $\partial \Omega''\cap \Omega$. Hence the distance function is constant over $\Omega$.

  \textbf{Case (iii-b)}: $\zeta$  encloses a simply connected domain $V$ on $U$. In particular, the distance function is strictly less than $M$ on $V$ and identically $M$ on $\partial V\cap U$. The arguments in case (i) and case (ii) still hold (see the notes at the end of the proof of case (i) and case (ii)). Hence, similarly to case (i) and case (ii), we see that the distance function is constant on $V$. 
  For any $r\geq M$, Let $\Omega_r$ be the component of $\{p\in \Omega: \dist(u(p),v(p))>r\}$ that admits a  complementary component, say $V_r$, which contains $V$. Let $\mathbf r$ be the supremum of $r\geq M$ such that $V_r$ is simply connected. Then for each $r<\mathbf r$, we see that the distance function is constant on $V_r$. By continuity, the distance function is constant on the closure $\overline{V_\mathbf{r}}$. On the other hand, the definition of $\mathbf{r}$ implies that 
 for each $r>\mathbf{r}$, the complementary component $V_r$, as it contains the non-simply connected $V_\mathbf{r}$ that is homotopic to the cylinder $U\supset \Omega$, is itself not simply connected.  Then the discussion in case (iii-a) implies that the distance function is constant on $\Omega_r$, so is also constant on the union $\cup_{r>\mathbf{r}} \Omega_r $. Note that 
 $\Omega= \cup_{r>\mathbf{r}} \Omega_r \cup (\overline{V_\mathbf{r}}\cap \Omega)$.
 In summary, the distance function is constant on $\Omega$. This completes the proof.
   \end{proof}
  
   \subsection{Generalized maximum principle}
 \label{subsec:generalized:maximum:principle} 
 We continue with the notations $X_u$, $X_v$, $f_u:X_u\to Y$, $f_v:X_v\to Y$,  $\Sing(u)$, $\Crit(u)$,  $\Reg(u)$, $\Singt(u)$, $\Critt(u)$,  $\Regt(u)$, $\Sing(v)$, $\Crit(v)$,  $\Reg(v)$, $\Singt(v)$, $\Critt(v)$, and $\Regt(v)$ introduced in the beginning of Section \ref{subsec:two:differentials}.

 Recall that
    \begin{itemize}
        \item 
     $\Regt(u)$ is the subset of points $p\in \widetilde{Y}$ such that $u(p)$ is not a vertex of $T$;
    \item $\Critt(u)$ is the subset of points $p\in \widetilde{Y}$ such that $u(p)$ is  a vertex of $T$ but admit a neighbourhood whose image under $u$ is a geodesic segment;
    \item $\Singt(u)$ is the subset of points $p\in \widetilde{Y}$ such that $u(p)$ is vertex and admit a neighbourhood  whose image under $u$ is a star with $m\geq 3$ edges/prongs. 
     \end{itemize}

  \begin{lemma}\label{lem:equi-angle}  Let $p\in Y$ be a point. 
  \begin{enumerate}[(a)]
      \item  If $p\in Y\setminus\Sing(u)$, then there exists a neighbourhood $\Omega$ of $p$ such that the component through $p$ of the projection of $u^{-1}(u(p))$ into $\Omega$  cuts $\Omega$ into two sectors, each of which bounds an angle of $\pi$ at $p$.
      
      \item If $p\in\Sing(u)$ is a singular point such that $u(p)$ is a vertex of valence $m\geq3$.  Then there is a neighbourhood $\Omega\subset Y$ of $p$ such that the component through $p$ of the projection of $u^{-1}(u(p))$ into $\Omega$
      cuts $\Omega$ into $m$ sectors, each of which bounds an angle of $2\pi/m$ at $p$.
  \end{enumerate}
  Similar conclusions also hold for $\Sing(v)$. 
  \end{lemma} 
  \begin{proof}
    Let $f_u:X_u\to Y$ be the harmonic map from the (closed) minimal graph to $Y$.     Consider the horizontal foliation of the Hopf differential $\Hopf(f_u)$ near $f_u^{-1}(p)$. 
    
    For  statement (a), there exists a neighbourhood $\Omega'\subset X_u$ of $f_u^{-1}(p)$ such that the component of $f_u^{-1}\circ u^{-1}(u(p))\cap \Omega'$ containing $f_u^{-1}(p)$ is a smooth curve crossing $\Omega'$. Since $f_u$ is a diffeomorphism, we see that  $\Omega:=f_u(\Omega')$ is cut out by  the component of $u^{-1}(u(p))\cap \Omega$ containing $p$  into two sectors, each of which bounds an angle of $\pi$ at $p$. 
    
    For statement (b), there exists a neighbourhood $\Omega'$ of $f_u^{-1}(p)$ cut out by    the component of $f_u^{-1}\circ u^{-1}(u(p))\cap \Omega'$ containing $f^{-1}_u(p)$  into $m$ sectors, each of which bounds an angle of $2\pi/m$ at   $f_u^{-1}(p)$ (\cite[Section 6]{Strebel1984}). Since $f_u$ is conformal  at  $f_u^{-1}(p)$ (because the Beltrami differential of $f_u$ is zero at $f_u^{-1}(p)$), it follows that $\Omega:=f_u(\Omega')$ is cut out by the component of $u^{-1}(u(p))\cap \Omega$ containing $p$ into $m$ sectors, each of which bounds an angle of $2\pi/m$ at $p$. 
  \end{proof}
  
  \begin{lemma}[{generalized maximum principle}]
  \label{lem:generalized:maximum:principle}
 Let $u: \widetilde{Y}\to T$ and $v: \widetilde{Y}\to T$ be two minimal graphs. If there exists a point  $p\in \widetilde Y$  which realizes the supremum of the distance function on $\widetilde{Y}$, i.e. 
   $$ \dist(u(p),v(p))=
   \sup_{q\in \widetilde Y}
   \dist(u(q),v(q))$$
   then $u=v$.
 \end{lemma}
 
 \begin{remark}
 We will apply this lemma in the case when the two maps $u$ and $v$ share an admissible partial boundary map, but, as we do not need that hypothesis in the proof, we state the result in a more general form.
 \end{remark}
 
\begin{proof}
 Let $q$ be an arbitrary point on $Y$. Consider the images of a neighbourhood of $q$ under $u$ and $v$. There are two possibilities.  
  \begin{enumerate}[(a)]
  \item for any neighbourhood $\Omega$ of $q$,  the convex hull of $u(\Omega)\cup v(\Omega)$ is not a geodesic.
      \item $q$ admits a neighbourhood $\Omega$ such that the convex hull of $u(\Omega)\cup v(\Omega)$ is a geodesic segment.
  \end{enumerate}
  Notice that any singular point $q\in \Sing(u)\cup\Sing(v)$  satisfies the condition (a).  Moreover, the subset of points satisfying the condition (b) is a proper open subset of $Y$, while the subset of points satisfying the condition (a) is a nonempty closed subset of $Y$.
  
  \vskip5pt
  \textbf{Case (a):}  for any neighbourhood $\Omega$ of $p$,  the convex hull of $u(\Omega)\cup v(\Omega)$ is not a geodesic. Let $\Omega$ be a small neighbourhood of $p$ such that both $u(\Omega) $ and $v(\Omega)$ are stable, meaning that smaller neighbourhoods  share the same image as $\Omega$, up to isotopies in $T$ fixing $u(p)$ and $v(p)$ respectively. Suppose that $u(\Omega)$ and $v(\Omega)$ are respectively stars with $m_u\geq 2$ and $m_v\geq 2$ edges, where we take respectively $u(p)$ and $v(p)$ as the vertices of the stars.     Let $\{\Omega_i(u):1\leq i\leq m_u\}$  be the set of complementary sectors in $\Omega$ of the connected locus of $u^{-1}(u(p))\cap \Omega$  containing $p$. Let  $\{\Omega_j(v):1\leq j\leq m_v\}$) be defined similarly.  There is a one-to-one correspondence between $\{\Omega_i(u)\}$ (resp. $\{\Omega_j(v)\}$ ) and the edges of $u(\Omega)$ (resp. $v(\Omega)$).

Suppose $u(p)\neq v(p)$ and $m_u=m_v=2$.  There are two subcases depending on whether the convex hull of $u(\Omega)\cup v(\Omega)$ is a tripod or not. If it is not a tripod, then any point $q$ in the nonempty set $\left(\Omega_1(u)\cup \Omega_2(u)\right)\cap \left(\Omega_1(v)\cup\Omega_2(v)\right)$ satisfies $d(u(q),v(q))>d(u(p),v(p))$, contradicting the assumption that $p$ attains the supremum of $\distfunction$ over $Y$.
  If the convex hull of $u(\Omega)\cup v(\Omega)$ is a tripod, then exactly one of $u(p)$ and $v(p)$  is a vertex of valence three of this tripod (since $v(p)\neq u(p)$). Without loss of generality, we may assume that $u(p)$ is the vertex. Then there exists one of  $\{\Omega_1(v),\Omega_2(v)\}$, say $\Omega_1(v)$, such that any point $q$ in $\Omega_1(v)\cap \left(\Omega_1(u)\cup\Omega_2(u)\right)$, which is nonempty by Lemma \ref{lem:equi-angle}, satisfies $d(u(q),v(q))>d(u(p),v(p))$, contradicting the assumption that $p$ attains the supremum of $\distfunction$ over $Y$.

Suppose $u(p)\neq v(p)$ and either $m_u$ or $m_v$ is bigger than 2. Without loss of generality, we may assume that $m_u>2$.   Consider the stars $u(\Omega)$ and $v(\Omega)$. There exists at most one (open) edge of $u(\Omega)$, say $u(\Omega_1)$, such that any point of this edge to $v(\Omega)$ is strictly less than $d(u(p),v(p))$. Similarly, there exists  at most one (open) edge of $v(\Omega)$, say $v(\Omega_1)$, such that the distance of any point of this edge  to $u(\Omega)$ is strictly less than $d(u(p),v(p))$.  In particular, if 
  \begin{equation}
      \label{eq:Omega:uv}
      \left(\cup_{2\leq i\leq m_u}   \Omega_i(u)\right)\cap \left(\cup_{2\leq j\leq m_v} \Omega_j(v) \right)\neq \emptyset,
  \end{equation}
   then any point in this set  satisfies  
  $d(u(q),v(q))>d(u(p),v(p))$, contradicting the assumption that $p$ attains the supremum of $\distfunction$ over $Y$. Next, we shall prove (\ref{eq:Omega:uv}).   By Lemma \ref{lem:equi-angle}, the set $\Omega_1(u)$ bounds an angle of $2\pi/m_u\leq 2\pi/3$ while $\Omega_1(v)$ bounds an angle of $2\pi/m_v\leq \pi$. It follows that the closures (in $\Omega$) $\overline{\cup_{2\leq i\leq m_u} \Omega_i(u)}$ and  $\overline{\cup_{2\leq i\leq m_v} \Omega_i(v)}$ are sectors based at $p$ of angles at least $\frac{4\pi}{3}$ and $\pi$, respectively. Consequently, the intersection set $$\overline{\cup_{2\leq i\leq m_u}   \Omega_i(u)}\cap\overline{ \cup_{2\leq j\leq m_v} \Omega_j(v)} $$ has positive measure. On the other hand, the complement of ${\cup_{2\leq i\leq m_u} \Omega_i(u)}$ in $\overline{\cup_{2\leq i\leq m_u} \Omega_i(u)}$ is contained in $\Crit(u)\cap \Omega$, and hence has measure zero. Similarly, the complement of ${\cup_{2\leq i\leq m_v} \Omega_i(v)}$ in $\overline{\cup_{2\leq i\leq m_v} \Omega_i(v)}$ also has measure zero.   
  Therefore, 
  $$\cup_{2\leq i\leq m_v}   \Omega_i(u)\cap \cup_{2\leq j\leq m_v} \Omega_j(v) \neq \emptyset,$$ completing the proof of (\ref{eq:Omega:uv}).
  
  The discussion above excludes the possibility that $u(p)\neq v(p)$. Therefore, $u(p)=v(p)$. In particular, $d(u(p),v(p))=0$. The assumption that $p$  attains the supremum of $\distfunction$ over $Y$ implies that $\distfunction\equiv0$  on $Y$. Hence $u=v$.

  \textbf{Case (b):}  $p$ admits a neighbourhood $\Omega$ such that the convex hull of $u(\Omega)\cup v(\Omega)$ is a geodesic segment. We notice that this geodesic admits two orientations; we pick one of these. This allows us to define two functions $\mathbf{u}$ and $\mathbf{v}$ on $\Omega$ as follows:
 \[\mathbf{u}(q):=u(q)-u(p),~~\mathbf{v}(q):=v(q)-v(p),~~\forall q\in \Omega.\] 
 In other words, both $\mathbf{u}$ and $\mathbf{v}$, restricted to $\Omega$, may be viewed as minimal graphs in $\H^2\times \R$ over a bounded simply connected domain. The  maximum principle then implies that the distance function $\distfunction$ is a constant on $\Omega$. Then, again by the  maximum principle, we see that the distance function $\distfunction$  is a constant over the whole component, say $\Theta$, of the subset of points satisfying condition (b) which contains $p$.  By continuity, the distance function is also a constant over the closure $\overline{\Theta}$ of $\Theta$. On the other hand, any point in the (interior) boundary $\partial\Theta\cap Y$, which is nonempty,  satisfies the condition (a). It then follows from \textbf{Case (a)} that $u=v$.  This completes the proof.
\end{proof} 
  
  \subsection{Finishing the proof of Theorem \ref{thm:minimalgraph:uniqueness}} 
   \begin{proof}[Proof of Theorem \ref{thm:minimalgraph:uniqueness}] The theorem now follows from  Lemma \ref{lem:uv:localconstant} and Lemma \ref{lem:generalized:maximum:principle}.
   \end{proof}

\section{Subconvergence of harmonic maps rays}
\label{sec:subconvergence:stretchrays}
The goal of this section is to 
begin the proof of Theorem~\ref{thm:HR:SR}. Roughly that theorem asserts the convergence of some harmonic map rays to a Thurston geodesic. In this section, we prove the {\it sub}convergence of that family (see Theorem \ref{thm:HR:limit:geodesic}). We will be left to prove the (full) convergence, whose proof will occupy sections \ref{sec:stretchline:construction} and \ref{sec:convergence:stretchrays}. Looking ahead, that proof of convergence will involve two steps:  in section~\ref{sec:convergence:stretchrays}, we will show that the harmonic map rays defined from our degenerating family $X_t$ converge to a harmonic map ray from a punctured surface $X_{\infty}$; This will be combined with the work of section~\ref{sec:stretchline:construction} that will show that such a harmonic map ray from $X_{\infty}$, whose range is naturally a family of crowned surfaces, extends to a Thurston geodesic through closed surfaces.

\subsection{Uniformly Lipschitz property of harmonic map rays}
We begin by showing that the parametrization of a harmonic map ray captures a Lipschitz bound. 
\begin{lemma}[Uniformly Lipschitz]\label{lem:1:lip}
 Let $X\in \T(S)$ and $\Phi\in Q(X)$. Let
 $$ \HR_{X,\Phi}(\cdot):[0,\infty) \longrightarrow \T(S) $$
 be the harmonic map ray (see Definition \ref{def:HR}).  Let $f_t: X\to \HR_{X,\Phi}(t)$ be the harmonic map at a fixed time $t$. Then the map $f_s\circ f_t^{-1}:\HR_{X,\Phi}(t)\to\HR_{X,\Phi}(s)$ is $\sqrt{s/t}$-Lipschitz for all $s\geq t>0$. 
\end{lemma}
\begin{proof}
 By \cite[Proposition 4.3]{Wolf1989}, it follows that for all $p\in M$ with $\Phi(p)\neq 0$, the function $|\nu(p,t)|$ is an increasing function of $t\in(0,\infty)$, where $\nu(p,t)$ is the Beltrami differential of  $f_t$.

  Let $\G(p,t)=\log (1/|\nu(p,t)|)$. Then $\G$  is a decreasing function of $t\in(0,\infty)$ at $p$ with $\Phi(p)\neq 0$.
  Let $z=x+\i y$ be a canonical coordinate chart of $\Phi$ near $p$. Then by (\ref{eq:pullback:metric:2}), the pullback metric of $Y_t:= \HR_{X,\Phi}(t)$ by $f_t$ on $X$ is:
  \begin{equation}%\label{eq:pullback}
    f_t^*Y_t=2t(\cosh \G(z,t)+1)\d x^2+2t(\cosh \G(z,t)- 1)\d y^2.
  \end{equation}
  Then for any $s>t>0$,  since $\G$  is a decreasing function of $t\in(0,\infty)$, we have
  $$ \mathrm{Lip} (f_s\circ (f_t)^{-1})|_{f_t(p)}\leq \max\left\{\sqrt{\frac{2s(\cosh \G(p,s)+1)}{2t(\cosh \G(p,t)+1)}}, \sqrt{\frac{2s(\cosh \G(p,s)-1)}{2t(\cosh \G(p,t)-1)}}\right\}< \sqrt{\frac{s}{t}}.  $$
  This implies that $f_s\circ (f_t)^{-1}$ is $\sqrt{(s/t)}$-Lipschitz outside the zero locus of $\Phi$. By continuity, $f_s\circ (f_t)^{-1}$  is $\sqrt{(s/t)}$-Lipschitz on the whole surface  $Y_t$. 
\end{proof}

 For $X,Y\in\T(S)$, let $\HR_{X,Y}:[1,+\infty)\to\T(S)$ be the harmonic map ray such that $\Hopf(X,\HR_{X,Y}(s))=s\Hopf(X,Y)$. In particular, $\HR_{X,Y}(1)=Y$. A direct consequence of Lemma \ref{lem:1:lip} is the following corollary.
 
\begin{corollary}\label{cor:normal:family} 
Let $Y$ be a fixed hyperbolic surface. For any $s>1$, the family of maps $\{\HR_{X,Y}:[1,s]\to \T(S)\}_{X\in \T(S)}$ is uniformly bounded and equi-continuous. 
\end{corollary}

\subsection{Subconvergence of harmonic map rays}
 We first show that harmonic map rays are almost geodesics with respect to the Thurston metric. (We refer to the notion of length of a foliation defined in subsection~\ref{subsec:ExtremalHyperbolicLength}.)
 \begin{proposition}\label{prop:HR:almost:geodesic}
   For for any $\epsilon>0$, there exists  $\mathbf T>0$ depending only on $\epsilon$ and the topology of $S$, such that for
     any harmonic map ray  $\HR_{X,\Phi}$ with $\|\Phi\|=1$ and for all $s>t\geq \mathbf{T}$,
     \begin{equation}\label{eq:length:ratio}
     \sqrt\frac{s}{t} \cdot(1-\epsilon)  \leq \frac{\ell_{Y_{s}}(\Hor(\Phi))}
    {\ell_{Y_t}(\Hor(\Phi))}\leq \sqrt\frac{s}{t}
     \end{equation}
     and
   \begin{equation}\label{eq:HR:almost:geodesic}
      \log\sqrt{\frac{s}{t}}-\epsilon \leq \d_{\mathrm{Th}}(Y_t,Y_s)\leq \log\sqrt{\frac{s}{t}},
   \end{equation}
   where $Y_t=\HR_{X,\Phi}(t)$ and $Y_s=\HR_{X,\Phi}(s)$.
 \end{proposition}
 \begin{proof} 

  Let $Y_s:=\HR_{X,\Phi}(s)$ and let $\lambda:=\Hor(\Phi)$ be  the horizontal measured foliation of $\Phi$.  Then the horizontal foliation of $s\Phi$ is $\sqrt{s}\lambda$.
   From \eqref{eq:phi:length} in Lemma \ref{lem:hyplength:quadraticnorm}, we see that
  \begin{equation*}
   -C\leq  \ell_{Y_s}(\sqrt{s}\lambda)- {2\|s\Phi\|} \leq C
  \end{equation*}
  for some constant $C$ depending only the topology of $S$. Combining with the assumption that $\|\Phi\|=1$, we infer that
  \begin{equation*}
     -C\leq  \ell_{Y_s}(\sqrt{s}\lambda)- 2s \leq C.
  \end{equation*} 
Then for any $\epsilon>0$,  there exists $\mathbf{T}>0$ depending only on $\epsilon$ and the topology of $S$, such that for all $s\geq \mathbf{T}$, we have 
  \begin{equation*}
    e^{-\epsilon}<  \frac{\ell_{Y_s}(\sqrt{s}\lambda)}{2s} <e^\epsilon.
  \end{equation*}
  In particular, for all $s>t\geq \mathbf T$,
\begin{eqnarray}
   \log\frac{\ell_{X_{s}}(\lambda)}
    {\ell_{X_t}(\lambda)}
    &=&  \log \left(
  \frac{\sqrt{t}}{\sqrt{s}}\cdot
  \frac{\ell_{X_{s}}(\sqrt{s}\lambda)}{\ell_{X_t}(\sqrt{t}\lambda)}
  \right) \nonumber\\
   &\geq& \log\left( \frac{\sqrt{t}}{\sqrt{s}}\cdot
   { \frac{2s\cdot e^{-\epsilon}}{2t\cdot e^\epsilon} }
   \right)\nonumber\\
   &\geq&\log\sqrt{\frac{s}{t}}-2\epsilon.
    \label{eq:almostmaximal}
  \end{eqnarray}
 Therefore
  \begin{eqnarray*}
    \d_{\text{Th}}(X_{t},X_{s})&=& \log \sup_{\mu\in\ML(S)}
    \frac{\ell_{X_s}(\mu)}{\ell_{X_t}(\mu)}
    \\ &\geq &\log\frac{\ell_{X_{s}}(\lambda)}
    {\ell_{X_t}(\lambda)}
    \\
   &\geq&\log\sqrt{\frac{s}{t}}-2\epsilon.
  \end{eqnarray*}

  On the other hand, by Lemma \ref{lem:1:lip}, we see that
   \begin{equation}\label{eq:dth:upperbd}
      \d_{\mathrm{Th}}(Y_t,Y_s)
      \leq \log \sqrt{s/t}, ~~\forall s>t\geq1.
   \end{equation}
 \end{proof}

 We now consider the compactness of a  family of harmonic map rays passing through a fixed hyperbolic surface. The first part of the following theorem is a restatement of Theorem \ref{thm:HR:limit:geodesic}.
 
\begin{theorem}\label{thm:HR:limit:geodesic2}
  For any fixed $Y\in\T(S)$, let $X_n\in\T(S)$ be any divergent sequence. Then the sequence of harmonic map rays $\HR_{X_n,Y}:[1,\infty)\to\T(S)$ contains a subsequence which converges to some (reparametrized) Thurston geodesic  locally uniformly.

  If moreover,  $\Hor(\Hopf(X_n,Y))$ converges to some $\lambda\in\PML(S)$ as $n\to\infty$, then $\lambda$ is a subset of the maximally stretched lamination of the limit geodesic. 
\end{theorem}

 \begin{proof}
 By Corollary \ref{cor:normal:family}, we see that $\HR_{X_n,Y}:[1,\infty)\to\T(S)$ contains a subsequence, still denoted  by $\HR_{X_n,Y}$ for simplicity, which converges locally uniformly to some continuous map:  $\mathbf{R}:[1,\infty)\to\T(S)$. 
  Let $\Phi_n:=\Hopf(X_n,Y)$, and  $Y_{n,s}:=\HR_{X_n,Y}(s)$.  Then
 $\mathbf{R}(s)=\lim\limits_{n\to\infty}Y_{n,s}$ and
 $\|\Phi_n\|\to\infty$, as $n\to\infty$.

 Set $\widehat\Phi_n:=\Phi_n/\|\Phi_n\|$ to be a quadratic differential on $X_n$ of unit norm. Let $\HR_{X_n,\widehat\Phi_n}:[\|\Phi_n\|,\infty)\to\T(S)$ be a reparametrization of  the harmonic map ray $\HR_{X_n,Y}:[1,\infty)\to\infty$ by setting $\HR_{X_n,\widehat\Phi_n}(s\|\Phi_n\|)=\HR_{X_n,Y}(s)$ for $s\geq1$.  In particular, $\HR_{X_n,\widehat\Phi_n}(s\|\Phi_n\|)=Y_{n,s}\to \mathbf R(s)$ as $n\to\infty$. Since $\|\Phi_n\|\to\infty$, it follows that for any $\epsilon>0$, there exists $N_\epsilon$ such that for all $n>N_\epsilon$ we have $\|\Phi_n\|>\mathbf T$, where $\mathbf T$ is the constant from Proposition \ref{prop:HR:almost:geodesic}.  Applying Proposition \ref{prop:HR:almost:geodesic} to $\HR_{X_n,\widehat\Phi_n}$, we see that for all $n>N_\epsilon$ and all $s>t>1$, we have 
 \begin{equation*}
 	 \log\sqrt{\frac{s}{t}}-\epsilon \leq \d_{\mathrm{Th}}(\HR_{X_n,\widehat\Phi_n}(t\|\Phi_n\|),\HR_{X_n,\widehat\Phi_n}(s\|\Phi_n\|))\leq \log\sqrt{\frac{s}{t}}.
 \end{equation*}
 By letting $n\to\infty$, we see that
 \begin{equation*}
   \log\sqrt{\frac{s}{t}}-\epsilon \leq \d_{\mathrm{Th}}(\mathbf{R}(t),\mathbf{R}(s))\leq \log\sqrt{\frac{s}{t}}
 \end{equation*}
   By the arbitrariness of $\epsilon$, we see that
   \begin{equation*}
   \d_{\mathrm{Th}}(\mathbf{R}(t),\mathbf{R}(s))= \log\sqrt{\frac{s}{t}}.
 \end{equation*}
This proves that $\mathbf{R}$ is a Thurston geodesic.

\bigskip
It remains to show that $\mathbf{R}$ maximally stretches $\lambda$. Let $\lambda_n:=\Hor(\Phi_n)$.  Notice that $\|\Phi_n\|\to\infty$ as $n\to\infty$. It then follows from Equation (\ref{eq:length:ratio}) that  for any $\epsilon>0$, there exists ${N}_\epsilon$, such that for all $n>N_\epsilon$ and all $s>t>1$,
 \begin{equation*}
   \sqrt{\frac{s}{t}}\cdot (1-\epsilon)\leq
   \frac{\ell_{X_{n,s}}(\lambda_n)}{\ell_{X_{n,t}}(\lambda_n)}
   \leq \sqrt{\frac{s}{t}}.
 \end{equation*}
 Letting $n\to\infty$, we see that 
 \begin{equation}\label{eq:length:maximizing}
  \sqrt{\frac{s}{t}}\cdot (1-\epsilon)\leq \frac{\ell_{X_{s}}(\lambda)}{\ell_{X_{t}}(\lambda)}
  \leq \sqrt{\frac{s}{t}}.
 \end{equation}
 where $X_s=\mathbf{R}(s)$ and $X_t=\mathbf{R}(t)$.
 The arbitrariness of $\epsilon$ then implies that
 \begin{equation*}
  \frac{\ell_{X_{s}}(\lambda)}{\ell_{X_{t}}(\lambda)}
  = \sqrt{\frac{s}{t}}.
 \end{equation*}
 \end{proof}

We end this subsection with a similar estimate for harmonic map dual rays.
 \begin{proposition}\label{prop:AHR:almost:geodesic}
 For  any $\epsilon>0$, there exists $\mathbf T>0$ depending only on $\epsilon$ and the topology of $S$, such that  for  any  harmonic map dual ray $\HDR_{Y,\lambda}$ with scaling $\ell_Y(\lambda)=1$, and for all $s>t\geq \mathbf{T}$,
        \begin{equation}\label{eq:AHR:almost:geodesic}
         \d_{{T}}(\HDR_{Y,\lambda}(t),\HDR_{Y,\lambda}(s))\geq
         \log\sqrt{\frac{s}{t}}-\epsilon.
        \end{equation}
 \end{proposition}

 \begin{proof}
Let $X_t:=\HDR_{Y,\lambda}(t)$ and 
   $\Psi_t:=\Hopf(X_t,Y)$. Then the horizontal measured foliation/lamination of $\Psi_t$ is  $t\lambda$. 
  Hence $\Ext_{X_t}(t\lambda)=\|\Psi_t\|$. Combining this with \eqref{eq:phi:length} in Lemma \ref{lem:hyplength:quadraticnorm} we see that
  \begin{equation*}
      -C\leq 2\Ext_{X_t}(t\lambda)-\ell_Y(t\lambda)\leq C
  \end{equation*}
  for some constant $C$ depending only on the topology of $S$. Inserting the assumption that $\ell_Y(\lambda)=1$ into the above equation yields,
  \begin{equation*}
      -C\leq 2\Ext_{X_t}(t\lambda)-t\leq C.
  \end{equation*}
  Then for any $\epsilon>0$, there exists  $\mathbf T>0$ depending only on $\epsilon$ and the topology of $S$, such that for all $t\geq \mathbf{T}$, we have 
  \begin{equation*}
    e^{-\epsilon}<  \frac{t}{2\Ext_{X_t}(t\lambda)} <e^\epsilon.
  \end{equation*}  
  In particular, for all $s>t\geq \mathbf T$,
  \begin{eqnarray*}
    \d_{{T}}(X_{s},X_{t})&=& \frac{1}{2}\log \sup_{\mu\in\ML(S)}
    \frac{\Ext_{X_{t}}(\mu)}{\Ext_{X_s}(\mu)}
    \\ &\geq&
   \frac{1}{2}\log \frac{\Ext_{X_{t}}(\lambda)}{\Ext_{X_s}(\lambda)} \\
   &=& \frac{1}{2}\log\left(\frac{s^2}{t^2}
    \frac{\Ext_{X_{t}}(t\lambda)}{\Ext_{X_s}(s\lambda)} \right) \\
   & \geq &
   \frac{1}{2}\log\left(\frac{s^2}{t^2}
  {\frac{te^{-\epsilon}/2}{se^\epsilon/2}} \right)\\
   &\geq&\log\sqrt{\frac{s}{t}}-\epsilon.
  \end{eqnarray*}
  
 \end{proof}

%===========================
\section{Construction of piecewise harmonic stretch maps }\label{sec:stretchline:construction}

The goal of this section is to prove Theorem \ref{thm:generalized:stretchmap}, which will be used in Section \ref{sec:convergence:stretchrays} to finish the proof of  Theorem \ref{thm:HR:SR}. 

It is perhaps worth taking a moment to recall how this passage will sit in our general theory.  In this section, we prove the existence of \enquote{piecewise harmonic \ss rays}. These are paths that generalize Thurston's construction of concatenation of stretch lines when the maximally stretched lamination $\lambda$ is not maximal. Here instead of concatenating stretch maps with maximal laminations that extend $\lambda$, we use harmonic maps to define the stretch line on the complementary regions that are not ideal triangles. 

On the other hand, we will later, in Section~\ref{sec:stretchray:uniqueness}, construct a family of maps that extend the Thurston theory for these general laminations in a different way. In that section, given hyperbolic surfaces $Y,Z \in \T(S)$, a harmonic stretch ray will be a limit of a family of harmonic map rays from base points $X_n \in \T(S)$ that all proceed through   $Y$ to $Z$,  as the base points $X_n$ degenerate in $\T(S)$. In some sense, we will pick out, from some possible ways of stretching from $Y$ to $Z$, a canonical path that minimizes a particular energy that we call the \emph{harmonic stretch segment} from $Y$ to $Z$.

With that context in mind, we now begin the discussion of the piecewise harmonic stretch maps.

 %======================
 \subsection{Thurston's construction of  stretch maps
 for maximal geodesic lamination}
 \label{sussec:outline}
Given a closed hyperbolic surface $Y$ and a maximal geodesic lamination $\lambda$, Thurston constructs stretch maps in two steps. 
\begin{itemize}
  \item[Step 1.] 
 The first step \cite[Proposition 2.2]{Thurston1998} is to define a change  in the hyperbolic structure of each ideal triangle complementary to $\lambda$ on $Y$ such that the sides of each triangle are expanded by a factor of $e^t$. The change is realized by a built-in Lipschitz map having Lipschitz constant $e^t$ on the boundary and at most $e^t$ in the interior. This is the only place where Thurston uses the assumption of $\lambda$ being maximal. These changes match up along $\lambda$ so that they change the arc length of the leaves of $\lambda$ by a factor of $e^t$ as well.
 \item[Step 2.]
  The second step \cite[Section 4]{Thurston1998} is to extend the (new) hyperbolic structures on the complement of $\lambda$ over a neighbourhood $\mathcal N$ of $\lambda$ by describing its developing map as an infinite product in the group of isometries of the hyperbolic plane, with the help of an induced measured foliation $F_\mathcal{N}(\lambda)$ on $\mathcal{N}$ transverse to $\lambda$. In this step, Thurston makes  assumptions neither about the existence of transverse measures for $\lambda$, nor about the topological or combinatorial type of the complement of $\lambda$. In other words, this second step works for all geodesic laminations. The hyperbolic structure on $\mathcal{N}$ is determined by the transverse measured foliation $F_\mathcal{N}(\lambda)$ and a \textit{sharpness function} in a neighbourhood of each ideal vertex (spike), which records the length of each horocycle leaf inside each spike of $Y\backslash\lambda$. The hyperbolic structure on $Y\backslash\lambda$ constructed in Step 1 provides a unique sharpness function. So the construction results in a unique hyperbolic structure on $Y$ with a built-in Lipschitz homeomorphism.  This step is carried out in detail in \cite[Section 3]{PapadopoulosTheret2007} and  \cite[Section 5]{Bonahon1996} (See also  \cite{CalderonFarre2021}).
 \end{itemize}
 To generalize the construction from the case of maximal laminations to non-maximal ones, the key point is to change the hyperbolic structure on each component of $Y\backslash\lambda$ in a suitable way. Here we provide an approach to doing this by considering harmonic maps from punctured surfaces to crowned hyperbolic surfaces (possibly with simple closed geodesic boundary components).

  \subsection{Harmonic maps from punctured surfaces to crowned hyperbolic surfaces}\label{subsec:parametrization}
 
   The change of hyperbolic metric on $Y\backslash\lambda$ is based on the following theorem.
 
 \begin{theorem}[\cite{Gupta2017} Theorem 1.2]
 \label{thm:parametrization}
 Let $X$ be a closed Riemann surface with a set of marked points $D=\left\{p_{1}, p_{2}, \ldots, p_{k}\right\}$ with fixed coordinate disks around them. For a collection of principal parts $\mathcal{P}=\left\{P_{1}, P_{2}, \ldots, P_{k}\right\}$ having poles of orders $n_{i} \geq 2$ for $i=1,2, \ldots k,$
  \begin{itemize}
    \item let
          $Q({X}, D, \mathcal{P})$ be the space of meromorphic quadratic differentials on $ X$ with principal part $P_{i}$ at $p_{i},$
          and
    \item $ \mathcal{T}(\mathcal{P})$ be the space of marked crowned hyperbolic surfaces homeomorphic to $ X \backslash D$ with $k$ crowns, each having $\left(n_{i}-2\right)$ boundary cusps, with metric residues (see Definition \ref{def:metric:residue}) twice the absolute value of the real part of the residues of the principal parts $P_{i}$.
    \end{itemize}
    Then we have a homeomorphism
$$
{\Psi}: \mathcal{T}(\mathcal{P}) \rightarrow Q( X, D, \mathcal{P})
$$ 
that assigns, to any marked crowned hyperbolic surface $ Y$, the Hopf differential of the unique harmonic diffeomorphism in the appropriate homotopy class from ${X}\backslash D$ to $ Y$ with prescribed principal parts.
  \end{theorem}
 
We will use Theorem \ref{thm:parametrization} several times. For reader's convenience, we provide a sketch of the proof.  Given a principal part $P$ of order $n$ on the punctured disc and a hyperbolic crown end with $n+2$ ideal geodesics whose metric 
residue is twice the absolute value of the real part of the residues of the principal parts $P$, Gupta \cite[Section 3]{Gupta2017} constructed an asymptotic model harmonic map from the punctured disc to the idea polygon whose Hopf differential has $P$ as the principal part at the puncture point $z=0$. Next, consider a punctured Riemann surface $X\backslash \{p_1,\cdots,p_k\}$ with a compact exhaustion $\{X_i\}$,   a meromorphic differential with principal part $\{P_1,\cdots,P_k\}$ at $p$ and a crowned hyperbolic surface $Y$ with metric 
residue at the crown end corresponding to $p_j$ twice the absolute value of the real part of the residues of the principal parts $P_j$. Then using the constructed asymptotic model map to give boundary values on $\partial X_i$ and solving the Dirichlet problem with these boundary values, we get a sequence of harmonic diffeomorphisms $f_i$ from $X_i$ to its image on $Y$.  Generalizing the estimate developed in \cite{Wolf1991b}, Gupta proved that $f_i$ converges to a unique harmonic diffeomorphism from $X\backslash \{p_1,\cdots,p_k\}$ to $Y$ whose Hopf differential has principal part $P_j$ near the puncture $p_j$. This defines the map  $
{\Psi}: \mathcal{T}(\mathcal{P}) \rightarrow Q( X, D, \mathcal{P})
$ mentioned in Theorem \ref{thm:parametrization}. Since  $\mathcal{T}(\mathcal{P})$ and $Q( X, D, \mathcal{P})$ are simply connected and have the same dimension, by Invariance of Domain, Gupta proved that $\Psi$ is a homeomorphism by showing that $\Psi$ is continuous, proper and injective.

\begin{remarks}\label{rmk:parametrization} 
(1)The original first statement in \cite[Theorem 1.2]{Gupta2017} is for  $n_i\geq3$, but the method also works for  $n_i=2$, so we state the result more generally.  (2) Given an ideal hyperbolic polygon, there is more than one harmonic diffeomorphism from the complex plane to the prescribed hyperbolic polygon (see Remark \ref{rmk:polygon:ADT}). The prescribed principal part at each puncture $p_i$ in the statement of Theorem \ref{thm:parametrization} is needed to ensure that the map $\Psi$ is well-defined. (3) Let $\Phi$ be a meromorphic differential on $ X\backslash D$ having poles of order at least two at each punctured/marked point. We take the principal parts of $\Phi$ as the prescribed principal part in Theorem \ref{thm:parametrization}. Then, by Theorem \ref{thm:parametrization}, there exists a unique (up to isometry homotopic to the identity) hyperbolic crowned surface $ Y$, and a unique harmonic diffeomorphism $f: X\backslash D\to  Y$ in the appropriate homotopy class whose Hopf differential is $\Phi$.
\end{remarks}

 \subsection{Deformation of crowned hyperbolic surfaces}\label{subsec:deformation:crownedsurfaces}
  We start with the following lemma, an analogue of the Lemma \ref{lem:1:lip}. 
  
 \begin{lemma}[Uniformly Lipschitz]\label{lem:bounded:lip}
 Let ${f}_t:{X}\backslash D\to {Y}_t$ be a path of  harmonic diffeomorphisms from a punctured Riemann surface ${X}\backslash D$ to crowned hyperbolic surfaces ${Y}_t$ with Hopf differential ${\Phi}_t=t{\Phi}$. Then for any $0<t<s$, the composition map $f_s\circ (f_t)^{-1}: {Y}_t\to {Y}_s$ is a  Lipschitz map having pointwise Lipschitz constant strictly less than $\sqrt{s/t}$.  Moreover, it extends to  an affine  map of factor $\sqrt{s/t}$ from $\partial{Y}_t$ to  $\partial{Y}_s$. 
 \end{lemma}

 We postpone the proof of Lemma \ref{lem:bounded:lip} to the end of this section.  
\begin{lemma}[Analyticity]\label{lem:analyticity1}
     Let ${f}_t:{X}\backslash D\to {Y}_t$ be a path of harmonic diffeomorphisms from a punctured Riemann surface ${X}\backslash D$ to crowned hyperbolic surfaces ${Y}_t$ with Hopf differential ${\Phi}_t=t{\Phi}$. Then ${Y}_t$, ${f_t}$, and the extended  homeomorphism  ${f_t}\circ ({f_1})^{-1}: \overline{Y_1}\to \overline{Y_t}$ (obtained by Lemma \ref{lem:bounded:lip}) are real analytic in $t>0$.
\end{lemma}

We postpone the proof of Lemma \ref{lem:analyticity1} to the end of Section \ref{sec:limit:HR:PSL}. 
[It is important to note here that we will use the other statements of Theorem~\ref{thm:generalized:stretchmap} in the arguments in the earlier part of Section~\ref{sec:limit:HR:PSL}.]

With the notations as in Lemma \ref{lem:analyticity1}, let $l$ be a proper and bi-infinite vertical leaf of $\Phi$. The image $f_t(l)$ is a proper and bi-infinite arc on $Y_t$. For each fixed $t>0$, by Theorem \ref{thm:Minsky:traintrack}, as the distance from zeros of $t\Phi$ goes to infinity, the images of horizontal leaf segments of $t\Phi$ converge smoothly to segments on the boundary $\partial Y_t$. Hence, both ends of $f_t(l)$ approach $\partial Y_t$ perpendicularly. 
  
  \subsection{Extending foliations across a geodesic lamination}\label{sec:extend:foliation}
  Lemma \ref{lem:bounded:lip} allows us to deform crowned hyperbolic surfaces in a controlled way. To deform a closed hyperbolic surface, we need to \enquote{glue} those deformed crowned hyperbolic surfaces to get a closed hyperbolic surface. This gluing process relies on extending a foliation across a geodesic lamination (Lemma \ref{lem:glue:foliation}). The following lemma is key to that extension.

 \begin{lemma}[Lipschitz line field]\label{lem:field:lip}
  	Let ${f}:{X}\backslash D\to {Y}$ be a harmonic diffeomorphism from a punctured Riemann surface ${X}\backslash D$ to a crowned hyperbolic surface ${Y}$ with Hopf differential ${\Phi}$. Then there exist constants $R_0>0$ and $L=L(R_0)$ such that, for all $R>R_0$, the line field tangent to the pushforward of the vertical leaves of $\Phi$ on $ Y\backslash  f(\mathscr P_R)$ is $L$-Lipschitz, where $\mathscr P_R$ is Minsky's polygonal region of $\Phi$ (see Theorem \ref{thm:Minsky:polygon}).
  \end{lemma}

  \begin{proof}
  	Let $\mathbf{h}$ and $\mathbf{v}$ be respectively the line field tangent to the horizontal and vertical foliations of $ \Phi$ on $ X\backslash \mathrm{Zero}(\Phi)$, where $\mathrm{Zero}(\Phi)$ represents the zero set of $\Phi$.   Let $p_0\in  X\backslash \mathscr P_R$.  Let $\widetilde X$ and $\widetilde Y$ be the universal cover of $ X$ and $ Y$ respectively. Let $\widetilde{\mathbf{h}}$, $\widetilde{\mathbf{v}}$, $\widetilde \Phi$, $\widetilde f$, and $\widetilde p_0$ be respectively lifts of $\mathbf{h}$, $\mathbf{v}$, $ \Phi$, $ f$, and $p_0$.	 Choosing orientation representatives for $\widetilde{\mathbf{h}}$ and $\widetilde{\mathbf{v}}$, we may assume that $\widetilde{\mathbf{h}}$ and $\widetilde{\mathbf{v}}$ are unit vector fields with respect to the flat  $|\widetilde \Phi|$-metric.  To simplify notations, we also denote by  $\widetilde{\mathbf{h}}$ and $\widetilde{\mathbf{v}}$ their pushforward to $\widetilde Y$ by $\widetilde f$.  Let  $|\widetilde{\mathbf{h}}|_{\widetilde Y}$ and $|\widetilde{\mathbf{v}}|_{\widetilde Y}$ be their pointwise length with respect to the hyperbolic metric $\widetilde Y$.  Finally, let $\widetilde \nabla^{\widetilde Y}$ be the covariant derivative with respect to $\widetilde Y$. To prove the lemma, it suffices to prove that  the gradient  $\widetilde \nabla^{\widetilde Y} \frac{\widetilde {\mathbf v}}{|\widetilde {\mathbf v}|_{\widetilde Y}}$  is bounded  on  the lift of $ X\backslash \mathscr P_R$. Since  ${ \frac{\widetilde {\mathbf h}}{|\widetilde {\mathbf h}|_{\widetilde Y}}}$ and $\frac{\widetilde {\mathbf v}}{|\widetilde {\mathbf v}|_{\widetilde Y}}$ are orthogonal to each other (and hence define a frame at each point $X\setminus\mathrm{Zero}(\Phi)$), it is equivalent to prove
  	\begin{itemize}
  		\item[(i)] $\left\|\widetilde \nabla^{\widetilde Y}_{ \frac{\widetilde {\mathbf v}}{|\widetilde {\mathbf v}|_{\widetilde Y}}} \frac{\widetilde {\mathbf v}}{|\widetilde {\mathbf v}|_{\widetilde Y}}\right\|_{\widetilde Y}$ is bounded from above, and 
  		\item[(ii)] $\left\|\widetilde \nabla^{\widetilde Y}_{ \frac{\widetilde {\mathbf h}}{|\widetilde {\mathbf h}|_{\widetilde Y}}} \frac{\widetilde {\mathbf v}}{|\widetilde {\mathbf v}|_{\widetilde Y}}\right\|_{\widetilde Y}$ is bounded from above.
  	\end{itemize}
  	
  	\bigskip
  	 Notice that $\widetilde \nabla^{\widetilde Y}_{ \frac{\widetilde {\mathbf v}}{|\widetilde {\mathbf v}|_{\widetilde Y}}} \frac{\widetilde {\mathbf v}}{|\widetilde {\mathbf v}|_{\widetilde Y}}$ is the geodesic curvature of vertical leaves with respect to the hyperbolic metric on $\widetilde Y$.  Choose the local coordinate near $\widetilde p_0$ so that $\widetilde\Phi =dz^2$. By \eqref{eq:pullback:metric:2}, the pullback of the hyperbolic metric on $\widetilde Y$ near $\widetilde p_0$ is: 
  	\begin{equation*}
  		2(\cosh \mathscr G(z)+1) dx^2 +2(\cosh \mathscr G(z)-1)dy^2,
  	\end{equation*} 
  	where  $\mathscr G(p)= \log (1/|\nu(p)|)$ and $\nu(p)$ is the Beltrami differential of $\widetilde f$.  Let 
  	\begin{equation*}
  		g_{11}= 2(\cosh \mathscr G(z)+1), \quad g_{12}=g_{21}=0, \quad g_{22}=2(\cosh \mathscr G(z)-1)
  	\end{equation*}
  	and
  		\begin{equation*}
  		g^{11}= \frac{1}{2(\cosh \mathscr G(z)+1)}, \quad g^{12}=g^{21}=0, \quad g^{22}=\frac{1}{2(\cosh \mathscr G(z)-1)}.
  	\end{equation*}	
  	Accordingly, the Christoffel symbols of the second type ($\Gamma_{ij}^m=\frac{1}{2}\sum_{l=1}^2 g^{lm}(\partial_ j g_{il}+\partial_i g_{jl}-\partial_l g_{ji})$) are:
  	\[\begin{array}{ll}
  		\Gamma_{11}^1=\frac{1}{2}g^{11}\partial_1 g_{11} =\frac{\sinh \mathscr G }{2(\cosh \mathscr G+1)}\frac{\partial \mathscr G}{\partial x}, &\Gamma_{11}^2=-\frac{1}{2}g^{22}\partial_2 g_{11}=  -\frac{\sinh \mathscr G }{2(\cosh \mathscr G-1)}\frac{\partial \mathscr G}{\partial y},\\
  		\Gamma_{12}^1=\frac{1}{2}g^{11}\partial_2 g_{11} =\frac{\sinh \mathscr G }{2(\cosh \mathscr G+1)}\frac{\partial \mathscr G}{\partial y}, &\Gamma_{12}^2=\frac{1}{2}g^{22}\partial_1 g_{22}=  \frac{\sinh \mathscr G }{2(\cosh \mathscr G-1)}\frac{\partial \mathscr G}{\partial x},\\
  		\Gamma_{21}^1=\Gamma_{12}^1, &\Gamma_{21}^2=\Gamma_{12}^2\\
  			\Gamma_{22}^1=-\frac{1}{2}g^{11}\partial_1 g_{22} =-\frac{\sinh \mathscr G }{2(\cosh \mathscr G+1)}\frac{\partial \mathscr G}{\partial x}, &\Gamma_{22}^2=\frac{1}{2}g^{22}\partial_2 g_{22}=  \frac{\sinh \mathscr G }{2(\cosh \mathscr G-1)}\frac{\partial \mathscr G}{\partial y}.
  	\end{array}\]
  
  	Note that with respect to the chosen coordinates, the vector fields are:
  	\begin{equation}\label{eq:vec:field}
  		 \frac{\widetilde {\mathbf h}}{|\widetilde {\mathbf h}|_{\widetilde Y}}=\frac{1}{\sqrt{g_{11}}}\partial_1, \quad \frac{\widetilde {\mathbf v}}{|\widetilde {\mathbf v}|_{\widetilde Y}}=\frac{1}{\sqrt{g_{22}}}\partial_2,
  	\end{equation}
  	where $\partial_1:=\frac{\partial}{\partial_x}$ and $\partial_2:=\frac{\partial}{\partial_y}$  are coordinate vector fields. 
  	Hence,
  	\begin{eqnarray*}
\widetilde \nabla^{\widetilde Y}_{ \frac{\widetilde {\mathbf v}}{|\widetilde {\mathbf v}|_{\widetilde Y}}} \frac{\widetilde {\mathbf v}}{|\widetilde {\mathbf v}|_{\widetilde Y}}&=& \widetilde \nabla^{\widetilde Y}_{\frac{1}{\sqrt{g_{22}}}\partial_2}\left(\frac{1}{\sqrt{g_{22}}}\partial_2\right)\\
&=& \frac{1}{{g_{22}}}\Gamma_{22}^1\partial_1+\frac{1}{\sqrt{g_{22}}}\left(\partial_2 \frac{1}{\sqrt{g_{22}}}+\frac{1}{\sqrt{g_{22}}}\Gamma_{22}^2\right)\partial_2\\
&=&  \frac{1}{{g_{22}}}\Gamma_{22}^1\partial_1.
  	\end{eqnarray*}
  	Therefore,   
  	\begin{eqnarray}
  \nonumber\left\|\widetilde \nabla^{\widetilde Y}_{ \frac{\widetilde {\mathbf v}}{|\widetilde {\mathbf v}|_{\widetilde Y}}} \frac{\widetilde {\mathbf v}}{|\widetilde {\mathbf v}|_{\widetilde Y}}\right\|_{\widetilde Y}&=&	\left\|	\frac{1}{{g_{22}}}\Gamma_{22}^1\partial_1\right\|_{\widetilde Y}\\
  \nonumber &=& {\frac{\sqrt{g_{11}}}{g_{22}}} |\Gamma_{22}^1|\\
   \nonumber	&=&  {\frac{\sqrt{2(\cosh \mathscr G+1)}}{2(\cosh \mathscr G-1)}}\frac{\sinh \mathscr G }{2(\cosh \mathscr G+1)}\left|\frac{\partial \mathscr G}{\partial x}\right|\\	
  &= & {\frac{1}{\sqrt{2(\cosh \mathscr G+1)}}}\frac{\sinh \mathscr G }{{2(\cosh \mathscr G-1)}}\left|\frac{\partial \mathscr G}{\partial x}\right|\\
  &\leq & \frac{1}{4\sinh (\mathscr G/2)}\left|\nabla \mathscr G\right|,\label{eq:grad:vert:}
  	\end{eqnarray}
  	where $\nabla \mathscr G$ represents the gradient of $\mathscr G$ with respect to the flat $|\Phi|$-metric.
  Similarly,
    	\begin{eqnarray*}
\widetilde \nabla^{\widetilde Y}_{ \frac{\widetilde {\mathbf h}}{|\widetilde {\mathbf h}|_{\widetilde Y}}} \frac{\widetilde {\mathbf v}}{|\widetilde {\mathbf v}|_{\widetilde Y}}&=& \widetilde \nabla^{\widetilde Y}_{\frac{1}{\sqrt{g_{11}}}\partial_1}\left(\frac{1}{\sqrt{g_{22}}}\partial_2\right)\\
&=& \frac{1}{\sqrt{g_{11}g_{22}}}\Gamma_{12}^1\partial_1+\frac{1}{\sqrt{g_{11}}}\left(\partial_1 \frac{1}{\sqrt{g_{22}}}+\frac{1}{\sqrt{g_{22}}}\Gamma_{12}^2\right)\partial_2\\
&=& \frac{1}{\sqrt{g_{11}g_{22}}}\Gamma_{12}^1\partial_1.
  	\end{eqnarray*}
and
  	\begin{eqnarray}
  \nonumber\left\|\widetilde \nabla^{\widetilde Y}_{ \frac{\widetilde {\mathbf h}}{|\widetilde {\mathbf h}|_{\widetilde Y}}} \frac{\widetilde {\mathbf v}}{|\widetilde {\mathbf v}|_{\widetilde Y}}\right\|_{\widetilde Y}&=&	\left\| \frac{1}{\sqrt{g_{11}g_{22}}}\Gamma_{12}^1\partial_1\right\|_{\widetilde Y}\\
  \nonumber &=& {\frac{1}{\sqrt{g_{22}}}} |\Gamma_{12}^1|\\
  	&= &  {\frac{1}{\sqrt{2(\cosh \mathscr G-1)}}}\frac{\sinh \mathscr G }{2(\cosh \mathscr G+1)}\left|\frac{\partial \mathscr G}{\partial y}\right|\\
    &\leq & \frac{1}{4\cosh(\mathscr G/2)} \left|{\nabla \mathscr G}\right|.\label{eq:grad:hor:}
  	\end{eqnarray}

  	In the remainder of the proof, we shall estimate ${|\nabla \mathscr G|}$ and $  \frac{|\nabla \mathscr G|}{\mathscr G}$.	
  By Lemma \ref{lem:Minsky:decay}, we see that 
  	\begin{equation}\label{eq:G:appendix}
  		 \mathscr{G}(\widetilde p_0)\leq 
   \frac{\sinh^{-1}(|\chi( X)|/(R-4)^2)}{\exp(R-4)}  	\end{equation}
   holds for every $\widetilde p_0\in \widetilde{ X\backslash \mathscr P_{R}}$.
  	Together with \eqref{eq:grad:vert:} and \eqref{eq:grad:hor:}, this implies that there exists $R_0>0$ depending only on the topology of $ X$ such that 
  	\begin{equation}\label{eq:grad:hv}
  		\left\|\widetilde \nabla^{\widetilde Y}_{ \frac{\widetilde {\mathbf v}}{|\widetilde {\mathbf v}|_{\widetilde Y}}} \frac{\widetilde {\mathbf v}}{|\widetilde {\mathbf v}|_{\widetilde Y}}\right\|_{\widetilde Y}\leq \frac{|\nabla \mathscr G|}{\mathscr G} \text{ and } \left\|\widetilde \nabla^{\widetilde Y}_{ \frac{\widetilde {\mathbf h}}{|\widetilde {\mathbf h}|_{\widetilde Y}}} \frac{\widetilde {\mathbf v}}{|\widetilde {\mathbf v}|_{\widetilde Y}}\right\|_{\widetilde Y}\leq  {|\nabla \mathscr G|},
  	\end{equation}
  	holds for every $\widetilde p_0\in \widetilde{ X\backslash \mathscr P_{R_0}}$.  	
  	Combining Lemma \ref{lem:grad:decay} and \eqref{eq:G:appendix}, we see that,  there exists positive constants $c_1=c_1(R_0)$ and $c_2=c_2(R_0)$ depending only on $R_0$ and the topology of $ X$ such that  
  	\begin{equation}\label{eq:grad:g}
  {|\nabla \mathscr G(\widetilde p_0)|}\leq c_1 \text{ and }	 \frac{|\nabla \mathscr G|}{\mathscr G}(\widetilde p_0)\leq c_2
  	\end{equation}
  	for every $\widetilde p_0\in \widetilde{ X\backslash \mathscr P_{R_0}}$.   	
  	It then follows from \eqref{eq:grad:hv} and \eqref{eq:grad:g} that, for any $R>R_0$, using that $\mathscr{P}_{R_0}\subset \mathscr{P}_{R}$, the line field $\frac{\widetilde {\mathbf v}}{|\widetilde {\mathbf v}|_{\widetilde Y}}$ is $\max\{c_1,c_2\}$-Lipschitz on the lift of $ {  f( X\backslash \mathscr P_R)}$ to the universal cover $\widetilde Y$.  Hence, the line field $\frac{ {\mathbf v}}{| {\mathbf v}|_{  Y}}$ is $\max\{c_1,c_2\}$-Lipschitz on $  f ({ X\backslash \mathscr P_R})= Y\backslash  f (\mathscr P_R)$.  This completes the proof.
  \end{proof}
\begin{remark}\label{rmk:curvature}
   Theorem \ref{thm:Minsky:traintrack} also holds in the current setting of punctured Riemann surfaces and crowned hyperbolic surfaces.
   \begin{itemize}
       \item  By \cite[Lemma 2.2 (v)]{Wolf1991a}, the geodesic curvature of the pushforward of the horizontal leaves of $\Phi$ in  $ Y\backslash  f(\mathscr P_R)$ goes to zero uniformly as $R\to\infty$. Using the above computation, the horizontal geodesic curvature is 
       \begin{equation}
           \left\|\widetilde \nabla^{\widetilde Y}_{\frac{\widetilde {\mathbf h}}{|\widetilde {\mathbf h}|_{\widetilde Y}}} \frac{\widetilde {\mathbf h}}{|\widetilde {\mathbf h}|_{\widetilde Y}}\right\|_{\widetilde Y}
      =\frac{\sqrt{g_{22}}}{g_{11}}|\Gamma_{11}^2|
  =\frac{\sinh \mathscr G}{2(\cosh\mathscr G+1)\sqrt{2(\cosh\mathscr G-1)}}\left|\frac{\partial \mathscr G}{\partial y}\right|
       \leq \frac{1}{4\cosh\frac{\mathscr G}{2}}\left|\nabla\mathscr G\right|
       \end{equation}
   which goes to zero uniformly as $R\to\infty$, by \eqref{eq:G:appendix} and \eqref{eq:grad:g}. 
   \item 
    By \eqref{eq:G:appendix} and \eqref{eq:pullback:metric:2}, the hyperbolic length of every vertical leaf of $\Phi$ on $ Y\backslash  f (\mathscr P_R)$ converges to zero exponentially as $R\to\infty$. In particular, for any $\epsilon>0$, there exists a constant $R$ such that $ Y\backslash f (\mathscr P_R)$ is contained in the $\epsilon$ neighborhood of the boundary of ${Y}$.
   \end{itemize}
\end{remark}
  \begin{lemma}[Gluing partial foliations]\label{lem:glue:foliation}
      Let $Y\in\T(S)$  be any closed hyperbolic surface, and let $\lambda$ be any geodesic lamination. Then, for any harmonic diffeomorphism $f:X\to Y\backslash \lambda$ from some (possibly disconnected) punctured surface $X$, the pushforward of the vertical measured foliation of the Hopf differential of $f$ on $Y\backslash\lambda$ extends to a unique measured foliation on $Y$.    \end{lemma}
\begin{proof}
Let $\Phi$ be the Hopf differential of $f$.
   Let $V$ be the pushforward of the vertical measured foliation of $\Phi$ on $Y\backslash\lambda$.  Let $\mathbf v$ be the line field tangent to $V_{Y\backslash\lambda}$, outside the images of zeros of $\Phi$ under $f$,  with respect to the hyperbolic metric on $Y$. By Lemma \ref{lem:field:lip}, we know that there exist  constants $R>0$ and $c>0$ such that $\mathbf v$ is $c$-Lipschitz  on $f(X\backslash \mathscr P_{R})=Y\backslash(\lambda\cup f(\mathscr P_{R}))$, where $\mathscr P_{R}$ is Minsky's polygonal region of $\Phi$. 

   Let $H$ be the pushforward of the horizontal measured foliation of $\Phi$ on $Y\backslash\lambda$. 
   By \cite[Lemma 2.2 (v)]{Wolf1991a} (see also Remark \ref{rmk:curvature}), as $R$ goes to infinity, the hyperbolic geodesic curvature of leaves of $H$ on  $f(X\backslash \mathscr P_{R_0})=Y\backslash(\lambda\cup f(\mathscr P_R))$ converges to zero uniformly. Combining this with Theorem \ref{thm:Minsky:traintrack} (see also the second item in Remark \ref{rmk:curvature}) and the assumption $f(X)=Y\backslash\lambda$, we see that the line field tangent to leaves of $H$ on $Y\backslash(\lambda\cup f(\mathscr P_R))$ extends continuously to $Y\backslash  f(\mathscr P_R)$ by assigning to each point $p$ in $\lambda$ the line of tangent direction to $\lambda$ at $p$. On the other hand, by \eqref{eq:pullback:metric}, (the pushforward of) leaves of $H$ and $V$ are orthogonal to each other. Hence, the restriction to $Y\backslash(\lambda\cup f(\mathscr P_R))$ of the line field $\mathbf{v}$ continuously extends to a line field $\mathbf{v}'$ on $Y\backslash f(\mathscr P_R)$ by assigning to each point $p$ in $\lambda$ the line of normal direction to $\lambda$.  Notice that the restriction to $\lambda$ of the extended line field $\mathbf{v}'$ is Lipschitz (with respect to the hyperbolic metric on $Y$ instead of the intrinsic metric on $\lambda$).  By Lemma \ref{lem:field:lip}, the restriction to $Y\backslash (\lambda\cup f(\mathscr P_R))$ of the extended line field $\mathbf{v}'$  is also Lipschitz. Hence, the extended line field $\mathbf{v}'$ is Lipschitz on $Y\backslash f(\mathscr P_{R})$. Therefore, we can continuously extend the partial foliation $V_{Y\backslash(\lambda\cup\mathscr P_R)}$ to a unique foliation  $V_{Y\backslash\mathscr P_R}$ on $Y\backslash\mathscr P_R$ by integrating the extended line field $\mathbf v'$. Together with the partial foliation $V_{f(\mathscr P_R)}$, this extends the foliation $V$, originally defined on $Y\backslash\lambda$, to a unique foliation $V'$ on the whole surface $Y$. Recall that $V$ admits a transverse measure arising from the vertical transverse measure of $\Phi$. Define a transverse measure for $V'$ near $\lambda$ as follows: for any arc $\eta$ near $\lambda$ and transverse to $\lambda$, we project $\eta$ to a segment $\hat k$ of $\lambda$ along leaves of $V'$ and define intersection number between $\eta$ and $V'$ to be half the hyperbolic length of the geodesic segment $\hat k$.     By Theorem \ref{thm:Minsky:traintrack}, this construction is well defined. Again by Theorem \ref{thm:Minsky:traintrack}, these two transverse measures coincide on $Y\backslash(\lambda\cup f(\mathscr P_R))$. Hence, we construct a transverse measure on $V'$ whose restriction to $Y\backslash\lambda$ is exactly the transverse measure of $V$. This completes the proof.
\end{proof}

For convenience of a later reference, we summarize some basic properties of the extended vertical foliation as follows, where the first conclusion is proved in the second paragraph of the proof of Lemma \ref{lem:glue:foliation} and the second conclusion follows from Theorem \ref{thm:Minsky:traintrack} (also follows from \eqref{eq:pullback:metric:2}, Lemma \ref{lem:Minsky:decay} and Remark \ref{rmk:curvature}).
\begin{lemma}\label{lem:intersection:extend}
7   Let $V'$ be the extended measured foliation in Lemma \ref{lem:glue:foliation}. Then $V'$ is transverse to $\lambda$ and intersects $\lambda$ orthogonally at every intersection point.   Furthermore, for any arc $\kappa$ in a small neighborhood of $\lambda$ that is transverse to $V'$, the intersection number $i(\kappa,V')$ equals half of the hyperbolic length of projection of $\kappa$ to a leaf of $\lambda$ along leaves of $V'$. 
\end{lemma}

  \subsection{Proof of Theorem \ref{thm:generalized:stretchmap}}\label{subsec:construction:PSL}
  We now construct piecewise harmonic stretch lines based on Lemma \ref{lem:bounded:lip}, Lemma \ref{lem:field:lip}, and Lemma \ref{lem:glue:foliation}, proving Theorem \ref{thm:generalized:stretchmap}.

\begin{proof}[{Proof of Theorem \ref{thm:generalized:stretchmap}}]
  The proof is similar to the proof of \cite[Corollary 4.2]{Thurston1998}, except that we use Lemma \ref{lem:bounded:lip} to deform the hyperbolic metrics on each component of $S\backslash\lambda$.  Let $Y\in\T(S)$ be a closed hyperbolic surface and $\lambda$ a geodesic lamination on $Y$. Let $Y^1,\cdots, Y^m$ be the connected components of $Y\backslash\lambda$, which are crowned hyperbolic surfaces. Suppose that $Y^i$ has $k_i$ crowns. For each $1\leq i\leq m$, let us choose a closed Riemann surface $X^i$ with marked points $D^i=\{p^i_1,\cdots,p^i_{k_i}\}$ and a harmonic diffeomorphism $f^i:X^i\to Y^i$ with the Hopf differential denoted by $\Phi^i$.  Note that $\Phi^i$ has a pole of order $r^i_j\geq2$  at each marked point $p^i_j$ of $X^i$. For each $1\leq i\leq m$ and each $t>0$,  let $\T_i(t)$ be the family of spaces of crowned hyperbolic surfaces homeomorphic to $X^i$  with $k_i$ crowns, each having $r^i_j-2$ boundary cusps, with metric residue being twice of the absolute value of the real part of the residues of the principal parts of $t\Phi^i$.  Applying Theorem \ref{thm:parametrization}  to the pair $(X^i,t\Phi^i)$ and taking the principal part of $t\Phi^i$ at each puncture of $X^i$ as the prescribed principal part of Theorem \ref{thm:parametrization},  we see that  for each $t>0$,   there exist a  unique crowned hyperbolic surface $Y^i_t$ in $\T_i(t)$ and a unique harmonic diffeomorphism $f^i_t:X^i\to Y^i_t$ with $Y^i_1=Y^i$ such that $\Hopf(f^i_t)=t\Hopf(f^i)=t\Phi^i$. This proves item (a) of Theorem \ref{thm:generalized:stretchmap}.   
  
  Next, we shall glue $Y^1_t,\cdots,Y^m_t$ to obtain a family of closed hyperbolic surfaces $\{Y_t\}$ as follows. For $t>0$, let $V^i_t$ be the pushforward on $Y_t^i$ of the vertical foliation of the Hopf differential of $\Phi^i$ via $f_t^i:X^i\to Y^i_t$. For $t=1$, by Lemma \ref{lem:glue:foliation}, there exists a unique measured foliation $V_1$ on $Y$ such that $V_1|_{Y^i}=V^i_1$. For $t>0$, by Lemma \ref{lem:bounded:lip},  the map $u^i_{1,t}:=f^i_t\circ (f^i_1)^{-1}:Y^i_1\to Y^i_t$ extends to an affine map of linear factor $\sqrt{t}$ on the boundary. Notice that there are many hyperbolic surfaces $Y_t'\in \T(S)$ whose restriction to the component corresponding $Y^i_1$ is $Y^i_t$. Any two such metrics differ by a \enquote{shearing} along the geodesic lamination $\lambda$. Constructing the desired hyperbolic metric is equivalent to specifying the \enquote{shearing} along $\lambda$. This is done through the measured foliation $V_1$ as follows. From the construction of $V^i_t$, we see that $(u^i_{1,t})^* V^i_t=\sqrt{t}V^i_1$. This allows us to glue $Y^i_t$ in such a way that the extended foliation on the resulting surface arising from $V^i_t,\cdots,V^m_t$ by Lemma \ref{lem:glue:foliation} is exactly $\sqrt{t}V_1$.  Let $Y_t$ be the resulting hyperbolic surface and $V_t$ be the extended measured foliation arising from $V^1_t,\cdots, V^m_t$ by Lemma \ref{lem:glue:foliation}. (Later in this proof we will compute the monodromy of $Y_t$ in terms of $V^i_t$ and $V_t=\sqrt{t}V_1$. In particular, this demonstrates the uniqueness of $Y_t$.) 

 Having defined the surface $Y_t$, we next wish to extend the maps $u^i_{1,t}$ to a map $u_{1,t}:Y_1 \to Y_t$. To do this note that, for any sufficiently small neighborhood $U$ of $\lambda$ on $Y_t$, two leaf segments of the restriction to $U$ of  $V^1_t,\cdots,  V^m_t$ are contained in a single leaf of the restriction to $U$ of $V_t$ if and only if their preimages under $u^1_{1,t},\cdots, u^m_{1,t}$ are contained in a single leaf of the restriction to $\lambda\cup (u^1_{1,t})^{-1}(U\backslash \lambda)\cup\cdots\cup (u^m_{1,t})^{-1}(U\backslash\lambda)$ of $V_1$. Therefore, the maps $u^1_{1,t},\cdots, u^m_{1,t}$ may be glued together to get a homeomorphism $u_{1,t}$ from $Y_1$ to $Y_t$ that sends leaves of $V_1$ to $V_t$; note that then $(u_{1,t})^*V_t=\sqrt{t}V_1$.
  
  By Lemma \ref{lem:bounded:lip}, we see that for any $0<s<t$,  the composition map $u^i_{s,t}:=f^i_t\circ (f^i_s)^{-1}: Y^i_s\to Y^i_t$ is a Lipschitz map having pointwise Lipschitz constant strictly less than $\sqrt{t/s}$ and extends to an affine map of factor $\sqrt{t/s}$ from  $\partial{Y}^i_s$ to  $\partial{Y}^i_t$. Note that $u^i_{s,t}=u^i_{1,t}\circ (u^i_{1,s})^{-1}$. Combining with the discussion in the preceding paragraph, we see that we may extend  $u^1_{s,t},\cdots, u^m_{s,t}$ to obtain a homeomorphism $u_{s,t}:=u_{1,t}\circ (u_{1,s})^{-1}:Y_s\to Y_t$ which is an affine map of factor $\sqrt{t/s}$ on $\lambda$ and which has pointwise Lipschitz constant strictly less than  $\sqrt{t/s}$ in $Y\backslash \lambda$. This proves item (b) of Theorem \ref{thm:generalized:stretchmap}.

\bigskip
We now compute the monodromy of $Y_t$ and prove the real analytic dependence on $t$ (assuming Lemma \ref{lem:analyticity1}), following the idea of Thurston \cite[Section 4]{Thurston1998} (see also \cite[Theorem A]{Bonahon1992} or \cite[Section 3.5 in Chapter 2]{PapadopoulosTheret2007}). We start with an overview of the argument. Notice that the length functions of simple closed geodesics are real analytic over $\T(S)$ and that each hyperbolic surface is uniquely determined by the lengths of finitely many simple closed geodesics (\cite{Schmutz1993,Okumura1996}), hence by the monodromy matrices of finitely many simple closed curves.  In view of this, to show that $Y_t$ varies real analytically on $t$, it suffices to show that the monodromy matrix of any simple closed curve varies real analytically in $t$. Let $\alpha$ be an arbitrary simple closed curve. In \cite{Thurston1998}, Thurston realizes $\alpha$ as a  concatenation of  geodesic subarcs of (isolated) leaves of $\lambda$  and  subleaves of the horocycle foliation transverse to $\lambda$  (assume that $\lambda$ is maximal), and then expresses the monodromy matrix of $\alpha$ as an infinite product of Mobius transformations $$\prod{N_i}.$$
Here $N_i$ is a Mobius transformation which translates along  the imaginary axis or the horocycle leaves. Here for each $i$, the translation $N_i$ varies real analytically in $t$.  Thurston proves the analyticity by showing that a specific sequence of finite product approximations of $\prod N_i$ converges uniformly (see \cite[Section 3.5 in Chapter 2]{PapadopoulosTheret2007}) for a detailed calculation).

In our setting, we replace the horocycle foliation by the foliation $\eta$ which is the extension to $Y$ of the pushforward of the vertical foliation of the Hopf differential of $f:X\to Y\backslash\lambda$ (Lemma \ref{lem:glue:foliation}), and express the monodromy of $\alpha$ as infinite product of Mobius transformations 
$$ \prod M_i' $$
where $M_i'$ is a Mobius transformation that translates along the imaginary axis or the leaves of $\eta$. Here Lemma \ref{lem:analyticity1} guarantees the analyticity of $M_i'$ in $t>0$.   Applying the same computation as in   \cite[Section 4]{Thurston1998} (see also \cite[Section 3.5 in Chapter 2]{PapadopoulosTheret2007}) allows us to conclude that  the sequence of finite product approximations of $\prod M_i'$ converges uniformly.  In particular, the limit $\prod M_i'$ also varies real analytically in $t$. 

\bigskip
We now elaborate on the details. Let $\eta$ be the extension to $Y$ of the pushforward of the vertical foliation of the Hopf differential of $f:X\to Y\backslash \lambda$ (Lemma \ref{lem:glue:foliation}).  Consider the $\epsilon$ neighbourhood of $\lambda$. There exists a realization $\alpha^*$ of the simple closed curve $\alpha$ which is the concatenation of finitely many subarcs of leaves bounding complementary regions of $\lambda$ and finitely many subarcs  of leaves of $\eta$.  Fix an orientation of $\alpha^*$. We may relabel those subarcs in a cyclic order as $\{\gamma_i\}$  according to the orientation of $\alpha^*$.    Next, we associate a matrix $M_i$ to each of   $\gamma_i$ ,  and   write the monodromy of $\alpha$ as a finite product  $\prod_{1\leq i\leq m} M_i$. In particular, let $V$ be a continuous unit vector field along $\alpha^*$ such that at every endpoint $p$ of $\gamma_i$, the vector $V(p)$ is tangent to the leaf of $\lambda$ that contains $p$.
 The orientation of $\alpha^*$ induces an orientation of each arc $\gamma_i$.  For $\gamma_i$, we lift the vector at the starting point to the upward vertical unit vector at $i\in\H^2$. Let $M_i$  be the Mobius transformation which moves the lifted vector at $i$ to the lifted vector at the other endpoint of the lift of $\gamma_i$.  The monodromy of $\alpha^*$ is the product $\prod_{1\leq i\leq m} M_i$.

To show that $ \prod_{1\leq i\leq m} M_i$ is real analytic in $t$, it suffices to show that each $M_i$ is real analytic in $t$.  If the subarc $\gamma_i$ corresponding to $M_i$ is contained in a (boundary) leaf of $\lambda$, then, by the construction/definition of $M_i$ in the proceeding paragraph, the Mobius transformation $M_i$ moves $i\in\mathbb H^2$ along the imaginary axis by the distance $t\ell_{Y}(\gamma_i)$, hence is real analytic.   We now consider the case where $\gamma_i$ is a subarc contained in some leaf of $\eta$.  As we mentioned earlier, we will express this $M_i$ as an infinite product of Mobius transformations corresponding to the subarcs of $\gamma_i\backslash\lambda$, show that each term in the product is real analytic in $t$, and prove that a specific sequence of finite approximations converges uniformly. 

\bigskip
Before that, we need some estimates on lengths of $\gamma_i\setminus \lambda$ leading to estimates on the Mobius transformations $M_i$ (here and in the next three paragraphs).
   Notice that for any fixed $\epsilon'>0$, the number of components of $\gamma_i\backslash\lambda$ with length at least $\epsilon'$ is finite. In particular, there are only finitely many components of $\gamma_i\backslash\lambda$ not contained in the cusp region of $Y_t\backslash\lambda$.  Let $\wedge$ be an arbitrary cusp region of $Y\backslash\lambda$ with $\lambda_0\subset \lambda$ being one of the bounding ideal geodesics of $\wedge$.  Let $q_0$ be a point on $\lambda_0$. Let $\{q_{i_n}\}$ be the endpoints of $\gamma_i\backslash\lambda$ on $\lambda_0$, labelled in an increasing way according to their distances to the base point $q_0$ along $\lambda_0$. We label the components of $\gamma_i\backslash\lambda$ in $\wedge$ by $\hat\gamma_{i_n}$ such that the endpoint of $\hat\gamma_{i_n}$ on $\lambda_0$ is $q_{i_n}$. 

Let $[a,b]\subset (0,\infty)$ be an arbitrary closed interval of finite length on which we will prove analyticity in $t$.  Let $r_0$ be the infimum of  the length of the shortest simple closed curves on $Y_t$ as $t$ varies in $[a,b]$. Since $Y_t$ varies continuously in $t$ on this compact set $[a,b]$, we see that $r_0>0$. After dividing $\gamma_i$ into smaller arcs if necessary, we may assume that the length of $\gamma_i$ is much smaller than $\frac{1}{2}r_0$. This implies there exists some constant $r>\frac{1}{2}r_0$ such that the hyperbolic distance on $Y_a$ between any pair of consecutive points of $\{q_{i_n}\}$ is at least $r$ (or else we would find a closed curve comprising a path in $\lambda_0$ and $\gamma_i$ of length less than $r_0$). To simplify the notation, since our attention is focused on a single arc $\gamma_i$, we can set $q_{i_n} = q_n$.

 Since the induced map $Y_a\to Y_t$ expands by a factor of $\sqrt{t/a}$ along $\lambda$,  the distance $d_t(q_0,q_n)$ along $\lambda_0$ on $Y_t$ satisfies: $$d_t(q_0,q_n)=\sqrt{t/a}\cdot d_a(q_0,q_n)\geq \sqrt{t/a}nr.$$   Combining this with  \eqref{eq:pullback:metric:2} and Lemma \ref{lem:Minsky:decay},   we see  that the preimage of $\hat\gamma_{i_n}$ on $X$ under the map $f_t:X\to Y_t\backslash\lambda$ has $|\Hopf(f_t)|$-distance at least $C \sqrt{t/a}nr$ from the zeros  of $\Hopf(f_t)$ for some constant $C$ which is independent of $t\in[a,b]$: here the point is this. The arc $\hat\gamma_{i_n}$ has pre-image which is a bi-infinite vertical leaf in two vertical half-planes (possibly with some bi-infinite strips of finite height in between) which are being mapped to the cusp -- and since $\hat\gamma_{i_n}$ is located at least at hyperbolic distance $C \sqrt{t/a}nr$ into the cusp, it has preimage at least at a proportional distance in the $|\Hopf(f_t)|$-distance from some base point, such as the zeroes of the differential. Again by \eqref{eq:pullback:metric:2} and Lemma \ref{lem:Minsky:decay}, we see that the hyperbolic length of $\hat\gamma_{i_n}$ on $Y_t$ satisfies:
\begin{equation}\label{eq:hat:gamma}
	\ell(\hat\gamma_{i_n})\leq C e^{-C\sqrt{t/a}nr} \end{equation} 
for some positive constant $C$ independent of $t\in[a,b]$.

Fix a base point $x_0\in\H$ in the hyperbolic plane.
Let $\mathbf{d}$ be  the complete, left invariant metric on $PSL(2,\R)$ defined as (cf. \cite{PapadopoulosTheret2007}):
\begin{equation*}
	\mathbf{d}(A,B):=\sup_{x\in\H}d_{\H}(Ax,Bx)e^{-d_\H(x,x_0)},\qquad\forall A,~B\in PSL(2,\R).
\end{equation*}
 For any $A\in PSL(2,\R)$, define  $\|A\|:=\mathbf{d}(Id,A)$. Let $\widehat M_{i_n}$ be the Mobius transformation corresponding to $\hat\gamma_{i_n}$ defined similarly as $M_i$.  By \eqref{eq:hat:gamma}, we see that there exists another constant $C'$, independent of $t\in[a,b]$, such that $\|\widehat M_{i_n}\|\leq C' e^{-C\sqrt{t/a}nr}$. In particular,  there exist two   positive constants $ C_1$ and $C_2$  such that for any $t\in [a,b]$ and any $n\geq1$, such that
\begin{equation}\label{eq:small:M}
	\|\widehat M_{i_n}\|\leq C_1 e^{-C_2 n}.
\end{equation}
From this we infer that the distance to $Id\in PSL(2,\R)$ from any finite product of matrices in $\{\widehat M_{i_n}\}$ is bounded from above by $C_1\sum_{n\geq1}e^{-C_2n}=C_1(1-e^{-C_2})^{-1}.$

\bigskip
We now return to the construction of $M_i$, the monodromy corresponding to $\gamma_i$ as an infinite product (recall that $\{\hat\gamma_{i_n}\}$ are the components of $\gamma_i\backslash\lambda$, here relaxing the notation to allow the components $\{\hat\gamma_{i_n}\}$ to be subsets of any one of the finitely many cusp regions. [because there are only finitely many cusp regions, the estimates of the previous paragraphs continue to hold]).
For each $m\geq1$, consider the (finite) subset $I_m$ of $\{\hat\gamma_{i_n}\}$ which comprises those $\hat\gamma_j$ whose length is at least $1/m$.  let $\mathbb M_m$ be  the product of those $\widehat M_{i_n}$ corresponding to $\hat\gamma_{i_n}\in I_m$, according to the order inherited from the orientation of $\gamma_i$ (recall that $\hat\gamma_{i_n}\subset \gamma_i$ ). In particular, $\mathbb M_{m+1}$ is obtained from $\mathbb M_m$ by inserting finitely many  $\widehat M_{i_n}$ corresponding to $\hat\gamma_{i_n}\in I_{m+1}$.  Applying a standard calculation using \eqref{eq:small:M} as in \cite[the proof of Proposition 3.15 in Page 161-162]{PapadopoulosTheret2007}, we see that $\mathbb M_m$ converges uniformly in $t\in[a,b]$, where the limit is exactly   the Mobius transformation $M_i$ corresponding to $\gamma_i$ (recall that $\{\hat\gamma_{i_n}\}$ are the components of $\gamma_i\backslash\lambda$). For each $\hat{\gamma}_{i_n}$, there is exactly one $j_n\in\{1,2,\cdots,k\}$ such that the segment $\hat{\gamma}_{i_n}$ is contained in the component  $Y^{j_n}_t$ of $Y_t\backslash\lambda$. 

By Lemma \ref{lem:analyticity1} both $Y^{j_n}_t$ and the realization $f^{j_n}_t\circ (f^{j_n}_1)^{-1}(\hat{\gamma}_{i_n})$ of $\hat{\gamma}_{i_n}$ on $Y^{j_n}_t$ are both real analytic in $t$. Therefore, the Mobius transformation $\widehat{M}_{i_n}$ corresponding to $\hat{\gamma}_{i_n}$ is real analytic in $t$. Combined with the uniform convergence, this proves that $M_i$, as the limit of $\mathbb M_m$, is analytic in $t\in(a,b)$. Consequently, the monodromy of $\alpha$, as a product of finitely many $M_i$, varies analytically in $t$.  Since $Y_t$ is determined by finitely many such $\alpha$, this implies that $Y_t$ is analytic in $t\in(a,b)$. The arbitrariness of $[a,b]\subset (0,\infty)$ then implies that $Y_t$ varies analytically in $t\in(0,\infty)$.
\end{proof}

\begin{definition}
    Given a closed hyperbolic surface $Y\in\T(S)$, a geodesic lamination $\lambda$ on $Y$, and a harmonic diffeomorphism $f:X\to Y\backslash\lambda$ from some (possibly disconnected) punctured Riemann surface. The path 
    \begin{equation*}
        \PSL_{Y,\lambda,f}:(0,\infty)\to\T(S)
    \end{equation*}
    constructed in Theorem \ref{thm:generalized:stretchmap} is called a {\it piecewise harmonic stretch line} defined by the triple $(Y,\lambda,f)$.  The restriction to $[1,\infty)$ is called a {\it piecewise harmonic stretch ray} defined by $(Y,\lambda,f)$, denoted by $\PSR_{Y,\lambda,f}:[1,\infty)\to\T(S)$.
\end{definition}
From the construction, we see that $\PSL_{Y,\lambda,f}(1)=Y$. In particular, the piecewise harmonic stretch ray $\PSR_{Y,\lambda,f}$ starts at $Y$.  Each piecewise harmonic stretch line/ray admits a canonical orientation coming from the orientation of the positive real ray $\{t>0\}$. In that orientation, a piecewise harmonic stretch line is a (reparametrized) geodesic in the Thurston metric (Theorem \ref{thm:generalized:stretchmap}(b)).

We end this subsection with the following identification.
\begin{lemma}\label{lem:id:PSR:SR} 
    Let $Y\in\T(S)$  and $\lambda$ be a maximal geodesic lamination $\lambda$. Then for any harmonic diffeomorphism $f:X\to Y\backslash\lambda$   from a punctured Riemann surface $X$, the piecewise harmonic stretch line $\PSL_{Y,\lambda,f}$ and the Thurston stretch line $\SL_{Y,\lambda}$ (see Section \ref{subsec:stretch:maps} for definition) coincide.
\end{lemma}
\begin{proof} 
   We begin with recalling the construction of the Thurston stretch line $\SL_{Y,\lambda}$. Let $Y^1,\cdots,Y^{4g-4}$ be the set of ideal triangles complementary to $\lambda$ on $Y$. Let $F$ be the horocyclic partial foliation of $\lambda$ on $Y$. On each component $Y^i$, the horocycle partial foliation $\hat F^i:=F|_{Y^i}$ foliates $Y^i$ except for a central region bounded by three pairwise tangent horocyclic arcs meeting each boundary edge of $Y^i$ perpendicularly.  Let $F^i$ be a measured foliation obtained from $\hat{F}^i$ by pushing the leaves of $\hat{F}^i$ toward the central region of $Y^i$ while keeping endpoints on $\partial Y$ invariant and keeping the resulting leaves orthogonal to $\partial Y^i$. Note that any two thus obtained $F^i$'s differ by an isotopy of $Y^i$ fixing $\partial Y^i$ pointwise. The transverse measure on $F^i$ is defined as follows. For any arc $\tau$ on $Y^i$ that is transverse to $F^i$, there exists exactly one boundary geodesic of $Y^i$, to which we may  project $\tau$ along leaves of $F^i$.  The transverse measure on $\tau$ is then defined to be the hyperbolic length of the projection of $\tau$ to that geodesic.

   Let $X^1,\cdots,X^{4g-4}$ be the corresponding components of $X$ and $f^i:X^i\to Y^i$ be the restriction of $f$ to $X^i$. Note that for each $i$, the Riemann surface $X^i$ is conformally the complex plane. That $Y^i$ is an ideal triangle implies that the Hopf differential $\Phi^i$ of $f^i$ is a linear polynomial, that is, $\Phi^i=a^i(z-b^i)dz^2$ for some pair of complex numbers $(a^i,b^i)$ (item (i) of Lemma \ref{lem:Hopf:differential}). The critical graph of the horizontal measured foliation of $\Phi^i$ is a tripod, whose complement consists of three half-planes. In particular, the horizontal foliation is symmetric about each of the three half-infinite critical leaves.  On the other hand, each ideal triangle is symmetric about any of the three half-infinite geodesic rays $\{h^{i1},h^{i2},h^{i3}\}$, where each starts perpendicularly from a boundary edge and converges to the opposite ideal point. Applying Theorem \ref{thm:parametrization} to the triple $(X^i,Y^i,\Phi^i)$ with respect the principal part given by $\Phi^i$, we see that the harmonic diffeomorphism $f^i:X^i\to Y^i$ is symmetric about each of $\{h^{i1},h^{i2},h^{i3}\}$. Accordingly, the pushforward $V^i$ of the vertical foliation of $\Phi^i$ is also symmetric about each of $\{h^{i1},h^{i2},h^{i3}\}$. By Lemma \ref{lem:intersection:extend}, the transverse measure on $V^i$  coincides with half of the projection to boundary geodesics of $Y^i$. More precisely, for any arc $\tau$ on $Y^i$ that is transverse to $V^i$, there exists exactly one boundary geodesic of $Y^i$, to which we may project $\tau$ along leaves of $V^i$. By Lemma \ref{lem:intersection:extend}, the transverse measure on $\tau$ equals half of the hyperbolic length of the projection of $\tau$ to that geodesic. 
   
   Combining the discussion in the preceding two paragraphs, we see that  $2V^i$ and $F^i$ 
   differ by an isotopy of $Y^i$ fixing $\partial Y^i$ pointwise. This also holds for $2tV^i$ and $tF_i$ for any $t>0$. Therefore, the double of the extended foliation $tV$ obtained from $tV^1,\cdots,tV^{4g-4}$ via Lemma \ref{lem:glue:foliation} also differs from $tF$ by an isotopy of $Y$  that fixes $\cup_i\partial Y^i$ pointwise. This implies that for each $t>0$, the surface $\PSL_{X,\lambda,f}(t)$ and $\SL_{Y,\lambda}(t)$ coincide since both of them are obtained from gluing ideal triangles with the same transverse measured foliation $tF=2t V$.
\end{proof}

%======================
%======================
\subsection{Limits of piecewise harmonic stretch lines}
Recall that every piecewise harmonic stretch line  is directed.  We next prove a proposition that will be useful when we define canonical harmonic stretch lines from given base points to points on the Thurston boundary of \tec space in Section~\ref{sec:exponential}. 
\begin{proposition} \label{prop:convergencestretchlines:PML}
 Let $Y\in\T(S)$ be a closed hyperbolic surface and $\lambda$ a geodesic lamination on $Y$. Let $f:X\to Y\setminus\lambda$ be a harmonic diffeomorphism from some punctured surface $X$. Let $\beta$ be the pushforward to $Y$ of the vertical foliation of $\Hopf(f)$. Then the piecewise harmonic stretch line determined by $(Y,\lambda,X,f)$ converges to $[\beta]\in\PML(S)$ in the Thurston compactification in the forward direction. 
\end{proposition}

\begin{proof}
 Let $\PSL$ be the piecewise harmonic stretch line
 determined by $(Y,\lambda,X,f)$. Let $\widehat{Y}_t$ be the (possibly disconnected) crowned hyperbolic surface determined by the pair $(X,t\Hopf(f))$ by Theorem \ref{thm:parametrization}.  Let $Y_t\in\PSL$ be the hyperbolic surface such that $Y_t\setminus\lambda=\widehat{Y}_t$.
 Let $\alpha$ be an arbitrary simple closed curve on $Y$. Consider the pushforward to $Y$ of the horizontal and vertical foliations of $\Hopf(f)$. By Lemma \ref{lem:glue:foliation}, the pushforward of the vertical foliation extends to a unique measured foliation on the whole surface $Y_t$. For any $\epsilon>0$, let $\alpha^*$  be a representative of $\alpha$ satisfying the following:
  \begin{itemize}
      \item $\alpha^*$ consists of segments of the pushforward of horizontal and vertical foliations of $\Hopf(f)$ alternatively, 
    \item $\alpha^*$  avoids the zeros  of $\Hopf(f)$, 
    \item $\alpha^*$ almost realizes the  (minimal) intersection number with $\beta$, the vertical foliation  of $\Hopf(f)$, namely 
     \begin{equation}\label{eq:almost:minimal} \left|\int_{\alpha^*}|\mathrm{Re}\sqrt{\Hopf(f)}|-i(\alpha,\beta)\right|<\epsilon. 
     \end{equation}
  \end{itemize}
  To see such a $\alpha^*$ exists, let $\alpha'$ be a representative of $\alpha$ which realizes the minimal intersection number with $\beta$, that is, $\int_{\alpha^*}|\mathrm{Re}\sqrt{\Hopf(f)}|=i(\alpha,\beta)$.  Next, we deform $\alpha'$ to get a new representative $\alpha''$ that consists of segments of the pushforward of horizontal and vertical foliations of $\Hopf(f)$ alternatively. Then $\alpha''$ also realizes the minimal intersection number with $\beta$. Finally, at each zero of $\Hopf(f)$ on $X$, we pick a polygon $Q$ consisting of  horizontal and vertical segments alternatively each of which is of distance $\epsilon'\ll\epsilon$ to the enclosed zero under the flat metric induced by the Hopf differential $\Hopf(f)$, and then push the subsegments of $\alpha''$ inside the image $f(Q)$ of this polygon to the boundary of $f(Q)$;  here the choice of $\epsilon'$ depends on $\alpha$. This gives the desired representative $\alpha^*$.
  
  Let $d>0$ be the distance of $f^{-1}(\alpha^*)$ to the zeros of $\Hopf(f)$. 
  Let $K_v$ be the number of vertical segments of $\alpha^*$. Then by (\ref{eq:pullback:metric:2}) and Lemma \ref{lem:Minsky:decay}, for  $t>0$ sufficiently large, the total length $ \mathrm{Length}_{Y_t}(^{v}\alpha^*)$ of vertical segments of $\alpha^*$ with respect to the hyperbolic metric $Y_t\in \PSL$ satisfies:
    \begin{eqnarray}
    \nonumber
  \mathrm{Length}_{Y_t}(^v \alpha^*)
  &\leq& 
  K_v \int_{td}^{\infty} 2 \sinh^{-1}(\chi(Y)/s^2) \exp(-s)ds\to0, \text{ as } t\to+\infty. 
  \end{eqnarray} 
  and the total length $ \mathrm{Length}_{Y_t}(^{h}\alpha^*)$ of horizontal segments of $\alpha^*$ with respect to the hyperbolic metric $Y_t\in \PSL$ satisfies:
  \begin{eqnarray*}
    \nonumber
  \mathrm{Length}_{Y_t}(^h \alpha^*)
  &\leq& 
 2 \int_{\alpha^*}|\mathrm{Re}\sqrt{t\Hopf(f)}|(1+{ \exp(-td)\sinh^{-1}(\chi(Y)/(td)^2) })\\
  &\to& 2\sqrt{t} \int_{\alpha^*}|\mathrm{Re}\sqrt{\Hopf(f)}| , \text{ as }t\to+\infty.
  \end{eqnarray*} 
  Therefore, the length $\mathrm{Length}_{Y_t}(\alpha^*)$ of $\alpha^*$ on $Y_t$ satisfies: \begin{eqnarray*}
  \limsup_{t\to+\infty}  \frac{\mathrm{Length}_{Y_t}(\alpha^*)}{\sqrt{t}}
  &=&\limsup_{t\to+\infty}  \frac{\mathrm{Length}_{Y_t}(^v\alpha^*)+\mathrm{Length}_{Y_t}(^h\alpha^*)}{\sqrt{t}}\\
  &\leq& 2\int_{\alpha^*}|\mathrm{Re}
  \sqrt{\Hopf(f)}|.
  \end{eqnarray*}
  Hence the length of the geodesic representative of $\alpha$ on $Y_t$ satisfies:
  \begin{equation}\label{eq:length:alpha*:upper}
   \limsup_{t\to+\infty}  \frac{\ell_{Y_t}(\alpha)}{\sqrt{t}}\leq 
     \limsup_{t\to+\infty}  \frac{\mathrm{Length}_{Y_t}(\alpha^*)}{\sqrt{t}}
  \leq 2\int_{\alpha^*}|\mathrm{Re}
  \sqrt{\Hopf(f)}|.
  \end{equation}
  On the other hand, by (\ref{eq:pullback:metric:2}),  \begin{eqnarray*}
  {\ell_{Y_t}(\alpha)}
  \geq2\int_{\alpha^*}|\mathrm{Re}\sqrt{t\Hopf(f)}|
  = 2\sqrt{t}\int_{\alpha^*}|\mathrm{Re}\sqrt{\Hopf(f)}|.
  \end{eqnarray*}
  Combined with (\ref{eq:length:alpha*:upper}), this implies that 
  \begin{eqnarray*}
  \lim_{t\to+\infty}  \frac{\ell_{Y_t}(\alpha)}{\sqrt{t}}
  =2\int_{\alpha^*}|\mathrm{Re}
  \sqrt{\Hopf(f)}|. 
  \end{eqnarray*}
  It then follows from (\ref{eq:almost:minimal})
 and the arbitrariness of $\epsilon$ that 
 \begin{eqnarray*}
  \lim_{t\to+\infty}  \frac{\ell_{Y_t}(\alpha)}{\sqrt{t}}
  =2i(\alpha,\beta).
  \end{eqnarray*}
  This proves that $Y_t$ converges to $[\beta]\in\PML(S)$ as $t\to+\infty$.
 \end{proof}

%=============================

  \subsection{Generalized maximum principle} 
  To prove Lemma \ref{lem:bounded:lip}, we need the following generalized maximum principle from \cite{ChengYau1975}.
  
\begin{theorem}[\cite{ChengYau1975}Theorem 3 and Theorem 8]\label{thm:generalized:maximum}
  Let $M$ be a complete noncompact manifold with Ricci curvature bounded from below by $-K$ for some constant $K\geq0$. Let $u$ be a $C^2$ function on $M$.
  \begin{enumerate}[(I)]
    \item If $u$ is bounded from above, then there  exists a sequence of points $p_k\in M$ such that
       \begin{eqnarray*}
         \lim_{k\to\infty}u(p_k)=\sup_{p\in M} u(p), &
         \lim_{k\to\infty}|\nabla u(p_k)|=0,&
         \limsup_{k\to\infty} \Delta u(p_k)\leq0.
       \end{eqnarray*}
    \item If  $u$ satisfies the differential inequality $\Delta u\geq f(u)$, where $f$ is a function on $\R$ with the property that there exists a continuous non-decreasing function $g$ positive on some interval $[a,\infty)$ such that:
        \begin{equation*}
          \liminf_{t\to\infty}\frac{f(t)}{g(t)}>0
        \end{equation*}
        and
        \begin{equation*}
          \int_{b}^{\infty}\left(\int_{a}^{t}
          g(\tau)d\tau\right)^{-1/2}dt<\infty ~~\text{ for some } b\geq a,
        \end{equation*}
        then $u$ is bounded from above. Furthermore, if $f$ is lower semi-continuous, $f(\sup u)\leq0$.
  \end{enumerate}
  \end{theorem}

\subsection{Geometry of harmonic map rays from degenerated surfaces} \label{sec:stretchline:construction:degenerate}
 To prove Lemma \ref{lem:bounded:lip}, we need a generalization of \cite[Proposition 4.3]{Wolf1989}.  Before that we need to recall some notation.   Let $(M,\sigma|dz|^2)$ be the hyperbolic plane or the complex plane and $(\H^2,\rho|dw|^2)$ the hyperbolic plane. Let $f_t:M\to (\H^2,\rho|dw|^2) $ be a harmonic diffeomorphism with Hopf differential $t\Phi$. Recall the conformal energy density $\mathcal{H}_t=\frac{\rho(f_t)}{\sigma}|\frac{\partial f_t}{\partial z}|^2$, the anti-conformal energy density $\mathcal{L}_t=\frac{\rho(f_t)}{\sigma}|\frac{\partial f_t}{\partial \bar{z}}|^2$, the Laplacian operator 
 $\Delta_\sigma=
    \frac{4}{\sigma}\frac{\partial^2}
    {\partial z\partial \bar{z}}$,  the Beltrami differential $\nu_t=\frac{(\partial f_t/\partial{\bar{z}})d\bar{z}}{(\partial f_t /\partial z)dz}$, and identities relating them \eqref{eq:density:H:L}--\eqref{eq:Phi}.  
   
\begin{proposition}\label{prop:monotonicity}
  Let $(M,\sigma |dz|^2)$ be the hyperbolic plane $\H^2$ or the complex plane $\C$. Let $\Phi dz^2$ be a holomorphic quadratic differential on $M$; in the case where $M$ is the complex plane, we assume $\Phi$ is a polynomial.
 Let $f_t: (M,\sigma |dz|^2)\to \H^2$ be a family of harmonic diffeomorphisms onto the corresponding images with Hopf differentials $\Phi_t=t\Phi$, where $t>0$.
   Then the following statements hold, where  we use primes to indicate differentiation with respect to $t$. 
    \begin{enumerate}[(i)]
      \item For any $p\in M$, ${\mh'_t}(p)\geq0$.
      \item For any $p\in M$ ${\ml'_t}(p)\geq0$.
      \item For any $p\in M$, $|\nu_t|'(p)\geq0$, and   ${|\nu_t|'}(p)>0$ if  $\Phi(p)\neq 0$.
      \item The function $\frac{|\nu_t|'}{|\nu_t|}(p)$ extends to a bounded positive analytic function on $M$. In particular, for any compact subset $N$ of $M$, there exists a constant $\epsilon>0$, such that $\frac{|\nu_t|'}{|\nu_t|}(p)>\epsilon$ for any $p\in N$.
    \end{enumerate}
    \end{proposition}

   \begin{proof}
      The idea is to use the maximum principle and the generalized maximum principle.

      Notice that for any fixed $t\geq 0$,
      \begin{equation}\label{eq:H:inf}
        \inf_{p\in M} \mh_t(p) >0.
      \end{equation}
      In fact,  if $M$ is the hyperbolic plane, then it follows from \cite[Theorem 12 and Proposition 10]{Wan1992a} that $\mh_t (p)\geq 1$ for all $t\geq0$ and all $p\in M$. If $M$ is the complex plane, since $\Phi$ is a nonconstant polynomial, then by (\ref{eq:quadratic:norm}), since $\mh_t(p)\geq\ml_t(p)$ (notice $f_t$ is orientation-preserving), we have that
      $\mh_t(p)\geq t|\Phi(p)|\to\infty$, as $p\to\infty$. In particular, there exists $R>0$, such that $\mh_t(p)>1$ for all $|p|>R$. On the other hand, on $\{p\in\C:|p|\leq R\}$, we have that $\mh_t(p)$ is positive and continuous.  Therefore, in both cases, (\ref{eq:H:inf}) holds.

   \vskip 5pt
     First, we  show that ${\mh'_t}(p)\geq0$ for all $p\in M$. To this end, it is equivalent to show that for any fixed $p\in M$,  $\mh_t(p)$ is an increasing function of $t>0$.  For any $0<t<s$, let $\mw_t=\frac{1}{2}\log \mh_t$, then by (\ref{eq:Bochner:H}), we have the Bochner-type equations: 
     \begin{eqnarray}
       \Delta_\sigma \mw_t&=& e^{2\mw_t} -\frac{t^2|\Phi|^2}{\sigma^2}  e^{-2\mw_t}+K(\sigma) \label{eq:Bochner:W},\\
         \Delta_\sigma \mw_s&=& e^{2\mw_s} -\frac{s^2|\Phi|^2}{\sigma^2}  e^{-2\mw_s} + K(\sigma)\notag.
     \end{eqnarray}
     By subtracting the two equations, we get
     \begin{equation*}
      \Delta_\sigma( \mw_t-\mw_s)= (e^{2\mw_t}-e^{2\mw_s}) -\frac{t^2|\Phi|^2}{\sigma^2}  e^{-2\mw_t}
      +\frac{s^2|\Phi|^2}{\sigma^2}  e^{-2\mw_s}.
     \end{equation*}
     Let $\tilde{\sigma} = \tilde{\sigma}(s,z):=\sigma e^{2\mw_s}|dz|^2$ and $\eta:=\mw_t-\mw_s$. Then
      \begin{eqnarray}
      \Delta_{\tilde{\sigma}}\eta&=& e^{2\eta}-1 -\frac{t^2}{s^2} |\nu_s|^2  e^{-2\eta}
      +|\nu_s|^2 \qquad (\text{ by }(\ref{eq:Phi}))
      \label{eq:laplacian:eta}\\
      &\geq&  e^{2\eta}- e^{-2\eta}-1\qquad \qquad(\text{ since } |\nu_s|\leq 1).\nonumber
     \end{eqnarray}
      Recall that $\tilde{\sigma} = \tilde{\sigma}(s)$ is a complete metric with curvature bounded from below: indeed, for fixed $s$,
     \begin{eqnarray}
       \nonumber K(\tilde{\sigma})&=&-2\Delta_{\tilde{\sigma}}\log \tilde{\sigma}\\
       \nonumber  &=&
      - \frac{2\sigma}{\tilde{\sigma}}
     \Delta_\sigma (\log\sigma+2\mw_s)\\
    \nonumber  &=&- e^{-2\mw_s}
     (2\Delta_\sigma \log\sigma+4\Delta_\sigma\mw_s)\\
   \nonumber   &{=}&-e^{-2\mw_s}\left(-K(\sigma)+ 4 e^{2\mw_s} -4\frac{s^2|\Phi|^2}{\sigma^2}  e^{-2\mw_s}+4K(\sigma)\right)\quad({\text{by }(\ref{eq:quadratic:norm},\ref{eq:Bochner:H})})\\
   \nonumber   &{=}&-3K(\sigma)e^{-2\mw_s}-4+4|\nu_s|^2 \qquad \qquad\qquad \qquad({\text{by }(\ref{eq:Phi})})\\
          &\geq&-4 .\qquad\qquad(\text{since } K(\sigma)\in\{0,-1\})
   % \nonumber   &\geq&-5e^{-2\mw_s}-4+4|\nu_s|^2
      %\qquad\qquad(\text{since } K(\sigma)\in\{0,-1\})\\
   %\nonumber   &=&-\frac{5}{\mh_s}-4+4|\nu_s|^2\\
  % \nonumber   &\geq&-\frac{5}{\mh_s}-4\\
  %   &\geq&-\frac{5}{\inf_{p\in M}\mh_s(p)} -4>-\infty. \qquad\qquad (\text{by (\ref{eq:H:inf})})
\label{eq:curv:tilde:sigma}
     \end{eqnarray}
     It then follows from Theorem \ref{thm:generalized:maximum} by  taking $f(t)=e^{2t}-e^{-2t}-1$ and $g(t)=e^{2t}$,  that $\bar{\eta}:=\sup_{p\in M} \eta(p)<+\infty$.  Combining this boundedness with item (I) of Theorem \ref{thm:generalized:maximum} and the fact that the curvature of $\tilde\sigma$ is bounded from below \eqref{eq:curv:tilde:sigma}, we see that there exists a sequence $p_k\in M$, such that
     \begin{equation*}
       \lim_{k\to\infty}\eta(p_k)=\bar{\eta}=\sup_{p\in M}\eta(p),\qquad
      \limsup_{k\to\infty}
       \Delta_{\tilde{\sigma}}\eta(p_k)\leq0.
     \end{equation*}
     By taking a subsequence if necessary, we can assume that
     \begin{equation*}
       \lim_{k\to\infty}|\nu_s(p_k)|^2=a\in[0,1].
     \end{equation*}
     It then follows from (\ref{eq:laplacian:eta}) that
     \begin{eqnarray*}
       0&\geq & \limsup_{k\to\infty}\Delta_{\tilde{\sigma}}\eta (p_k) \\
       &=& \limsup_{k\to\infty}\left( e^{2\eta(p_k)}-1 -\frac{t^2}{s^2} |\nu_s(p_k)|^2  e^{-2\eta(p_k)}
      +|\nu_s(p_k)|^2 \right)\\
       &=& e^{2\bar{\eta}}-\frac{t^2}{s^2}a
       e^{-2\bar{\eta}}-1+a\\
       &\geq &  e^{2\bar{\eta}}-a
       e^{-2\bar{\eta}}-1+a \qquad (\text{ since } 0<t<s)
     \end{eqnarray*}
     which implies that $\bar{\eta}\leq0$.
     Therefore,
     \begin{equation*}
       \mw_t(p)-\mw_s(p)=\eta(p)\leq \bar{\eta}\leq0
     \end{equation*}
     Hence, $\mh_t(p)\leq\mh_s(p)$ for all $p\in M$. Equivalently, ${\mh'_t}(p)\geq 0$ for all $p\in M$.

\vskip 5pt
      Next, we show that ${\ml'_t}(p)\geq 0$ for all $p\in M$.  Notice that $t^2|\Phi|^2/\sigma^2 =\mh_t\ml_t$ and $s^2|\Phi|^2/\sigma^2 =\mh_s\ml_s$
    Then \begin{equation}
       \label{eq:HL:st}    \frac{\ml_t}{\ml_s}(p)=\frac{t^2}{s^2} \frac{\mh_s}{\mh_t}(p),\qquad \Phi(p)\neq0
         \end{equation}
         extends to  a well-defined analytic function on the whole surface $M$, still denoted by $\frac{\ml_t}{\ml_s}(p)$ for simplicity.
    Let $\delta(p)=\frac{1}{2}\log\frac{\ml_t}{\ml_s}(p)$. The  equation above implies that $\delta(p)=\log{\frac{t}{s}}+\frac{1}{2}\log\frac{\mh_s}{\mh_t}(p)$. It then follows that
      \begin{eqnarray}
      \Delta_{{\sigma}}\delta&=&\frac{1}{2} \Delta_\sigma \log\mh_s -\frac{1}{2} \Delta_\sigma \log\mh_t\\
      &=& \mh_s-\mh_t-\frac{s^2|\Phi|^2}{\sigma^2}\mh_s^{-1}+
      \frac{t^2|\Phi|^2}{\sigma^2}\mh_t^{-1}.
     \end{eqnarray}
     Let  $\hat{\sigma}:=\sigma {\mh_t}|dz|^2$, then
     \begin{eqnarray*}
       \Delta_{\hat{\sigma}}\delta 
       &=& \frac{\mh_s}{\mh_t}-1-\frac{s^2|\Phi|^2}{\sigma^2\mh_s\mh_t}
       +\frac{t^2|\Phi|^2}{\sigma^2\mh_t^2}\\
         &=& \frac{\mh_s}{\mh_t}-1-\frac{t^2|\Phi|^2}{\sigma^2\mh_t^2}
       \frac{s^2}{t^2}\frac{\mh_t}{\mh_s}
       +\frac{t^2|\Phi|^2}{\sigma^2\mh_t^2}\\
       &{=}& \frac{\mh_s}{\mh_t}-1-|\nu_t|^2
       \frac{s^2}{t^2}\frac{\mh_t}{\mh_s}+|\nu_t|^2 \qquad({\text{by } (\ref{eq:Phi})})\\
       &{=}& \frac{s^2}{t^2}e^{2\delta}-1-|\nu_t|^2 e^{-2\delta}+|\nu_t|^2\qquad({\text{by }(\ref{eq:HL:st})})\\
       &\geq&\frac{s^2}{t^2}e^{2\delta}-1- e^{-2\delta}.
     \end{eqnarray*}
    By applying an argument similar to the one for the case of $\eta$, we see that $\sup_{p\in M}\delta<0$. In particular, $\ml_t<\ml_s$ for all $0\leq t<s$ and all $p\in M$ with $\Phi(p)\neq0$. This implies that $\ml_t'(p)\geq0$ for all $t\geq0$ and all $p\in M$ with $\Phi(p)\neq0$. It then follows from continuity of $\ml_t'$ that  $\ml_t'(p)\geq0$ for all $t\geq0$ and all $p\in M$.

     It remains to show item (iii) and item (iv).    Suppose that $\Phi(p)\neq0$.  Taking derivatives respect to $t$ for (\ref{eq:quadratic:norm}) and (\ref{eq:Phi}), we see that 
     \begin{eqnarray}
         \label{eq:sum:HL} \frac{{\mh}'_t}{\mh_t}(p)+ \frac{{\ml}'_t}{\ml_t}(p)&=&\frac{2}{t},\\
         \label{eq:sum:HV} \frac{{\mh}'_t}{\mh_t}(p)+ \frac{{|\nu|_t}'}{|\nu_t|}(p)&=&\frac{1}{t}.
        \end{eqnarray}
        Therefore, for any fixed $t>0$, both $ \frac{{\ml}'_t}{\ml_t}$ and $\frac{{|\nu|}'_t}{|\nu|_t}$ extend to be bounded non-negative analytic functions on the whole of $M$. We next improve this bound on $ \frac{{\ml}'_t}{\ml_t}$.
        
         Applying  Theorem \ref{thm:generalized:maximum} to $ \frac{{\mh}'_t}{\mh_t}$, which is a bounded non-negative analytic function by (\ref{eq:sum:HL}), we see that there exists a sequence of points $p_k\in M$ with $p_k\to\infty$, such that
        \begin{equation}\label{eq:log:H:maximum}
          \lim_{k\to\infty}\frac{{\mh}'_t}{\mh_t}(p_k)
          =\sup_{p\in M}\frac{{\mh}'_t}{\mh_t}(p),
          ~~\limsup_{k\to\infty}
          \Delta_\sigma\frac{{\mh}'_t}{\mh_t}(p_k)\leq0.
        \end{equation}
       Therefore,
        \begin{eqnarray*}
        &&2\sup_{p\in M}\frac{{\mh}'_t}{\mh_t}(p)-\frac{2}{t}\\
        &=&
         \lim_{k\to\infty} \left( 2\frac{{\mh}'_t}{\mh'_t}(p_k)-
          \frac{2}{t} \right)\\ &=&
          \lim_{k\to\infty} \left( \frac{{\mh}_t'}{\mh_t}(p_k)-
         \frac{{\ml}_t'}{\ml_t}(p_k) \right) \\
          &\leq&
          \lim_{k\to\infty} \left( \frac{{\mh}_t'}{\mh_t}(p_k)-
         \frac{{\ml}_t'}{\mh_t}(p_k) \right)  \qquad (\text{since } \mh_t>\ml_t\geq0 \text{ and } \ml_t'\geq0)\\
          &=&\lim_{k\to\infty}
         \frac{\frac{1}{2}\Delta_\sigma  \frac{{\mh}_t'}{\mh_t}(p_k)}{\mh_t(p_k)}
          \qquad ( \text{by } (\ref{eq:Bochner:H}))
          \\
          &\leq &0. \qquad\qquad\qquad(\text{by } (\ref{eq:H:inf}),(\ref{eq:log:H:maximum}))
        \end{eqnarray*}
        Consequently,
        \begin{equation}\label{eq:H:up}
          \sup_{p\in M} \frac{{\mh}'_t}{\mh_t}(p)\leq
          \frac{1}{t}.
        \end{equation}
        Next, we claim that 
        \begin{equation}
          \label{eq:H:strictub} \frac{{\mh}'_t}{\mh_t}(p)<
          \frac{1}{t}
        \end{equation}
        holds for any $p\in M$. Suppose to the contrary that there exists some $p_0\in M$ such that 
          $ \frac{{\mh}'_t}{\mh_t}(p_0)=\frac{1}{t}$. Combined with (\ref{eq:H:up}), this implies that 
           \begin{equation*}
      \Delta_\sigma \frac{{\mh}'_t}{\mh_t}(p_0)
      \leq0.
        \end{equation*}
    It then follows from (\ref{eq:Bochner:H})  that
        \begin{equation*}
     \mh_t'(p_0)-\ml'_t(p_0)= \Delta_\sigma \frac{{\mh}'_t}{\mh_t}(p_0)
      \leq0,
        \end{equation*}
        which implies that $\mh_t'(p_0)\leq\ml_t'(p_0)$.  Combined with (\ref{eq:sum:HL}) and the fact that $\mh_t(p_0)>\ml_t(p_0)$, this yields that 
        \begin{equation*}
          \frac{{\mh}'_t}{\mh_t}(p_0)<
          \frac{1}{2}\left(\frac{{\mh}'_t}{\mh_t}(p_0)
          +\frac{{\ml}'_t}{\ml_t}(p_0)\right)=\frac{1}{t}.
        \end{equation*}
        This contradicts the assumption that $\frac{\mh_t'}{\mh_t}(p_0)=\frac{1}{t}$.  Hence 
        \begin{equation*}
     \frac{{\mh}'_t}{\mh_t}(p)<
          \frac{1}{t}, \forall p\in M \text{ and }\forall t>0.
        \end{equation*}
        It then follows from (\ref{eq:sum:HV}) that
        \begin{equation*}
          \frac{|\nu_t|'}{|\nu_t|}(p)=\frac{1}{t}- \frac{{\mh}'_t}{\mh_t}(p)>0, \qquad \forall p\in M \text{ and } \forall t>0.
        \end{equation*} 
        Combined with the continuity of $        \frac{|\nu_t|'}{|\nu_t|}$, this  proves (iii) and (iv).   
   \end{proof}

\begin{proof}[Proof of Lemma \ref{lem:bounded:lip}]
 Let $p\in {X}\backslash D$ be an arbitrary point with $\Phi(p)\neq0$. 
 By the third item of Proposition \ref{prop:monotonicity}, we see that $|\nu_t|_p$ is a strictly  increasing function of $t\in(0,\infty)$. Therefore,  the function  $\G(p,t):=\log (1/|\nu_t(p)|)$  is a strictly decreasing function of $t\in[0,\infty)$ at $p$ with $\Phi(p)\neq 0$. Recall that with respect to the canonical coordinate charts  $z=x+\i y$ of $\Phi$ near $p$,  the pullback metric of ${Y}_t$ by $f_t$  is:
  \begin{equation}\label{eq:pullback}
    2t(\cosh \G(p,t)+1)\d x^2+2t(\cosh \G(p,t)- 1)\d y^2.
  \end{equation}
   Then for any $s>t>0$, the Lipschitz constant of $f_s\circ (f_t)^{-1}$ at $f_t(p)$ satisfies 
  $$ \mathrm{Lip} (f_s\circ (f_t)^{-1})|_{f_t(p)}\leq \max\left\{\sqrt{\frac{2s(\cosh \G(p,s)+1)}{2t(\cosh \G(p,t)+1)}}, \sqrt{\frac{2s(\cosh \G(p,s)-1)}{2t(\cosh \G(p,t)-1)}}\right\}
  < \sqrt{\frac{s}{t}}.  $$
  Combined with  Lemma \ref{lem:Minsky:decay}, this implies  that as $p$ approaches $D$, the set of punctures of $  {X\backslash D}$, we have that
   $ \mathrm{Lip} (f_s\circ (f_t)^{-1})|_{f_t(p)}$  converges uniformly to $\sqrt{\frac{s}{t}}$.
   
  We now  consider Lipschitz constants of $f_s\circ f_t^{-1}$ at zeros of $\Phi$. Let $p_0$ be an arbitrary zero of $\Phi$. Then by (\ref{eq:Beltrami}), the value $\nu_t(p_0)=0$ for all $t\geq 0$.  Hence 
\begin{equation}
    \label{eq:G:diverge}
    \G(p_0,t)=+\infty \text{ and }\lim_{p\to p_0}\G(p,t)=+\infty.  
\end{equation}
  Let $U$ be a small neighbourhood of $p_0$. By item (iv) of Proposition \ref{prop:monotonicity}, for any $0<t<s$, there exists $\epsilon>0$, such that $\frac{|\nu_r|'}{|\nu_r|}(p)>\epsilon$ for any $p\in U$ and any $r\in[t,s]$. It then follows by integrating this inequality that
  \begin{equation}
      \label{eq:nu:st}
     \frac{|\nu_s|}{|\nu_t|}(p)=\exp(\log|\nu_s|(p)-\log|\nu_t|(p))>\exp(\epsilon(s-t))
  \end{equation} 
  holds for any $p\in U$.
  Therefore, the Lipschitz constant of $f_s\circ (f_t)^{-1}$ at $f_t(p_0)$ satisfies
  \begin{eqnarray*}
  &&\mathrm{Lip} (f_s\circ (f_t)^{-1})|_{f_t(p_0)}\\
   &=&\lim_{p\to p_0}\mathrm{Lip} (f_s\circ (f_t)^{-1})|_{f_t(p)} \\
   &\leq& \lim_{p\to p_0} \max\left\{\sqrt{\frac{2s(\cosh \G(p,s)+1)}{2t(\cosh \G(p,t)+1)}}, \sqrt{\frac{2s(\cosh \G(p,s)-1)}{2t(\cosh \G(p,t)-1)}}\right\}\\
   &=&\lim_{p\to p_0}\sqrt{\frac{s}{t} \frac{\frac{1}{|\nu_s|(p)}}{\frac{1}{|\nu_t|(p)}}} \qquad (\text{ by } (\ref{eq:G:diverge}))
   \\
   &=&\lim_{p\to p_0}\sqrt{\frac{s}{t} \frac{|\nu_t|(p)}{|\nu_s|(p)}}\\
   &\leq & \sqrt{\frac{s}{t} \exp(-\epsilon(s-t))} \qquad (\text{ by } (\ref{eq:nu:st})) \\
   &<& \sqrt{\frac{s}{t}}.
  \end{eqnarray*}
  This completes the proof.
\end{proof}

%===================
%===================
 \section{Convergence of harmonic maps rays}\label{sec:convergence:stretchrays}

 In this section, we complete the proof of Theorem \ref{thm:HR:SR}.
 To prove the theorem,  it will suffice by Theorem \ref{thm:generalized:stretchmap}  to prove that the family of harmonic maps $f_t:X_t\to Y$ converges in the sense of Section \ref{sec:limit:varingdomains}, as $X_t$ degenerates along the harmonic map dual ray $\HDR_{Y,\lambda}$. To this end, by Lemma \ref{lem:precompactness},  it suffices to show that every convergent sequence of the family $\{f_t:X_t\to Y\}_{t\geq 0}$ shares the same limit.  Let $f_{t_n}:X_{t_n}\to Y$ be an arbitrary convergent sequence with limit harmonic diffeomorphism $f_\infty:X_\infty\to Y\backslash\lambda_\infty$ for some geodesic lamination $\lambda_\infty$. 
 
\subsection{Identifying limits of harmonic map rays with piecewise harmonic stretch lines}\label{sec:limit:HR:PSL} By Theorem \ref{thm:HR:limit:geodesic2}, any (divergent) sequence of harmonic map rays through a fixed point subconverges to a Thurston geodesic ray/line. The goal of this subsection is to identify that limit ray with some piecewise harmonic stretch ray/line. 
Let $Y\in \T(S)$ be a fixed hyperbolic surface. Let $X_n\in\T(S)$ be a divergent sequence of Riemann surfaces such that the harmonic map rays $\HR_{X_n,Y}:[0,\infty)\to \T(S)$ converge to some Thurston geodesic $\mathbf{GR}:(0,\infty)\to \T(S)$. Note that our notation for the harmonic map rays sets $\HR_{X_n,Y}(1)\equiv Y$. Let $f_{n,t}:X_n\to \HR_{X_n,Y}(t)$ be the harmonic diffeomorphism in the given homotopy class.
\begin{proposition}\label{prop:limit:HR:PSL}
  With the notation and assumption as above, there exists a subsequence of $X_n$, still denoted by $X_n$ for simplicity, with the following properties.
    \begin{enumerate}[(i)]
        \item The sequence $f_{n,1}: X_n\to Y$ converges to a harmonic diffeomorphism $f_{\infty,1}:X_\infty\to Y\backslash\lambda$ for some punctured Riemann surface $X_\infty$ and some chain recurrent geodesic lamination $\lambda\subset Y$.
        \item The limit geodesic $\mathbf{GR}:(0,\infty)\to \T(S)$ is the piecewise harmonic stretch line $\PSL_{Y,\lambda,f_{\infty,1}}:(0,\infty)\to\T(S)$ defined by the triple $(Y,\lambda,f_{\infty,1})$.
        \item For any $s>0$, the sequence $f_{n,s}: X_n\to \HR_{X_n,Y}(s)$ converges to a harmonic diffeomorphism $f_{\infty,s}:X_\infty\to \mathbf{GR}(t)\backslash\lambda$ such that the Hopf differentials satisfy $\Hopf(f_{\infty,s})=s\Hopf (f_{\infty,1})$.
        \item For any $0<s<t$, the sequence $f_{n,t}\circ (f_{n,s})^{-1}: \HR_{X_n,Y}(s)\to \HR_{X_n,Y}(t)$  converges to a  $\sqrt{t/s}$-Lipschitz homeomorphism $\mathbf{GR}(s)\to \mathbf{GR}(t)$ that  extends $f_{\infty,t}\circ (f_{\infty,s})^{-1}:\mathbf{GR}(s)\backslash\lambda\to \mathbf{GR}(t)\backslash\lambda$.
    \end{enumerate}
\end{proposition}
\begin{proof}
   We begin with an outline of the argument. By Lemma \ref{lem:precompactness}, the sequence $f_{n,1}: X_n\to Y$ (sub)converges to a harmonic diffeomorphism $f_{\infty,1}:X_\infty\to Y\backslash\lambda$ for some punctured Riemann surface $X_\infty$ and some chain recurrent geodesic lamination $\lambda\subset Y$.
   We inherit a complication from this construction: the surface $X_\infty$ is a punctured surface and the limiting map $f_{\infty,1}$ of the sequence $f_{n,1}:X_{n}\to Y$ of maps, takes values in the open surface $Y\setminus \lambda$; in particular, the lamination $\lambda$ is not in the image, so there is no {\it a priori} extension of the image $Y\setminus \lambda$ to the full surface $Y$.

 We observe that the identical issue arises elsewhere along the Thurston geodesic $\mathbf{GR}$:  for $Y_{n,s}=\HR_{X_n,Y}(s)\in\HR_{X_{n},Y}$ along the ray $\HR_{X_{n},Y}$, we have 
$\lim\limits_{n\to\infty} Y_{n,s}=\mathbf{GR}_s$ (again because of the overall subconvergence of rays to a geodesic), but the limiting map
$f_{n,s}:X_{n}\to Y_{n,s}$ has limit 
$f_{\infty,s}:X_\infty\to \mathbf{GR}_s\backslash\lambda$ (with limit Hopf differential $s \Hopf(f_\infty)$) -- here we note again that the image is a subdomain $\mathbf{GR}_s\backslash\lambda$ of $\mathbf{GR}_s$.

Now,  the Thurston geodesic $\mathbf{GR}_s$ has (noncanonical) corresponding  homeomorphisms  $\mathbf{L}_s:\mathbf{GR}_1 \to \mathbf{GR}_s$ which stretch exactly $\sqrt{s}$ along $\lambda$. (Here, we shall choose $\mathbf{L}_s$ as the limit of $f_{n,s}\circ(f_{n,1})^{-1}:Y\to Y_{n,s}$.)
Separately, because  $\mathbf{GR}_s$ is a compact surface, there is also an obvious extension $\mathbf{Ext}_s: \mathbf{GR}_s\setminus \lambda \to \mathbf{GR}_s$ (by continuity) of $\mathbf{GR}_s\setminus \lambda$ across $\lambda$ to $\mathbf{GR}_s$.  (When we carry out this plan, it will be convenient to realize this extension by taking limits along the images of the leaves of the vertical foliation of $\Hopf(f_{\infty,s})$. Equivalently, we will describe the extended foliation on $\mathbf{GR}_s$ of the pushforward of the vertical foliation of $\Hopf(f_{\infty,s})$ on $\mathbf{GR}_s\backslash\lambda$.)

Analogously, in terms of the harmonic maps, on the one hand, there is also a well-defined deformation of $f_{\infty,s}$ from $f_{\infty}$ obtained by scaling the Hopf differentials, i.e. $\Hopf(f_{\infty,s})=s\Hopf(f_{\infty})$.
In contrast to the case in the previous paragraph of the Thurston geodesic, there is however, not an immediately apparent extension of either $f_{\infty,s}$ or $f_{\infty}$ across any punctures.

 Thus, to prove (ii) in the proposition, evidently what we must show is that 
 \begin{equation}\label{eq:commutation1}
     \mathbf{L}_s(Y_1)= \mathbf{L}_s \circ \mathbf{Ext}_1 \circ f_{\infty,1} (X_\infty) = \mathbf{Ext}_s\circ\mathbf{L}_s \circ f_{\infty, 1}(X_\infty);
 \end{equation}
  here we note that in the right hand side, we have that $f_{\infty, 1}(X_\infty)$ has image $\mathbf{GR}_1 \setminus \lambda$, so that $\mathbf{L}_s$ is acting on the complement of the lamination $\lambda$, while on the left-hand side, we will have $\mathbf{L}_s$ acting on the extension of $f_{\infty, 1}(X_\infty)$, and so acts as a homeomorphism of a compact surface. This commutation relation \eqref{eq:commutation1} allows us to identify the limit geodesic ray $\mathbf{GR}$ with the piecewise harmonic stretch ray defined by $(Y,\lambda,X_\infty, f_\infty)$ and complete the proof.

As a final remark, we note that we may expand the last displayed commutation relation between $\mathbf{L}$ and $\mathbf{Ext}$ and see it as an interchange of limits.  In particular, the map $\mathbf{L}_s$ may be realized as the limit of the maps $f_{n,s} \circ [f_{n,1}]^{-1}: Y \to Y_{n,s}$.  On the complement of the lamination $\lambda$, that composition limits on $f_{\infty, s} \circ [f_{\infty,1}]^{-1}: Y \setminus \lambda \to Y_s \setminus \lambda$.  This then reduces what we must prove to a statement about the extension map $\mathbf{Ext}$. 

This relates to the correspondence, across the puncture(s), between the ends of the vertical trajectories of $\Hopf(f_{\infty,s})$ (itself the limit of $\Hopf(f_{n,s})$): we must show that the extension map $\mathbf{Ext}$, obtained as above by continuity across the lamination $\lambda$ is also obtained as the (well-defined) limit of the image under $f_{\infty,s}$ of corresponding ends of the vertical trajectories we just described.  [Put another way as an interchange of \enquote{temporal} and \enquote{spatial} limits of vertical leaves of Hopf differentials, because the vertical trajectories of $\Hopf(f_{n, s})$ realize the extension across $\lambda$ for the image $f_{n,s}(X_{n})$, we need to show that that limiting (in $n$) extension map across $\lambda$ of the image under $f_{n,s}(X_{n})$ of a vertical leaf of $\Hopf(f_{n,s})$ agrees with the (spatial) limit along the (disconnected) image of the end of that limiting corresponding  vertical leaf of $\Hopf(f_{n, s})$.]

\bigskip
We now fill in the details of this sketch.

 By Lemma \ref{lem:precompactness}, any sequence of maps $f_{n}:X_{n}\to Y$ contains a convergent subsequence, still denoted by the same notation for simplicity, which converges to a harmonic diffeomorphism $f_{\infty,1}:X_\infty\to Y\backslash\lambda$ for some punctured Riemann surface $X_\infty$ and some chain recurrent geodesic lamination $\lambda\subset Y$. This proves the conclusion item (i).

  We next address item (iii). By Theorem \ref{thm:generalized:stretchmap}\footnote{In this argument, we will use 
  statements (a) and (b) of Theorem~\ref{thm:generalized:stretchmap} that we proved in Section~\ref{sec:stretchline:construction}; we will not use the claim on analyticity in Theorem~\ref{thm:generalized:stretchmap} 
  relying on Lemma~\ref{lem:analyticity1}, whose proof we deferred until just after the present proof.}, the harmonic diffeomorphism  $f_{\infty,1}:X_\infty\to Y\setminus\lambda$ defines a piecewise harmonic stretch line $\PSL_{Y,\lambda,f_{\infty,1}}:(0,\infty)\to\T(S)$ such that the Hopf differential $\Hopf(f_{\infty,s}')$ of the induced harmonic diffeomorphism $f'_{\infty,s}:X_\infty\to \PSL_{Y,\lambda,f_{\infty,1}}(s)\backslash\lambda$ satisfies
\begin{equation}\label{eq:Hopf:infty:s:prime1}
    \Hopf(f'_{\infty,s})=s \Hopf(f_{\infty,1}).
\end{equation}
Note that $\PSL_{Y,\lambda,f_{\infty,1}}(1)=Y$  and $f_{\infty,1}'=f_{\infty,1}$.  

Now, by Lemma \ref{lem:glue:foliation}, we note the crucial point that the pushforward of the vertical foliation of $\Hopf(f_{\infty,s}')$ on $\PSL_{Y,\lambda,f_{\infty,1}}(s)\backslash\lambda$ extends to a unique measured foliation $V_s$ on $\PSL_{Y,\lambda,f_{\infty,1}}(s)$.  From the construction of  the piecewise harmonic stretch lines (cf. the second paragraph in the proof of Theorem \ref{thm:generalized:stretchmap} in Section \ref{subsec:construction:PSL}), we see that
\begin{equation}\label{eq:vert:s11}
	V_s=\sqrt{s}V_1.
\end{equation}

 By assumption, the family of harmonic map rays $\HR_{X_n,Y}:[0,\infty)\to \T(S)$  locally uniformly converges to the Thurston geodesic ray $\mathbf{GR}:(0,\infty)\to\T(S)$. To simplify notations, we set $\mathbf{GR}_s=\mathbf{GR}(s)$. Note that $\mathbf{GR}_1=Y$. For $s>0$, let $Y_{n,s}=\HR_{X_n,Y}(s)\in\HR_{X_{n},Y}$ be the hyperbolic surface such that the Hopf differential of the harmonic diffeomorphism $f_{n,s}:X_{n}\to Y_{n,s}$ is $s\Hopf (f_{n,1})$. Then $\mathbf{GR}_s=\lim\limits_{n\to\infty} Y_{n,s}$. Furthermore,  combining the fact that $f_{n,1}:X_{n}\to Y$ converges to $f_{\infty,1}:X_\infty\to Y\backslash\lambda$, the fact that $\Hopf(f_{n,s})=s \Hopf(f_{n,1}) $ for all $s>0$, and the assumption that $\HR_{X_{n},Y}:[0,\infty)\to \T(S)$ converges to $\mathbf{GR}:(0,\infty)\to\T(S)$ locally uniformly, we see that $f_{n,s}:X_{n}\to Y_{n,s}$ \emph{subconverges} to a (not necessarily surjective) harmonic map $f_{\infty,s}:X_\infty \to \mathbf {GR}_s$ with $\Hopf(f_{\infty,s})=s\Hopf(f_{\infty,1})$. Note that by the second bullet property in Lemma \ref{lem:precompactness}, the geodesic lamination $\lambda$ satisfies: 
   \begin{equation}
       \lambda =\lim\limits _{R\to\infty}\lim\limits_{n\to\infty}\overline{f_{n,1} ( X_{n}\backslash\mathscr{P}_{R}(\Hopf(f_{n,1})))},
   \end{equation} where  $\mathscr P_R(\Hopf(f_{n,1}))$ is the Minsky polygonal region of the Hopf differential $\Hopf(f_{n,1})$. Combining Theorem \ref{thm:Minsky:traintrack} and the fact (Lemma \ref{lem:1:lip}) that $f_{n,s}\circ (f_{n,1})^{-1}$ is a $\sqrt{s}$ Lipschitz homeomorphism homotopic to the identity, we see that 
   \begin{equation}
       \lambda =\lim\limits _{R\to\infty}\lim\limits_{n\to\infty}\overline{f_{n,s} ( X_{n}\backslash\mathscr{P}_{R}(\Hopf(f_{n,1})))}.
   \end{equation}
   
   Applying a similar argument as in the proof of Lemma \ref{lem:precompactness} shows that the image  $f_{\infty,s}(X_\infty)$ is exactly $\mathbf{GR}_s\backslash\lambda$ and that $f_{\infty,s}$ is injective. In summary, the harmonic map $f_{\infty,s}:X_\infty\to\mathbf{GR}_s\backslash\lambda$ is a diffeomorphism with Hopf differential
   \begin{equation}
       \label{eq:Hopf:infty:prime:f}\Hopf(f_{\infty,s})=s\Hopf(f_{\infty,1}).
   \end{equation} 
   Combined with \eqref{eq:Hopf:infty:s:prime1} and Theorem \ref{thm:parametrization} (with respect to the principal part of $s\Phi$ at the punctures of $X_\infty$), this implies that $\mathbf{GR}_s\backslash\lambda=\PSL_{Y,\lambda,f_{\infty,1}}(s)\backslash \lambda$ and that $f_{\infty,s}=f_{\infty,s}'$. In particular, for any fixed $s$, the map $f_{\infty,s}$ obtained as a limit of a subsequence in $n$ of $f_{n,s}$, is independent of the choice of a convergent subsequence (in $n$) of $f_{n,s}:X_n\to Y_{n,s}$. Hence,  the sequence $f_{n,s}:X_n\to Y_{n,s}$ converges to $f_{\infty,s}:X_\infty \to \mathbf{GR}_s\backslash\lambda$.  This proves the conclusion item (iii) of the proposition. 
   
    To finish the proof of the identification $\mathbf{GR}_s=\PSL_{Y,\lambda,f_{\infty,1}}(s)$, it remains to show that $\mathbf{GR}_s\backslash\lambda$ extends to $\mathbf{GR}_{s}$ in the same way as $\PSL_{Y,\lambda,f_{\infty,1}}(s)\backslash\lambda$ extends to $\PSL_{Y,\lambda,f_{\infty,1}}(s)$. In other words, we need to show that the extended foliation $\mathbf{V}_s$ on $\mathbf{GR}_s$ of the pushforward of the vertical foliation of $\Hopf(f_{\infty,s})$ on $\mathbf{GR}_s\backslash\lambda$ is identical to extended foliation ${V}_s$ on $\PSL_{Y,\lambda,f_{\infty,1}}(s)$ of the pushforward of the vertical foliation of $\Hopf(f_{\infty,s}')$ on $\PSL_{Y,\lambda,f_{\infty,1}}\backslash\lambda$. Recall that by \eqref{eq:vert:s11}, we have $V_s=\sqrt{s}V_1$.  Note further that $\mathbf{GR}_1=Y=\PSL_{Y,\lambda,f_{\infty,1}}(1)$, so $\mathbf{V}_1=V_1$. In view of these two identities, to prove $\mathbf{V}_s=V_s$, it suffices to show that $\mathbf{V}_s=\sqrt{s}\mathbf{V}_1$.   

For $s\geq 1$, from the construction of harmonic map rays, we see that the composition map $f_{n,s}\circ (f_{n,1})^{-1}:Y\to Y_{n,s}$, which is a $\sqrt{s}$ Lipschitz homeomorphism,  sends leaves on $Y$ of the pushforward of the vertical foliation of the Hopf differential $\Hopf(f_{n,1})$ by $f_{n,1}$ to corresponding leaves on $Y_{n,s}$ of the pushforward of the vertical foliation of the Hopf differential $\Hopf(f_{n,s})=s \Hopf(f_{n,1})$ by $f_{n,s}$. Recall that $\lim_{n\to\infty}Y_{n,s}=\mathbf{GR}_s$. On the one hand, being $\sqrt{s}$ Lipschitz uniformly implies that $f_{n,s}\circ (f_{n,1})^{-1}:Y\to Y_{n,s}$ {\it subconverges} to a $\sqrt{s}$ Lipschitz map from $Y$ to $\mathbf{GR}_s$.  On the other hand, by the conclusion item (iii) proved above, for any $s\geq1$, we have $f_{n,s}$ converges to $f_{\infty,s}$ as $n\to\infty$. Hence the restriction $f_{n,s}\circ (f_{n,1})^{-1}|_{Y\backslash\lambda} $ converges to $f_{\infty,s}\circ (f_{\infty,1})^{-1}:Y\backslash\lambda\to \mathbf{GR}_{s}\backslash\lambda$. Furthermore, by Theorem \ref{thm:HR:limit:geodesic2}, any limit of $f_{n,s}\circ (f_{n,1})^{-1}:Y\to Y_{n,s}$ stretches  $\lambda$ by the factor $\sqrt{s}$.  Therefore, the limit map of any convergent subsequence of the sequence $f_{n,s}\circ (f_{n,1})^{-1}:Y\to Y_{n,s}$ is the $\sqrt{s}$-Lipschitz homeomorphism $\mathbf{L}_s: Y\to \mathbf{GR}_s$ that  extends the limit composition map $f_{\infty,s}\circ (f_{\infty,1})^{-1}:Y\backslash\lambda\to \mathbf{GR}_{s}\backslash\lambda$, proving that $f_{n,s}\circ (f_{n,1})^{-1}$ converges to $\mathbf{L}_s$ and establishing the conclusion item (iv). 
Now, we also conclude from these considerations that
the limit homeomorphism $\mathbf{L}_s: Y\to \mathbf{GR}_s$ sends leaves on $Y$ 
   of the extended foliation $\mathbf V_1$ (of the pushforward of the vertical foliation of the Hopf differential $\Hopf(f_{\infty,1})$ by $f_{\infty,1}$) to corresponding leaves on $\mathbf{GR}_s$ of the extended foliation $\mathbf V_s$ (of the pushforward of the vertical foliation of the Hopf differential $\Hopf(f_{\infty,s})=s \Hopf(f_{\infty,1})$ by $f_{\infty,s}$). (In other words, this establishes the commutation relation \eqref{eq:commutation1}.) This implies that \begin{equation}\label{eq:vert:s1:bf1}
	\mathbf V_s=\sqrt{s}\mathbf V_1,
\end{equation}
 where the factor $\sqrt{s}$ is due to the identity $\Hopf(f_{\infty,s})=s\Hopf(f_{\infty,1})$.  This completes the proof that $\mathbf{GR}_s=\PSL_{Y,\lambda,f_{\infty,1}}(s)$ for $s\geq 1$.  The proof for $0<s<1$ is similar. Therefore, $\mathbf{GR}_s=\PSL_{Y,\lambda,f_{\infty,1}}(s)$ for $s>0$, establishing the conclusion item (ii). The proof is now complete.    
\end{proof}

  \begin{proof}[Proof of Lemma \ref{lem:analyticity1}] Let $\eta$ be  the horizontal foliation of $\Phi$. Note that the  pushforward ${f_t}_* (t\eta)$ of $t\eta$ via  ${f_t}:{X}\backslash D\to {Y}_t$ is an admissible measured foliation on ${Y}_t$. (See Definition \ref{def:admissibility:general} of \enquote{admissible measured foliations}.) Consider $t=1$.
      By Proposition \ref{prop:appendix:convergence}, there exist
      \begin{itemize}
          \item[(i)] a closed hyperbolic surface $W$ and a chain recurrent geodesic lamination $\lambda$ on $W$ such that $W\backslash\lambda$ is the union of two or four isometric copies of (the crowned surface) ${Y}$,
          \item[(ii)] a sequence of Riemann surfaces $X_n\in\T(W)$ and harmonic diffeomorphisms $h_n:X_n\to W$ homotopic to the identity such that $h_n$ converges to a harmonic diffeomorphism $h_\infty:X_\infty\to W\backslash\lambda$ (in the sense of Definition \ref{defn:limitharmonicmap}) from some punctured Riemann surface $X_\infty$, such that
           on each component $X_\infty'$ of $X_\infty$,  the horizontal measured foliation of the Hopf differential of the restriction to $h_\infty|_{X_\infty'}:X_\infty'\to  Y$ is exactly $\eta$. 
      \end{itemize} 
      Combining the second property with Theorem \ref{thm:minimalgraph:uniqueness}, we see that 
      \begin{equation}
        \label{eq:limit:f:inf}X'_\infty={X}\backslash D \text{ and }  h_\infty|_{X_\infty'}={f_1}.
      \end{equation}

      Now consider the sequence of harmonic map rays $\HR_{X_n,W}:[1,\infty)\to\T(W)$.  In particular, $\HR_{X_n,W}(1)=W$ for all $n$. By Theorem \ref{thm:HR:limit:geodesic2}, there exists a convergent subsequence of $\HR_{X_n,W}$, still denoted by $\HR_{X_n,W}$ for simplicity, that locally uniformly converges to some Thurston geodesic. By Proposition \ref{prop:limit:HR:PSL}, that limit geodesic is the piecewise harmonic stretch line $\PSL_{W,\lambda,h_{\infty}}:(0,\infty)\to\T(W)$ defined by the triple $(W,\lambda,h_{\infty})$.  Let $g_{n,t}:X_n\to \HR_{X_n,W}(t)$ be the harmonic diffeomorphism in the given homotopy class.   Then by Lemma \ref{lem:1:lip}, for each $t>s>0$, the composition map 
      \begin{equation*}
          g_{n,t}\circ (g_{n,s})^{-1}:\HR_{X_n,W}(s)\to \HR_{X_n,W}(t)
      \end{equation*}
      is a $\sqrt{t/s}$-Lipschitz homeomorphism. 

      By Lemma \ref{lem:analyticity2}, for each $n$, the harmonic map ray $\HR_{X_n,W}:[1,\infty)\to \T(W)$ is real analytic. Therefore, as the (locally uniform) limit of $\HR_{X_n,W}$, the piecewise harmonic stretch line $\PSL_{W,\lambda,h_\infty}:(0,\infty)\to\T(W)$ is also real analytic in $t>0$.  Note that the proof of Lemma \ref{lem:analyticity2} proves that the pullback metric (see \eqref{eq:pullback:t})  of the hyperbolic metric of $\HR_{X_n,W}(t)$ via the harmonic diffeomorphism $g_{n,t}:X_n\to \HR_{X_n,W}(t)$ is real analytic in $t>0$.
     According to \cite[Theorem 3.1]{EL81} (see also \cite[Proposition 3.3]{Slegers}), we see that the family of harmonic diffeomorphisms $g_{n,t}$ is also real analytic in $t>0$. In particular, the composition map $ g_{n,t}\circ (g_{n,s})^{-1}:\HR_{X_n,W}(s)\to \HR_{X_n,W}(t)$ is real analytic in $t\in(s,\infty)$. By Proposition \ref{prop:limit:HR:PSL} (iii), we see that the sequence of harmonic diffeomorphisms $g_{n,t}:X_n\to \HR_{X_n,W}(t)$ converges to a harmonic diffeomorphism $g_{\infty,t}:X_\infty\to \PSL_{W,\lambda,h_\infty}(t)\backslash\lambda$ such that the Hopf differentials  satisfy 
     \begin{equation}\label{eq:Hopf:infty:gh}
         \Hopf(g_{\infty,t})=t \Hopf (h_\infty).
     \end{equation}
      By Proposition \ref{prop:limit:HR:PSL} (iv), the sequence of $\sqrt{t/s}$-Lipschitz homeomorphisms 
     $g_{n,t}\circ (g_{n,s})^{-1}$ converges to a $\sqrt{t/s}$-Lipschitz homeomorphism $$\mathbf{L}_{s,t}:\PSL_{W,\lambda,h_\infty}(s)\to \PSL_{W,\lambda,h_\infty}(t)$$ that 
     extends $g_{\infty,t}\circ (g_{\infty,s})^{-1}:\PSL_{W,\lambda,h_\infty}(s)\backslash\lambda\to \PSL_{W,\lambda,h_\infty}(t)\backslash\lambda$. 
     Combined with the fact that $\HR_{X_n,W}:(0,\infty)\to\T(S)$ locally uniformly converges to $\PSL_{W,\lambda,h_\infty}:(0,\infty)\to \T(W)$, this implies that the limit lipschitz homeomorphism $\mathbf{L}_{s,t}:\PSL_{W,\lambda,h_\infty}(s)\to \PSL_{W,\lambda,h_\infty}(t)$ is also real analytic in $t\in(s,\infty)$. In particular, the restriction map $g_{\infty,t}\circ (g_{\infty,s})^{-1}:\PSL_{W,\lambda,f_\infty}(s)\backslash\lambda\to \PSL_{W,\lambda,f_\infty}(t)\backslash\lambda$ is also real analytic in $t\in(s,\infty)$. 
 Combining \eqref{eq:limit:f:inf}, \eqref{eq:Hopf:infty:gh}, and  Theorem \ref{thm:parametrization}, we see that $\PSL_{W,\lambda,f_\infty}(t)\backslash\lambda$ is the union of two or four isometric copies of ${Y}_t$, and that on each component $X_\infty'$ ($={X}\backslash D$ by \eqref{eq:limit:f:inf}) of $X_\infty$, we have $g_{\infty,t}|_{X_\infty'}={f_t}$. In particular, on each component of $\PSL_{W,\lambda,h_\infty}(1)\backslash\lambda$, which is also rewritten as $W\backslash\lambda$, the restriction of the limit map $\mathbf{L}_{s,t}: \PSL_{W,\lambda,h_\infty}(s)\to \PSL_{W,\lambda,h_\infty}(t)$ coincides with the composition map ${f_t}\circ ({f_s})^{-1}:{Y}_s\to {Y}_t$. Therefore, 
 the analyticity of $\mathbf{L}_{s,t}$ and the ray $\PSL_{W,\lambda,f_\infty}:(0,\infty)\to \T(W)$ implies that the surface ${Y}_t$  and the map ${f_t}\circ ({f_s})^{-1}:{Y}_s\to {Y}_t$ along with its Lipschitz homeomorphic extension from the closure of ${Y_s}$ to the closure of ${Y_t}$ (either obtained from  Lemma \ref{lem:bounded:lip} or induced from $\mathbf{L}_{s,t}$) are real analytic in $t\in(s,\infty)$. The analyticity of ${f_t}$ follows by precomposing ${f_t}\circ ({f_s})^{-1}$  by ${f_s}$. This completes the proof. 
  \end{proof}

The proof of Theorem \ref{thm:generalized:stretchmap} is now complete.

%========
\subsection{Identifying the limit horizontal foliation}
 We now turn to the proof of Theorem \ref{thm:HR:SR}.  To begin, let $X_{t_n}\in\HDR_{Y,\lambda}$ and $f_{t_n}:X_{t_n}\to Y$ be an arbitrary convergent sequence with limit  harmonic diffeomorphism $f_\infty:X_\infty\to Y\backslash\lambda_\infty$ for some geodesic lamination $\lambda_\infty$ (Lemma \ref{lem:precompactness}).

 \begin{lemma}\label{lem:topology:limit} With the notations introduced as above, we have $\lambda_\infty=\lambda$.  
 \end{lemma}
 
 \begin{proof}   
   Let $\mathscr P_{R_n}(\Phi_{t_n})$ be the Minsky's polygonal region on $X_{t_n}$ with  $t_n\to\infty$ and $R_n:=t_n^{1/4}\to\infty$ as $n\to\infty$. Then by Theorem \ref{thm:Minsky:traintrack} the complement $f_{t_n}(X_{t_n}-\mathscr P_{R_n}(\Phi_{t_n}))$ is contained in an $\epsilon_n$ neighbourhood of $\lambda$ on $Y$, where $\epsilon_n\to0$ as $n\to\infty$. Hence $\lambda_\infty\subset\lambda$ and $f_\infty(X_\infty)\supset Y-\lambda$.   For the other direction, consider the  minimal components of $\lambda$, denoted by $\lambda^i,\cdots,\lambda^k$. Note that the horizontal foliation of $\Phi_{t_n}$ is $t_n\lambda$. Let $X^1_{t_n},\cdots, X^k_{t_n}$ be the complementary components of the critical graph (.cf Definition \ref{def:critical:graph}) of the horizontal foliation of $\Phi_{t_n}$, labelled in such a way that   the restriction to $X^i_{t_n}$ of the horizontal foliation of $\Phi_{t_n}$ corresponds to the minimal component $t_n\lambda^i$ of $t_n\lambda$.   By \eqref{eq:phi:length:component}, we see that $2\|\Phi_{t_n}|_{X^i_{t_n}}\|\geq \ell_Y(t_n\lambda^i)-C$ for some constant $C$.
    Recall that $R_n=t_n^{1/4}$.  By the item (v) of Theorem \ref{thm:Minsky:polygon}, the area of $\mathscr P_{R_n}(\Phi_{t_n})$ is at most $c\sqrt{t_n}$, which is much less than $\|\Phi_{t_n}|_{X^i_{t_n}}\|$. Hence,  the subsurface ${X^i_{t_n}}$ is not contained in $\mathscr P_{R_{t_n}}(\Phi_{t_n})$. In particular, from the item (vi) of Theorem \ref{thm:Minsky:polygon}, we see that there exists a  horizontal leaf segment $l^i_{t_n}$  of $ \Phi_{t_n}|_{X^i_{t_n}}$ in the boundary of  $\mathscr P_{R_{t_n}}(\Phi_{t_n}))$ with length $|l^i_{t_n}|\geq K_3R_n=K_3 t_n^{1/4}$ for some constant $K_3$. By Theorem \ref{thm:Minsky:traintrack}, the image $f_{t_n}(l^i_{t_n})$ is nearly a geodesic segment contained in the $\epsilon_n$ neighbourhood of the minimal component $t_n\lambda^i$. That $\lambda^i$ is minimal and that $t_n\to\infty$ imply that  $f_{t_n}(l^i_{t_n})$ converges to the whole component $\lambda^i$ as $n\to\infty$.  Therefore, the limit lamination $\lambda_\infty$, which is the limit of $f_{t_n}(X_{t_n}-\mathscr P_{R_{t_n}}(\Phi_{t_n})) $ as $n\to\infty$, contains $\lambda$. Combined with the discussion in the beginning of the proof, this proves that $\lambda_\infty=\lambda$. 
 \end{proof}

 A direct consequence of the above lemma is that $X_\infty$ and $Y\backslash\lambda$ are homeomorphic. Let $Y^1,\cdots, Y^m$ be the components of $Y\backslash\lambda$. Let  $X_\infty=X^1\cup X^2\cup\cdots\cup X^m$ with $X^i$ homeomorphic to $Y^i$. Set $f^i:=f_\infty |_{X^i}$.

 \begin{lemma}\label{lem:hopf:limit}
 The critical graph (Definition \ref{def:critical:graph}) of the Hopf differential $\Hopf(f^i)$ is connected. Moreover, the complement of the critical graph of $\Hopf(f^i)$  consists of half-infinite cylinders corresponding to the closed geodesic boundary components of $Y^i$ and half-planes corresponding to the ideal geodesic boundary components  of $Y^i$.
 \end{lemma}
 \begin{proof}
 Let $\lambda^i$ be the horizontal measured foliation of $\Hopf(f^i)$.

 \textbf{Claim 1: If $\lambda$ contains a simple closed curve $\alpha$, then the corresponding (maximal) cylinder component $A_{t_n}$ of the horizontal foliation of $\Hopf(f_{t_n})$ limits on two half-infinite cylinders on $\Hopf(f_\infty)$.}
  As $n\to\infty$, the height of $A_{t_n}$ on $X_{t_n}$ goes to infinity. Combining with the Minsky's estimate (Theorem \ref{thm:Minsky:traintrack}), we see that the circumference of $A_{t_n}$ converges to half the hyperbolic length of the geodesic representative of $\alpha$ on $Y$, which is both finite and positive. Therefore the cylinder $A_{t_n}$ converges to two half-infinite cylinders on $X_\infty$ with circumferences equal to half the hyperbolic length of the geodesic representative of $\alpha$ on $Y$. This proves Claim 1.

 \textbf{Claim 2: for each $i$, the measured foliation $\lambda^i$  has no compact components} (see section \ref{subsec:meromorphic:leaves} for definition). Otherwise, suppose that $\lambda^i$ contains a compact component $\lambda^i_0$. Since the horizontal foliation of $\Hopf(f_t)$ is simply $t\lambda$,  we see that  $\lambda^i_0$ has to be a component of $t\lambda$. Hence
 the scaling $\frac{1}{t}\lambda^i_0$ is a component of $\lambda$, which goes to zero as $t$ goes to infinity. As a consequence, $\lambda^i_0$ is not a component of $\lambda$, which contradicts that $\lambda^i_0$ is a component of $t\lambda$. This proves Claim 2.

 By claim 2, we see that each component of the complement of the critical graph of $\Hopf(f^i)$ in  $X^i$  is either a half-infinite cylinder, a half-plane, or an infinite strip.

 \textbf{Claim 3: the complement of the critical graph of $\Hopf(f^i)$ contains no infinite-strip components} (see section \ref{subsec:meromorphic:leaves} for definition). Suppose to the contrary that the complement of the critical graph of $\Hopf(f^i)$ contains a component $C$ which is an infinite strip.  Let $\beta$ be a saddle connection in $C$ which connects two zeros $z^+,z^-$  of $\Hopf(f^i)$ belonging to different boundary components of $C$. Let $U$ be a $\delta$-neighbourhood of $\beta$ on $X^i$ under the flat metric $|\Hopf(f^i)|$ which contains no zeros of $\Hopf(f^i)$ other than $z^\pm$. Since the flat surface $(X_{t_n},\Hopf(f_{t_n}))$ converges to $(X_\infty,\Hopf(f_\infty))$, there exists a sequence of approximating maps $\eta_{t_n}:U\to X_{t_n}$ homeomorphic to the corresponding images such that
 \begin{itemize}
     \item $\eta_{t_n}(z^\pm)$ are zeros of $\Hopf(f_{t_n})$, and
    \item
      the pullback metrics $\eta_{t_n}^*(|\Hopf(f_{t_n})|$ converge to the restriction of $|\Hopf(f^i)|$ to $U$.
 \end{itemize}
 (Here note that there may be many choices for $\eta_{t_n}(z^\pm)$, say if multiple zeroes of $\Hopf(f_{t_n})$ converge to $z^\pm$: the precise choice will not be important in what follows.)
  The assumption that $U$ contains no zeros of $\Hopf(f^i)$ other than $z^\pm$ implies that as $n$ goes to infinity, every zero contained in the image $\eta_{t_n}(U)$ would collapse to $z^+$ or $z^-$. Since $\beta$ is a geodesic segment which contains no zeros of $\Hopf(f_\infty)$ in the interior, we may modify $\eta_{t_n}$ so that  the geodesic representative of $\eta_{t_n}(\beta)$  is also a saddle connection for large enough $n$. Let us fix such a large enough $\mathbf{\bar n}$.  It then follows from the recurrence property of leaves of measured foliations on closed surfaces  that there exists a closed curve $\gamma$ which is a concatenation of a horizontal leaf segment $\zeta_{t_\mathbf{\bar n}}$ and  a subset segment $\xi_{t_\mathbf{\bar n}}$ of $\eta_{t_\mathbf{\bar n}}(\beta)$ centered near the midpoint of $\eta_{t_{\mathbf{\bar n}}}(\beta)$ and with length $0<|\xi_{t_\mathbf{\bar n}}|<|\beta|/2$ (that $0<|\xi_{t_\mathbf{\bar n}}|$ follows from claim 1). In particular, the intersection number of $\gamma$ and the horizontal measured foliation of $\Hopf(f_{t_\mathbf{\bar n}})$, which is equivalent to $t_\mathbf{\bar n}\lambda$, satisfies:
  $$ 0<i(\gamma,{t_\mathbf{\bar n}}\lambda)<|\xi_{t_\mathbf{\bar n}}|<|\beta|/2. $$
  Therefore,
  \begin{equation}\label{eq:gamma:intersection:number}
      i(\gamma,t_n\lambda)=\frac{t_n}{{t_\mathbf{\bar n}}}i(\gamma,{t_\mathbf{\bar n}}\lambda)\to\infty, \text{ as } n\to\infty {\text{ because } \frac{t_n}{{t_\mathbf{\bar n}}} \to \infty }.
  \end{equation}
  Notice that for each $n>\mathbf{\bar n}$ the (topological) horizontal foliation of $\Hopf(f_{t_n})$ can be obtained from that of $\Hopf(f_{t_\mathbf{\bar n}})$ by a homeomorphism followed by a sequence of Whitehead moves. This implies that $\gamma$ can also be realized on $X_{t_n}$ as a concatenation of a horizontal leaf segment $\zeta_{t_n}$ and a subsegment $\xi_{t_{ n}}$ of $\eta_{t_{ n}}(\beta)$, which gives 
  $$ i(\gamma,t_n\lambda)<|\eta_{t_n}(\beta)|\to |\beta|, \text{ as }n\to\infty. $$
  This contradicts (\ref{eq:gamma:intersection:number}), completing the proof of Claim 3.

  The lemma then follows from claim  1, claim 2, and claim 3.
 \end{proof}
  As a direct consequence of Lemma \ref{lem:hopf:limit}, we have

 \begin{lemma}\label{lem:dual:tree}
    The dual tree of the lift of the horizontal measured foliation of the Hopf differential $\Hopf(f^i)$ to the universal cover consists of exactly one vertex and countably many half-infinite edges corresponding to the boundary of $\widetilde{Y^i}$.
 \end{lemma}

 Consider the family of harmonic maps $f_t:X_t\to Y$ with Hopf differential $\Hopf(f_t)=\Phi_t$.  The image of $X_{t}\backslash\mathscr P_R(\Phi_t)$ under the map $f_t$ is a (thickened) train track we denote by $\tau_{t,R}$. An {\it $\epsilon$ train track} that carries $\lambda$ is a train track that carries $\lambda$ and is contained in the $\epsilon$ neighbourhood of $\lambda$.  
  Combining Lemma \ref{lem:hopf:limit} and  Theorem \ref{thm:Minsky:traintrack}, we see that
  
   \begin{lemma}[Train-track approximation]\label{lem:TT:approximation}
   For any $\epsilon>0$, there exists $T_0>0$  such that for any $t>T_0$, the train track $\tau_{t,t^{1/4}}$ in $Y$  corresponding to $X_t\backslash\mathscr{P}_{t^{1/4}}(\Phi_t)$ is an $\epsilon$ train track which carries $\lambda$.
    \end{lemma}
 \begin{proof} 
 For any $\epsilon>0$, by Theorem \ref{thm:Minsky:traintrack}, there exists $T_0>0$ such that for any $t>T_0$ the thickened train track $\tau_{t,t^{1/4}}$,  being the image of $X_t\backslash \mathscr{P}_{t^{1/4}}(\Phi_t)$ under the harmonic diffeomorphism $f_t:X_t\to Y$, is contained in the $\epsilon$ neighbourhood of $\lambda$. It remains to show that $\tau_{t,t^{1/4}}$ carries $\lambda$.
 
   By Lemma \ref{lem:hopf:limit}, there exists $T_0>0$ such that every non-critical horizontal leaf in $\mathscr{P}_{t^{1/4}}(\Phi_t)$ is
   \begin{itemize}
     \item either a closed leaf homotopic to one of the boundary components of $\mathscr{P}_{t^{1/4}}(\Phi_t)$, or
     \item contained in a maximal leaf which is homotopic \textit{rel} $\partial_v \mathscr{P}_{t^{1/4}}(\Phi_t)$ to a horizontal boundary leaf of $\mathscr{P}_{t^{1/4}}(\Phi_t)$, where $\partial_v \mathscr{P}_{t^{1/4}}(\Phi_t)$ comprises the vertical boundary segments of $\partial\mathscr{P}_{t^{1/4}}(\Phi_t)$, or
     \item contained in a horizontal strip of finite height.
   \end{itemize}
   For the third case, since all these strips are strips of $t\lambda$, the topological type of these strips stabilizes as $t\to\infty$.   Let $\mathrm{Strip}_t$ be such a strip. By Lemma \ref{lem:in:annulus}, each maximal leaf in $\mathrm{Strip}_t$ has length at least $K_3'R_t=K_3't^{1/4}$. By Theorem \ref{thm:Minsky:polygon}(v), it follows that $\mathrm{Strip}_t$ has height  $d_t\leq C't^{1/4}$.  Let $\beta_t$ be a saddle connection in $\mathrm{Strip}_t$ connecting the two distinct boundary  components of $\mathrm{Strip}_t$. Then, by the recurrence property of horizontal leaves of $\Phi_t$, there exists a closed curve $\gamma$ which is a concatenation of a subsegment $\beta'_t$ of $\beta_t$ and a leaf segment $\zeta_t$ of $t\lambda$ which is the horizontal foliation of $\Phi_t$. In particular, we have $0<i(\gamma,t\lambda)<d_t$. Equivalently, $0<i(\gamma,\lambda)<d_t/t$. Letting $t\to\infty$ yields $0<i(\gamma,\lambda)= 0$, which is impossible. So the third case can not happen. For the first two cases, the image of every non-critical horizontal leaf of $\mathscr{P}_{t^{1/4}}(\Phi_t)$ is carried by the train-track $\tau_{r,t^{1/4}}$ corresponding to $X_t\backslash\mathscr{P}_{t^{1/4}}(\Phi_t)$. Therefore, the train track $\tau_{t,t^{1/4}}$ carries $\lambda$ for all $t>T_0$.
 \end{proof}
  \remarks Note that for Hopf differentials and choices of size $R$ of Minsky regions in less restrictive settings than we are considering here, the train track $\tau$ in $Y$  corresponding to the complement of $\mathscr{P}_R$ is not sufficient to carry $\lambda$, because it may not carry those bi-infinite leaves which are contained in  $\mathscr{P}_R$.

  We will use this lemma in Section~\ref{sec:convergence:Stretch-EarthquakeDisk}.

%========
 \subsection{Proof of Theorem \ref{thm:HR:SR} }
\label{sec:proof:thm1}

 \begin{proof}[Proof of Theorem \ref{thm:HR:SR}] The proof will be divided into two steps. In the first step, we will establish the convergence of harmonic map rays $\HR_{X_t,Y}$ when $X_t$ diverges along harmonic map dual rays. In the second step, we will identify the limit geodesic with the Thurston stretch line provided that the lamination $\lambda$ is maximal.  

\textbf{Step 1: convergence when $X_t$ diverges along the harmonic map dual ray $\HDR_{Y,\lambda}$.}

 Let $X_t=\HDR_{Y,\lambda}(t)$. Consider the family of harmonic maps $f_t:X_t \to Y$. By Lemma \ref{lem:precompactness}, any sequence of maps $f_{t_n}:X_{t_n}\to Y$ contains a convergent subsequence. For each component $Y^i$ of $Y\backslash\lambda$,   by Lemma \ref{lem:dual:tree}, the dual tree of the lift on $\widetilde{Y^i}$ of the pushforward of the horizontal foliation of any limit Hopf differential is the regular tree with one vertex and countably many half-infinite edges corresponding to the boundary edges of $\widetilde{Y^i}$.
 By Theorem \ref{thm:minimalgraph:uniqueness},
  the limit harmonic maps of all convergent sequences of $f_t:X_t\to Y$ are the same, say $f_\infty:X_\infty\to Y\setminus\lambda$. By Theorem \ref{thm:HR:limit:geodesic2}, any sequence of the family of harmonic map rays $\HR_{X_t,Y}:[1,\infty)\to \T(S)$ contains a subsequence which locally uniformly converges to some Thurston geodesic ray.
  However, by Proposition \ref{prop:limit:HR:PSL}, that limit Thurston geodesic ray of any convergent subsequence of $\HR_{X_t,Y}:[1,\infty)\to \T(S)$ is the piecewise harmonic stretch ray $\PSR_{Y,\lambda,f_\infty}:[1,\infty)\to\T(S)$ defined by $(Y,\lambda,f_\infty)$. (Note that $\HR_{X_t,Y}:[1,\infty)\to \T(S)$ is the subray of $\HR_{X_t,Y}:[0,\infty)\to\T(S)$  starting from $Y$ and that $\PSR_{Y,\lambda,f_\infty}$ is the subray of $\PSL_{Y,\lambda,f_\infty}$ starting from $Y$.)
  Hence, the family of harmonic map rays $\HR_{X_t,Y}$ converges to the piecewise harmonic stretch ray $\PSR_{Y,\lambda,f_\infty}$ (independent of the choice of subsequences). That the convergence is locally uniform follows from Theorem \ref{thm:HR:limit:geodesic}.

   \bigskip\textbf{Step 2: identifying the limiting geodesic with the Thurston stretch line for maximal geodesic laminations.} This follows from Lemma \ref{lem:id:PSR:SR} and completes the proof.
 \end{proof}

For convenience of a later reference, we summarize the following result from the discussion in step 1 of the proof of Theorem \ref{thm:HR:SR}. Let $Y\in\T(S)$ be a hyperbolic surface and let $\lambda$ be a measured geodesic lamination on $Y$. Let $\eta$ be the admissible measured foliation (Definition \ref{def:admissibility:general}) on $Y\backslash\lambda$ that comprises half-infinite cylinders foliated by closed leaves parallel to closed geodesic boundary components of $Y\backslash\lambda$ and half-planes foliated by bi-infinite leaves parallel to ideal boundary geodesics of $Y\backslash\lambda$. 
\begin{theorem}\label{thm:unique:eta}
Let $Y,\lambda,\eta$ be as above. 
    Then there exists a unique (possibly disconnected) punctured Riemann surface $X_\infty$, homeomorphic to $Y\backslash\lambda$, and a unique harmonic diffeomorphism $f_\infty:X_\infty\to Y\backslash\lambda$, homotopic to the identity map, such that the pushforward to $Y\backslash\lambda$ of the horizontal foliation of the Hopf differential of $f_\infty$ is measure equivalent to $\eta$.
\end{theorem}

We close this section with a discussion when $X_t$ diverges along a Teichm\"uller geodesic. 
\begin{theorem}\label{thm:convergence:TR}
    Let $Y\in\T(S)$ and let $\lambda$ be a maximal measured geodesic lamination on $Y$. Let $X_t\in\TR_{Y,\lambda}$. Then the path of harmonic map rays $\HR_{X_t,Y}$ converges to the Thurston stretch line $\SL_{Y,\lambda}$.
\end{theorem}
\begin{proof} 
 By Theorem \ref{thm:HR:limit:geodesic2}, for any divergent sequence $X_{t_n}\in \TR_{Y,\lambda}$, the sequence  $\HR_{X_{t_n},Y}:[0,\infty)\to \T(S)$  of harmonic map rays contains a subsequence that converges to a Thurston geodesic. By Proposition \ref{prop:limit:HR:PSL}, that limit Thurston geodesic is a piecewise harmonic stretch line $\PSL_{Y,\lambda',f_\infty}:(0,\infty)\to \T(S)$ defined by a chain recurrent geodesic lamination $\lambda'$ and a harmonic diffeomorphism $f_\infty:X_\infty \to Y\backslash\lambda'$  from some (possibly disconnected) punctured Riemann surface $X_\infty$. For simplicity, we assume that $\HR_{X_{t_n},Y}$ itself converges. By Theorem \ref{thm:generalized:stretchmap}, the piecewise harmonic stretch line $\PSL_{Y,\lambda',f_\infty}$ maximally stretches exactly $\lambda'$.  By \cite[Theorem 8.1]{Minsky1992}, we see that, up to considering a subsequence if necessary, the horizontal measured foliation of the harmonic map $h_{t_n}:X_{t_n}\to Y$ homotopic to the identity projectively converges to a measured foliation supported on $\lambda$. Then Theorem \ref{thm:HR:limit:geodesic2} implies $\lambda\subset\lambda'$.  Since $\lambda$ is maximal, it follows that $\lambda'=\lambda$. By Lemma \ref{lem:id:PSR:SR}, we see that the limit geodesic $\PSL_{Y,\lambda',f_\infty}$ is the Thurston stretch line $\SR_{Y,\lambda}$. The theorem then follows from the arbitrariness of the convergence subsequence of $\HR_{X_t,Y}$.
\end{proof}

\section{Convergence to Teichmuller rays} \label{sec:convergence:TeichRays}

Let $X\in\T(S)$ be a fixed Riemann surface and $\Phi$ be a holomorphic quadratic differential on $X$. Let $\HR_{X,\Phi}$ be the harmonic map  ray defined by $X$ and $\Phi$. For $s\geq0$, let $Y_s=\HR_{X,\Phi}(s)\in \HR_{X,\Phi}$ be the hyperbolic surface such that the Hopf differential of the harmonic map from $X$ to $Y_s$ is $s\Phi$. Let $T_h$ be the $\R$-tree dual to the horizontal measured foliation of $\tilde{\Phi}$, the lift of $\Phi$ to the universal cover $\tilde{X}$. Then $\tilde{X}$ is the minimal surface in $\tilde{Y}_s\times (T_h,2s^{1/2}d)$, also a rescaled minimal surface  in $s^{-1/2}\tilde{Y}_s\times (T_h,2d)$. As $s\to\infty$, $\tilde{X}$ is exactly the minimal surface in $T_v\times T_h$, where $T_v$ is the $\R$-tree dual to the vertical measured foliation of $\widetilde\Phi$. Let $X_{s,t}\in \HDR_{Y_s,\sqrt{s}\lambda}$ be such that $\Hor (\Hopf(X_{s,t}\to Y_s)=t\sqrt{s}\lambda$.
 Then $X_{s,1}\equiv X$  for all $s>0$. The goal of this section is to prove the following:

 \begin{theorem}\label{thm:HR:HDR:limit}
   For any fixed $X$ and $\lambda$, as $s\to\infty$, the family of harmonic map dual rays 
  \[\HDR_{Y_s,\sqrt{s}\lambda}:[1,\infty)\to\T(S)\]
  \color{black}
    locally uniformly converge to the Teichm\"uller geodesic ray $\TR_{X,\Phi}:[1,\infty)\to \T(S)$ with $\Hor(\Phi)=\lambda$.
 \end{theorem}
  {\it Convention.} In the remainder of this section, to simplify the notation, we will denote the dual tree $(T_\eta,2d)$  by $T_\eta$. 
  
 \subsection{Minimal surfaces in the product of trees} 
 Let $\alpha$ and $\beta$ be a pair of transverse measured foliations on the closed surface $S$, that is, for any homotopically nontrivial simple closed curve $\gamma$ we have $i(\alpha,\gamma)+i(\beta,\gamma)>0$.  By \cite[Theorem 5.1]{GardinerMasur1991}, there exists a unique holomorphic quadratic differential $\Psi$ whose  horizontal and vertical measured foliations are $\alpha$ and $\beta$  respectively.  Let $T_\alpha$ and $T_\beta$ be respectively the dual trees of $\alpha$ and $\beta$. Define the energy map  $E(\bullet,T_\alpha\times T_\beta):\T(S)\longrightarrow \R$ which associates to $Z\in\T(S)$ the (equivariant) energy of the equivariant harmonic map $\tilde{Z}\to T_\alpha\times  T_\beta$: the integral of the energy density over a fundamental domain of $\widetilde Z$ under the Fuchsian group defining $Z$ (see \cite[Page 111]{Wolf1996}). 
 
 \begin{lemma}\label{lem:proper:E:alpha:beta}
 The function 
  $E(\bullet,T_\alpha\times T_\beta):\T(S)\longrightarrow \R$ is  proper. Moreover,
  \[E(Z,T_\alpha\times T_\beta)\geq 4\|\Psi\|,\]
     where the equality holds if and only if $Z$ is the underlying Riemann surface of the quadratic differential $\Psi$.
    \end{lemma}
    \begin{proof}
    Recall  that
      $E(Z,T_\alpha\times T_\beta)=E(Z,T_\alpha)+E(Z,T_\beta)=2\Ext_Z(\alpha)+2\Ext_Z(\beta)$, where the second equation follows from \eqref{eq:energy:tree}.
      The properness  then follows from  the fact that
      $\Ext_Z(\alpha)+\Ext_Z(\beta)$ is a proper function over $\T(S)$.
      Moreover,
     \begin{eqnarray*}
       &&\Ext_Z(\alpha)+\Ext_Z(\beta)\\
        &\geq& 2 \sqrt{ \Ext_Z(\alpha)\cdot\Ext_Z(\beta)} \\
        &\geq &2i(\alpha,\beta)\qquad\qquad\qquad (\text{by \cite[Theorem 5.1]{GardinerMasur1991}})\\
        &=&2\|\Psi\|,
     \end{eqnarray*}
     where the equalities hold if and only if
     \begin{itemize}
      \item $\Ext_Z(\alpha)=\Ext_Z(\beta)$,
       \item $Z$ is the underlying Riemann surface of $\Psi$.
     \end{itemize}
 \end{proof}

   \subsection{Estimating the energy of harmonic maps to  trees}
 For any $Z\in\T(S)$ and each measured foliation $\mu$ on $Z$,  we denote by $E(Z,T_\mu)$ the equivariant energy of the harmonic map from $\tilde{Z}$ to the dual tree of the lift of $\mu$. By \eqref{eq:energy:tree}, we have $E(Z,T_\mu)=2\Ext_{Z}(\mu)$.  
  \begin{lemma}\label{lem:energy:tree}
    Let $Y_s\in \HR_{X,\Phi}$. Let $\bar{\lambda}$ be the vertical measured lamination of $\Phi$. Then for any $Z\in\T(S)$,
  we have $E(Z,T_{\bar{\lambda}})\leq E(Z,s^{-1}Y_s)$, where $E(Z,s^{-1}Y_s)$ is the energy of the harmonic diffeomorphism from $Z$ to $s^{-1}Y_s$ (see Section \ref{subsec:harmonic:map}).
  \end{lemma}
  
  \begin{proof}
     We define a family of equivariant projection maps:
 \[j_s: s^{-1}\tilde{Y}_s\longrightarrow T_{\bar{\lambda}}\]
  as follows. Recall that on the natural coordinates of $\Phi$ on $X$, the hyperbolic metric on $Y_s$ can be written as:
  \[ f^*_s Y_s= 2s(\cosh \G(z,s)+1)\d x^2+2s(\cosh \G(z,s)- 1)\d y^2,\]
  where $f_s: X\to Y_s$ is the unique harmonic map.
  Let $j_s: s^{-1}\tilde{Y}_s\longrightarrow T_{\bar{\lambda}}$ be the projection map along the vertical leaves of $\tilde{\Phi}$, the lift of $\Phi$ to $\tilde{X}$.
  Then $j_s$ collapses the vertical leaves while it scales the horizontal leaves by a factor of $\sqrt{\frac{2}{\cosh\G(z,s)+1}}<1$ at $f_s(z)\in \tilde{Y}_s$ (with respect to the hyperbolic metric on $\tilde{Y}_s$).  

  Let $F_s:\tilde Z\to s^{-1}\tilde Y_s$ be the harmonic map, then
  $$ E(Z,T_{\bar\lambda})\leq E(j_s\circ F_s) \leq E(F_s) =
  E(Z,s^{-1}Y_s).$$
  \end{proof}

  \subsection{Proof of Theorem \ref{thm:HR:HDR:limit}}
  \begin{proof}[Proof of Theorem \ref{thm:HR:HDR:limit}]
  Let $\bar\lambda$ be the vertical measured lamination of $\Phi$.
  Let $\Phi_{s,t}$ be the Hopf differential of $X_{s,t}\to Y_s$. Then the Hopf differential of $\widetilde{X}_{s,t}\to T_{\sqrt{s}t\lambda}$ is $-\widetilde{\Phi}_{s,t}$, the lift of $-\Phi_{s,t}$. Consider the equivariant harmonic map $\tilde{X}_{s,t}\to T_{\bar{\lambda}} \times tT_\lambda$. By Lemma \ref{lem:energy:tree}, we have
    \begin{eqnarray}
     \nonumber    &&E(X_{s,t},T_{\bar{\lambda}}\times T_{t\lambda})\\
     \nonumber
    &=& E(X_{s,t},T_{\bar{\lambda}})+
    E(X_{s,t},T_{{t\lambda}})\\
   \nonumber  &{\leq}&
     E(X_{s,t},s^{-1}Y_s)+
    E(X_{s,t},T_{{t\lambda}})
    \qquad\qquad ({{\text{by Lemma } \ref{lem:energy:tree}}})
    \\
   \nonumber  &=&
    s^{-1}E(X_{s,t},Y_s) +s^{-1}E(X_{s,t},T_{\sqrt{s}t\lambda})\\
    \nonumber   &=& s^{-1}E(X_{s,t},Y_s)+2s^{-1}\|\Phi_{s,t}\|\\
   \nonumber   &\leq & s^{-1}(\ell_{Y_s}(\sqrt{s}t\lambda)+C)+s^{-1}(\ell_{Y_s}(\sqrt{s}t\lambda)+C)\qquad
   (\text{by Lemma }\ref{lem:hyplength:quadraticnorm})\\
  \nonumber   &{=}&  2s^{-1}\left(t\ell_{Y_s}(\sqrt{s}\lambda)+C\right) \qquad\qquad\qquad \\
  \nonumber   &{\leq}&  2s^{-1}  (2t\|s\Phi\|+tC)+C)
  \qquad\qquad\qquad ({{\text{by Lemma } \ref{lem:hyplength:quadraticnorm}}})\\
    &=&4t  \|\Phi\|+2s^{-1}C(t+1). \label{eq:energy:Xst}
     \end{eqnarray}
     
    Combined with Lemma \ref{lem:proper:E:alpha:beta}, this implies that for any fixed $t$, $E(X_{s,t},T_{\bar{\lambda}}\times T_{t\lambda})\to 4t\|\Phi\|=4i(\bar{\lambda},t\lambda)$, as $s\to \infty$.
     It then follows from Lemma \ref{lem:proper:E:alpha:beta}  that
     $X_{s,t}\to \TR_{\Phi}(t)$, as $s\to\infty.$

     \vskip 5pt
     To see the locally uniform convergence, consider the function $\psi:[0,\infty)\to [0,\infty)$ defined as following. Let $X_{\infty,t}\in \TR_{X,\Phi}$ be the Riemann surface underlying the quadratic differential whose horizontal and vertical measured foliations are $t\lambda$ and $\bar{\lambda}$, respectively. Then for any  $Z\in\T(S)$, we have
      $E(Z,T_{\bar{\lambda}}\times T_{t\lambda})
     \geq E(X_{\infty,t},T_{\bar{\lambda}}\times T_{t\lambda})$. 
     For any $t_0\geq 1$, let
     \begin{equation}\label{eq:comparison:psi}
       \psi_{t_0}(r)=\max_{1\leq t\leq t_0}\max\{d_T(Z,X_{\infty,t}): E(Z,T_{\bar{\lambda}}\times T_{t\lambda})
     -E(X_{\infty,t},T_{\bar{\lambda}}\times T_{t\lambda})\leq r\}.
     \end{equation}
     Then $\psi$ is continuous and increasing in $r$.
     By (\ref{eq:energy:Xst}), we see that,
     \begin{eqnarray*}
      &&  E(X_{s,t},T_{\bar{\lambda}}\times T_{t\lambda})
     -E(X_{\infty,t},T_{\bar{\lambda}}\times T_{t\lambda})\\ &\leq &
     2s^{-1}C(t_0+1)\to 0, \text{ as } s\to\infty.
     \end{eqnarray*}
     Then
     \[d_T(X_{s,t},X_{\infty,t})\leq \psi_{t_0} \left(2s^{-1}C(t_0+1)\right)\to 0\]
     uniformly in $t\in[1,t_0]$ as $s\to \infty$. 
  \end{proof}

\section{Convergence to Teichmuller disks}\label{sec:convergence:TeichDisks}
In this section, we introduce two models of harmonic dual disks, both of which will converge to Teichm\"uller disks.

Let $Y_{s,\theta}$ be the hyperbolic surface such that $\Hopf(X,Y_{s,\theta})=se^{2i\theta}\Phi$. Let $\lambda_\theta$ be the horizontal foliation of $e^{2i\theta}\Phi$.
\begin{enumerate}
 \item[\textbf{M1.}] Let
      $$\HDD_{X,\Phi,s}=\bigcup_{0\leq\theta\leq\pi} \HDR_{Y_{s,\theta},\lambda_{\theta}}$$
      denote the (variable target) harmonic map dual disk.
  \item[\textbf{M2.}]   Let
      $$ \widehat{\HDD}_{X,\Phi,s}=\bigcup_{-\pi/2< \theta< \pi/2}\HDR_{Y_s,\lambda_\theta}$$
      denote the (fixed target) harmonic map dual disk.
\end{enumerate}
In both versions, there is a dependence of the lamination $\lambda_\theta$ on $\theta$.
In the first version $\bf{M1}$, the family of terminal points $Y_{s,\theta}$ also changes with $\theta$, while in the second version ${\bf M2}$, there is a single target surface $Y_s$, and the dependence on $\theta$ is only in the lamination.

 Let $\TD_{X,\Phi}$ be the Teich\"muller disk determined by $X$ and $\Phi$; i.e. the complex line in $\T(S)$ comprising surfaces whose Teich\"muller differentials from a given base point $X$ are all complex multiples of
 $\Phi$ on $X$.
 Let $ h_a:\mathcal{QT}(S)\to \mathcal{QT}(S)$  be the horocyclic flow on the quadratic differential bundle $\mathcal{QT}(S) \to \mathcal{T}(S)$ corresponding to the lower triang{ular} matrix 
 $$\left\{\begin{pmatrix}
     1 & 0 \\
     a & 1
   \end{pmatrix}:a\in\R\right\}.$$
   The horocycle flow acts in the standard way on $\R^2$ once we define charts from the Riemann surface to $\R^2$ by using natural coordinates. Given a Riemann surface $X$ and a holomorphic quadratic differential $\Phi$ on $X$, we choose local coordinates $z=x+iy$ so that $\Phi=dz^2$. With this coordinate,  the Teichm\"uller ray $\TR_{X,\Phi}$ may be interpreted as a family of Riemann surfaces underlying the quadratic differentials $\Phi_s:=(s^{1/2}dx+s^{-1/2}dy)^2$. For $a\in\R$, we denote by $h_a(\TR_{X,\Phi})$ the family of Riemann surfaces underlying $h_a(\Phi_s)$ for $s\geq1$.
   Let $$\TD^{h}_{X,\Phi}:=\bigcup_{a\in\R}h_a(\TR_{X,\Phi})$$ 
   be  the Teichm\"uller horodisk, which is the union of the horocyclic translates of $\TR_{X,\Phi}$.

\begin{theorem}\label{thm:convergence:disks}
   $\HDD_{X,\Phi,s}$ locally uniformly converges to the Teichm\"uller disk $\TD_{X,\Phi}$. $\widehat{\HDD}_{X,\Phi,s}$ locally uniformly converges to the Teichm\"uller horodisk $\TD^{h}_{X,\Phi}$.
\end{theorem}

Before proving Theorem \ref{thm:convergence:disks}, we state a lemma we will need, deferring its proof to the end of this section.

\begin{lemma}\label{lem:length:Ys:theta}
  Let $Y_s\in \HR_{X,\Phi}$ and let $\lambda_\theta$ be the horizontal foliation of $e^{2i\theta}\Phi$. Let $\delta>0$ be a constant which  is smaller than half of the shortest distance between zeros of $\Phi$. Then there exists ${s_0}>0$, which depends only on $\delta$, such that for $s > s_0$, we have
  $$ \ell_{Y_s}(\lambda_\theta)\leq 2\sqrt{s}i(\lambda_{\pi/2},\lambda_\theta)+
  96(g-1)\delta^2\sqrt{s}+\|\Phi\|.$$
\end{lemma}

\begin{proof}[Proof of Theorem \ref{thm:convergence:disks}]
  The proof is almost the same as the proof of Theorem \ref{thm:HR:HDR:limit}.

  (1) Convergence of $\HDD_{X,\Phi,s}$.
  Let $X_{s,t,\theta}$ be the Riemann surface such that the maximal stretch lamination of the harmonic map $X_{s,t,\theta}\to Y_{s,\theta}$ is $s^{1/2}t \lambda_\theta$. Let $\Phi_{s,t,\theta}$ be the Hopf differential of $X_{s,t,\theta}\to Y_{s,\theta}$. Then the Hopf differential of $\widetilde{X}_{s,t,\theta}\to T_{\sqrt{s}t,\theta}$ is $-\widetilde{\Phi}_{s,t,\theta}$, the lift of $-\Phi_{s,t,\theta}$.
  Therefore,
  \begin{eqnarray*}
    &&E(X_{s,t,\theta}, T_{\lambda_{\theta+\pi/2}}\times T_{t\lambda_\theta})\\
    &=& E(X_{s,t,\theta},T_{{\lambda_{\theta+\pi/2}}})+
    E(X_{s,t,\theta},T_{{t\lambda_\theta}})\\
    &\leq&
     E(X_{s,t,\theta},s^{-1}Y_{s,\theta})+
    E(X_{s,t,\theta},T_{{t\lambda_\theta}})
    \qquad (\text{by Lemma \ref{lem:energy:tree}})\\
    &=& s^{-1}E(X_{s,t,\theta},Y_{s,\theta})+
    s^{-1}E(X_{s,t,\theta},T_{{\sqrt{s}t\lambda_\theta}})\\
    &=& s^{-1}E(X_{s,t,\theta},Y_{s,\theta})+
    2s^{-1} \| \Phi_{s,t,\theta} \|\\
    &\leq & s^{-1}(\ell_{Y_{s,\theta}}(\sqrt{s}t\lambda_{\theta})+C)
    +s^{-1}(\ell_{Y_{s,\theta}}(\sqrt{s}t\lambda_{\theta})+C)\qquad
   (\text{by Lemma }\ref{lem:hyplength:quadraticnorm})\\
    &=& 2s^{-1} t \ell_{Y_{s,\theta}}(\sqrt{s}\lambda_{\theta})+2s^{-1} C\\
    &\leq & 2s^{-1} t(2 \|\Phi_{1,s,\theta}\|+C)+2s^{-1}C \qquad\qquad\text{by Lemma \ref{lem:hyplength:quadraticnorm}}  \\
    &=&4t\|\Phi_{1,1,\theta}\|+ 2s^{-1}C(t+1) \\
    &=&4t\|\Phi\|+ 2s^{-1}C(t+1).
  \end{eqnarray*}
  Therefore, $ \lim\limits_{s\to\infty}E(X_{s,t,\theta}, T_{\lambda_{\theta+\pi/2}}\times T_{t\lambda_\theta})\leq 4t\|\Phi\|=4\i(\lambda_{\theta+\pi/2},t\lambda_\theta)$.
  It then follows from Lemma \ref{lem:proper:E:alpha:beta} that $ X_{s,t,\theta}$ converges to the Riemann surface
  $ X_{\infty,t,\theta}\in\TD_{X,\Phi}$ underlying the quadratic differential $\Psi_{t,\theta}$ whose horizontal and vertical measured foliations are $t\lambda_\theta$ and $\lambda_{\theta+\pi/2}$.

  By considering a function similar to the one defined in (\ref{eq:comparison:psi}), we see that for any fixed $t>0$, the convergence is uniform in $(t,\theta)\in [1,t_0]\times [0,2\pi]$.

  \bigskip
  (2) Convergence of  $\widehat{\HDD}_{X,\Phi,s}$.
   Let $\widehat{X}_{s,t,\theta}$ be the Riemann surface such that the maximal stretch lamination of the harmonic map $\widehat X_{s,t,\theta}\to Y_{s}$ is $s^{1/2}t \lambda_\theta$.  Let $\widehat\Phi_{s,t,\theta}$ be the Hopf differential of $X_{s,t,\theta}\to Y_{s}$.
  Then
  \begin{eqnarray*}
    &&E(\widehat X_{s,t,\theta}, T_{\lambda_{\pi/2}}\times T_{t\lambda_\theta}) \\
    &=& E(\widehat X_{s,t,\theta},T_{{\lambda_{\pi/2}}})+
    E(\widehat X_{s,t,\theta},T_{{t\lambda_\theta}})\\
    &\leq&
     E(\widehat X_{s,t,\theta},s^{-1}Y_{s})+
    E(X_{s,t,\theta},T_{{t\lambda_\theta}})\qquad
   (\text{by Lemma }\ref{lem:energy:tree})\\
     &=& s^{-1}E(X_{s,t,\theta},Y_{s})+
    s^{-1}E(X_{s,t,\theta},T_{{\sqrt{s}t\lambda_\theta}})\\
    &=& s^{-1}E(X_{s,t,\theta},Y_{s})+
    2s^{-1} \| \widehat \Phi_{s,t,\theta} \|\\
    &\leq & s^{-1}(\ell_{Y_{s}}(\sqrt{s}t\lambda_{\theta})+C)
    +s^{-1}(\ell_{Y_{s}}(\sqrt{s}t\lambda_{\theta})+C)\qquad
   (\text{by Lemma }\ref{lem:hyplength:quadraticnorm})\\
   &=&2s^{-1/2}t\ell_{Y_{s}}(\lambda_\theta)+2s^{-1}C(t+1).
  \end{eqnarray*}
  Combined with  Lemma \ref{lem:length:Ys:theta}, this implies that for any $\delta$ which is smaller than half of the shortest distance between zeros of $\Phi$, there exists ${s_0(\delta)}>0$, which depends  only on $\delta$, such that for $s> s_0$, we have
  \[  E(\widehat X_{s,t,\theta}, T_{\lambda_{\pi/2}}\times T_{t\lambda_\theta})\leq
  4i(\lambda_{\pi/2},t\lambda_\theta)+192(g-1)\delta^2t
  +2ts^{-1/2}\|\Phi\|+2s^{-1}C(t+1).\]
Combining Lemma \ref{lem:proper:E:alpha:beta} and the arbitrariness of $\delta$, we see that as $s\to\infty$, we have that
   $E(\widehat X_{s,t,\theta}, T_{\lambda_{\pi/2}}\times T_{t\lambda_\theta})$ converges to $4i(\lambda_{\pi/2},t\lambda_\theta)$  uniformly in $(t,\theta)\in[1,t_0]\times(-\pi/2,\pi/2)$.
   Let $\widehat X_{\infty,t,\theta}\in{\TD^h_{X,\Phi}}$ be
  the Riemann surface  underlying the quadratic differential $\widehat{\Phi}_{t,\theta}$ whose horizontal and vertical measured foliations are $t\lambda_\theta$ and $\lambda_{\pi/2}$.
  It then follows from Lemma \ref{lem:proper:E:alpha:beta} that $\widehat X_{s,t,\theta}$ converge to $\widehat X_{\infty,t,\theta}\in\TD^h_{X,\Phi}$.
   By considering a function similar to the one defined in (\ref{eq:comparison:psi}), we see that for any fixed $t_0>0$, the convergence is uniform in $(t,\theta)\in [1,t_0]\times (-\pi/2,2\pi/2)$.
\end{proof}

\begin{proof}[Proof of Lemma \ref{lem:length:Ys:theta}]
  Let $\delta>0$ be a fixed constant which is smaller than half of the shortest distance between zeros of $\Phi$.
  For a zero $z_j$ of $e^{2i\theta}\Phi$ which is of order $n_j$, let $\mathscr{P}_j$ be a horizontal-vertical $(2n_j)$-polygon of $e^{2i\theta}\Phi$ around $z_j$ such that   every horizontal  segment and every vertical segment  of $\partial \mathscr{P_j}$ have the same $|\Phi|$-length $2\delta$. In particular, the $|\Phi|$-distance from $z_j$ to $\partial\mathscr P_j$ is $\delta$.
  Next, we decompose the horizontal foliation $F_\theta$ of $e^{2i\theta}\Phi$ outside the union $\cup_{j}\mathscr{P}_j$ into rectangles $\mathscr{R}_i$. Therefore,
  $$ F_\theta=\left(\cup_i (F_\theta\cap \mathscr{R}_i)\right)\bigcup \left(\cup_j (F_\theta\cap \mathscr{P}_j)\right). $$

  The hyperbolic length of leaves of
    $(F_{\theta}\cap \mathscr P_j)$ on $Y_s$ is not convenient to estimate. To overcome this, we need a modification. Notice that  the critical  leaves of  $F_\theta\cap\mathscr P_j$ decompose $\mathscr P_j$ into several rectangles $\{R_{jk}\}_{1\leq k\leq n_j}$, where $n_j$ is the order of the zero $z_j$ of $\Phi$.    We homotope (relative to its endpoints) each non-critical leaf $L$ of  $F_\theta\cap\mathscr P_j$ to a curve $L'$ which is contained in the boundary  $\partial R_{jk}\cap\partial \mathscr P_j$. Equivalently, $L\cup L'$ is the boundary of   the component of ${R}_{jk}\backslash L$ whose closure does not contain the zero $z_j$ of $\Phi$.
   Then the lengths of $L'$ and $\partial R_{ji}\cap \partial{\mathscr P}_j$  satisfy
   $$ \mathrm{Length}_{Y_s}(L')\leq \mathrm{Length}_{Y_s}(\partial R_{jk}\cap \partial{\mathscr P}_j)$$
   where we use the fact that $L'\subset \partial R_{ji}\cap \partial\mathscr{P}_i$.
     Let $\mathrm{Length}_{Y_s}(F_{\theta}\cap \mathscr{R}_{jk})$ be the hyperbolic length of the leaves of the foliation $F_{\theta}$ restricted to $\mathscr{R}_{jk}$, integrated over the induced transverse measure of the family of leaves. Then we have
    \begin{equation}\label{eq:Pj:Rij}
        \mathrm{Length}_{Y_s}(F_{\theta}\cap \mathscr{R}_{jk})\leq \delta \cdot \mathrm{Length}_{Y_s}(\partial R_{jk}\cap \partial{\mathscr P}_j).
    \end{equation}
    Let $\mathrm{Length}_{Y_s}(F_{\theta}\cap \mathscr{R}_i)$ be the hyperbolic length of the leaves of the foliation $F_{\theta}$ restricted to $\mathscr{R}_i$, integrated over the induced transverse measure of the family of leaves.
    By definition, $\ell_{Y_s}(\lambda_\theta)$ is the $Y_s$-length of geodesic representatives of leaves of $F_\theta$, integrated over the induced  transverse measure of family of leaves.
   Hence, 
   \begin{eqnarray}
   \nonumber
   \ell_{Y_s}(\lambda_\theta)
   &\leq& \sum_i \mathrm{Length}_{Y_s}(F_\theta\cap \mathscr{R}_i)+ \sum_j  \mathrm{Length}_{Y_s}(F_\theta\cap \mathscr{P}_j)\\
     \nonumber
   &=& \sum_i \mathrm{Length}_{Y_s}(F_\theta\cap \mathscr{R}_i)+ \sum_j  \sum_{1\leq k\leq n_j} \mathrm{Length}_{Y_s}(F_\theta\cap R_{jk})\\
     \nonumber
   &\overset{(\text{by \eqref{eq:Pj:Rij}})}{\leq}& \sum_i \mathrm{Length}_{Y_s}(F_\theta\cap \mathscr{R}_i)+ \sum_j  \delta \sum_{1\leq k\leq n_j} \mathrm{Length}_{Y_s}(\partial R_{jk}\cap \partial\mathscr{P}_j)
   \\
   &=&  \sum_i \mathrm{Length}_{Y_s}(F_\theta\cap \mathscr{R}_i)+\delta \sum_j \mathrm{Length}_{Y_s}(\partial\mathscr{P}_j). \label{eq:length:lambda:theta}
   \end{eqnarray}

  Recall that with respect to the canonical coordinate of $\Phi$, the hyperbolic metric on $Y_s$ can be expressed as:
  \begin{eqnarray}
    \nonumber ds_Y^2&=&2s(\cosh \G(z,s)+1)d x^2+2s(\cosh \G(z,s)- 1)d y^2\\
    &=&
    4sd x^2+2s(\cosh \G(z,s)- 1)(d x^2+d y^2),\label{eq:metric:Ys}
  \end{eqnarray}
  where $\G(z,s)=\log (1/|\nu(z,s)|)$ and $\nu(z,s)$ is the Beltrami differential of the harmonic map $X\to Y_s$.
  By Lemma \ref{lem:Minsky:decay}, we see that there exists $\mathbf{s_0}>1$, such that for every $s>\mathbf{s_0}$ and every $z\in \overline{X\backslash(\cup_j \mathscr P_j)}$, we have
  \begin{equation}\label{eq:s:cosh:G}
      2s(\cosh\G(z,s)-1)<1.
  \end{equation}

  Let $w_i$ and $h_i$ be respectively the horizontal width and vertical height of $\mathscr{R}_i$ with respect to $e^{2i\theta}\Phi$.
  Combining  (\ref{eq:metric:Ys}) and (\ref{eq:s:cosh:G}), we see that
  \begin{eqnarray}
  \nonumber\mathrm{Length}_{Y_s}(F_\theta\cap \mathscr{R}_i)
  &\leq& \left(2\sqrt{s}w_i \cos\theta + w_i\right) h_i\\
  \nonumber
  &=& 2\sqrt{s}i(F_{\pi/2}\cap \mathscr{R}_i, F_\theta\cap \mathscr{R}_i) + w_ih_i\\
  &=&2\sqrt{s}i(F_{\pi/2}\cap \mathscr{R}_i, F_\theta\cap \mathscr{R}_i) + \|\Phi\|_{\mathscr{R}_i}.
  \label{eq:length:Ri}
  \end{eqnarray}

  Recall that $\partial\mathscr{P}_j$ consists of $n_j$  horizontal arcs and $n_j$ vertical arcs with respect to $e^{2i\theta}\Phi$. Moreover, every such arc has the same $|\Phi|$-length $\delta$. Hence by  (\ref{eq:metric:Ys}) and (\ref{eq:s:cosh:G}),
  \begin{eqnarray}
  \nonumber && \mathrm{Length}_{Y_s}(\partial\mathscr{P}_j)\\
  \nonumber  &\leq& n_j \left(2\sqrt{s}\delta\cos\theta+
  2\delta \right)+
  n_j \left(2\sqrt{s}\delta\sin\theta+
  2\delta\right)\\
  &\leq& 8n_j\sqrt{s}\delta. \label{eq:length:Pj}
  \end{eqnarray}

   Combining (\ref{eq:length:lambda:theta}), (\ref{eq:length:Ri}), and (\ref{eq:length:Pj}), we have
 \begin{eqnarray*}
    \ell_{Y_s}(\lambda_\theta)&\leq &2\sqrt{s}i(F_{\pi/2}, F_\theta) + \| \Phi\| +\sum_j 8n_j\delta^2\sqrt{s}
   \\
   &\leq& 2\sqrt{s}i(F_{\pi/2}, F_\theta) +\| \Phi\| + 96(g-1)\delta^2\sqrt{s},
  \end{eqnarray*}
  where the last inequality follows from the fact $\sum_j n_j\leq 12(g-1)$.
\end{proof}

\section{Convergence to stretch-earthquake disks} \label{sec:convergence:Stretch-EarthquakeDisk}

We recall some basic facts about earthquake deformations. For more details, we refer to \cite[Section II]{Kerckhoff1983} (see also  \cite{Kerckhoff1983b} and \cite{mirzakhani2008ergodic}). Let 
 $\mu\alpha$ be a weighted simple closed geodesic on $X$. The \emph{time $s$ earthquake} $\mathcal{E}^s_{\mu\alpha}(X)$ of $X$ is defined to be the hyperbolic surface obtained from $X$ by twisting left by a distance $s\mu{\color{red}\ell_X(\alpha)}$ along $\alpha$. The weight $\mu$ determines the speed of twisting. This construction extends to general measured laminations as follows. For a measured geodesic lamination $\lambda$, let $\mu_i\alpha_i$ be a sequence of weighted simple closed geodesics that converges to $\lambda$ in $\ML(S)$. The \emph{time $s$ earthquake} $\mathcal{E}^s_{\lambda}(X)$ of $X$ is defined as the limit surface of $\mathcal{E}^s_{\mu_i\alpha_i}(X)$. It is a non-trivial fact that the limit exists and is independent of the choice of the approximating sequence $\mu_i\alpha_i$ (\cite[Corollary 2.5]{Kerckhoff1983}). 

   \begin{definition}[Stretch-earthquake disk]
   Let $Y\in\T(S)$  be a hyperbolic surface and $\lambda\in \ML(S)$ a measured foliation/lamination. Let $\PSL_{Y,\lambda}$ be the piecewise harmonic stretch line obtained from Theorem \ref{thm:HR:SR} as the limit of harmonic map rays $\HR_{X_t,Y}$ where $X_t$ diverges along the harmonic map dual ray $\HDR_{Y,\lambda}$.  Let $\mathcal{E}^s_{\lambda}(Y)$ be the surface obtained from $Y$ by acting by a time $s$ (left) earthquake along $\lambda$. Define the stretch-earthquake disk $\mathbf{SED}(Y,\lambda)$ of $(Y,\lambda)$ to be the set:
   \begin{equation*}
     \bigcup_{-\infty<s<+\infty}\PSL_{\mathcal{E}^s_{\lambda}(Y),\lambda}(0,+\infty).
   \end{equation*}
 \end{definition}

 \begin{definition}[Hopf differential disk]
   Define the Hopf differential disk $\mathbf{HDD}(X,\Phi)$ of $(X,\Phi)$ to be the set
   \begin{equation*}
     \bigcup_{-\pi\leq \theta\leq \pi} \HR_{X,e^{\i\theta}\Phi}(0,\infty).
   \end{equation*}
 \end{definition}

 Let $X_t\in\HDR_{Y,\lambda}$ be the Riemann surface such that the horizontal foliation of  $\Hopf(X_t,Y)=\Phi_t$ is $t\lambda$.
 Let $Y(t,r,s)$ be the hyperbolic surface such that $\Hopf(X_t,Y(t,r,s))=re^{i\frac{s}{2t}}\Phi_t$ and $Y_r=\PSL_{Y,\lambda}(r)\in \PSL_{Y,\lambda}$.
 
The goal of this section is to prove Theorem~\ref{thm:Hopfdisk:limit1}, which we restate here.

 \begin{theorem}\label{thm:Hopfdisk:limit}
    Let $Y\in\T(S)$  be a hyperbolic surface and $\lambda\in \ML(S)$ a measured foliation/lamination. Let $X_t\in\HDR_{Y,\lambda}$ be as described above.

    Then the family of Hopf differential disks $(\mathbf{HDD}(X_t,\Phi_t),Y)$ with base point $Y$ locally uniformly converges to the stretch-earthquake disk $(\mathbf{SED}(Y,\lambda),Y)$ with base point $Y$. Namely, for any prescribed $\mathbf{s}>0$ and $0<\mathbf{r}<\mathbf{r}'$, the point
    $ Y(t,r,s)$  converges to $\mathcal{E}^s_{\lambda}(Y_r)$, uniformly in $(r,s)\in[\mathbf{r},\mathbf{r}']\times [-\mathbf{s},\mathbf{s}]$, as $t\to\infty$.
 \end{theorem}

 Recall that the family of harmonic maps $f_t:X_t\to Y$ converges to the harmonic map $f_\infty:X_\infty\to Y$ with $\Hopf(f_\infty)=\Phi_\infty$, where $\Phi_\infty$ is the union of half-infinite cylinders and half-planes (Lemma \ref{lem:hopf:limit}). 

   The proof of  Theorem  \ref{thm:Hopfdisk:limit} relies on a generalization of shearing coordinates of $\T(S)$,  namely the shear-shape coordinates of $\T(S)$ developed by Calderon-Farre (\cite{CalderonFarre2021}).

 \subsection{Shearing}
\label{sec:hss}
In order to glue crowned hyperbolic surfaces, one needs a foliation transverse to the boundary of crowned surfaces. Thurston \cite{Thurston1998} uses the horocycle foliation (where the crowned hyperbolic surfaces are ideal triangles), Calderon and Farre \cite{CalderonFarre2021} use the orthogeodesic foliation. Here we use the extended foliation of the vertical foliation of Hopf differentials as in the construction of piecewise harmonic stretch lines.

In this subsection, we will construct \emph{shear-shape cocycles} of $\lambda$ using the vertical foliation of Hopf differentials, following the ideas and constructions of Calderon-Farre. 

\subsubsection{Extended vertical foliations}\label{sec:extended:vertical:foliation:11}

Let $Y\in\T(S)$ be an arbitrary hyperbolic surface and let $\lambda$ be a measured geodesic lamination on $Y$. Let $\eta$ be the admissible measured foliation (Definition \ref{def:admissibility:general}) on $Y\backslash\lambda$ that comprises half-infinite cylinders foliated by closed leaves parallel to closed geodesic boundary components of $Y\backslash\lambda$ and half-planes foliated by bi-infinite leaves parallel to ideal boundary geodesics of $Y\backslash\lambda$. Note that $\eta$ is determined by $\lambda$. By Theorem \ref{thm:unique:eta}, there exists a unique (possibly disconnected) punctured Riemann surface $X_\infty$, homeomorphic to $Y\backslash\lambda$, and a unique harmonic diffeomorphism $f_\infty:X_\infty\to Y\backslash\lambda$, homotopic to the identity map, such that the pushforward to $Y\backslash\lambda$ of the horizontal measured foliation of the Hopf differential of $f_\infty$ is measure equivalent to $\lambda$. By Lemma \ref{lem:glue:foliation} and Lemma \ref{lem:intersection:extend}, the pushforward to $Y\backslash\lambda$ of the vertical foliation of the Hopf differential of $f_\infty$ extends to a unique measured foliation $\beta_Y$ on $Y$ which is transverse to $\eta$ and $\lambda$ and intersects $\lambda$ orthogonally at every intersection point.  

\begin{quote}
  \it   In the remainder of this section, for any $Z\in\T(S)$, whenever we mention $\eta$ and $\beta$ on $Z\backslash\lambda$, we mean (unless otherwise stated) their realization on $Z\backslash\lambda$ as the pushforward by $f_\infty$ of the horizontal foliation and vertical foliation of the Hopf differential of $f_\infty$ (where $X_\infty$, $f_{\infty}$, and $\Hopf{f_{\infty}}$ refer to $Z$ and $\lambda$).
\end{quote}
 
Let $\T(S\backslash \lambda)$ be the \tec space of crowned hyperbolic surfaces homeomorphic to $S\backslash\lambda$.  Let $\MF(\lambda)\subset \MF(S)$ be the  subset of measured foliations transverse to $\lambda$. 
The construction above  defines a map 
\begin{equation}\label{eq:pi:vert}
    \pi_{\lambda} :\T(S)\to \T(S\backslash \lambda)\times \MF(\lambda)
\end{equation}
 which sends $Z\in\T(S)$ to the pair $(Z\backslash\lambda,\beta_Z)$.  Notice that $\beta_Z$ prescribes identifications for gluing both the leaves of $\beta_Z|_{Z\setminus \lambda}$ as well as the components of $Z\backslash \lambda$ across $\lambda$.    
To describe the gluing process, we introduce an analogue of the \textit{shear-shape cocycle} developed by \cite{CalderonFarre2021} using pointed geodesics.

\subsubsection{Dual arc systems} We start with the definition of dual arc systems.  
Let $Y,\lambda$ and $\eta$ be as above. Let $\pi_{\lambda}(Y)=(Y\backslash\lambda,\beta_Y)$. For simplicity, we denote $\beta_Y$ by $\beta$. Recall that every leaf of $\beta|_{Y\backslash\lambda}$ approaching $\partial Y\subset\lambda$ hits $\partial {Y\backslash\lambda}$ orthogonally (Lemma \ref{lem:intersection:extend}). For every leaf segment $e$ of $\eta$ that connects a pair of singular points of $\eta$ (a \emph{\enquote{saddle connection of $\eta$}}), there exists a strip $\mathrm{Strip}_e$ of $\beta|_{Y\backslash\lambda}$ foliated by regular leaves of $\beta|_{Y\backslash\lambda}$ that intersect $e$ orthogonally and approaches a pair of boundary geodesics of $\partial ({Y\backslash\lambda})$ orthogonally. Let $c_e$ be the height of $\mathrm{Strip}_e$. Then by Lemma \ref{lem:intersection:extend} we see that $c_e$ is half the hyperbolic length of the projection of $\mathrm{Strip}_e$ to a boundary geodesic of ${Y\backslash\lambda}$ along leaves of $\beta|_{Y\backslash\lambda}$.  Note that any two leaves of $\beta|_{Y\backslash\lambda}$ in $\mathrm{Strip}_e$ are homotopic rel $\partial ({Y\backslash\lambda})$. Let $\alpha_e$ be the regular leaf of $\beta|_{Y\backslash\lambda}$ that cuts $\mathrm{Strip}_e$ into two substrips of equal height. The \textit{dual arc system} $\underline{\alpha}({Y\backslash\lambda})$ of ${Y\backslash\lambda}$ is defined to be the union $\cup_e\alpha_e$ where $e$ ranges over all saddle connections of $\eta$. Define the \emph{weighted dual arc system} as the formal sum:
  \[\underline{A}({Y\backslash\lambda}):=\sum_{e} c_e \alpha_e,\]
  where $e$ ranges over all saddle connections of $\eta$. Notice that $\eta$ has finitely many saddle connections.
    The dual arc system cuts ${Y\backslash\lambda}$ into pieces, each of which contains exactly one singularity of $\eta$ (here $\eta$ and $\beta|_{Y\backslash\lambda}$ have the same set of singularities). 
  Following Caldron-Farre, we call each such piece a \emph{hexagon}, no matter its shape. Let $\mathbf{H}$ be the set of hexagons of ${Y\backslash\lambda}$.

\subsubsection{Pointed geodesics and shear}\label{sec:pointed:geodesics}  Let $\widetilde{Z}$, $\widetilde{Y\backslash\lambda}$, $\widetilde{\beta}$, $\widetilde{\eta}$, and $\widetilde{\mathbf{H}}$ be respectively the lifts of $Z$, $Y\backslash\lambda, \beta,\eta,\mathbf{H}$. 
  We now define a family of base points associated to the boundary leaves of $\widetilde{\lambda}$. Let $H_v\in\widetilde{\mathbf H}$ be  a hexagon with a singular point $v$ of $\widetilde{\beta|_{Y\backslash\lambda}}$. For a boundary leaf  $g_v$ of $\widetilde{\lambda}$ intersecting $\partial H_v$, define $p_v$ to be the  projection of $v$ to $g$ along the half-infinite critical leaf of $\widetilde{\beta|_{Y\backslash\lambda}}$ that starts at $v$ and approaches $g$. The pair $(g_v,p_v)$ is called a \textit{pointed geodesic}.
  For any pair of hexagons $H_v,H_w\in\widetilde{\mathbf H}$ belonging to distinct components of $\widetilde{Y\backslash\lambda}$, there is a unique geodesic $g_v$ intersecting $\partial H_v$ that separates $v$ from $w$. Symmetrically, there is a unique geodesic $g_w$ intersecting $\partial H_w$ that separates $w$ from $v$. The pair $(g_v,p_v)$ and $(g_w,p_w)$ is said to be \emph{simple} if  $g_v$ and $g_w$ are connected by a regular leaf segment $\xi_v^w$ (not unique) of $\widetilde{\beta}$ such that for each component $\Omega$ (intersecting $\xi_v^w$) of $\widetilde{Y\backslash\lambda}$ , the segment $\xi_v^w$ cuts out a cusp region from $\Omega$. Let $q_v$ (resp. $q_w)$ be the intersection point between $g_v$ (resp. $g_w$) and $\xi_v^w$. The \emph{shear} $\sigma_\lambda(Y)(v,w)$  among $H_v$ and $H_w$ is defined to be the sum of the signed distance from $p_v\in g_v$ to $q_v\in g_v$ along $g_v$ and the signed distance from $p_w\in g_w$ to $q_w\in g_w$ along $g_w$, with respect to the orientations of $g_v$ and $g_w$ induced from the orientation of $\widetilde{Y\backslash\lambda}$.  

\subsubsection{Weighted system}\label{sec:weighted:system}
   Let $\tau$ be an \emph{$\epsilon$ train track} of $\lambda$, that is, a train track which carries $\lambda$ and which is contained in the $\epsilon $ neighbourhood of $\lambda$ such that $Y\backslash\lambda$ and $Y\backslash\tau$ have the same topological type. The orientation of $Y$ induces an orientation for each component of $Y\backslash\lambda$. A {\it standard smoothing}  train track $\tau_{\underline{\alpha}}$ of $\tau\cup\underline{\alpha}({Y\backslash\lambda})$ is a smoothing at each intersection point $\tau\cap \underline{\alpha}({Y\backslash\lambda})$ so that the incoming tangent vector corresponding to  $\underline{\alpha}({Y\backslash\lambda})$ points in the positive direction with respect to the boundary orientation of $Y\backslash \lambda$ (see \cite[Figure 12 on page 2049]{CalderonFarre2021}).  
 The shears among simple pairs of pointed geodesics defined as above define a \emph{weighted system} $\mathbf{w}(Y)$ of $\tau_{\underline{\alpha}}$ as follows.
 \begin{itemize}
   \item To each branch corresponding to a component $\alpha$ of $\underline{\alpha}$, assign the weight to be the height of the strip of $\beta|_{Y\backslash\lambda}$ foliated by leaves parallel to $\alpha$.
   \item To each branch $b$ which is not a component of $\underline{\alpha}$, choose a lift $\widetilde{b}$ of $b$. Let $H_v,H_w\in \widetilde{\mathbf H}$ be the pair of hexagons adjacent to $\widetilde{b}$ and set the weight to be $\sigma_\lambda(Y)(v,w)$.
 \end{itemize}
  The proof that $\mathbf{w}(Y)$ satisfies the switch condition of $\tau_{\underline{\alpha}}$ is identical to the proof of \cite[Lemma 13.6]{CalderonFarre2021}.  The advantage of working with weighted systems is that each weighted system is determined by the weight of branches, whose number is uniformly bounded in terms of the topology of $Y$. This finiteness will be used later in the proof of locally uniform convergence.

  \subsubsection{Shear-shape cocycle}\label{sec:ss:cocycle}
By \cite[Proposition 7.13 and Proposition 9.5]{CalderonFarre2021}, each weighted system of $\tau_{\underline{\alpha}}$ defines a shear-shape cocycle $(\underline{A}({Y\backslash\lambda}),\sigma)$ in the sense of \cite[Definition 7.11]{CalderonFarre2021}, where $\sigma$ is a function which assigns to every arc $k$ transverse to $\lambda$ and disjoint from $\underline{\alpha}:=\cup \alpha_i$  a real number $\sigma(k)$, satisfying the following axioms:
    \begin{enumerate}[(1)]
      \item (support): If $k$ does not intersect $\lambda$, then  $\sigma(k)=0$.
      \item (transverse invariance): If $k$ and $k'$ are isotopic transverse to $\lambda$  and disjoint from $\underline{\alpha}$ then $\sigma(k)=\sigma(k')$.
      \item (finite additivity): If $k=k_1\cup k_2$ where $k_i$ have disjoint interiors then $\sigma(k)=\sigma(k_1)+\sigma(k_2)$.
      \item ($\underline{A}$-compatibility): Suppose that $k$ is transverse to $\lambda$ and isotopic rel endpoints to some arc which may be written as $t_i\cup l$, where $t_i$ is an arc which meets $\alpha_i$ exactly once and is disjoint from $\lambda\cup \underline{\alpha}\backslash\{\alpha_i\}$,  and $l$ is an arc transverse to $\lambda$ and disjoint from $\underline{\alpha}$. Then the loop $k\cup t_i\cup l$ encircles a unique point $p$ of $\lambda\cap \underline{\alpha}$, and $$ \sigma(k)=\sigma(l)+\epsilon c_i $$
          where $\epsilon$ denotes the winding number of $k\cup t_i\cup l$ about $p$, where the loop is oriented so that the edges are traversed $k$ then $t_i$ then $l$ (see \cite[Figure 9 on page 2038]{CalderonFarre2021}).
    \end{enumerate}
    
    Applying verbatim the argument in the proof of \cite[Lemma 13.9]{CalderonFarre2021}, we see that the resulting  shear-shape cocycle is independent of the choice of $\tau_{\underline{\alpha}}$.  This defines a map 
    \begin{align*}
        \Pi_\lambda:\T(S)&\to \T(S\backslash\lambda)\times\mathcal{SH}(\lambda)\\
        Y &\mapsto (Y\setminus \lambda, (\underline{A}(Y\backslash\lambda),\sigma_\lambda(Y)))
    \end{align*}
    where $\mathcal{SH}(\lambda)$ is the space of shear-shape cocycles of $\lambda$. 
    \begin{proposition}\label{prop:cocycle:injective}
        The map $\Pi_\lambda:\T(S)\to \T(S\backslash\lambda)\times\mathcal{SH}(\lambda)$ is homeomorphic onto its image.
    \end{proposition}
 \begin{proof}
 We first show that $\Pi_\lambda$ is injective.
    Suppose that there are $Z,Z'\in\T(S)$ such that $\Pi_{\lambda}(Z')=\Pi_\lambda(Z)$, then by the computation of monodromy in the proof of Theorem \ref{thm:generalized:stretchmap}, we see that $Z$ and $Z'$ have the same monodromy group, so they are equal. Another proof is to construct an equivariant isometry from the universal cover $\widetilde{Z}$ to $\widetilde{Z'}$ by following the argument as in the proof of \cite[Lemma 11]{Bonahon1996} or \cite[Proposition 13.12]{CalderonFarre2021} (except for the second paragraph there because we include $\T(S\backslash\lambda)$ in the codomain of the map $\sigma_\lambda$).

    Next, note that $\Pi_\lambda$ is clearly continuous. 
    For the inverse map $\Pi_\lambda^{-1}:\Pi_\lambda(\T(S))\to \T(S)$ from the image $\Pi_\lambda(\T(S))$, we reconstruct the hyperbolic surface $Z$ from the components $Z\backslash\lambda\in\T(S\backslash\lambda)$ of the complement of $\lambda$ and the shear-shape cocycle $(\underline{A}(Z\backslash\lambda),\sigma_\lambda(Z))\in\mathcal{SH}(\lambda)$ via a canonical construction: as mentioned in Section \ref{sec:extended:vertical:foliation:11}, the crowned surface $Z\backslash\lambda$ and the admissible foliation $\eta$ (determined by $\lambda$) on $Z\backslash\lambda$ gives a unique vertical foliation $\beta|_{Z\backslash\lambda}$ on $Z\backslash\lambda$. Then the function $\sigma_\lambda(Z)$ provides the data for gluing the components of $\beta|_{Z\backslash\lambda}$ across $\lambda$ to get a unique measured foliation $\beta_Z$ on the whole surface; using that extension of the foliation $\beta_Z$, we then glue together components of $Z\backslash\lambda$ in a way similar to the method we used in the construction of a piecewise harmonic stretch line (cf. subsection~\ref{subsec:construction:PSL}).  The computation of monodromy in the proof of Theorem \ref{thm:generalized:stretchmap} then implies that the inverse map $\Pi_\lambda^{-1}:\Pi_\lambda(\T(S))\to \T(S)$  is also continuous. Therefore, the map $\Pi_\lambda$ is homeomorphic onto its image $\Pi_\lambda(\T(S))$.
 \end{proof}
  Let
  \begin{eqnarray*}
      \hat{\sigma}_\lambda:\T(S)&\to&\mathcal{SH}(\lambda) \\
      Y&\mapsto & (\underline{A}(Y\backslash\lambda),\sigma_\lambda(Y))
  \end{eqnarray*}
 be the composed map $\T(S)\to \T(S\backslash\lambda)\times\mathcal{SH}(\lambda)\to\mathcal{SH}(\lambda)$, where the second map is the projection map.

 \subsection{
Proof of Theorem \ref{thm:Hopfdisk:limit} } 
In view of the map $\Pi_{\lambda}:\T(S)\to\T(S\backslash\lambda)\times \mathcal{SH}(\lambda)$,  any limit of $Y(t,r,s)$, as $t\to\infty$,  is described by a crowned surface and a shear-shape cocycle of $\lambda$ which encodes the shearing between the components of the limit crowned surface.

  \begin{lemma}[shape convergence] \label{lem:shape:convergence} $Y(t,r,s)\backslash\lambda$ converges to $Y_r\backslash \lambda$, locally uniformly in $(r,s)\in[\mathbf{r},\mathbf{r}']
 \times[-\mathbf{s},\mathbf{s}]$, as $t\to\infty$.
  \end{lemma}

\begin{proof}

Recall that $\Hopf(X_t,Y(t,r,s))=e^{i\frac{s}{2t}}\Hopf (X_t,Y(t,r,0))=re^{i\frac{s}{2t}}\Phi_t$, where $\Phi_t=\Hopf(X_t,Y)$. By (\ref{eq:quadratic:norm},\ref{eq:Bochner:H}), we see that the energy density $\mathbf{e}$, which is defined to be $\mathcal{H}+\mathcal{L}$ in the notations introduced in (\ref{eq:quadratic:norm},\ref{eq:Bochner:H}), satisfies:
\[\mathbf{e}(X_t,Y(t,r,s))=\mathbf{e}(X_t,Y(t,r,0)).\]
The pullback of the hyperbolic metric $\rho(t,r,s)$ of $Y(t,r,s)$  to $X_t$ via the  harmonic map $f_{t,r,s}:X_t\to Y(t,r,s)$ is:
\[(f_{t,r,s})^*\rho(t,r,s)=re^{i\frac{s}{2t}}\Phi_t +\overline{re^{i\frac{s}{2t}}\Phi_t }+\sigma_t\mathbf{e}(X_t,Y(t,r,s)),\]
where $\sigma_t$ is the hyperbolic metric on $X_t$.
This gives
\[(f_{t,r,s})^*\rho(t,r,s)-(f_{t,r,0})^*\rho(t,r,0)=
r(e^{i\frac{s}{2t}}-1)\Phi_t +\overline{r(e^{i\frac{s}{2t}}-1)\Phi_t }.\]
Combining this with the fact that  $\Phi_t$ converges to $\Phi_\infty$,  uniformly on compact subsets of $X_\infty\overset{homeo}{\sim}{Y\backslash\lambda}$,   as $t\to\infty$, we see that on any compact   subset of $Y\backslash\lambda$, we have $Y(t,r,s)$ converges to $Y_r$ uniformly in $(r,s)\in[\mathbf{r},\mathbf{r}']
 \times[-\mathbf{s},\mathbf{s}]$ as $t\to\infty$. Since any crowned hyperbolic surface is determined by its truncation (i.e. removing a small neighbourhood for each cusp) and each truncation is of arbitrary small distance to compact subsets of the crowned surface, so $Y(t,r,s)\backslash\lambda$ itself converges to $Y_r\backslash\lambda$ uniformly in $(r,s)\in[\mathbf{r},\mathbf{r}']
 \times[-\mathbf{s},\mathbf{s}]$ as $t\to\infty$.
 \end{proof}

  \begin{lemma}[shear-shape convergence]\label{lem:shear:convergence}
  The shear-shape cocycle $\hat \sigma_\lambda({Y(t,r,s)})$ converges to $\hat \sigma_\lambda(\mathcal{E}^t_\lambda (Y_r))$ uniformly in $(r,s)\in[\mathbf{r},\mathbf{r}']
 \times[-\mathbf{s},\mathbf{s}]$ as $t\to\infty$. 
   \end{lemma}

 \begin{proof}
  Notice that the set of dual arc systems
  $\{\underline{\alpha}( Y_r\backslash\lambda): r\in[\mathbf{r},\mathbf{r}']\}$ is a finite set. The locally uniform convergence of $Y(t,r,s)\backslash\lambda$ to $Y_r\backslash\lambda$ in $(r,s)\in [\mathbf{r},\mathbf{r}']\times [-\mathbf{s},\mathbf{s}]$ as $t\to \infty$ then implies that
  $$\{\underline{\alpha}(Y(t,r,s)\backslash\lambda):
  t>T_0,r\in[\mathbf{r},\mathbf{r}'],
  s\in [-\mathbf{s},\mathbf{s}] \}$$ is a  also finite and that the weighted dual arc system $\underline{A}(Y(t,r,s)\backslash\lambda)$ uniformly converges to the weighted dual arc system $\underline{A}(Y_r\backslash\lambda)$ for $r\in[\mathbf{r},\mathbf{r}']$. To prove the lemma, we may assume that the set above contains exactly one element, say the dual arc system $\underline{\alpha}$.

  Next, we calculate the transverse cocycle  $\sigma_\lambda(Y(t,r,s))$ of $Y(t,r,s)$  with respect to $\lambda$, using the train-track coordinates of transverse cocycles described in Section \ref{sec:weighted:system}.  We begin with a preliminary discussion in the setting where $s=0$, and then extend the analysis to the case where $s \neq 0$.

For the analysis, we fix a convenient choice of scale $R_t$ for the Minsky region $\mathscr{P}_{R_t}(\Phi_t)$. In particular, let $R_t:=t^{1/4}$. By Lemma \ref{lem:TT:approximation}, for any $\epsilon>0$, there exists $T_0>0$ such that for any $t>T_0$, the train track $\tau_{t,R_t}$ in $Y$  corresponding to $X_t\backslash\mathscr{P}_{R_t}(\Phi_t)$ carries $\lambda$ and is contained in the $\epsilon$ neighbourhood of $\lambda$. Lift everything to the universal cover.  Together with the dual arc system $\underline{\alpha}$, this train-track gives a finite coordinatization of $\sigma_\lambda(Y(t,r,s))$, which consists of a finite collection of simple pairs of (sufficiently close) pointed geodesics, denoted by $\mathscr G$. Let $((g_v,p_v),(g_w,p_w))\in\mathscr G$ be such a simple pair of pointed geodesics with $p_v\in g_v$ and $p_w\in g_w$ being the base points.
 Let $\beta$ be the extended vertical foliation on $Y$ of $\Phi_\infty$ (which is also the measured foliation in $\MF_{\eta}(\lambda)$ given by the map $\pi_{\lambda}(Y)$ in \eqref{eq:pi:vert}). 
  Since $(g_v,p_v)$ and $(g_w,p_w)$ are (sufficiently close) simple pairs,  there exists a curve $ k_{vw}$ on $\widetilde{Y}$ that is a concatenation of a (horizontal) leaf segment $ k_{vw}^1$ of $g_v$ from $p_v$ to some point $q_v$, a (vertical) leaf segment $\xi_{v}^w$ of $\beta$ from $q_v$ to some point $q_w\in g_w$, and a (horizontal) leaf segment $k_{vw}^2$ of $g_w$ from $q_w$ to $p_w$. Hence, by definition, we have
  \begin{eqnarray*}
    &&\sigma_{\lambda}(Y_r)((g_v,p_v),(g_w,p_w))\\
    &=&\sigma_\lambda(Y_r)(k_{vw})
   \\ 
  &=&  2\sqrt{r}  \sum_{j=1}^2\epsilon_j\cdot i(k_{vw}^j,\beta)
  \end{eqnarray*}
  where $\epsilon_j\in\{\pm1\}$ depends on the oriented distance along $g_v$ from $p_v$ to $q_v$ (resp. along $g_w$ from $p_w$ to $q_w$).
  
  Now, by Theorem \ref{thm:Minsky:traintrack} and Lemma \ref{lem:TT:approximation}, we could exponentially closely approximate
  $(g_v,p_v)$ and $(g_w,p_w)$ by a family of pointed $\Phi_t$-horizontal leaves $(g_v(t),p_v(t))$ and $(g_w(t),p_w(t))$, and exponentially closely approximate 
  $k_{vw}$ by the image of a family of $\Phi_t$-polygonal curves $k_{vw}(t)$ on $\widetilde{Y}$ that is a concatenation of a horizontal leaf segment $k_{vw}^1(t)$ of $g_v(t)$ from $p_v(t)$ to some point $q_v(t)$, a vertical leaf segment $\xi_{v}^w(t)$  from $q_v(t)$ to some point $g_w(t)$, and a horizontal leaf segment $k_{vw}^2(t)$ of $g_w(t)$ from $q_w(t)$ to $p_w(t)$; here we may choose $k_{vw}(t)$ in the complement of $\mathscr{P}_{R_t}(X_t)$ as $\lambda$ is in the complement of the image of $\mathscr{P}_{R_t}(X_t)$. \color{black} Then, from Theorem~\ref{thm:Minsky:traintrack}, because the leaves of $\lambda$ are well-approximated by images of $\Phi_t$-horizontal arcs, with distances along $\lambda$ well-approximated by $4\Phi_t$-horizontal lengths (outside of the polygonal region $\mathscr{P}_{R_t}(X_t)$), we may estimate that 
  \begin{eqnarray*}
    &&\sigma_{\lambda}(Y(t,r,0))((g_v,p_v),(g_w,p_w))\\
    &=&\sigma_\lambda(Y(t,r,0)(k_{vw}) 
   \\ &=&2\sqrt{r}  \sum_{j=1}^2\epsilon_j\cdot {i(k_{vw}^j(t),\mathrm{Vert(\Phi_t)})}(1+O(\exp(-brR_t)))
  \\
  &\to &  2\sqrt{r}  \sum_{j=1}^2\epsilon_j\cdot {i(k_{vw}^j,\beta)} \qquad (\text{as }t\to\infty)\\
  &=&\sigma_{\lambda}(Y_r)((g_v,p_v),(g_w,p_w))
  \end{eqnarray*}
   where  $b$ in the third line above is the constant from Theorem \ref{thm:Minsky:traintrack} and  $\epsilon_j\in\{\pm1\}$ depends on the oriented distance along $g_v$ from $p_v$ to $q_v$ (resp. along $g_w$ from $p_w$ to $q_w$).
   The final equality is by substituting the results of the previous displayed computation.

We now consider the effect of a non-zero \enquote{rotation factor} $s\neq 0$. 
By Lemma \ref{lem:TT:approximation}, for $t$ sufficiently large, the image of $X_t\backslash\mathscr{P}_{2R_t}(\Phi_t)$ under the harmonic map $X_t\to Y(t,r,s)$ is a thickened train track 
that carries $\lambda$ and is contained in the $\epsilon$ neighbourhood of $\lambda$ (here we use the fact $\mathrm{Hor}(\Phi_t)=t\lambda$).
  Notice that for any fixed $\mathbf{s}>0$ (and as usual, for $t$ sufficiently large), the complement $X_t\backslash\mathscr{P}_{R_t}(e^{is/2t}\Phi_t)$ contains  $X_t\backslash\mathscr{P}_{2R_t}(\Phi_t)$ for $-\mathbf{s}\leq s\leq \mathbf{s}$. 
  Hence, the image of $X_t\backslash\mathscr{P}_{R_t}(e^{is/2t}\Phi_t)$  under the harmonic map $X_t\to Y(t,r,s)$, which is contained in the $\epsilon$ neighbourhood of the geodesic lamination $\mathrm{Hor}(e^{is/2t}\Phi_t)$ by Theorem \ref{thm:Minsky:traintrack},   is also a thickened train track carrying $\lambda$. Therefore,
  \begin{eqnarray*}
   && \sigma_{\lambda}(Y(t,r,s))((g_v,p_v),(g_w,p_w)))\\&=&
  \sigma_\lambda(Y(t,r,s))(k_{vw})\\  &=&2\sqrt{r}\left(\sum_{j=1}^2\epsilon_j\cdot i(k_{vw}^j(t),\Vert(\Phi_t))\cos\frac{s}{2t} +  i(k_{vw}(t),\mathrm{Hor}(\Phi_t))\sin\frac{s}{2t}\right)(1+O(\exp(-brR_t)))\\
  &\to&  2\sqrt{r}  \sum_{j=1}^2\epsilon_j\cdot i(k_{vw}^j,{\beta})+s\sqrt{r} i(k_{vw},{\lambda}) \quad (\text{as }t\to\infty)\\&=& \sigma_\lambda(Y_r)((g_v,p_v),(g_w,p_w)) +s\sqrt{r} i(k_{vw},{\lambda}).
  \end{eqnarray*}
  This means that $$\sigma_\lambda(Y(t,r,s)\to \sigma_\lambda(\mathcal{E}^s_{\lambda}( Y_r))$$ as $t\to\infty$. The finiteness of the pairs of oriented geodesics then implies that the convergence is locally uniform in $(r,s)\in[\mathbf{r},\mathbf{r}']\times [-\mathbf{s},\mathbf{s}]$.
 \end{proof}

 \begin{proof}[Proof of Theorem \ref{thm:Hopfdisk:limit}] The theorem now follows from Proposition \ref{prop:cocycle:injective},  Lemma \ref{lem:shape:convergence} and Lemma \ref{lem:shear:convergence}.
 \end{proof}

\section{Harmonic-stretch lines between hyperbolic surfaces}
\label{sec:stretchray:uniqueness}

 The goal of this section is, in outline, to find for every ordered pair $(Y,Z)$ of distinct hyperbolic surfaces in $\T(S)$, a unique harmonic stretch line proceeding from $Y$ to $Z$, i.e. a unique Thurston geodesic proceeding from $Y$ to $Z$  determined only by $Y$ and $Z$, together with an extra side condition resulting from a requirement of minimizing a variational quantity. The main result is stated in Theorem~\ref{thm:unique:geodesic} which characterizes those stretch lines in terms of \enquote{admissible triples}.  We give that definition in the first subsection before proceeding to state and prove Theorem~\ref{thm:unique:geodesic} in the following subsections.
  %===========================
 
%---------------------------
%===========================
 \subsection{Harmonic stretch lines and admissible triples }\label{subsec:admissible:triples}
 
 Recall that (Theorem \ref{thm:generalized:stretchmap}) a piecewise harmonic stretch line is constructed  from a  closed hyperbolic surface $U$, a  geodesic lamination $\lambda$ on $U$, and a harmonic diffeomorphism $X\to U\setminus \lambda$  from some (possibly disconnected) punctured Riemann surface.
 
 \begin{definition} \label{defn:harmonic:stretch:line}
 We say that a piecewise harmonic stretch line is a {\it harmonic stretch line}, if it is a limit of a sequence of harmonic map rays. We refer to the restriction, of a harmonic stretch line to a half-infinite subsegment that begins at $Y$, a {\it harmonic stretch ray} from $Y$. We adopt the notations $\HSL$ and $\HSR$, respectively, often decorated with their defining data.
  \end{definition}
  
 Let $Y$ and $Z$ be two fixed hyperbolic surfaces in $\T(S)$.  Let $\lambda$ be the maximally stretched chain-recurrent lamination from $Y$ to $Z$ (Section \ref{subsec:maximally:stretched:lamination}). Let $L$ be the least Lipschitz constant for Lipschitz maps from $Y$ to $Z$ in the prescribed homotopy class. We will write, for example, $L^{-2}Z$ to indicate the constant curvature metric obtained by scaling the hyperbolic metric $Z$ by the constant $L^{-2}$ (so that the line element is scaled by $L^{-1}$). 
 \begin{definition}\label{def:ad:triples}   A triple $(X,f,h)$ is said to be \textit{admissible} if the following hold:
  \begin{enumerate}[(i)]
       \item $X$ is a punctured Riemann surface (possibly disconnected);
       \item $f:X\to Y\backslash\lambda$ and $h:X\to L^{-2}Z\backslash\lambda$ are harmonic diffeomorphisms with meromorphic Hopf differentials satisfying $\Hopf(f)=\Hopf(h)$;
       \item  the piecewise harmonic stretch line defined by $f:X\to Y\setminus\lambda$ passes through $Z$;
       \item the piecewise harmonic stretch line defined by $f:X\to Y\setminus\lambda$ is a harmonic stretch line.
   \end{enumerate}
 \end{definition}
 \remark \label{rmk:hsr:YZ} 
 We emphasize that the final condition in the definition requires that (.cf Definition~\ref{defn:harmonic:stretch:line}) that the piecewise harmonic stretch line is a {\it limit} of  harmonic map rays.  Note the nuance in the final two conditions:  we require the piecewise harmonic stretch line to pass through $Y$ and $Z$, we require that piecewise harmonic stretch line to be a limit of harmonic map rays, but we do not require the harmonic map ray approximates implied in condition (iv) to pass through both $Y$ and $Z$. 
 
 For the existence of admissible triples, we refer forward to Subsection \ref{subsec:approximate}. 
 By item (b) of Theorem \ref{thm:generalized:stretchmap}, we see that $h\circ f^{-1}$ extends to a 1-Lipschitz homeomorphsim from $Y$ to $L^{-2}Z$.

 In section~\ref{sec:exponential}, we will extend this definition to the case where the terminal surface $Z$ is replaced by an $\R$-tree. (Naturally, we will have to lift all of the definitions to spaces on which $\pi_1(S)$ acts equivariantly by automorphisms.) See Definition~\ref{def:ad:triples:trees} for that analogous construction.
 
 \subsection{Harmonic Stretch lines between points in \tec space.}
 
  With these definitions in hand, we may now state our main goal for this section, a restatement of Theorem~\ref{thm:unique:geodesic}.

\begin{theorem}
\label{thm:unique:HSR}
    For any two distinct hyperbolic surfaces $Y, Z\in \T(S)$, there is a unique  harmonic stretch line proceeding from $Y$ through $Z$. 
  \end{theorem}

 The remainder of this section is devoted to proving Theorem \ref{thm:unique:HSR}.

 \subsection{Existence of  admissible triples}
\label{subsec:approximate}
Let $Y,Z\in\T(S)$ be two hyperbolic surfaces with the optimal Lipschitz constant $L:=\exp(d_{Th}(Y,Z))$.  Then for any $X\in\T(S)$, we have
\begin{equation}\label{eq:energy:YZ}
    E(X,Z)\leq L^2 E(X,Y).
\end{equation}
(This is because the energy-minimizing map from $X$ to $Z$ will have less energy than any candidate map, in particular the composition of the energy-minimizer from $X$ to $Y$ with an optimal ($L$-)Lipschitz map from $Y$ to $Z$: that composition would have energy bounded from above by $L^2 E(X,Y)$.)
For  any $0\leq r < L^{-2}$, consider the energy difference
\begin{eqnarray*}
F_r:\T(S)&\longrightarrow &\R \\
X&\longmapsto& E(X,Y)- E(X,rZ).
\end{eqnarray*}
By Equation (\ref{eq:energy:YZ}), we see that for any $r\in[0,L^{-2})$,
$$F_r(X)\geq E(X,Y)-r L^2 E(X,Y)= \left(1-r L^2\right) E(X,Y).$$

In particular, $F_r$ is positive and proper. Therefore, $F_r$ has at least one critical point. Let $X_r$ be such a critical point. Then $0=d F_r\left|\right._{X_r}=\Hopf(X_r,Y)-r\Hopf(X_r,Z)$ (Using quadratic differentials to represent the differential of the energy functional over $\T(S)$ is a classic result, see \cite{Jost,Tromba,Wolf1998,Wentworth} for example). Hence $\Hopf(X_r,Y)=r\Hopf(X_r,Z)$.
 \begin{theorem}[Tholozan \cite{Tholozan2017}]\label{thm:Tholozan}
 For any $r\in[0,L^{-2})$, there exists a unique Riemann surface $X_r$ in $\T(S)$ such that $\Hopf(X_r,Y)=r\Hopf(X_r,Z)$. Moreover, $X_r\to \infty$ as $r\to L^{-2}$.
 \end{theorem}
 That $X_r$  in Theorem \ref{thm:Tholozan} diverges follows from \cite[Corollary 2.3]{Tholozan2017} or Lemma \ref{lem:1:lip} in the current paper.
  
Using the family of harmonic map rays $\HR_{X_r,Y}$, we have the following existence result.
  \begin{proposition}\label{prop:ad:existence}
    Let $Y$  and $Z$ be two distinct hyperbolic surfaces in $\T(S)$. Then there exists a  harmonic stretch line passing through $Y$ to $Z$. Consequently, there exists an admissible triple of $(Y,Z)$.
  \end{proposition}
  \begin{proof}
    For any $r\in [0,L^{-2})$, let $X_r$ be the Riemann surface obtained from Theorem \ref{thm:Tholozan}. Consider the family of harmonic map rays $\HR_{r}$ which starts at $X_r$ and passes through $Y$ and $Z$. By Theorem \ref{thm:HR:limit:geodesic}, there exist a sequence $r_n\to L^{-2}$ such that $\HR_{r_n}$ converges to a Thurston geodesic, say $\HSR$. Then $\HSR$ is a harmonic stretch line passing through $Y$ to $Z$.  This proves the proposition.   \end{proof}

\subsection{Equivalent admissible triples}  
 We begin the proof of the uniqueness part of Theorem \ref{thm:unique:HSR} with the following definition. 
 \begin{definition}
   Two admissible triples $(X,f,h)$ and $(\hat X,\hat f, \hat h)$  of $(Y,Z)$ are said to be \textit{equivalent} if there exists a conformal map $\eta: {X}\to \hat X$ such that ${f}=\hat f\circ\eta$.
 \end{definition}

It might seem unnatural that the above definition does not involve a condition on $h$ and $\hat{h}$. However, if two admissible triples $(X,f,h)$ and $(\hat X,\hat f, \hat h)$  of $(Y,Z)$ are equivalent, then, because of the identity between Hopf differentials of $f$ and $h$ (resp, $\hat{f}$ and $\hat{h}$), we also have $h=\hat{h}\circ\eta$ by Theorem \ref{thm:parametrization}. Indeed, we see the effect of the definition as stated in the consequence described in the next lemma.
 \begin{lemma}\label{lem:ad:stretch}
 Let $(X,f,h)$ and $(\hat{X},\hat{f},\hat{h})$ be two equivalent admissible triples of $(Y,Z)$. Then the harmonic stretch lines defined by them coincide.
 \end{lemma}
 
 \begin{proof}
   Let $\HSL:(0,\infty)\to\T(S)$ be the harmonic stretch line defined by $(X,f,h)$ such that the harmonic map $f_t:X\to \HSL(t)\backslash \lambda$ has Hopf differential $\Hopf(f_t) =t \Hopf(f)$. In particular $\HSL(1)=Y$. Define $\widehat{\HSL}(t)$ and  $\hat{f}_t$ similarly. 
  
   The assumption that $(X,f,h)$ and $(\hat{X},\hat{f},\hat{h})$ are equivalent implies that the composition map $\eta:=\hat{f}^{-1}\circ f$ from $X$ to $\hat{X}$ is conformal. Hence $\hat{f}_t\circ \eta:X\to \widehat{\HSL}(t)\backslash\lambda
   $ is also harmonic.
    Moreover,
   \begin{eqnarray*}
     \Hopf(\hat{f}_t\circ \eta)&=&\eta^*( \Hopf(\hat{f}_t))=\eta^*(t \Hopf(\hat{f}))
     \\ &=&t \Hopf(\hat{f}\circ \eta)=
     t \Hopf(f) =    \Hopf(f_t).
   \end{eqnarray*}
   By Theorem \ref{thm:parametrization}, there is a unique crowned hyperbolic surface defined via a prescribed Hopf differential on $X$. In particular, $\widehat{\HSL}(t)=\HSL(t)$ holds for all $t$. 
   This completes the proof.
 \end{proof}
 
The following uniqueness result is key to the proof of Theorem \ref{thm:unique:HSR}.
   \begin{proposition}\label{prop:ad:uniqueness}
   Let $(Y,Z)\in\T(S)\times \T(S)$ be a pair of distinct hyperbolic surfaces.   Then all admissible triples of $(Y,Z)$ are equivalent. 
  \end{proposition} 
  
  We prove this proposition in subsection~\ref{subsec:uniqueness}.
 
 \subsection{Energy difference}
Note that both $f:X{\to } Y$ and $h:X{\to } L^{-2}Z $ have infinite energy.  Nevertheless, we are able to talk about the energy difference of admissible triples in the following sense. Let $(X,f,h)$ be an admissible triple of $(Y,Z)$.  Let $e(f)$ and $e(h)$ be respectively the energy densities of $f$ and $h$. It then follows from Lemma \ref{lem:1:lip} that we have the pointwise estimate $e(h)\leq e(f)$. Combined with  Lemma \ref{lem:energy:difference:infinity:2}, this implies that for any compact exhaustion  $\{\mathcal{K}_j\}$  of $X$, the limit  $\lim_{j\to\infty}E(f|_{\mathcal{K}_j})-
    E(h|_{\mathcal{K}_j})$ exists and is always a (finite) non-negative real number. Moreover, again by Lemma \ref{lem:energy:difference:infinity:2}, the limit  is independent of the choice of the compact exhaustion. Set $E(f)-E(h):= \lim_{j\to\infty}E(f\left|\right._{\mathcal{K}_j})-
    E(h\left|\right._{\mathcal{K}_j})$. It is clear that $E(f)-E(h)\geq 0.$

 \begin{definition}
   The energy difference $E(X,f,h)$ of an admissible triple $(X,f,h)$ is defined by setting $E(X,f,h):=E(f)-E(h)$.
 \end{definition} 
 
 As we will see, the energy difference plays a key role in establishing the uniqueness of admissible triples.

%===========================
%==========================

  \subsection{Uniqueness of admissible triples}
  \label{subsec:uniqueness}
  In this subsection, we shall prove Proposition \ref{prop:ad:uniqueness}.
 The idea is to show that all admissible triples of $(Y,Z)$ have the same energy difference (Lemma \ref{lem:ed:equal}) and that equality of energy difference implies that equivalence of admissible triples (Lemma \ref{lem:equal:enery}). The key is to compare an arbitrary admissible triple $(\hat{X},\hat{f},\hat{h})$ to our favourite admissible triple $(X,f,h)$ that comes from Proposition \ref{prop:ad:existence}.   For Lemma \ref{lem:ed:equal}, we use the continuity and properness of the energy function over $\T(S)$. Lemma \ref{lem:equal:enery} is more subtle. To compare the energy differences $E(f)-E(h)$ and $E(\hat{f})-E(\hat{h})$, both of which can be approximated by energy differences on  Minsky regions $\mathscr{P}_R$, we first use Tholozan's local computation about energy densities to conclude that an auxiliary energy difference  $E(f_{\mathscr{P}_R})-E(\hat h\circ \hat f^{-1}\circ f|_{\mathscr{P}_R})$ is at least $E(\hat f|_{\mathscr{P}_R})-E(\hat h|_{\mathscr{P}_R})$ with equality holding if and only if $\hat{f}^{-1}\circ f$ is conformal. We  
  then combine Minsky's estimate outside $\mathscr{P}_R$ and the assumption that our favourite admissible triple $(X,f,h)$ arises as a limit of minimizers of energy difference between \emph{closed surfaces} to conclude that $\lim\limits_{R\to\infty}(E(\hat h\circ \hat f^{-1}\circ f|_{\mathscr{P}_R})-E(h_{\mathscr{P}_R}))\geq0$ (Lemma \ref{lem:positivity}). All together, this proves Lemma \ref{lem:equal:enery}. We now fill in the details.  

  We begin by asserting a comparison for maps restricted to the Minsky region for our favourite map $h: X \to Z\setminus \lambda$. In that setting, that map $h$ has no more energy than that of a competitor defined on a different domain, under a (naturally defined) map that changes domains.

  \begin{lemma}\label{lem:positivity}
   Let $(X,f,h)$ be an admissible triple of $(Y,Z)$ obtained from Proposition \ref{prop:ad:existence}, and $(\hat X,\hat{f},\hat{h})$ an arbitrary admissible triple of $(Y,Z)$. Then $$\lim_{R\to\infty} E(\hat h\circ \hat f^{-1}\circ {f}\left|_{\mathscr{P}_R}\right.)-
   E({h}\left|_{\mathscr{P}_R}\right.)\geq0,$$
   where ${\mathscr{P}_R}$ is the Minsky polygonal region of $\Hopf(f)$.
  \end{lemma}

 \begin{proof}
 
 We define the elements used in the construction of $(X,f,h)$ (described in Proposition~\ref{prop:ad:existence}) as follows.
 
 Let $0<r_n<L^{-2}$ be a sequence in which $r_n$ converges to $L^{-2}$ as $n\to\infty$. Let $X_n\in\T(S)$ be the unique Riemann surface  satisfying (Theorem \ref{thm:Tholozan}):
\begin{equation*}
   E(X_n,Y)-r_nE(X_n, Z)= \min_{X\in\T(S)}
     (E(X,Y)-r_n E(X,Z)).
\end{equation*} 
 Let $f_n:X_n\to Y$ and $h_n:X_n\to r_n Z$ be the corresponding harmonic maps. Then $\Hopf(f_n)=\Hopf(h_n)$. We rewrite the above equation as:
 \begin{equation}
    \label{eq:Xn:minimum}
   E(f_n)-E(h_n)= \min_{X\in\T(S)}
     (E(X,Y)-r_n E(X,Z)).
\end{equation} 
Finally, we obtain $(X,f,h)$ from $(X_n, f_n, h_n)$ as in Proposition~\ref{prop:ad:existence}. Then from Remark~\ref{rmk:hsr:YZ}, we notice that  $\hat h\circ\hat f^{-1}:Y\setminus \lambda \to L^{-2} Z\setminus \lambda$ extends to a 1-Lipschitz map from $Y$ to $ L^{-2}Z$. 
  Let $\Phi_{n}=\Hopf(f_n)=\Hopf(h_n)$.
  Then using the common Hopf differentials to provide for a common Minsky domain on which to integrate, we find
    \begin{eqnarray}
    \label{eq:estimate:positivity}
  &&\int_{X_n\backslash\mathscr{P}_R(\Phi_n)} (e(\hat h\circ\hat f^{-1}\circ f_n)-e(h_n))dA_{\Phi_n}\\
  \nonumber
  &=&\int_{X_n\backslash\mathscr{P}_R(\Phi_n)} (e(\hat h\circ\hat f^{-1}\circ f_n)-2) +\int_{X_n\backslash\mathscr{P}_R(\Phi_n)} (2-e(h_n))dA_{\Phi_n} \\
  \nonumber
  &\leq&  \int_{X_n\backslash\mathscr{P}_R(\Phi_n)} (e( f_n)-2)dA_{\Phi_n}+  \int_{X_n\backslash\mathscr{P}_R(\Phi_n)}(2-e(h_n))dA_{\Phi_n} \\
  \nonumber
  &&\qquad\qquad \text{ (since $\hat{h}\circ \hat{f}^{-1}$ is 1-Lipschitz)}\\
  \nonumber
  &\leq&  \int_{X_n\backslash\mathscr{P}_R(\Phi_n)} |e( f_n)-2|dA_{\Phi_n}+  \int_{X_n\backslash\mathscr{P}_R(\Phi_n)}|2-e(h_n)|dA_{\Phi_n} 
  \\
  \nonumber
  &\leq&2C|\chi(S)|\exp({-R/2}). \qquad\qquad\text{        (by Lemma \ref{lem:energy:difference:infinity}) }
  \end{eqnarray}

  Since $h_n:X_{n}\to r_nZ$ is an (energy-minimizing) harmonic map between closed surfaces and that $\hat h\circ \hat f^{-1}\circ f_n:X_n\to L^{-2}Z$ and $h_n:X_n\to rZ$ are homotopic, it follows that 
      $rL^2E(\hat h\circ \hat f^{-1}\circ f_n)\geq E(h_n).$
  The assumption that  $0<r<L^{-2}$ then implies that 
    $
      E(\hat h\circ \hat f^{-1}\circ f_n)\geq E(h_n).
 $ Therefore, using that $h_n: X_n \to rZ$ minimizes energy between the {\it closed} surfaces $X_n$ and $rZ$, we find
  $$ \int_{X_n} (e(\hat h\circ \hat f^{-1}\circ f_n)-e(h_n) )dA_{\Phi_n}=
  E(\hat h\circ \hat f^{-1}\circ f_n)- E(h_n)\geq 0.$$
  Combined with (\ref{eq:estimate:positivity}), this implies that
  $$ \int_{\mathscr{P}_R(\Phi_n)} \left(e( \hat h\circ \hat f^{-1}\circ f_n)-e(h_n)\right)dA_{\Phi_n}\geq -2C|\chi(S)|\exp({-R/2}).$$
  Letting $n\to\infty$ gives
   $$ \int_{\mathscr{P}_R(\Phi)}\left(e( \hat h\circ \hat f^{-1}\circ f)-e(h) \right)dA_{\Phi}\geq -2C|\chi(S)|\exp({-R/2}),
   $$
   where $\Phi$ is the Hopf differential of both $f$ and $h$. 
   Hence 
   \begin{eqnarray*} &&\lim_{R\to\infty}E(\hat h\circ \hat f^{-1}\circ {f}\left|_{\mathscr{P}_R(\Phi)}\right.)-
   E({h}\left|_{\mathscr{P}_R(\Phi)}\right.)\\
   &=&\lim_{R\to\infty}\int_{\mathscr{P}_R(\Phi)}\left(e( \hat h\circ \hat f^{-1}\circ f)-e(h) \right)dA_{\Phi}\\
   &\geq& -\lim_{R\to\infty}2C|\chi(S)|\exp({-R/2})=0,
   \end{eqnarray*}
   which completes the proof.
  \end{proof}
    
    We next observe that the energy difference for a solution to the problem of Proposition~\ref{prop:ad:existence} will have at least the energy difference as a solution to a similar problem with a less stringent constraint, i.e. that does not require the approximating lines to pass through exactly $Y$ and $Z$.

 \begin{lemma}\label{lem:equal:enery}
   Let $(X,f,h)$ be an admissible triple of $(Y,Z)$ obtained from Proposition \ref{prop:ad:existence}, and $(\hat X,\hat{f},\hat{h})$ an arbitrary admissible triple of $(Y,Z)$.  Then $$E(f)-E(h)\geq E(\hat{f})-E(\hat{h}),$$ 
   where the equality holds if and only if $\hat{f}^{-1}\circ f:X\to \hat{X}$ is conformal.
 \end{lemma}
 \begin{proof}
 By the definition of admissible triples, we have
     \begin{equation}\label{eq:equal:hopf}
       \Hopf(f)=\Hopf(h),~~
       \Hopf(\hat{f})=\Hopf(\hat{h}).
     \end{equation}
     Let $\sigma|dz|^2$ be a conformal metric on $\widehat{X}$,  let $\eta:=\hat{f}^{-1}\circ f:X\to \hat X$, let $\mathbf{e}:=e(\eta^{-1})$, and  let $\Psi:=\Hopf(\eta^{-1})$ (not necessary holomorphic).
    Now applying Tholozan's argument (\cite[Lemma 2.5]{Tholozan2017}) to the harmonic maps
    \begin{equation*}
        f:X\to Y,\qquad  h:X\to L^{-2}Z,
    \end{equation*}
    and
    \begin{equation*}
         \hat{f}:\hat{X}\to Y,\qquad \hat{h}:\hat{X} \to L^{-2}Z,
    \end{equation*}
    we see that for any compact subset $\widehat{K}$ of $\widehat{X}$, 
    \begin{eqnarray}\label{eq:decreasing}
    \nonumber
     && E(f|_{\eta^{-1}(\widehat{K})})-E(\hat h\circ \eta|_{\eta^{-1}(\widehat{K})})\\\nonumber
     &=& 
      \int_{\widehat{K}}
      \frac{1}{\sqrt{1-4|\Psi|^2/\sigma^2 \mathbf{e}^2}}(e(\hat{f}|_{\widehat{K}})-e(\hat{h}|_{\widehat{K}}))d\sigma
      \\ &\geq&
    E(\hat f|_{\widehat{K}})-E(\hat h|_{\widehat{K}}),
    \end{eqnarray}
      where the  equality holds if and only if $\Psi\equiv0$, i.e.  $\eta$ is conformal. 

     Let $\mathscr P_R$ be the Minsky's polygonal region of the Hopf differential of $f$. Then
     \begin{eqnarray*}
       && E({f})-E({h})
       \\
       &=& \lim_{R\to\infty}\left(  E({f}|_{\mathscr P_R})-E({h}|_{\mathscr P_R})\right)    \qquad\qquad(\text{by definition})\\
       &=&\lim_{R\to\infty}\left( E({f}\left|_{\mathscr P_R}\right.)- E(\hat h\circ \eta\left|_{\mathscr P_R}\right.)\right)\\
       &&~~~~+~\lim_{R\to\infty}\left(E(\hat h\circ \eta\left|_{\mathscr P_R}\right.)-E({h}|_{\mathscr P_R})\right)\\
       &\geq& \lim_{R\to\infty}  \left( E({f}\left|_{\mathscr P_R}\right.)- E(\hat h\circ \eta\left|_{\mathscr P_R}\right.)\right) ~~~~(\text{by Lemma \ref{lem:positivity}})\\
       &\geq& \lim_{R\to\infty}\left(
    E(\hat f|_{\eta(\mathscr P_R)})-E(\hat h|_{\eta(\mathscr P_R)})\right)\qquad(\text{by }(\ref{eq:decreasing}))\\
    &=& E(\hat{f})-E(\hat{h}),
    \qquad\qquad   (\text{by definition})
     \end{eqnarray*}
     where the equality holds if and only if $\eta$ is conformal. 
 \end{proof}
 
 Next, we use the energy difference minimization property of those admissible triples from Proposition \ref{prop:ad:existence} to show that all admissible triples have the same energy difference.
 \begin{lemma}
  \label{lem:ed:equal}
  Let $(X,f,h)$ be an admissible triple of $(Y,Z)$ from Proposition \ref{prop:ad:existence}, and $(\hat{X},\hat{f},\hat{h})$ an arbitrary admissible triple of $(Y,Z)$. Then 
  \begin{equation}
      \label{eq:ed:equal}
       E(\hat{f})-E(\hat{h})= E(f)-E(h).
  \end{equation}
 \end{lemma}
 \begin{proof}
 By definitions of admissible triples and harmonic stretch lines, there exists a  sequence of harmonic map rays $\HR_n$ converging to the (piecewise) harmonic stretch line determined by $(\hat{X},\hat{f},\hat{h})$. Let  $ X_n\in\T(S)$ be the initial point of $\HR_n$. 
 
 Let $0<r_m<L^{-2}$ be a sequence which converges to $L^{-2}$ as $m\to\infty$, such that 
$f_{r_m}:X_{r_m}\to Y$ and $h_{r_m}:X_{r_m}\to r_mZ$ converges to $f$ and $h$ respectively, where $X_{r_m}$ is the Riemann surface obtained from Theorem \ref{thm:Tholozan}. By Lemma \ref{lem:1:lip}, $h_{r_m}\circ (f_{r_m})^{-1}: Y\to r_m Z$ is 1-Lipschitz. Combined with Lemma \ref{lem:energy:difference:infinity}, this implies that 
 $\lim_{m\to\infty}(E(X_{r_m},Y)-r_m E(X_{r_m},Z))=E(f)-E(h)$. 
 Let $\epsilon>0$ be an arbitrary positive real number.
 Since $E(X_{r_m},Y)-r_m E(X_{r_m},Z)=\min_{X\in\T(S)}
     (E(X,Y)-r_m E(X,Z))$, 
  we may assume that, up to a subsequence, 
\begin{equation}
    \label{eq:1:over:m}
   \left|\min_{X\in\T(S)}
     (E(X,Y)-r_m E(X,Z))-(E(f)-E(h))\right|<{\epsilon}.
\end{equation}
 Notice that for any fixed $0<r<L^{-2}$, 
 \begin{equation*}
     \lim_{\tiny{Y'\to Y \atop Z'\to Z}}\min_{X\in\T(S)}
     (E(X,Y')-r E(X,Z'))=\min_{X\in\T(S)}
     (E(X,Y)-r E(X,Z)).
 \end{equation*}
 For each $m$, choose $X_{n_m}$, $\HR_{n_m}$,  $Y_{n_m}\in \HR_{n_m}$, and $Z_{n_m}\in \HR_{n_m}$ such that
 \begin{enumerate}[(a)]
  \item $Y_{n_m}\to Y$ and $Z_{n_m}\to Z$ as $m\to\infty$;
 \item $\Hopf(X_{n_m},Z_{n_m})=\frac{1}{r_m}\Hopf(X_{n_m},Y_{n_m})$, i.e. $X_{n_m}$ realizes the minimum $$\min_{X\in\T(S)}
     (E(X,Y_{n_m})-r_m E(X,Z_{n_m}));$$
     \item $r_m<\Lip(Y_{n_m},Z_{n_m})^{-2};$
    \item $$\left|\min_{X\in\T(S)}
     (E(X,Y_{n_m})-r_m E(X,Z_{n_m}))-\min_{X\in\T(S)}
     (E(X,Y)-r_m E(X,Z))\right|<{\epsilon}.$$
 \end{enumerate}
 Combining the items (c), (d)  with (\ref{eq:1:over:m}), we see that 
 \begin{equation}
     \label{eq:YZnm}
     \left|(E(X_{n_m},Y_{n_m})-r_m E(X_{n_m},Z_{n_m}))-(E(f)- E(h))\right|<{2\epsilon,}
 \end{equation}
 where we use the assumption (b) that $X_{n_m}$ solves the energy difference minimization problem.

  Let $\hat f_{n_m}:X_{n_m}\to Y_{n_m}$ and $\hat h_{n_m}:X_{n_m}\to r_mZ_{n_m}$ be the corresponding harmonic maps.  We may then rewrite equation (\ref{eq:YZnm}) as 
   \begin{equation}
     \label{eq:YZnm2}
     \left|(E(\hat f_{n_m})- E(\hat h_{n_m}))-(E(f)- E(h))\right|<{2\epsilon.}
 \end{equation}
Now, by Lemma \ref{lem:1:lip}, $\hat h_{n_m}\circ (\hat f_{n_m})^{-1}$ is 1-Lipschitz.    Combining this  with Lemma \ref{lem:energy:difference:infinity}, we see that
\begin{equation*}
   \left|\lim_{m\to\infty} \left(E(\hat f_{n_m})-E(\hat h_{n_m})\right)-\left(E(\hat{f})-E(\hat{h})\right)\right|\leq 2\epsilon.
 \end{equation*}
 The arbitrariness of $\epsilon$ then implies that 
 \begin{equation*}
   \lim_{m\to\infty} \left(E(\hat f_{n_m})-E(\hat h_{n_m})\right)=E(\hat{f})-E(\hat{h}).
 \end{equation*}
 Combined with (\ref{eq:YZnm2}), this yields 
 $$E(\hat{f})-E(\hat{h})=E(f)-E(h). $$
 \end{proof}

  \begin{proof}[Proof of Proposition \ref{prop:ad:uniqueness}]
  The proposition follows directly from Lemma \ref{lem:ed:equal} and Lemma \ref{lem:equal:enery}.
  \end{proof}

  \begin{proof}[Proof of Theorem \ref{thm:unique:HSR}]
 The theorem now follows directly from Proposition \ref{prop:ad:uniqueness} and Proposition  \ref{prop:ad:existence}.
 \end{proof}

\subsection{Continuity of harmonic stretch lines}\label{subsec:continuity}
By Theorem \ref{thm:unique:HSR} we see that for any two distinct hyperbolic surfaces $Y, Z\in \T(S)$, there is a unique  harmonic stretch line $\HSR_{Y,Z}$ proceeding from $Y$ through $Z$. A natural question one might ask is how stretch lines depend on the prescribed points $Y$ and $Z$. In this regard, we have the following continuity result.
\begin{proposition}
\label{prop:continuity}
Let $Y$ and $Z$ be two distinct points in $\T(S)$. Assume that $Y_n,Z_n\in\T(S)$ such that $\lim\limits_{n\to\infty}Y_n=Y$ and $\lim\limits_{n\to\infty}Z_n=Z$. Then 
$\HSR_{Y_n,Z_n}$ converges to $\HSR_{Y,Z}$ locally uniformly as $n\to\infty$.
 \end{proposition}
  \begin{proof}
 By the definition of harmonic stretch lines, we see that for every fixed $n$, there exists a sequence of harmonic map rays $\HR_{n,m}$ which converges to $\HSR_{Y_n,Z_n}$ locally uniformly as $m\to\infty$. Let $$r_n:=\max\{d_{Th}(Y,Y_n),d_{Th}(Y_n,Y),d_{Th}(Z,Z_n),d_{Th}(Z_n,Z)\}.$$ By assumption, we have $\lim\limits_{n\to\infty}r_n=0$.
 Now for each $n$, we choose a harmonic map ray $\HR_{n,m_n}$ whose $2r_n$-neighbourhood contains both $Y$ and $Z$.   This implies that the limit ray of any convergent subsequence of  $\{\HR_{n,m_n}\}_{n\geq1}$ proceeds from $Y$ to $Z$. Moreover, by Lemma \ref{lem:1:lip} and Definition \ref{defn:harmonic:stretch:line}, any subsequential limit is a harmonic stretch line proceeding from $Y$ to $Z$.  It then follows from the uniqueness part of Theorem \ref{thm:unique:HSR}  that $\HR_{n,m_n}$ converges to $\HSR_{Y,Z}$ locally uniformly as $n\to\infty$. This implies that $\HSR_{Y_n,Z_n}$ converges to $\HSR_{Y,Z}$ locally uniformly as $n\to\infty$.
 \end{proof}

\subsection{A characterization of admissible triples}
\label{subsec:characterizat} 

 We end this section with a characterization of harmonic stretch lines among piecewise harmonic stretch lines, in terms of their defining harmonic maps. We say that a triple  $(\hat{X},\hat{f},\hat{h})$ is \emph{quasi-admissible} if it satisfies  the  assumptions (i), (ii), (iii) in  the definition of admissible triples.  
 \begin{proposition}
 \label{prop:ad:characterization}
  A quasi-admissible triple of $(Y,Z)$ is admissible if and only if  it has the maximal energy difference among all quasi-admissible triples of $(Y,Z)$.
 \end{proposition}
 \begin{proof}
 Let $(X,f,h)$ be an admissible triple of $(Y,Z)$ obtained from Proposition \ref{prop:ad:existence}, and $(\hat X,\hat{f},\hat{h})$ an arbitrary admissible triple of $(Y,Z)$. Looking at the proof of
  Lemma \ref{lem:positivity} and Lemma \ref{lem:equal:enery}, we see that the only assumptions about $(\hat X,\hat{f},\hat{h})$ that we use are that $(\hat X,\hat{f},\hat{h})$  is a quasi-admissible triple. In particular, Lemma \ref{lem:positivity} and Lemma \ref{lem:equal:enery} still hold if $(X,f,h)$ is an admissible triple obtained from Proposition \ref{prop:ad:existence} while $(\hat X,\hat{f},\hat{h})$ is a  quasi-admissible triple.
 \end{proof}

The computation at the end of the proof of Lemma~\ref{lem:equal:enery} and this final result shows that the energy difference of an admissible triple is at least the energy difference of a quasi-admissible triple, with equality only if the quasi-admissible triple is actually (or conformal to) the admissible triple. Returning to the definitions of quasi-admissible and admissible triples, in the end we see that (non-admissible) quasi-admissible triples do not arise as limits of harmonic maps (or else these maps would have the same energy difference as an admissible triple and hence be admissible, acquiring the condition (iv) that the corresponding piecewise harmonic stretch line is in fact a harmonic stretch line).  It may be worth noting that a distinction between the two criteria is that an admissible triple will necessarily have identical residues of the Hopf differentials at the paired punctures at a node (where the differential has a second order pole). 

Note that while the energy difference of an admissible triple is {\it at least as large} as that for any quasi-admissible triple, every admissible triple is equivalent to the those arising as the subsequential limit of harmonic maps to $Y$ and $Z$, for which the energy difference $E(\cdot, Y) - E(\cdot, rZ)$ is a {\it minimum} for each choice of $r$ and declines as $r$ tends to $L^{-2}$ (where $L$ is the optimal Lipschitz constant from $Y$ to $Z$).

\section{An \enquote{Exponential map} for Thurston's metric} \label{sec:exponential}

 The goal of this section is to consider the \enquote{visual boundary} of the Thurston metric (Theorem  \ref{thm:exponentialmap}), and define two distinct versions of a Thurston geodesic flow.
 
 Recall that a harmonic stretch line is a limit of a sequence of harmonic map rays.  By Theorem \ref{thm:unique:HSR}, for any two distinct hyperbolic surfaces $Y, Z\in \T(S)$, there is a unique  harmonic stretch line proceeding from $Y$ to $Z$. Here we extend this to the following, where a {\it harmonic stretch ray} is a  ray contained in some harmonic stretch line with the induced orientation (Recall that a harmonic stretch line admits a canonical orientation).
 \begin{theorem}
 \label{thm:uniqueness:HSR:ray}
 
 For any  hyperbolic surface $Y\in \T(S)$ and any projective measured lamination $[\beta]\in\PML(S)$, there is a unique  harmonic stretch ray from $Y$ which converges to $[\beta]$ in the Thurston compactification.
 \end{theorem}
 
  {\it Convention.} In the remainder of this section, to simplify the notation, we will denote the dual tree $(T_\beta,2d)$  by $T_\beta$. 
 
  %==================
  %===================
 \subsection{Optimal equivariant Lipschitz maps to trees} 
 \label{sec:stretchline:end}
 Let $\beta$ be a representative of $[\beta]\in\PMF(S)=\PML(S)$. Let $T_\beta$ be the dual tree of the lift to the universal cover of $\beta$.
  Let $L$ be the least Lipschitz constant of equivariant (surjective) maps from the universal cover $\widetilde{Y}$ to $T_\beta$. Then  
 \begin{equation}\label{eq:lower:lip}
     L\geq \sup_{\mu\in\ML(S)}\frac{2i(\beta,\mu)}{\ell_Y(\mu)},
 \end{equation} 
 where the convention on the metric on the (dual) tree $T_{\beta}$ implies that lengths on the tree are measured by the intersection number $2i(\beta,\mu)$.
 (Later we will see that the two quantities are actually equal.)  For each $0<t<L^{-1}$, consider the energy difference function $E(\cdot,Y)-t^2 E(\cdot,T_\beta)$ on $\T(S)$.   Tholozan's argument (cf. the proof of Lemma~\ref{lem:equal:enery}) gives a unique minimizer $X_t \in\T(S)$.  Moreover, for each $0<r<L^{-1}$, the vertical foliation of $\Hopf(X_r,Y)$ is exactly $r\beta$.  Letting $r\to L^{-1}$ and using Lemma \ref{lem:precompactness} , we obtain a convergent subsequence $f_n: X_{r_n}\to Y$ which converges to a harmonic diffeomorphism $f_\infty:X_\infty \to Y\setminus \lambda$ for some chain-recurrent lamination $\lambda$. Correspondingly, the push forward of the vertical foliation of $\Hopf(f_\infty)$ via $f_\infty$  extends to a measured foliation on $Y$ which is exactly $L^{-1}\beta$, viewed  as the limit of ${r_n}\beta$ as $n\to\infty$. Therefore, by Proposition \ref{prop:convergencestretchlines:PML}, the piecewise harmonic stretch line determined by $f_\infty$ converges to $[\beta]\in\PML(S)$.

 Recall that in terms of the natural coordinates of $\Hopf(f_\infty)$ ($\Hopf(f_\infty)=dz^2$),  the pullback of the hyperbolic metric on $Y\setminus\lambda$ via $f_\infty$ to $X_\infty$ is 
 $$  2(\cosh\G+1)dx^2+2(\cosh\G-1)dy^2$$
 where $\G=\log(1/|\nu|)$ and $\nu$  is the Beltrami differential of $f_\infty$ (see (\ref{eq:pullback:metric:2})). Let $\pi:\widetilde{X_\infty}\to L^{-1} T_\beta$ be the projection map along leaves of vertical foliations of $\Hopf(\widetilde{f_\infty})$.  By the definition of $T_\beta$,  the pullback metric of $T_\beta$  on $\widetilde{X_\infty}$ via $\pi$ is exactly $4dx^2$.
 Then the composition map $\pi\circ \widetilde{f_\infty}^{-1}:\widetilde{Y}\setminus\widetilde{\lambda}\to L^{-1}T_\beta$ is a Lipschitz map with (pointwise) Lipschitz constant $\sqrt{2/({\cosh\G+1})}<1$. Moreover, by Lemma \ref{lem:Minsky:decay}, the (pointwise) Lipschitz constant of $\pi\circ \widetilde{f_\infty}^{-1}$ along horizontal leaves tends to $1$ as the distance to zeros of $\Hopf(\widetilde{f_\infty})$ goes to infinity.  Therefore, $\pi\circ \widetilde{f_\infty}^{-1}:\widetilde{Y}\setminus\widetilde{\lambda}\to L^{-1} T_\beta$ extends to 
  a Lipschitz map $\widetilde{Y}\to T_\beta$ whose restriction to  the  geodesic lamination $\widetilde{\lambda}$ is an affine map of factor $L$ and which has (pointwise) Lipschitz constant strictly less than $L$ outside $\widetilde{\lambda}$. 
 Since $\lambda$, the projection of $\widetilde{\lambda}$ to $Y$, is chain-recurrent, there exists a sequence of multicurves $\mu_n$ whose support converges to $\lambda$ in the Hausdorff topology. Consequently,  
 $$\lim_{n\to\infty}\frac{2i(\beta,\mu_n)}{\ell_Y(\mu_n)}=L.$$ 
 Combining this with (\ref{eq:lower:lip}), we see that 
 $$ L=\sup_{\mu\in\ML(S)}\frac{2i(\beta,\mu)}{\ell_Y(\mu)}. $$ 
  Moreover, this implies that every optimal Lipschitz map from $\widetilde{Y} $ to $T_\beta$ would maximally stretch $\widetilde{\lambda}$.  We define $\widetilde{\lambda}$ to be the \textit{maximally stretched lamination} from $Y$ to $T_\beta$, denoted by $\Lambda(Y,T_\beta)$.  
  
  Using an argument analogous to the proof of Proposition \ref{prop:ad:existence}, we see that the family of harmonic map rays $\HR_{X_r,Y}$ determined by the harmonic  map $f_r:X_r\to Y$ contains a sequence which converges to a  harmonic stretch line, as $r$ approaches $L^{-2}$.  
  
 In the remainder of this section, we set $L:=\sup\limits_{\mu\in\MF(S)}\frac{2i(\beta,\mu)}{\ell_Y(\mu)}$.
 We summarize the above discussion in the first two statements in the following, while the third statement follows (as discussed) from Proposition \ref{prop:convergencestretchlines:PML}.
 \begin{proposition}\label{prop:optimallip:trees}
 Let $Y\in\T(S)$ and $\beta\in\MF(S)$. Let $T_\beta$ be the tree dual to the lift of $\beta$ to the universal cover $\widetilde{S}$. Then we have the following.
 \begin{itemize}
     \item There exists a  harmonic stretch ray $\HSR_{Y,\beta}$ which starts at $Y$ and which maximally stretches along the maximally stretched lamination $\Lambda(Y,T_\beta)$.
    \item There exists an equivariant  Lipschitz map $f:\widetilde{Y}\to T_\beta$ whose restriction to the maximally stretched lamination $\Lambda(Y,T_\beta)$ is an affine map of factor $L$ and which has (pointwise) Lipschitz constant strictly less than $L$ outside $\Lambda(Y,T)$. 
    \item The ray $\HSR_{Y,\beta}$ converges to $[\beta]$ in Thurston's compactification.
 \end{itemize} 
 \end{proposition}

 \remark In \cite{Tabak1985}, Tabak proved that the push forward of $\Hopf(f_t:X_t\to Y)$ is a subsonic $\rho$-holomorphic quadratic differential on $Y$ with $\rho:Y\times[0,1/4)\to \R$ defined by $\rho(y,r)=(1-4r)^{-1/2}$.  In \cite{Sibner1970}, Sibner-Sibner proved a nonliear Hodge-De Rham theorem which states that for any measured foliation $\beta\in\MF(S)$, there exists a threshold value $t_0>0$ such that for all $0<t<t_0$, there exists a unique subsonic $\rho$-holomorphic quadratic  differential whose vertical foliation is $t\beta$. The discussion above gives a different proof of Sibner-Sibner's result for the very specific $\rho$ defined above  and  describes explicitly $L^{-1}$ as the threshold value of $\beta$.

 %=====================
 %=====================
 \subsection{Admissible triples of $(Y,T_\beta)$}
 Having established the maximally stretched lamination $\widetilde{\lambda}$ from $\widetilde{Y}$ to $T_\beta$, we are now in a position to consider admissible triples for $\widetilde{Y}$ and $T_\beta$, in the same way as in Section \ref{subsec:admissible:triples}.  Combining the construction of harmonic stretch lines and the discussion in Section \ref{sec:stretchline:end}, we see that every harmonic stretch ray which starts at $Y$ and which converges to $[\beta]\in\PML(S)$ in Thurston's compactification gives an optimal equivariant Lipschitz map $f:\widetilde{Y}\to T_\beta$ whose restriction to  $\widetilde{\lambda}$ is an affine map of factor $L$ and which has (pointwise) Lipschitz constant strictly less than $L$ outside $\widetilde{\lambda}$. 
 
 \begin{definition}\label{def:ad:triples:trees}
  A triple $(X,\tilde{f},\tilde{h})$ is said to be \textit{admissible} if the following hold:
  \begin{enumerate}[(i)]
       \item $X$ is a punctured Riemann surface (possibly disconnected);
       \item  $\tilde{f}:\widetilde{X}\to \widetilde{Y}\backslash\widetilde{\lambda}$  is an equivariant harmonic diffeomorphism   and $\tilde{h}:\widetilde{X}\to L^{-1} T_\beta$ is an equivariant harmonic map with meromorphic Hopf differentials satisfying $\Hopf(\tilde{f})=\Hopf(\tilde{h})$;
       \item the piecewise harmonic stretch line defined by the quotient map $f:X\to Y\setminus\lambda$ of $\tilde{f}:\widetilde{X}\to \widetilde{Y}\backslash\widetilde{\lambda}$ is a harmonic stretch line. 
    \end{enumerate}
\end{definition}
 \remark By Proposition \ref{prop:convergencestretchlines:PML}, item (ii) implies  that the harmonic stretch line defined by $f:X\to Y\setminus \lambda$ converges to $[\beta]\in\PML(S)$. Moreover, the composition map $\tilde{h}\circ \tilde{f}^{-1}$ extends to an equivariant 1-Lipschitz  map from $\widetilde{Y}$ to $T_\beta$. 
 
  \begin{definition}
   Two admissible triples $(\widetilde{X},\tilde{f},\tilde{h})$ and $(\widetilde{X}',\tilde{f}',\tilde{h}')$  of $(Y,T_\beta)$ are said to be \textit{equivalent} if there exists a conformal map $\eta: \widetilde{X}\to \widetilde{X}'$ such that $\tilde{f}=\tilde{f}'\circ\eta$. 
 \end{definition}
  By Proposition \ref{prop:optimallip:trees}, there exists at least one admissible triple. Applying the same  argument as in Section \ref{subsec:uniqueness}, we have
 \begin{proposition}\label{prop:ad:uniqueness:tree}
  Let $Y\in\T(S)$ and $\beta\in\MF(S)$. Let $T_\beta$ be the tree dual to the lift of $\beta$ to the universal cover $\widetilde{S}$. Then all admissible triples of $(\widetilde{Y},T_\beta)$
  are equivalent.
  \end{proposition} 
 \begin{proof}
 The proof is similar to that of Proposition \ref{prop:ad:uniqueness}, with $(Y,Z)$ replaced by $(\widetilde{Y},T_\beta)$.
 \end{proof}
 
  \begin{proof}[Proof of Theorem \ref{thm:uniqueness:HSR:ray}]
 The existence part follows from Proposition \ref{prop:optimallip:trees}. The uniqueness part follows from adapting the proof of Lemma \ref{lem:ad:stretch} and Proposition \ref{prop:ad:uniqueness:tree}.
 \end{proof}

  \begin{proof}[Proof of Theorem \ref{thm:exponentialmap}]
   The first part follows from Theorem \ref{thm:uniqueness:HSR:ray}.
   
   For the second part, we need to show that these harmonic stretch rays are either disjoint (away from their origin $Y$) or coincide, and we must also show that the union of the harmonic stretch rays covers the \tec space. That the union of the harmonic stretch rays covers $\T(S)$ follows basically from Theorem \ref{thm:unique:geodesic}: 
   given a point $Z \in \T(S)$, we take the harmonic stretch segment from $Y$ to $Z$ and extend it to a proper ray in $\T(S)$ by scaling the Hopf differential on the corresponding domain $X_{\infty}$ (.cf as found via an admissible triple). This harmonic stretch ray converges to a unique point on $\PML(S)$ by Proposition~\ref{prop:convergencestretchlines:PML}, and hence is a leaf of the foliation. That the harmonic stretch rays from $Y$ either coincide or are disjoint away from $Y$ is the content of Theorem~\ref{thm:unique:geodesic}, since if two harmonic stretch lines from $Y$ intersected at a point $Z \in \T(S)$, they would coincide on the harmonic stretch segment $[Y,Z]$ from $Y$ to $Z$  and hence extend beyond $Z$ identically.
   
   Finally, the third statement, that the harmonic stretch rays terminating at a point $[\eta] \in \PML(s)$ also foliate if we let the initial point $Y$ vary in $\T(S)$, follows easily from the disjointness and surjectivity arguments of the previous paragraph.
  \end{proof}

  %%%%======================
  %%%%======================
 \subsection{\enquote{Exponential map} rays}\label{sec:cont:exp}

 \begin{definition}
  Given $Y\in\T(S)$ and $[\beta]\in \PMF(S)$, the harmonic map ray that starts at $Y$ and converges to $[\beta]$ in the Thurston compactification (.cf Theorem \ref{thm:uniqueness:HSR:ray}) is called an \emph{exponential map ray}, denoted by $\ESR_{Y,[\beta]}$.  
 \end{definition}

\begin{remark}
The standard usage of the term \enquote{exponential map} in Riemannian geometry refers to a map from the tangent space at a base point to a neighborhood of that point which solves the initial value problem for geodesics with the data from the tangent space at that base point. The result is a ray structure on that neighborhood. Here we use the term \enquote{exponential map} in terms of that resulting ray structure, but not in the sense of a map from rays in a tangent space to arcs from a base point. 
\end{remark}
 
 Using the same argument as in the proof of Proposition \ref{prop:continuity}, we have  the following:
  \begin{proposition}
\label{prop:continuity2}
Let $Y\in\T(S)$ and $[\eta]\in\PML(S)$. Assume that $Y_n\in\T(S)$ and $[\eta_n]\in\PML(S)$ are such that $\lim\limits_{n\to\infty}Y_n=Y$ and $\lim\limits_{n\to\infty}[\eta_n]=[\eta]$. Then the exponential map ray 
$\ESR_{Y_n,[\eta_n]}$ converges to the exponential map ray $\ESR_{Y,[\eta]}$ locally uniformly as $n\to\infty$.
 \end{proposition}
 
 As a direct consequence, we obtain the following analog of Theorem \ref{thm:maximally:stretched:lamination}:
 \begin{proposition}
 \label{prop:maximally:stretched:lamination}
 Let $Y\in\T(S)$ and $\eta\in\ML(S)$. Assume that $Y_n\in\T(S)$ and $\eta_n\in\ML(S)$ such that $\lim\limits_{n\to\infty}Y_n=Y$ and $\lim\limits_{n\to\infty}\eta_n=\eta$.  Let $\Lambda(Y_n,\eta_n)$  (resp. $\Lambda(Y,\eta)$) be the maximally stretched lamination from $Y_n$ to $T_{\eta_n}$ (resp. from $Y$ to $T_\eta$). Then $\Lambda(Y,T_\eta)$ contains any geodesic lamination in the limit set (in the Hausdorff topology) of $\Lambda(Y_n,T_{\eta_n})$.
 \end{proposition}
 \begin{proof}
 Let $Z_n$ be a point in the exponential map ray $\ESR_{Y_n,[\eta_n]}$ such that $d_{Th}(Y_n,Z_n)=1$ and that $\ESR_{Y_n,[\eta_n]}$  proceeds from $Y_n$ to $Z_n$.  Let $Z$ be a point in the exponential map ray $\ESR_{Y,[\eta]}$ such that $d_{Th}(Y,Z)=1$ and that $\ESR_{Y,[\eta]}$  proceeds from $Y$ to $Z$. By Proposition \ref{prop:continuity2}, we see that $Z_n$ converges to $Z$ as $n\to\infty$.  Notice that $\Lambda(Y_n,T_{\eta_n})=\Lambda(Y_n,Z_n)$ and $\Lambda(Y,T_\eta)=\Lambda(Y,Z)$. The proposition now follows from Theorem \ref{thm:maximally:stretched:lamination}.
 \end{proof}

  \subsection{A comment on the infinitesimal exponential map}
We prove
 \begin{proposition} \label{prop:existence-exponential-lines}
 	Let $v \in T_Y\T(S)$ be a tangent vector to \tec space at $Y\in \T(S)$. Then there is a harmonic stretch line tangent to $v$ through $Y$.
 \end{proposition}
 
 The reader may recall that we had referred to this result in Remark~\ref{rmk:exp}.
 
 \begin{proof}
 
 	It is well-known (see, e.g. \cite{SchoenYau1979} or \cite{Wolf1998}) the the total energy function $E(\cdot)= E(\cdot, Y): \T(S) \to \R$, which records the total energy $E(X,Y)$ from a surface $X$ to the surface $Y$, is a proper function on the \tec space $\T(S)$; this function has a unique critical point at $Y\in \T(S)$, where the Hessian $\Hess(E)$ is positive definite.  Indeed, Tromba \cite{Tromba87}, \cite[Theorem 3.1.3]{Tromba} shows that $\Hess E$ at $Y\in \T(S)$ is a multiple of the Weil-Petersson metric on $\T(S)$, say $\Hess E = c_0 \Hess d_{WP}^2$ for some (universal) $c_0>0$, and thus $E = c_0 d_{WP}^2 + o(d_{WP}^2)$ in a small neighborhood of $Y \in \T(S)$, where $d_{WP}$ is the Weil-Petersson metric.
 	
 	We consider a small \enquote{level sphere} $S(\epsilon)$ in $\T(S)$ on which the Energy $E(X,Y)$ has value $E(Y,Y) + \epsilon$. There is a standard pairing between tangent and cotangent spaces to \tec space, given by integration of Beltrami differentials and holomorphic quadratic differentials; in terms of that pairing, the Hopf differential $\Hopf(X,Y)$ then has kernel tangent to $S(\epsilon)$, as $dE_{|X}$ is a multiple of the Hopf differential $\Hopf(X,Y)$. We focus on the harmonic map rays from such points $X \in S(\epsilon)$ through $Y$.
 	
 These rays are of course tangent to (the duals of) the Hopf differentials through $X$, but there is another distinguished curve from $X$ and passing through $Y$: this is the Weil-Petersson geodesic from $X$ to $Y$.  As $S(\epsilon)$ links $Y$ in $\T(s)$, the tangent vectors to these geodesics fill the unit tangent sphere $T^1_Y\T(S)$ at $Y\in \T(s)$ to the \tec space $\T(S)$: this is because the exponential map from $Y$ surjects onto any hypersurface linking $Y$. Indeed, because $c_0 d_{WP}^2$ and $E$ agree to order $o(\epsilon^2)$, we see that $S(\epsilon)$ is Weil-Petersson convex for small $\epsilon$, and hence the map from $S(\epsilon)$ to $T^1_Y\T(S)$ which records tangent vectors to the Weil-Petersson   geodesics to $Y$ has degree one.
 
 We wish to relate the harmonic map rays to the Weil-Petersson  geodesics: here we recall \cite{Wolf1989} that, analogously to the Tromba result just described, the harmonic map rays from $X$ are  Weil-Petersson  geodesics at $X$ so that the two rays agree in $C^1$ up to an error of $o(\epsilon)$.  Thus the harmonic map rays from $X \in S(\epsilon)$ to $Y$ agree in $C^1$ with the  Weil-Petersson   geodesics from $X \in S(\epsilon)$ to $Y$ to an error of $o(\epsilon)$. Thus the differential of the map from the normal bundle to $S(\epsilon)$ to $T^1_Y\T(S)$ is the identity, up to an error of $o(1)$: we conclude that the map, to $T^1_Y\T(S)$, of tangents to harmonic map rays from $S(\epsilon)$  has the same degree as the map (described just above) to the tangent vectors to the Weil-Petersson  geodesics from $S(\epsilon)$ to $Y$. In particular, that map from normal bundle to $S(\epsilon)$ to $T^1_Y\T(S)$ has degree one.
 
 Finally, we exhaust $\T(S)$ by energy level spheres.  These form a continuous family, so the rays from any of these level spheres have the same degree, as a map from the normal bundle to the sphere to $T^1_Y\T(S)$.  So, choose $v \in T^1_Y\T(S)$: we can find a diverging family $X_t \in \T(S)$ so that the harmonic map ray from $X_t$ through $Y$ is tangent to $v \in T^1_Y\T(S)$. Then Theorem~\ref{thm:HR:limit:geodesic} provides for a subsequential limit of these rays which is a Thurston geodesic and, as a locally uniform limit of smooth curves, converges in $C^1$, and hence passes through $v \in T^1_Y\T(S)$.
  \end{proof}
 \color{black}

 %===================
  %===================
  \subsection{Two versions of geodesic flow of Thurston metric}
  \label{sec:geodesic:flow}
  In this subsection, we define two versions of the geodesic flow of the Thurston metric. The first version is defined using the exponential map obtained in this section. The second version is defined using Theorem \ref{thm:HR:SR}.

  The exponential map we obtained in this section allows us to define a Thurston geodesic flow  $$\psi_t: \T(S)\times \PML(S)\longrightarrow \T(S)\times\PML(S)$$
  as follows.
   For each pair $(Y,[\lambda])\in \T(S)\times\PML(S)$, let $\ESR_{Y}^{[\lambda]}:[1,\infty)\to\T(S)$ be the harmonic stretch ray which starts at $Y$ and converges to $[\lambda]\in\PML(S)$. Then we define the flow $\psi_t$ by setting $\psi_t(Y,[\lambda]):=(\ESR_{Y}^{[\lambda]}(e^{2t}), [\lambda])$: that this flow is well-defined follows from the second paragraph in Theorem~\ref{thm:exponentialmap}. In particular, every $\psi$-orbit is a harmonic stretch line. Moreover, it follows from  Theorem \ref{thm:exponentialmap} that every harmonic stretch line appears as a $\psi$-orbit.

 We next define a second version of the Thurston geodesic flow  
 $$\phi_t:\T(S)\times\PML(S)\to  \T(S)\times \PML(S).$$
 For each pair $(Y,[\lambda])\in\T(S)\times\PML(S)$, let $\HSR_{Y,[\lambda]}:[1,\infty)\to\T(S)$ be the harmonic stretch ray obtained as the limit of harmonic map rays $\HR_{X_t,Y}$ where $X_t$ degenerates along the harmonic map dual ray $\HDR_{Y,\lambda}$(Theorem \ref{thm:HR:SR}). We then define the flow $\phi_t$ by setting $\phi_t(Y,[\lambda])=(\HSR_{Y,[\lambda]}(e^{2t}),[\lambda])$. 
 In particular,   every $\phi$-orbit is a harmonic stretch line. However, $\phi$-orbits are \enquote{rare} in the following sense.   For any non-uniquely ergodic lamination $\lambda$, the projection to $\T(S)$ of the orbit of $(Y, [\lambda])$ under $\phi_t$ is independent of the transverse measure of $\lambda$ (by Lemma \ref{lem:dual:tree} and Theorem \ref{thm:minimalgraph:uniqueness}). Moreover, by \cite[Theorem 10.7]{Thurston1998}, for any $Y$ and any simple {\it closed} curve $\lambda$, the set $\mathscr{Z}_{Y,\lambda}$, which  consists of surfaces $Z\in\T(S)$ such that the maximally stretched lamination from $Y$ to $Z$ is $\lambda$, is an open subset of $\T(S)$. Consequently, for any fixed $Y\in\T(S)$, the union of projection of the $\phi$ orbits of $(Y,[\lambda])$, as $[\lambda]$ varies in $\PML(S)$, is a proper subset of $\T(S)$. 

  Concerning the continuity of  both flows, we have
 \begin{proposition}\label{prop:flow:continuity} 	The $\psi_t$-flow is continuous while  the $\phi_t$-flow is not.
 \end{proposition}
 \begin{proof}
 	That the $\psi_t$-flow is continuous  follows from Proposition \ref{prop:continuity2}.  To see that $\phi_t$ is not continuous, consider two disjoint simple closed curves $\alpha$ and $\beta$. Let  $\lambda_n=\alpha+\frac{1}{n}\beta$. Then $\lambda_n\to \alpha$ in $\ML(S)$. Notice that in the forward direction,  the orbit path $\{\phi_t(X,\alpha):t\in\R\}$ maximally stretches exactly along $\alpha$,  while   the orbit path $\{\phi_t(X,\lambda_n):t\in\R\}_n$ maximally stretches exactly along $\alpha\cup\beta$  for every $n\geq1$. By Theorem \ref{thm:maximally:stretched:lamination}, as $n\to\infty$,  the orbit paths $\{\phi_t(X,\lambda_n):t\in\R\}_n$ does not converge to the orbit path $\{\phi_t(X,\alpha):t\in\R\}$. In particular, $\phi_t$ is not continuous at $(X,\alpha)$.
 \end{proof}

 How does the earthquake flow interact with the second version $\phi_t$ of the Thurston geodesic flow? Do they define an action of the upper triangular subgroup of $SL(2,\R)$?
 From the construction of piecewise harmonic stretch lines in Section \ref{sec:stretchline:construction}, we know that the translates of harmonic stretch lines by the earthquake flow are piecewise harmonic stretch lines.  Here we show that they are also harmonic stretch lines.

 \begin{proposition}\label{prop:compatibility}
 Let $\mathbf{R}$ be a harmonic stretch line in $\T(S)$ which maximally stretches along a measured geodesic lamination $\lambda$. Let $\mathcal{E}_{\lambda}(\mathbf{R})$ be a translate of $\mathbf{R}$ by an earthquake directed by $\lambda$. Then $\mathcal{E}_{\lambda}(\mathbf{R})$ is also a harmonic stretch line.
 \end{proposition}
 \begin{proof}
 Let $Y,Z\in \mathbf{R}$ be two hyperbolic surfaces such that $\mathbf{R}$ proceeds from $Y$ to $Z$. Then by Theorem \ref{thm:unique:geodesic},
  there exists a sequence of harmonic maps $f_n:X_n\to Y$ with $X_n\in\T(S)$ which converges to a harmonic map $f:X\to Y\setminus \lambda$ from a punctured surface $X$ and that $(X,f)$ defines $\mathbf{R}$ in the sense of Theorem \ref{thm:generalized:stretchmap}. Let $\lambda_n$ be the horizontal measured foliation of $\Hopf(f_n)$. Let $X_n'$ be the unique Riemann surface in $\T(S)$ such that the Hopf differential of the harmonic map $f_n':X_n'\to \mathcal{E}_{\lambda}(Y) $ is also $\lambda_n$ (\cite{Wolf1998}).  By Lemma \ref{lem:precompactness}, the sequence of maps   $\{f_n'\}$ contains a  convergent subsequence, still denoted by  $\{f_n'\}$ for simplicity. Let $f':X'\to \mathcal{E}_{\lambda}(Y)\setminus\lambda'$ be the limit harmonic diffeomorphism. Then, since $\lambda_n$ limits on both $\lambda$ and $\lambda'$, we have that $\lambda=\lambda'$ (as geodesic laminations). Moreover, the horizontal foliations of $\Hopf(f)$ and $\Hopf(f')$ are the same. It then follows from Theorem \ref{thm:minimalgraph:uniqueness} --- using the identification $Y\backslash\lambda=\mathcal{E}_\lambda(Y)\backslash\lambda'=\mathcal{E}_\lambda(Y)\backslash\lambda$ --- that $X'=X$ and $f'=f$. In particular, $\mathcal{E}_\lambda(\mathbf{R})$ is a  harmonic stretch line defined by $f'$.
 \end{proof}
 \begin{remark}
  As a direct consequence, the earthquake flow and the \enquote{Thurston geodesic flow} $\phi_t$ are compatible, and hence define an action of the upper triangular subgroup of $SL(2,\R)$. More precisely,  let 
  $$a_t=\begin{pmatrix}
     e^t & 0 \\
      0& e^{-t}
   \end{pmatrix}, 
   ~h_r=\begin{pmatrix} 
     1 & r \\
      0& 1
   \end{pmatrix}.  $$
   We define $a_t(\lambda,Y):=\phi_t(\lambda,Y)$ and $h_r(\lambda,Y):=(\lambda,\mathcal{E}_{r\lambda}(Y))$.
    \end{remark}
\section{Concluding Remarks}\label{sec:concluding}

\subsection{Constructing geodesics in the Teichm\"uller space of hyperbolic surfaces with geodesic boundary}\label{subsec:surfacewithboundary} 

Let $S_{g,b}$ be an orientable surface of genus $g$ with $b$  boundary components. Let $\T(S_{g,b})$ be the Teichm\"uller space of hyperbolic surfaces with $b$ geodesic boundary components. There are presently three versions of a Thurston-type metric on $\T(S_{g,b})$ defined as follows. The first one is the so-called \textit{arc metric/distance} introduced by  Liu-Papadopoulos-Su-Th\' eret\cite{LPST2010}. Let $\mathcal{C}$ be the set of isotopy classes of simple closed cuves on $S_{g,b}$ and $\mathcal{A}$ the set of isotopy classes (rel $\partial S_{g,b}$) of (essential) simple  arcs on $S_{g,b}$  with endpoints on $\partial S_{g,b}$. The \textit{arc distance} is defined as: 
$$d_A(X,Y):=\log\sup_{\alpha\in\mathcal{C}\cup\mathcal{A}}\frac{\ell_Y(\alpha)}{\ell_X(\alpha)}.$$
The other two versions, introduced in \cite{AlessandriniDisarlo2019}, are defined via Lipschitz maps. For $X,Y\in\T(S_{g,b})$, let $\mathscr{L}(X,Y)$ be the set of Lipschitz maps from $X$ to $Y$ that commute with the marking up to homotopy. Let $\Lip(\phi)$ be the Lipschitz constant of $\phi\in\mathscr{L}(X,Y) $
Define
\begin{eqnarray*}
d_{L\partial}(X,Y)&:=&\log\inf\{\Lip(\phi):\phi\in\mathscr{L}(X,Y) \text{ with } \phi(\partial X)\subset \partial Y\},
\\
d_{Lh}(X,Y)&:=&\log\inf\left\{\Lip(\phi):
\begin{array}{cc}
    \phi\in\mathscr{L}(X,Y) \text{ with } \phi \\ \text{ a homeomorphism} 
\end{array} \right\}.
\end{eqnarray*}
Alessandrini-Disarlo\cite{AlessandriniDisarlo2019} showed that $d_A=d_{L\partial}$ on $\T(S_{g,b})$ and conjectured that \cite[Conjecture 1.8]{AlessandriniDisarlo2019} that $d_A=d_{Lh}$. Using  harmonic stretch lines, we verify this conjecture in the case of unpunctured surfaces.
\begin{theorem}
\label{thm:arc:Lip}
With notations as above, for any $X,Y\in\T(S_{g,b})$, we have
$d_A=d_{Lh}$. Moreover, the optimal Lipschitz constant from $X$ to $Y$ is always realized by a homeomorphism.
\end{theorem}
\begin{proof}
It is clear that $d_A\leq d_{Lh}$. To prove the theorem, it suffices to show that the optimal Lipschitz constant from $X$ to $Y$ is realized by a homeomorphsim from $X$ to $Y$.

 Let $X^d$ (resp. $Y^d$) be the  double of $X$  (resp. $Y$), obtained by gluing respectively the orientation-reversing isometric copy of $X$ (resp. $Y$) to $X$ along $\partial X$ (resp. to $Y$ along $\partial Y$). Consider the harmonic stretch line $[X^d,Y^d]$ in the Teichm\"uller space of the double of $S_{g,b}$.
 The doubling process induces an involution, denoted by $\iota$, on both  $X^d$ and $Y^d$.
 The uniqueness of harmonic stretch segment (Theorem \ref{thm:unique:HSR}) then implies the Lipschitz map  $\phi^d$ from $X^d$ to $Y^d$ induced by any admissible triple of $(X^d,Y^d)$ is symmetric about this  involution.  Hence $\phi^d$ descends to a Lipschitz homeomorphism $\phi$ from $X$ to $Y$ with the same Lipschitz constant as $\phi^d$. This implies that $d_{Lh}(X,Y)\leq \log \Lip (\phi)=  d_{Th}(X^d,Y^d)$. On the other hand, the double of any $L$-Lipschitz homeomorphism from $X$ to $Y$ gives an $L$-Lipschitz homeomorphism from $X^d$ to $Y^d$. Hence, using that the doubling operation provides candidate maps for the minimization problem for closed surfaces from candidates for the minimization problem for the surfaces-with-boundary $X$ and $Y$, we see that $d_{Lh}(X,Y)\geq d_{Th}(X^d,Y^d) = \log \Lip( \phi)$. Consequently, $d_{Lh}(X,Y)=d_{Th}(X^d,Y^d) = \log \Lip( \phi)$. Combining with the fact $d_A(X,Y)=d_{Th}(X^d,Y^d)$ proved by
   Liu-Papadopoulos-Su-Th\'eret \cite{LPST2010},  we know that $d_A(X,Y)=d_{Lh}(X,Y) = \log \Lip( \phi)$. 
\end{proof}

\subsection{Relation to orthogonal foliation introduced by Choi-Dumas-Rafi and Calderon-Farre}
 Let $Y\in\T(S)$ be a hyperbolic surface and  $\lambda$ a measured geodesic lamination. The harmonic stretch ray obtained in Theorem \ref{thm:HR:SR} is determined by a (possibly disconnected) punctured Riemann surface $X$ and a harmonic diffeomorphism $f:X\to Y\backslash\lambda$. The pushforward of the vertical foliation of the Hopf differential $\Hopf(f)$ extends to a measured foliation on $Y$ which is transverse to $\lambda$. On the other hand, there is an orthogonal measured foliation associated to $\lambda$ constructed in \cite{CDR2012,CalderonFarre2021}. If $\lambda$ is maximal, then these two measured foliations coincide.   A natural question is to consider the relationship between these two types of transverse measured foliations of $\lambda$ on $Y$ for non-maximal $\lambda$.

%======================
%======================
 \subsection{Optimal Lipschitz map from hyperbolic surface $Y$ to negatively curve surface $Z$}

In this paper, we have considered the sequence of minimizers of the energy difference function: $E(X,Y)-t E(X,Z)$ for constant curvature surfaces $Y$ and $Z$. It is natural to reflect on how these results might generalize to the case of negatively curved surfaces $Y$ and $Z$ for $0<t<\mathrm{Lip}(Y,Z)^{-2}$.

\appendix
  \section{
 Existence for the  generalized Jenkins-Serrin problem}\label{sec:appendix}
  
  A principal tool in this paper was an extension of the Jenkins-Serrin theory for minimal graphs in Euclidean three-space with asymptotic boundary values to minimal graphs over hyperbolic surfaces with boundary which took values in real trees, also with asymptotic boundary values.  The uniqueness theory of such graphs was described in Theorem~\ref{thm:minimalgraph:uniqueness} and used throughout the paper.
  In this appendix, we prove the existence part of this general Jenkins-Serrin problem; this plays a more minor role in our arguments --- only appearing in the proof of Lemma~\ref{lem:analyticity1}---and since the proof is rather lengthy, we relegate it to an appendix. The proof is somewhat involved, but as it turns on the structures of the foliations of the Hopf differentials, it also provides for some further development of the more technical themes of this paper.  We begin with the statement of the main result; recall Definition~\ref{def:admissibleTree:general} of \enquote{admissible dual tree} and the subsequent construction of \enquote{admissible partial boundary map}  from Section~\ref{sec:GeneralizedJenkinsSerrin}.
   \begin{theorem}\label{thm:minimalgraph:existence}
    Let $Y$ be a crowned hyperbolic surface.  Let ${T}$ be an admissible dual tree. Then there exists a $\pi_1(Y)$-equivariant minimal graph in $\widetilde{Y}\times {T}$ with a prescribed admissible partial boundary map.
    \end{theorem}

    The idea of the proof is to first construct a sequence of harmonic maps $X_n \to W$ to a fixed closed surface $W$ consisting of copies of $Y$ which have been glued together, with  controlled horizontal foliations of the associated Hopf differentials.  The domains $X_n$ for the approximating harmonic maps are found by demanding that the Hopf differentials have approximating maximal stretch laminations.  We then prove that the limit of any convergent subsequence  gives an equivariant minimal graph in $\widetilde{Y}\times T$ (and in fact a unique one by Theorem \ref{thm:minimalgraph:uniqueness}).  The proof will be divided into two cases:
    \begin{enumerate}
        \item[Case I:] $Y$ has no crown ends, and the measured foliation defining $T$ consists of half-infinite cylinders corresponding to  boundary components of $Y$ and compactly supported subfoliations;
        \item[Case II:] the general case. 
    \end{enumerate}
     
     Here the division is entirely for expositional reasons: the basic structure of the argument will be apparent in the technically simpler Case I. Case II extends the technique to a more complicated setting.
     
  Now, the foliations on surfaces that arise in (lifts of) minimal surfaces over (lifts of) crowned surfaces have several qualitative types: (i) they can be closed curves, or (ii) they can be part of a half-plane adjacent to a boundary leaf of a crown, or (iii) they can be a strip of bounded width with an end tending to an ideal point in a crown or an end spiralling around a boundary curve, or (iv) the leaves can be none of these: in the latter case, the leaves may be collecting into a compact portion of the surface. Of course, the technical heart of the proof comprises checking that the limiting Hopf differential on the limiting surface $X_\infty$ has the desired trajectory structure. One crucial estimate is a bound on the diameter of the \enquote{compact part}, but one also needs to make sure that the leaves, in each end, limit in the expected way.
   
 \subsection{Case I: Geodesic boundary}
 Let $F$ be the measured foliation on $Y$ whose lift to the universal cover defines $T$. Suppose that $F$ 
 consists of $a$ half-infinite cylinders $\{A_i\}_{1\leq i\leq a}$ with core curves $\{\alpha_i\}_{1\leq i\leq a}$ -- which are also the boundary components of $Y$ -- and $b$ compactly supported subfoliations $\{B_i\}_{1\leq i\leq b}$.  (In this particular model case, we may take $b=1$, but we retain the notation for the more general case.)
 Let $\beta_i$ be the measured lamination corresponding to $B_i$. 
 
   Let $W$ be a closed hyperbolic surface obtained by gluing an isometric copy ${Y}'$ of $Y$ to $Y$ along $\partial Y$ in a way that preserves the orientation of the two copies.  (This choice of orientation is not important for this case, but will be important in Case II, so we introduce it here.)
 Let $\beta'$ be the copy of $\beta\in\MF(Y)$ on $Y'$.  Consider the measured foliation  $\mu_n:=2n\sum_i \alpha_i+\beta+{\beta}'\in \MF(W)$.  Let $X_n\in \T(W)$ be the Riemann surface such that the horizontal foliation of the Hopf differential of the harmonic map $f_n:X_n\to W$ is $\mu_n$; as usual, this is guaranteed by \cite{Wolf1998}. Equivalently, the universal cover $\widetilde{X_n}$ of $X_n$ is the unique minimal graph in $\widetilde{W}\times T_n$ where $T_n$ is the tree dual to $\widetilde{\mu_n}$.

 Next we decompose $X_n$ as $A_{1,n}\cup \cdots \cup A_{k,n}\cup B_{n}\cup {B}'_{n}$ according to the horizontal foliation of $\Phi_n:=\Hopf(f_n)$, where $A_{i,n}$ is the maximal flat cylinder corresponding to the curve $\alpha_i$ and $B_{n}$ (resp. ${B}'_{n}$) is the (precompact) subsurface corresponding to $\beta$ (resp. ${\beta}'$). In particular, the width of $A_{i,n}$ is $2n$ for all $i$.  By Lemma \ref{lem:hyplength:quadraticnorm}, we see that for each $A_{i,n}$ and $B_n$, we have 
 \begin{eqnarray}\label{eq:area:B}
 2n\ell_{W}(\alpha_i)-C\leq &2\left\| \Phi_n\left|_{A_{i,n}}\right.\right\|& \leq 2n\ell_{W}(\alpha_i)+C,\\\notag
   \ell_{W}(\beta)-C\leq &2\left\| \Phi_n\left|_{B_{n}}\right.\right\|& \leq \ell_{W}(\beta)+C,\\\notag
   \ell_{W}(\beta')-C\leq &2\left\| \Phi_n\left|_{{B}'_{n}}\right.\right\|& \leq \ell_{W}(\beta')+C,
 \end{eqnarray}
 where $C$ is a constant depending on the topology of $W$.

 Let  $\overline{B_{n}}\subset X_n$ be the closure of $B_{n}$.
 \begin{lemma}\label{lem:diameter:bound}
   There exists a constant $D>0$ depending on $W$ such that for any $n>1$, the diameters of  $\overline{B_{n}}$ and $\overline{B_{n}'}$ are at most $D$ with respect to the $|\Phi_n|$-metric.
 \end{lemma}
 \begin{proof}
   We demonstrate the proof for $\{\overline{B_n}\}$. The proof of $\{\overline{B_n}'\}$ is similar.
   Suppose to the contrary that there exists a subsequence of $\{\overline{B_{n}}\}_{n\geq1}$ whose diameter goes to infinity. Without loss of generality, we may assume this subsequence is $\overline{B_{n}}$ itself.

   The area bound (\ref{eq:area:B}) of $B_{n}$ implies that the injectivity radius of every point of $B_{n}$ is at most $D_1$ for some constant $D_1$ depending on $W$. Recall that $\Phi_n$ has at most $2|\chi(W)|$ zeros. The assumption that the diameter of  $\overline{B_{n}}$ goes to infinity implies that there exists $p_n\in B_{n}$ such that
   $$ d_n:=\mathrm{dist}(p_n,\mathbf{Sing}(\Phi_n))\to\infty, $$
   where $\mathbf{Sing}(\Phi_n)$ is the zero set of $\Phi_n$ and $\mathrm{dist}(\cdot,\cdot)$ represents the $|\Phi_n|$-distance. By
   \cite[Lemma 5.1]{MasurSmillie1991}, the point $p_n$ is contained in some smooth closed regular geodesic $\eta_n$ with length $\ell(\eta_n)<2D_1$.  Then any point of $\eta_n$ has distance at least $d_n-2D_1$ from $\mathbf{Sing}(\Phi_n)$. By \cite[Lemma 5.2]{MasurSmillie1991}, we see that the (maximal) cylinder $C'_n$ containing $\eta_n$ has width at least $2d_n-4D_1$. Therefore, there exists a subcylinder $C_n\subset C_n'$, with $p_n \in C_n$, of width $2d_n-4D_1-2R$ such that every point of $C_n$ is of distance at least $R$ to $\mathbf{Sing}(\Phi_n^d)$. (Of course, the constant $R$ here will refer to the constant in Minsky's region $\mathscr{P}_R$.)  

   Next, we claim that $C_n\subset B_{n}$.  Otherwise, suppose to the contrary that $C_n$ is not contained in $B_n$, then $C_n$  would intersect some $A_{i,n}$ because $\cup_ i A_{i,n}$ separates $B_n$ and $B_n'$.  Since $C_n$ intersects both $B_n$ and $A_{i,n}$, it follows that $C_n$ is not horizontal.  Let $q\in A_{i,n}\cap C_n$. Then the closed geodesic $\gamma_q\subset A_{i,n}$ which contains $q$ (and is a core curve of $A_{i,n}$) would cross $C_n$ because $C_n$ is not horizontal. This implies the $|\Phi_n|$-length $\gamma_q$ is at least the width of $C_n$ which goes to infinity. Thus, the circumference of $A_{i,n}$ also goes to infinity. On the other hand, recall that, by the construction of the measured foliation $\mu_n$, the width of $A_{i,n}$ is $2n$.  Let $\delta_{i,n}$ be the core curve of $A_{i,n}$ which is of distance $n$ from the boundary of $A_{i,n}$. Then by Minsky's estimate together with the fact that $A_{i,n}$ is horizontal, we see that the $|\Phi_n|$-length of $\delta_{i,n}$ converges to $\ell_{W}(\alpha_i)/2$ as $n\to\infty$. This implies that the circumference of $A_{i,n}$ also converges to $\ell_{W}(\alpha_i)/2$ as $n\to\infty$, which yields a contradiction.

   The extremal length of $\eta_n$ on $B_{n}$ satisfies:
   $$ \Ext_{B_{n}}(\eta_n)\leq \frac{1}{\mathrm{Mod}(C_n)}\leq \frac{2D_1}{2d_n-4D_1-2R}. $$
   By Theorem \ref{thm:Minsky:energy}, we see that
   \begin{eqnarray*}
     E(f_n\left|\right._{B_{n}})&\geq & \frac{1}{2}\frac{\ell_{W}^2(\eta_n)}
   {\Ext_{B_{n}}(\eta_n)}\\
   &\geq & \frac{(\mathrm{Syst}(W))^2(2d_n-4D_1-2R)}{4D_1},
   \end{eqnarray*}
   which goes to infinity as $n\to\infty$. This contradicts (\ref{eq:area:B}) because $ E(f_n\left|\right._{B_{n}})\leq 2\|\Phi_n|_{B_{1,n}}\|+2\pi|\chi(W)|$.
 \end{proof}

 \begin{proof}[Proof of Theorem \ref{thm:minimalgraph:existence} Case I]
     Let $z_n$ be a zero of $\Phi_n$ in $B_n$. By Lemma \ref{lem:diameter:bound}, the subset $B_n$ is contained in the Minsky polygonal region $P_R$ of $X_n$ for $R>D$.
     By  Lemma \ref{lem:inj:lowerbound}, we see that the sequence of pointed flat surfaces $(X_n,\Phi_n, z_n)$ contains a subsequence which converges to pointed singular flat surface $(X_\infty,\Phi_\infty,z_\infty)$; here the flat metric $|\Phi_\infty|$ is induced by a meromorphic quadratic differential $\Phi_\infty$ which has a pole of order two at the pinching locus of $\alpha_i$, for each $i$. Moreover, by Lemma \ref{lem:diameter:bound},  the horizontal foliation of $\Phi_\infty$ consists of $a$ half-infinite cylinders $\{A_{i,\infty}\}_{1\leq i\leq a}$ with core curve $\{\alpha_i\}_{1\leq i\leq a}$, and the compactly supported measured foliation $\beta$. 
    Namely, the horizontal foliation of $\Phi_\infty$ (on the limit surface $X_{\infty}$) is equivalent to the measured foliation $F$, whose lift to the universal cover defines $T$.

     Up to a subsequence, we see that $\widetilde{f}_n:(\widetilde{X_n},\widetilde{p_n})\to \widetilde{W}$ converges to a harmonic map $\widetilde{f_\infty}:(\widetilde{X_\infty},
     \widetilde{p})\to \widetilde{W}$.     Notice that the image of a core curve of  $A_i$ under $f_\infty:X\to W$ arbitrarily closely approximates the geodesic representative of $\alpha_i$ on $W$ (as the distance of the core curve from the compact part grows large). Hence the image $\widetilde{f_\infty}(\widetilde{X_\infty})$ is exactly the lift $\widetilde Y$ of $Y$. Namely, $\widetilde{X_\infty}$ is an equivariant minimal graph in the product $\widetilde{Y}\times T$. 
 \end{proof}

     \subsection{Case II: general case}
     Let $F$ be the measured foliation on $Y$ whose universal cover defines $T$ by duality.    Suppose that $F$ comprises $a$ half-infinite cylinders $\{A_i\}_{1\leq i\leq a}$ foliated by closed leaves,
     $k$ half-planes $\{G_j\}_{1\leq j\leq k}$, $m$  bi-infinite strips $\{V_j\}_{1\leq j\leq m}$, and  $b$ compactly support subfoliations $\{B_i\}_{1\leq i\leq b}$.
  Notice that a bi-infinite strip may spiral around some half-infinite cylinder of $F$ that is not foliated by closed leaves (this will correspond to a second order pole, with non-real residue, of the limit meromorphic differential $\Phi_\infty$ to be constructed at the end of the proof). Let $\{A_{a+i}\}_{1\leq i\leq s}$ be the set of half-infinite cylinders of $F$ that are not foliated by closed leaves.
  Correspondingly, $Y$ contains $a+s$ geodesic boundary components $\alpha_i$ corresponding to half-infinite cylinders $A_i$, and $k$ ideal geodesic boundary arcs $\gamma_j$ corresponding to half-planes $G_j$, as well as $m$ ideal geodesic arcs $\xi_j$ corresponding to $V_j$, and $b$ compactly supported measured laminations $\beta_i$ corresponding to $B_i$.

    We start with the construction of a closed hyperbolic surface which consists of four copies of $Y$.   We first take a copy $Y'$ of $Y$ and then glue it back to $Y$ with a "shear" of amount $ t \neq 0$ along every ideal geodesic arc $\gamma_i$. This yields a hyperbolic surface $\widetilde{Y}$ with $2a+2s+k$ boundary components $\alpha_1,\alpha_1',\cdots,\alpha_{a+s},\alpha_{a+s}'$ and $\delta_1,\cdots,\delta_k$, where $\alpha_i'$ is the geodesic boundary component of the copy ${Y}'$ corresponding to $\alpha_i$, and where $\delta_i$ is the boundary component obtained by the shearing along the ideal geodesic arc $\gamma_i$ of $Y\subset \widetilde{Y}$.   Moreover, the $  t$-shearing along each ideal geodesic yields that all of $\delta_i$ have the same length $ |2t|$ (all that is important here is the length of each boundary component of $\widetilde{Y}$ is positive).  We then get a closed hyperbolic surface $W$ that is obtained by gluing an isometric copy $\widetilde{Y}'$ of $\widetilde Y$ to $\widetilde Y$ in an orientation-preserving way along each geodesic boundary component.  Notice that $W$ consists of four copies of $Y$. We denote these four copies by $Y,Y',Y'',Y'''$ such that the projection map $W\to \widetilde{Y}$ sends $Y$ and $Y''$ (resp. $Y'$ and $Y''')$ to $Y$ (resp. $Y'$).  Each of $\alpha_i,\gamma_i,\xi_i, \beta_i$ gives a copy on $Y',Y'',Y'''$, denoted by $\alpha_i',\alpha_i'',
    \alpha_i'''$; $\gamma_i',\gamma_i'',\gamma_i'''$; $\xi_i',\xi_i'',\xi_i'''$ and $\beta_i',\beta_i'',\beta_i'''$. Because of the gluing process, each of the pairs $\{\gamma_i,\gamma_i'\}$, $\{\gamma_i'',\gamma_i'''\}$, $\{\alpha_i,\alpha_i''\}$, $\{\alpha_i',\alpha_i'''\}$ are identified on $W$. Let $\gamma_i,\widehat{\gamma}_i,\alpha_i,{\alpha}'_i$ be the resulting geodesics on $W$.    Consider the geodesic lamination $\mu$:
    \begin{eqnarray*}
    \mu:=&&\left( \bigcup_{1\leq i\leq a+s}(\alpha_i\cup{\alpha}'_i)\right) \bigcup\left(\bigcup_{1\leq i\leq k}(\gamma_i\cup\widehat\gamma_i)\right)\bigcup\left(\bigcup_{1\leq i\leq k}\delta_i\right)\\
    &&\bigcup\left(\bigcup_{1\leq j\leq m}(\xi_j\cup\xi_j'\cup\xi_j''\cup\xi_j''')\right)
    \bigcup\left(\bigcup_{1\leq j\leq b}(\beta_j\cup\beta_j'\cup\beta_j''\cup\beta_j''')\right).
    \end{eqnarray*}
    The sublamination
    $$    \left( \bigcup_{1\leq i\leq a+s}(\alpha_i\cup{\alpha}'_i)\right) \bigcup\left(\bigcup_{1\leq i\leq k}(\gamma_i\cup\widehat\gamma_i)\right)\bigcup\left(\bigcup_{1\leq i\leq k}\delta_i\right)$$ cuts $W$ into four components which are exactly $Y,Y',Y'',Y'''$.

    Now we construct a sequence of measured laminations on $W$ whose supports converge to that of $\mu$.
     Notice that both $\gamma_i$  and $\widehat{\gamma}_i$ spiral around two of $\{\delta_1,\cdots,\delta_k\}$, say $\delta_{i1}$ and $\delta_{i2}$. Let $^\perp\gamma_i$ (resp. $^\perp\widehat\gamma_i$)  be the (open) geodesic arc on $W\backslash(\delta_1\cup\cdots\cup\delta_k)$ whose closure is orthogonal to both $\delta_{i1}$ and $\delta_{i2}$ and which is freely homotopic to $\gamma_i$ (resp. $\widehat{\gamma}_i$). (Here, in the universal cover, we have that a lift of $^\perp\gamma_i$ connects lifts of $\delta_{i1}$ and $\delta_{i2}$ while a lift of $\gamma_i$ might be the common asymptotic to those lifts of $\delta_{i1}$ and $\delta_{i2}$.)  We then close up $^\perp\gamma_i$ and $^\perp\widehat\gamma_i$ using subarcs of $\delta_{i1}$ and $\delta_{i2}$. Denote the resulting simple closed curve by $\gamma_{i,0}$.  Applying the same operation to  the   pairs $\{\xi_i,\xi_i''\}$ and $\{\xi_i',\xi_i'''\}$, we get simple closed curves $\xi_{i,0}$ and $\xi_{i,0}'$. The closing up process is chosen so that the geodesic representatives of $\{\gamma_{i,0},\xi_{j,0},\xi_{j,0}'\}_{1\leq i\leq k,1\leq j\leq m}$ are pairwise disjoint.

    For each of  $\delta_l$,  $\alpha_{a+r}$, and ${\alpha}'_{a+r}$,  let $T_{\delta_l}^n$,  $ T^n_{\alpha_{a+r}}$  and $ {T}^n_{{\alpha}'_{a+r}}$  be respectively the $n$ times right or left Dehn twists about  $\delta_l$,  $\alpha_{a+r}$, and ${\alpha}'_{a+r}$,  where the direction right or left is chosen according to the spiralling direction of $\{\gamma_i,\xi_{j}\}_{1\leq i\leq k,1\leq j\leq m}$ around $\delta_l$.
    Let
    \begin{eqnarray*}
     \gamma_{i,n}&:=&{T}^{n}_{{\alpha}'_{a+s}}\circ\cdots \circ {T}^{n}_{{\alpha}'_{a+1}}\circ {T}^{n}_{{\alpha}_{a+s}}\circ\cdots \circ {T}^{n}_{{\alpha}_{a+1}}\circ {T}^{n}_{\delta_k}\circ\cdots \circ {T}^{n}_{\delta_1}(\gamma_{i,0}),\\
     \xi_{j,n}&:=&{T}^{n}_{{\alpha}'_{a+s}}\circ\cdots \circ {T}^{n}_{{\alpha}'_{a+1}}\circ {T}^{n}_{{\alpha}_{a+s}}\circ\cdots \circ {T}^{n}_{{\alpha}_{a+1}}\circ {T}^{n}_{\delta_k}\circ\cdots \circ {T}^{n}_{\delta_1}(\xi_{j,0}),\\
     \xi'_{j,n}&:=&{T}^{n}_{\alpha'_{a+s}}\circ\cdots \circ {T}^{n}_{{\alpha}'_{a+1}}\circ {T}^{n}_{{\alpha}_{a+s}}\circ\cdots \circ {T}^{n}_{{\alpha}_{a+1}}\circ {T}^{n}_{\delta_k}\circ\cdots \circ {T}^{n}_{\delta_1}(\xi'_{j,0}).
     \end{eqnarray*}
     Then as $n\to\infty$, the geodesic representatives of $\gamma_{i,n}$, $\xi_{j,n}$, and $\xi_{j,n}'$ converge respectively to the closures of $\gamma_i\cup\widehat{\gamma}_i$, $\xi_j\cup\xi_j''$, and $\xi_j'\cup\xi_j'''$. (This is elementary hyperbolic geometry, as the lifts of the curves $\gamma_{i,n}$, $\xi_{j,n}$, and $\xi_{j,n}'$ to the hyperbolic plane have unique limits, hence the ones specified.)
     In particular, as $n\to\infty$, the geodesic representative of $\sum_{i=1}^{k}\gamma_{i,n}+
     \sum_{j=1}^{m}(\xi_{j,n}+\xi_{j.n}')$ converges to
     \begin{eqnarray*}
     &&\left( \bigcup_{1\leq i\leq s}(\alpha_{a+i}\cup{\alpha}'_{a+i})\right) \bigcup\left(\bigcup_{1\leq i\leq k}(\gamma_i\cup\widehat\gamma_i)\right)\bigcup\left(\bigcup_{1\leq i\leq k}\delta_i\right)\\
    &&\bigcup\left(\bigcup_{1\leq j\leq m}(\xi_j\cup\xi_j'\cup\xi_j''\cup\xi_j''')\right).
    \end{eqnarray*}

    The desired sequence of measured laminations $\mu_n$ on $W$ is defined as:
      \begin{eqnarray*}
    \mu_n:=&& n\sum_{1\leq i\leq a}(\alpha_i+{\alpha}_i') +n\sum_{1\leq i\leq k}\gamma_{i,n}\\
    &&+\textbf{w}_j\sum_{1\leq j\leq m}(\xi_{j,n}+\xi_{j,n}')
    +\sum_{1\leq j\leq b}(\beta_j+\beta_j'+\beta_j''+\beta_j'''),
    \end{eqnarray*}
     where $\mathbf{w}_j$ is the width of the horizontal strip $I_j$ of $\Phi$ corresponding to $\xi_j$. Compare the constructions  of $\mu$  and $\mu_n$. It is clear that as $n\to\infty$, the support of $\mu_n$ converges to that of $\mu$.

     Let $X_n\in\T(W)$ be the Riemann surface such that the horizontal measured foliation of the Hopf differential $\Phi_n$ of the harmonic map $f_n:X_n\to W$ is equivalent to $\mu_n$. (Again, this follows from \cite{Wolf1998}, .cf subsection~\ref{subsec:minimalsuspension}.) 
     We decompose $X_n$ according to the realization of  components of $\mu_n$ by  the horizontal foliation of $\Phi_n$  as in Table \ref{tab:corresondence}, where we list the subsurfaces in the second row corresponding to the subfoliations in the first row.
      \begin{table}[h]
       \centering
       \begin{tabular}{|c|c|c|c|c|}
         \hline
         subfoliations of $\mu_n$& $\alpha_i,{\alpha}'_i$&
         $\gamma_{i,n}$ &$\xi_{j,n},\xi_{j,n}'$
         &$\beta_j,\beta_j',\beta_j'',\beta_j'''$ \\
         \hline
         subsurfaces of $X_n$&$A_{i,n},{A}'_{i,n}$&
         $G_{i,n}$&$V_{j,n},V_{j,n}'$& $B_{j,n},B_{j,n}',B_{j,n}'',B_{j,n}'''$ \\
         \hline
       \end{tabular}
       \caption{\small{Correspondence between subfoliations of $\mu_n$ and subsurfaces of $X_n$.}}\label{tab:corresondence}
     \end{table}

     The remainder of this subsection considers the convergence of the family of harmonic maps $f_n:X_n\to W$, in a certain sense. We start with the following lemma, which is an analogue of Lemma \ref{lem:diameter:bound}.

      \begin{lemma}\label{lem:diameter:bound:general}
   There exists a constant $D>0$ depending on $W$ such that for any $n>1$ and any $1\leq j\leq b$, the diameter of each of $\overline{B_{j,n}}$, $\overline{B_{j,n}'}$, $\overline{B_{j,n}''}$,$\overline{B_{j,n}'''}$ is at most $D$ with respect to the $|\Phi_n|$-metric.
 \end{lemma}

    \begin{proof}
       We demonstrate the proof for the family $\{\overline{B_{j,n}}\}_{n\geq1}$. The proofs for the other three families are exactly the same. Suppose to the contrary that there exists a subsequence of $\{\overline{B_{j,n}}\}_{n\geq1}$ whose diameter goes  to infinity, say $\{\overline{B_{1,n}}\}_{n\geq1}$. Without loss of generality, we may assume this subsequence is  $\overline{B_{1,n}}$ itself. As in the proof of Lemma \ref{lem:diameter:bound}, there exists a sequence of  flat cylinders $C_n$ (not necessary maximal) meeting $\overline{B_{1,n}}$ such that
       \begin{itemize}
         \item the circumference of $C_n$  is less than some constant $D_1$;
         \item the width of $C_n$ is bigger than $2d_n-4D_1-2R$;
         \item every point of $C_n$ is of distance at least $R$ from the zero set of $\Phi_n$,
       \end{itemize}
       where $\{d_n\}$ is a divergent positive sequence and $R$ is a sufficiently large constant.

       Next, we claim that $C_n\subset B_{1,n}$ for sufficiently large $n$. Otherwise the boundary of $B_{1,n}$, being a union of saddle connections, would cross $C_n$. Then $\partial B_{1,n}\cap C_n$ has $\Phi_n$-length at least the width of $C_n$ which is bigger than $2d_n-4D_1 - 2R$. On the other hand, by the third property of $C_n$ mentioned above, every point of $\partial B_{1,n}\cap C_n$ is of distance at least $R$ from the zero set of $\Phi$. By Minsky's estimate, the image of  $\partial B_{1,n}\cap C_n$ is very close to the geodesic representative of the horizontal foliation $\mu_n$ of $\Phi_n$ on $W$. Hence by (\ref{eq:pullback:metric}),
       $$ \ell_W([\partial B_{1,n}])\geq \left|\partial B_{1,n}\cap C_n\right|_{|\Phi_n|}\geq 2d_n-4D_1 -2R$$
       which goes to infinity as $n\to\infty$.  Combining \ref{thm:Minsky:traintrack} with the first and the third property of $C_n$ mentioned above, we see that the image of the core curve $\eta_n$  of $C_n$ on $W$ has length bounded from above.  Since there are only finitely many such curves, we may assume, up to a subsequence, that $\eta_n$ belongs to a fixed homotopy class.  On the other hand, note that the image on $W$ of the core curve $\eta_n$  of $C_n$ has length bounded from below. Hence,  there exists a constant $\theta\in [0,\pi/2)$ such that each closed leaf of $C_n$ is of angle at most $\theta$ with the horizontal direction.  So the image of each closed leaf of $C_n$ is contained in the $\epsilon$ neighbourhood of  the geodesic representative of $\eta_n$, and the image of $\partial B_{1,n}\cap C_n$ spirals around the geodesic representative of $\eta_n$.  (Recall that $\beta_n$ belongs a fixed homotopy class as mentioned above.) Since  the support of $\mu_n$ converges to the support of $\mu$, this implies that the core curve $\eta_n$ of $C_n$ is contained in $\mu$, meaning that $C_n$ is a horizontal cylinder.  But each horizontal cylinder of $\Phi_n$ is either disjoint from $B_{1,n}$ or identical to $B_{1,n}$.  Therefore, $C_n\subset B_{1,n}$.

       The remaining part of  the proof is exactly the same as that of Lemma \ref{lem:diameter:bound}.
    \end{proof}

     \begin{proof}[Proof of Theorem \ref{thm:minimalgraph:existence} Case II] 
      Let $z_{i,n}$ be a zero of $\Phi_n$ in $B_{i,n}$. By Lemma \ref{lem:diameter:bound:general}, the compact set $B_{i,n}$ is contained in the polygonal region $P_R$ of $X_n$ for $R>D$.
     By  Lemma \ref{lem:inj:lowerbound} and \cite[Theorem A.3.1]{McMullen1989}, we see that the sequence of the pointed flat surfaces $(X_n,\Phi_n, z_{i,n})$  contains a subsequence which converges to a singular flat surface $(X_{i,\infty},\Phi_{i,\infty}, z_{i,\infty})$  induced by some meromorphic quadratic differential $\Phi_{i,\infty}$.
     Moreover, the sequence of harmonic maps $f_n:X_n\to W$ (sub)converges to a limit harmonic map $f_{i,\infty}:X_{i,\infty}\to W$ with Hopf differential $\Phi_{i,\infty}$.

     \vskip 5pt
     \textbf{Step 1. Classification of possible contributions  to $(X_{i,\infty},\Phi_{i,\infty})$ from $A_{j,n},A_{j,n}'$, $G_{j,n}$, $V_{j,n}$, $V_{j,n}'$. }
     Since the height of $A_{j,n}$ goes to infinity for each $j$ as $n\to\infty$, the circumference of $A_{j,n}$ converges to $\ell_W(\alpha_j)/2$ by Minsky's estimate. This means that the \text{possible contribution} of  $A_{j,n}$ to $(X_{i,\infty},\Phi_{i,\infty})$ is a half-infinite horizontal cylinder. Similarly, the \text{possible contribution} of each of  $A_{j,n}'$ to $(X_{i,\infty},\Phi_{i,\infty})$ is also a half-infinite horizontal cylinder.

     Notice that since the height of $G_{j,n}$ also goes to infinity as $n\to \infty$, the image of the core curve $\gamma_{j,n}$ under $f_n:X_n\to W$ is approximately the geodesic representative of $\gamma_{j,n}$ on $W$. Recall that as $n\to\infty$, the curves $\gamma_{j,n}$ converge to the union of $\gamma_j \cup \widehat\gamma_j$ and the two closed curves from $\{\delta_i\}$ to which they spiral. Therefore,  as $n\to\infty$,  the circumference of $G_{j,n}$ for each $j$ converges to $\ell_W(\gamma_{j})/2$ which is infinite.  This means that the \text{possible contribution} of $G_{j,n}$  to $(X_{i,\infty},\Phi_{i,\infty})$ is a half-plane.

    As for $V_{j,n}$, by Lemma \ref{lem:hyplength:quadraticnorm}, we get
     $$2\|\Phi_n|_{V_{j,n}}\|\geq \ell_W(\xi_{j,n}) -C, $$
     which diverges as $n\to\infty$.
     Combining this with the fact that the height of $V_{j,n}$  is always $\mathbf{w}_j$, we see that for each $i$, the circumference  $\ell(V_{j,n})$ of $V_{j,n}$ diverges as $n\to\infty$. This implies that  the \text{possible contribution} of $V_{j,n}$  to $(X_{i,\infty},\Phi_{i,\infty})$ is an infinite  strip  with height $\mathbf{w}_j$. Similarly,  the \text{possible contribution} of $V_{j,n}'$  to $(X_{i,\infty},\Phi_{i,\infty})$ is also an infinite  strip of  height $\mathbf{w}_j$.

      \vskip 5pt
      \textbf{Step 2. Showing that  $f_{i,\infty}(X_{i,\infty})$ is disjoint from any of $$\{\alpha_l,\alpha_l',\gamma_j,
      \widehat{\gamma}_j:1\leq l\leq a+s, 1\leq j\leq k\}.$$ }
      
      We first show that $\alpha_l\cap f_{i,\infty}(X_{i,\infty})=\emptyset$ for $1\leq l\leq a$. (As a reminder, these first $a$ geodesics come from purely horizontal cylinders -- there is no spiraling of the horizontal foliation.) Let $L_{l,n}$ be the central core curve of $A_{l,n}$, namely the core curve whose distance to $\partial A_{i,n}$ is $n/2$. Then by Theorem \ref{thm:Minsky:traintrack}, we see that $f_n(L_{l,n})$ converges to $\alpha_l$ as $n\to\infty$.  Now suppose that $\alpha_l\cap f_{i,\infty}(X_{i,\infty})\neq \emptyset$. Let $p $  be a point in the intersection set. 
      Then there exists a neighbourhood $U$ of $p$ with $U\subset f_{i,\infty}(X_{i,\infty})$. Since $f_n:X_n\to W$ converges to $f_{i,\infty}:X_{i,\infty}\to W$, it follows that there exists small neighbourhoods $U_n\subset X_n$ with $f_n(U_n)\subset U$ and also with $U_n$ approximating some fixed region in $X_{i,\infty}$ and hence at a uniformly bounded distance from the zeroes of $\Phi_n$.   On the other hand, the discussion above implies that there exists a sequence of points $p_n\in L_{l,n}$  whose distance to the zeros of $\Phi_n$ diverges such that $f_n(p_n)\to p\in\alpha_l$. In particular, $f_n(p_n)\in f_n(U_n)$ but $p_n\notin U_n$.  This contradicts the fact that $f_n$ is a homeomorphism, proving that $\alpha_l\cap f_{i,\infty}(i,\infty)=\emptyset$ for $1\leq l\leq a$.   Taking similar core curves of $A_{l,n}'$ and $G_{j,n}$ we see that $\alpha_l'\cap f_{i,\infty}(i,\infty)=\emptyset$ for $1\leq l\leq a$ and $(\gamma_j\cup \widehat{\gamma}_j)\cap f_{i,\infty}(i,\infty)=\emptyset $ for $1\leq j\leq k$. 
      
      It remains to show that $(\alpha_{a+l}\cup\alpha_{a+l}')\cap f_{i,\infty}(X_{i,\infty})=\emptyset$ for $1\leq l\leq s$.
      Recall that, by the way we constructed the curves $\xi_{j,n}$ and $\xi'_{j,n}$, the Hausdorff limit of $\bigcup_{j=1}^m (\xi_{j,n}\cup\xi'_{j,n})$ contains $\bigcup_{l=1}^{s}(\alpha_{a+l}\cup{\alpha}'_{a+l})$.
      This implies that for each $\alpha_{a+l}$ and ${\alpha}'_{a+l}$  with $1\leq l\leq s$, there exists $1\leq j_l\leq k$ such that  
      $ \alpha_{a+l} $ and $\alpha_{a+l}'$ are respectively contained in the 
      Hausdorff limit of $\xi_{j_l,n}$ and $\xi_{j_l,n}'$ as $n\to\infty$.  Consider the corresponding cylinders $V_{j_l,n}$ and $V_{j_l,n}'$. By Step 1 we know that both of their circumferences  diverge to infinity as $n\to\infty$. Now, combining Theorem \ref{thm:Minsky:polygon}, Lemma \ref{lem:inj:lowerbound} and the fact that $V_{j_l,n}$ and $V_{j_l,n}'$ have fixed width $\mathbf{w}_{j_l}$, we see that for any $R>R_0$, there exists $c>0$  such that each of $V_{j_l,n}$ and $V_{j_l,n}'$ crosses the Minsky's polygonal region $\mathscr{P}_R$ at most $c$ times. 
   Consider the intersection $\mathscr{P}_R\cap V_{j_l,n}$ (resp. $\mathscr{P}_R\cap V_{j_l,n}'$). By Theorem \ref{thm:Minsky:polygon}, the area of $\mathscr{P}_R$ is at most $CR^2$, and the horizontal segments of $\partial \mathscr{P}_R$ have length at most $K_1R$. Combining with the fact that both  $ V_{j_l,n}$ and $ V_{j_l,n}'$ has fixed width $\mathbf{w}_{j_l}$, we see that  the lengths of horizontal leaves of every component of $\mathscr{P}_R\cap V_{j_l,n}$ (resp. $\mathscr{P}_R\cap V_{j_l,n}'$) are bounded above by some constant $C_1$. (Here we apply the area bound to substrips of $V_{j_l,n}$ that are contained in the interior of $\mathscr{P}_R$ and the boundary estimate to substrips of $V_{j_l,n}$ that only meet the boundary $\partial\mathscr{P}_R$; the two bounds together yield the uniform bound $C_1$.)  Therefore, the central core curve of $V_{j_l,n}$ ($V_{j_l,n}'$)  has total length at least $\ell(V_{j_l,n})-cC_1$ (resp. $\ell(V_{j_l,n}')-cC_1$) outside $\mathscr{P}_R$, where $\ell(V_{j_l,n})$ and $\ell(V_{j_l,n}')$ are respectively the $\Phi_n$-circumferences of $V_{j_l,n}$ and $V_{j_l,n}'$.  
      This then implies that for sufficiently large $n$, the central core curve of $V_{j_l,n}$ (reps.  $V_{j_l,n}'$) contains a subsegment   of length at least $\frac{\ell(V_{j_l,n})}{c}-C_1$ (resp. $\frac{\ell(V_{j_l,n}')}{c}-C_1$) which is contained in the complement $X_n\setminus\mathscr{P}_R$ of $\mathscr{P}_R$.
      By Theorem \ref{thm:Minsky:traintrack}, the images of these segments are contained in the $\epsilon_R$-neighbourhood of, and are nearly parallel to, $\xi_{j_l,n}\cup \xi_{j_l,n}'$ which converges to $\xi_{j_l}\cup \xi_{j_l}'\cup\xi_{j_l}''\cup \xi_{j_l}'''\cup \alpha_{a+l}\cup \alpha_{a+l}'$ as $n\to\infty$: here we have chosen $R$ so that $\epsilon_R$ is sufficiently small.   Since both $\ell(V_{j_l,n})$ and $\ell(V_{j_l,n}')$ diverge, it follows that $\alpha_{a+l}\cup \alpha_{a+l}'$ is contained in the Hausdorff limit of the images of these segments. 
      Now suppose that $(\alpha_{a+l}\cup\alpha_{a+l}')\cap f_{i,\infty}(X_{i,\infty})\neq \emptyset$. Similarly to our previous argument, let $p$  be a point in the intersection set.
      Then there exists a neighbourhood $U$ of $p$ with $U\subset f_{i,\infty}(X_{i,\infty})$. Since $f_n:X_n\to W$ converges to $f_{i,\infty}:X_{i,\infty}\to W$, it follows that there exists small neighbourhoods $U_n\subset X_n$ with $f_n(U_n)\subset U$ but with $U_n$ at uniformly bounded distance from the zeroes of $\Phi_n$.   On the other hand, the discussion above implies that there exists a sequence of points $p_n$ in the union of the central core curves of $V_{j_l,n}$ and $V_{j_l,n}'$ whose distance to the zeros of $\Phi_n$ diverges such that $f_n(p_n)\to p\in\alpha_l\cup\alpha_l'$. In particular, $f_n(p_n)\in f_n(U_n)$ but $p_n\notin U_n$.  This contradicts the fact that $f_n$ is a homeomorphism, proving that $(\alpha_{a+l}\cup\alpha_{a+l}')\cap f_{i,\infty}(X_{i,\infty})=\emptyset$ for $1\leq l\leq s$.

 Similarly, we see that $f_{i,\infty}(X_{i,\infty})\cap \gamma_j=\emptyset$ and $f_{i,\infty}(X_{i,\infty})\cap \widehat{\gamma}_j=\emptyset$ for $1\leq j\leq k$.
     
           \vskip 5pt
     \textbf{Step 3. Showing that the image $f_{i,\infty}(X_{i,\infty})$  is exactly one of  $Y,Y',Y''$ , and $Y'''$.} Recall that $\Phi_{i,\infty}$ is a meromorphic differential with infinite area. Let $p$ be a pole of $\Phi_{i,\infty}$ of order at least two.

      If $p$ has order at least three, then there exists a neighbourhood $U(p)$ near $p$ such that  $(X_{i,\infty},\Phi_{i,\infty})$ is realized as a union of half-planes and strips. It follows from Step 1 that every such half-plane is a limit of  $\{G_{j,n}\}_{g\geq1,1\leq j\leq k}$  and that every strip is a limit of $\{V_{j,n}, V_{j,n}'\}_{n\geq 1, 1\leq j\leq m}$.   By Theorem \ref{thm:Minsky:traintrack}, the neighbourhood  $U(p)$ is mapped by $f_{i,\infty}$ to a crown end whose ideal geodesic arcs are a subset of $\{\gamma_j\}_{1\leq j\leq m}$.
   
     If $p$ is a second order pole with real residue, then a neighborhood of $p$ defines a horizontal  half-infinite cylinder  $C(p)$  near $p$. By the analysis in Step 1, we see that this cylinder is a limit of    $\{A_{j,n}, A_{j,n}'\}_{n\geq1,1\leq j\leq a}$. By Theorem \ref{thm:Minsky:traintrack}, the  image $f_{i,\infty}(C(p))$  is an one-sided neighbourhood of some curve in $\{\alpha_l,\alpha_l'\}_{1\leq l\leq a}$.

    If $p$ is a second order pole with non-real complex residue, this provides for a non-horizontal half-infinite cylinder $C$ near $p$. Let $\omega_d\subset C$ be the core curve whose distance to the compact boundary of $C$ is $d$.  If $C$ is vertical, then by Theorem \ref{thm:Minsky:traintrack} we see that the length of $f_{i,\infty}(\omega_d)\subset W$ converges to zero as $d\to\infty$. This contradicts the fact that $W$ is a closed hyperbolic surface with the shortest closed geodesic having positive length.

       We are left with the case that $C$ is neither horizontal nor vertical.  Let $s$ be the slope of $C$. Then $0<|s|<\infty$. Since $C$ is not horizontal, there exists some horizontal infinite strip crossing it. By step 1, every such strip is the limit of  $\{V_{j,n},V_{j,n}'\}_{n\geq1}$. Without loss of generality, we assume that it is the limit of $\{V_{j,n}\}_{n\geq1}$.  Let $\ell(C)$ be the circumference of $C$. Recall that $(X_{i,\infty},\Phi_{i,\infty},z_{i,\infty})$ is a limit of $(X_n,\Phi_n,z_n)$. There exists a non-horizontal cylinder $C_n$ on $X_n$  satisfying the following:
      \begin{itemize}
          \item The slope $s_n$ of $C_n$ converges to $s$ as $n\to\infty$.
          \item The circumference $\ell(C_n)$ of $C_n$ converges to $\ell(C)$ as $n\to\infty$.
          \item The width $w_n$ goes to infinity as $n\to\infty$.
          \item The core curve of $C_n$ is homotopic to that of $C$ (because $f_n$ is homotopic to the identity).
          \item The cylinder $V_{j,n}$  crosses $C_n$ for sufficiently large $n$ (because the limit of $V_{j,n}$ crosses $C$).
      \end{itemize}
       Combining the five properties of $C_n$ mentioned above and Theorem \ref{thm:Minsky:traintrack}, we see that 
             for large $n$, the image $f_n(C_n)$ on $W$ meets a neighbourhood of the closed geodesic on $W$ homotopic to the core curve of $C$. Moreover, the image $f_n(\xi_{j,n})$ spirals around this geodesic nearly $\frac{w_n }{\ell(C_n)/|s_n| }$ times (recall that $\xi_{j,n}$ is the core curve of $V_{j,n}$). As $n\to\infty$, the limit of $f_n(\xi_{j,n})$ spirals infinitely many times around the geodesic homotopic to the core curve of $C$. Since $f_n$ is homotopic to the identity, the limit of $\xi_{j,n}$ also spirals infinitely many times around the geodesic homotopic to the core curve of $C$. On the other hand, by the construction of $\xi_{j,n}$ we know that the only curves around which $\xi_{j,n}$ spirals infinitely many times is some curve in $\{\alpha_{a+l},\alpha_{a+l}',\delta_j\}_{1\leq l\leq s, 1\leq j\leq k}$. Therefore, the core curve of $C$ is homotopic to some curve in $\{\alpha_{a+l},\alpha_{a+l}',\delta_j\}_{1\leq l\leq s, 1\leq j\leq k}$. It then follows from Theorem \ref{thm:Minsky:traintrack} that $f_{i,\infty}(C)$ is a one-sided neighbourhood of this curve.

      Summarizing  the above discussion in this step, we see that the image $f_{i,\infty}(X_{i,\infty})$ is a crowned surface bounded by some curves from $\{\alpha_l,\alpha_l',\gamma_j,\widehat{\gamma}_j,\delta_j:1\leq l\leq a+s,1\leq j\leq k\}$. Combining this with Step 2, we see that  $f_{i,\infty}(X_{i,\infty})$ is one of  $Y,Y'$ $Y''$ and $Y'''$, the four components of the complement of 
     $$    \left( \bigcup_{1\leq i\leq a+s}(\alpha_i\cup{\alpha}'_i)\right) \bigcup\left(\bigcup_{1\leq i\leq k}(\gamma_i\cup\widehat\gamma_i)\right)\bigcup\left(\bigcup_{1\leq i\leq k}\delta_i\right).$$

    \vskip 5pt
    \textbf{Step 4. Determining the horizontal foliation of $\Phi_\infty$. }
    For any $1\leq i\leq b$, by Step 1 and Step 3, we see that   $\Phi_{i,\infty}$ contains a compactly supported foliation $B_i$ corresponding to the lamination $\beta_i$ and that the image of  $f_{i,\infty}:(X_{i,\infty},\Phi_{i,\infty})\to W$ is $Y$ since $\beta_i\subset Y$ (and not one of the copies $Y'$, $Y''$ or $Y'''$). 
    Since  $\beta_i\subset Y$ (and, again, not one of the copies $Y'$, $Y''$ or $Y'''$) for all $i=1,2,\cdots,b$, 
    we then see that we may apply a diagonal argument and choose a sequence of points $z_n \in \cup_{1\leq i\leq b}\overline{B_{i,n}} \subset X_{n}$ so that there exists a subsequence of $f_n:(X_n,\Phi_n,z_n)\to W$ which converges to a harmonic map $f_\infty:(X_\infty,\Phi_\infty,z_\infty) \to Y$ such that the  horizontal foliation of the Hopf differential $\Phi_\infty$ contains all of the laminations $\beta_1,\cdots,\beta_b$.
    By Step 1, we see that the horizontal foliation of $\Phi_\infty$ contains $a$ half-infinite cylinders corresponding to $\{A_{i,n}\}_{1\leq i\leq a}$, $k$ half-planes corresponding to $\{G_{j,n}\}_{1\leq j\leq k}$, and $b$ compactly supported foliations corresponding to $\{B_{j,n}\}_{1\leq j\leq b,}$. It remains to
    consider the contributions from $\{V_{j,n},V_{j,n}'\}_{1\leq j\leq m}$ and $\{B_{j,n}',B_{j,n}'',B_{j,n}'''\}_{1\leq j\leq b}$.

    By the construction of $\xi_{j,n}$ for each $1\leq j\leq m$,  (the geodesic representative of) the image of the core curve $\xi_{j,n}$ of $V_{j,n}$ converges to the union of $\xi_j\cup \xi_j''$ and the closed geodesics to which they spiral; this becomes clear when one lifts the families to $\H^2$. Combined with the fact that $\xi_j\subset Y$,  this means that some portion of $V_{j,n}$ survives in  $\Phi_\infty$. By Step 1, the contribution from $\{V_{j,n}\}_{n\geq1}$ is an infinite strip. Consequently, the horizontal foliation of $\Phi_\infty$ contains $m$ strips corresponding to $\{V_{j,n}\}_{1\leq j\leq m}$. On the other hand, if $\Phi_\infty$ contains some contribution from $\{V_{j,n}'\}_{n\geq 1}$, then by Step 1 this contribution is an infinite strip whose image on $W$ contains curves homotopic to  $\gamma_i'$ or $\gamma_i'''$; recall that we have identified $\gamma_i'''$ and $\gamma_i''$.  But neither $\gamma_i'$  nor $\gamma_i'''$ is homotopic to some curve  contained in $Y$. Hence the horizontal foliation of $\Phi_\infty$ contains no contribution from any of  $\{V_{j,n}'\}_{1\leq j\leq m,n\geq 1}$.

    Finally, suppose that the horizontal foliation of $\Phi_\infty$ contains some contribution from  $\{\beta_j',\beta_j'',\beta_j'''\}_{1\leq j\leq b,n\geq1},$ say $\beta_j'$. Then by Lemma \ref{lem:diameter:bound:general}, it contains the whole of $\beta_j'$. Correspondingly, $Y=f_\infty(X_\infty)$ contains $f_\infty(\beta_j')$, which is homotopic to $\beta_j'$ (on $Y'$). 
    This contradicts the fact that $\beta_j'$ is contained in $Y'$ instead of $Y$. Therefore, the horizontal foliation of $\Phi_\infty$ does not contain any contribution from $\{\beta_j',\beta_j'',\beta_j'''\}_{1\leq j\leq b,n\geq1}$.
   
   Summarizing the discussion above, we see that the horizontal foliation of $\Phi_{i,\infty}$ consists of $a$ half-infinite cylinders corresponding to $\{\alpha_i\}_{1\leq i\leq a}$, $k$ half-planes corresponding to $\{\gamma_j\}_{1\leq j\leq k}$, $m$ bi-infinite strips corresponding to $\{\xi_j\}_{1\leq j\leq m}$ and $b$ compactly supported foliations corresponding to $\{\beta_j\}_{1\leq j\leq b}$.
   
      \vskip5pt
     \textbf{Step 5. Constructing a minimal graph.} By Step 4, we know  that the horizontal foliation of $\Phi_{\infty}$ consists of $a$ half-infinite cylinders corresponding to $\{\alpha_i\}_{1\leq i\leq a}$, $k$ half-planes corresponding to $\{\gamma_j\}_{1\leq j\leq k}$, $m$ bi-infinite strips corresponding to $\{\xi_j\}_{1\leq j\leq m}$ and $b$ compactly supported foliations corresponding to $\{\beta_j\}_{1\leq j\leq b}$. In other words,
    the horizontal foliation of $\Phi_\infty$  is equivalent to the measured foliation $F$, whose lift to the universal cover defines $T$. Therefore, $X_\infty$ is a minimal graph in $\widetilde{Y}\times T$. This completes the proof.
     \end{proof}

    \remark In the above proof, we implicitly assume that the measured foliation $F$ has non-trivial compact components.   If the measured foliation $F$ has no compact component (i.e. the foliation comprises only half-infinite cylinders, half-planes, and infinite strips), then we simply take $z_n$ to be an arbitrary zero of $\Phi_n$. 

For the convenience of reference, we also summarize the construction in this appendix as follows:

\begin{proposition}\label{prop:appendix:convergence}
     For any crowned hyperbolic surface $Y$ and any admissuble measured foliation $F$ on $Y$, there exist
      \begin{itemize}
          \item a closed hyperbolic surface $W$ and a chain recurrent geodesic lamination $\lambda$ on $W$ such that $W\backslash\lambda$ is the union of two or four isometric copies of ${Y}$, and
          \item a sequence of Riemann surfaces $X_n\in\T(W)$ and harmonic diffeomorphisms $f_n:X_n\to W$ homotopic to the identity which converges to a harmonic diffeomorphism $f_\infty:X_\infty\to W\backslash\lambda$ (in the sense of Definition \ref{defn:limitharmonicmap}) from some punctured Riemann surface $X_\infty$, such that
         on each component $X_\infty'$ of $X_\infty$, the horizontal measured foliation of the Hopf differential of $f_\infty|_{X_\infty'}:X_\infty'\to Y$  is equivalent to $F$.
      \end{itemize}
\end{proposition}

   The ideas in this section also lead to the following characterization of Thurston stretch lines.
   \begin{corollary}
    \label{cor:Thurston:harmonic:stretch}
    Let $Y\in\T(S)$ be a hyperbolic surface and $\lambda$ be a maximal geodesic lamination. 
    Let $\SR_{Y,\lambda}$ be the Thurston stretch line determined by $Y$ and $\lambda$. Then  $\SR_{Y,\lambda}$ is a harmonic stretch line if and only if $\lambda$ is chain-recurrent.
   \end{corollary}
    \begin{proof}
     Suppose that $\SR_{Y,\lambda}$ is a harmonic stretch line. Then by definition, there exists a sequence of harmonic map rays $\HR_{X_n,Y}$ for some $X_n\in\T(S)$ that converges to $\SR_{Y,\lambda}$ as $n\to\infty$. 
     
     By Lemma  \ref{lem:precompactness}, the sequence of harmonic maps $f_n:X_n\to Y$ (sub)converges to a harmonic diffeomorphism $f_\infty:X_\infty\to Y\setminus \mu$ for some chain-recurrent lamination $\mu\subset \lambda'$.
     Since by assumption $\SR_{Y,\lambda}$ is the limit of $\HR_{X_n,Y}$, we see by Proposition \ref{prop:limit:HR:PSL}(ii) that $\SR_{Y,\lambda}$ is the  piecewise harmonic stretch line constructed from $f_\infty:X_\infty\to Y\setminus \mu$ in Theorem \ref{thm:generalized:stretchmap}. In particular, the harmonic stretch line $\SR_{Y,\lambda}$ maximally stretches exactly $\mu$. Therefore, $\mu=\lambda$. This implies that $\lambda$ is chain-recurrent.

     Now we turn to the other direction. Suppose that $\lambda$ is chain-recurrent. Then there exists a sequence of multicurves $\alpha_n$ which converges to $\lambda$ in the Hausdorff topology (on the set of geodesic laminations on $Y$), as $n\to\infty$. 
     Let $X_n$ be the Riemann surface such that the horizontal foliation of the Hopf differential of the harmonic map $X_n\to Y$ is $n\alpha_n$. Let $\HR_{X_n,Y}$ be the corresponding harmonic map ray. By Lemma \ref{lem:precompactness}, we see that $\{\HR_{X_n,Y}\}_{n\geq1}$ contains a subsequence which converges to some harmonic stretch line $\HSR$. 
     Applying the idea of the proof of Theorem  \ref{thm:minimalgraph:existence}, we conclude  that the limiting harmonic stretch line $\HSR$ maximally stretches along $\lambda$. The assumption that $\lambda$ is maximal then implies that $\HSR=\SR_{Y,\lambda}$ (cf. Lemma \ref{lem:id:PSR:SR}).
    \end{proof}

%%==============================
 % BibTeX References
\bibliographystyle{alpha}

\bibliography{main}   % name my BibTeX data bases
%\nocite{*} %List no cited references

\end{document}